\documentclass[reqno,tbtags,a4paper,12pt]{amsart}

\usepackage[in]{fullpage}
\usepackage[matrix,arrow,curve]{xy}
\usepackage{tikz,tikz-cd}
\usetikzlibrary{decorations.pathreplacing}
\usetikzlibrary{decorations.markings}

\usepackage{amsmath,amssymb,amsthm,amsfonts,amscd,array,longtable,multirow,afterpage}
\usepackage{graphics,graphicx}

\usepackage{mathrsfs}
\renewcommand{\mathscr}{\mathcal}

\newcommand{\bS}{\mathbb{S}}

\newcommand{\Mod}{\mathop{\mathrm{Mod}}\nolimits}
\newcommand{\PMod}{\mathop{\mathrm{PMod}}\nolimits}

\newcommand{\Ann}{\mathop{\rm Ann}\nolimits}
\newcommand{\Homeo}{\mathop{\mathrm{Homeo}}\nolimits}
\newcommand{\Sp}{\mathrm{Sp}}
\newcommand{\SL}{\mathrm{SL}}

\newcommand{\Act}{\mathop{\mathrm{Act}}\nolimits}

\newcommand{\Out}{\mathop{\mathrm{Out}}\nolimits}

\newcommand{\G}{\mathcal{G}}
\newcommand{\tG}{\widetilde{\mathcal{G}}}
\newcommand{\I}{\mathcal{I}}
\newcommand{\K}{\mathcal{K}}
\newcommand{\CL}{\mathcal{L}}
\newcommand{\D}{\mathcal{D}}
\newcommand{\Tot}{\mathop{\mathrm{Tot}}\nolimits}
\newcommand{\tot}{\mathrm{tot}}

\newcommand{\sgn}{\mathop{\mathrm{sgn}}\nolimits}

\newcommand{\Spl}{\mathop{\mathrm{Split}}\nolimits}
\newcommand{\Z}{\mathbb{Z}}
\newcommand{\Q}{\mathbb{Q}}

\newcommand{\F}{\mathbb{F}}

\newcommand{\R}{\mathbb{R}}
\newcommand{\bA}{\mathbf{A}}
\newcommand{\fA}{\mathfrak{A}}

\newcommand{\CA}{\mathcal{A}}
\newcommand{\CB}{\mathcal{B}}
\newcommand{\tCB}{\widetilde{\CB}}
\newcommand{\CC}{\mathcal{C}}
\newcommand{\oCA}{\bar{\mathcal{A}}}
\newcommand{\CM}{\mathcal{M}}
\newcommand{\CE}{\mathcal{E}}
\newcommand{\CF}{\mathcal{F}}
\newcommand{\CG}{\mathcal{G}}
\newcommand{\CH}{\mathcal{H}}
\newcommand{\CN}{\mathcal{N}}
\newcommand{\CR}{\mathcal{R}}
\newcommand{\CP}{\mathcal{P}}

\newcommand{\CS}{\mathcal{S}}
\newcommand{\CT}{\mathcal{T}}
\newcommand{\CU}{\mathcal{U}}
\newcommand{\CV}{\mathcal{V}}
\newcommand{\CW}{\mathcal{W}}
\newcommand{\CX}{\mathcal{X}}
\newcommand{\CY}{\mathcal{Y}}
\newcommand{\CZ}{\mathcal{Z}}
\newcommand{\tCC}{\widetilde{\mathcal{C}}}
\newcommand{\tD}{\widetilde{\mathcal{D}}}
\newcommand{\tCE}{\widetilde{\mathcal{E}}}
\newcommand{\tCF}{\widetilde{\mathcal{F}}}
\newcommand{\tCK}{\widetilde{\mathcal{K}}}
\newcommand{\tCM}{\widetilde{\mathcal{M}}}
\newcommand{\tCN}{\widetilde{\mathcal{N}}}
\newcommand{\tCP}{\widetilde{\mathcal{P}}}
\newcommand{\tCQ}{\widetilde{\mathcal{Q}}}
\newcommand{\tCY}{\widetilde{\mathcal{Y}}}
\newcommand{\tCZ}{\widetilde{\mathcal{Z}}}
\newcommand{\tP}{\widetilde{P}}

\newcommand{\bbeta}{\boldsymbol{\beta}}
\newcommand{\bgamma}{\boldsymbol{\gamma}}

\newcommand{\talpha}{\widetilde{\alpha}}

\newcommand{\tgamma}{\widetilde{\gamma}}

\newcommand{\tDelta}{\widetilde{\Delta}}

\newcommand{\ttau}{\widetilde{\tau}}

\newcommand{\tphi}{\widetilde{\varphi}}

\newcommand{\oS}{\overline{S}}
\newcommand{\hS}{\widehat{S}}

\newcommand{\fX}{\mathfrak{X}}

\newcommand{\tE}{\widetilde{E}}

\newcommand{\td}{\tilde{d}}

\newcommand{\tq}{\widetilde{q}}

\newcommand{\tz}{\widetilde{z}}

\newcommand{\Yequiv}{\mathrel{\mathop{\equiv}\limits_{Y}}}
\newcommand{\Yequivs}{\mathrel{\smash{\mathop{\equiv}\limits_{Y}}}}
\newcommand{\Kequiv}{\mathrel{\mathop{\equiv}\limits_K}}
\newcommand{\Kequivs}{\mathrel{\smash{\mathop{\equiv}\limits_K}}}
\newcommand{\Lequiv}{\mathrel{\mathop{\equiv}\limits_L}}
\newcommand{\Lequivs}{\mathrel{\smash{\mathop{\equiv}\limits_L}}}
\newcommand{\Iequiv}{\mathrel{\mathop{\equiv}\limits_{I}}}

\newcommand{\ebr}{[\ \ ]}
\newcommand{\fC}{\mathfrak{C}}
\newcommand{\fD}{\mathfrak{D}}

\newtheorem{thm}{Theorem}
\newtheorem{corintro}[thm]{Corollary}

\newtheorem{theorem}{Theorem}[section]
\newtheorem{propos}[theorem] {Proposition}
\newtheorem{cor}[theorem] {Corollary}
\newtheorem{lem}[theorem]{Lemma}
\newtheorem{fact}[theorem]{Fact}
\theoremstyle{definition}

\newtheorem{defin}[theorem]{Definition}
\newtheorem{remark}[theorem]{Remark}

\newtheorem{problem}[theorem]{Problem}

\numberwithin{equation}{section}

\author{Alexander A. Gaifullin}

\thanks{This work was supported by the Russian Science Foundation under grant no.~26-11-00052,\\ \texttt{https:/\!/rscf.ru/en/project/26-11-00052/}.}

\keywords{Mapping class group, Torelli group, finite generation, spectral sequence, complex of minimizing cycles, boundedly generated group, abelian cycle}

\title[Finite generation of~$H_2(\I_3)$]{Finite generation of the second homology group of the genus three Torelli group}
\date{}

\address{Steklov Mathematical Institute of Russian Academy of Sciences, Moscow, Russia}
\address{Lomonosov Moscow State University, Russia}
\email{agaif@mi-ras.ru}

\subjclass{20F34 (Primary); 11H56; 20J05; 57M07 (Secondary)}

\begin{document}

\begin{abstract}
 We prove that the group~$H_2(\I_3;\Z)$ is finitely generated, where $\I_3$ is the Torelli subgroup of the mapping class group of a closed genus-three surface~$\Sigma_3$. A key role in the proof  is played by establishing the following property: for each symplectic transvection $t_x\in\Sp_{6}(\Z)$, we have
 $$
 \bigl(t^{120}_x-1\bigr)^{57}\bigl(t_x^{24}-1\bigr)(t_x-1)^6H_2(\I_3;\Z)=0.
 $$
 The proof is based on the study of the spectral sequences associated with the action of~$\I_3$ on the Bestvina--Bux--Margalit complex of minimizing cycles and with the actions on certain auxiliary complexes of the stabilizers of two- and three-component multicurves on~$\Sigma_3$. The remaining ingredients are Tavgen's theorem, which  asserts that $\Sp_{6}(\Z)$ has bounded elementary generation, and the result of Kassabov and Putman, which asserts that $H_2(\I_3;\Z)$ is finitely generated as a $\Z[\Sp_6(\Z)]$-module.
\end{abstract}

\maketitle

\section{Introduction}

\subsection{Main result}
Let $\Sigma_g$ be a closed oriented surface of genus~$g$, and let $\Homeo^+(\Sigma_g)$ be the group of orientation preserving homeomorphisms of~$\Sigma_g$. The \textit{mapping  class group} is the group
$$
\Mod(\Sigma_g)=\pi_0\Homeo^+(\Sigma_g).
$$
Set $H=H_1(\Sigma_g;\Z)$. Then $H\cong\Z^{2g}$, and $H$ is equipped with the intersection form, a unimodular skew-symmetric bilinear form. Let $\Sp(H)\cong \Sp_{2g}(\Z)$ be the group of automorphisms of~$H$ preserving the intersection form. The action of~$\Mod(\Sigma_g)$  on~$H$ gives a homomorphism $\kappa\colon\Mod(\Sigma_g)\to\Sp(H)$, which is well known to be surjective. By definition, the \textit{Torelli group} is the kernel of this homomorphism; then we have a short exact sequence
\begin{equation}\label{eq_ses_I}
 1\longrightarrow\I_g\longrightarrow\Mod(\Sigma_g)\stackrel{\kappa}{\longrightarrow}\Sp(H)\longrightarrow 1.
\end{equation}

It is a classical result by Dehn that $\Mod(\Sigma_1)\cong\SL_2(\Z)$ and $\I_1=1$. Mess~\cite{Mes92} proved that $\I_2$ is an infinitely generated free group. However, Johnson~\cite{Joh83} showed that $\I_g$ is finitely generated for $g\ge 3$. An old and important open problem is whether the groups~$\I_g$ with $g\ge3$ are finitely presented; see~\cite[Appendix, Problem~29]{Bir74}, \cite[Problem~2.9]{Kir97}, \cite[Problem~2.1]{Mor99}. A related question is whether the second homology groups $H_2(\I_g;\Z)$ are finitely generated. The latter question has recently been solved for $g\ge4$. First, Minahan and Putman~\cite{MiPu25} obtained a partial result: they proved that the rational homology group~$H_2(\I_g;\Q)$ is finite-dimensional for $g\ge5$. Then the author~\cite{Gai26} proved that the group~$H_k(\I_g;\Z)$ is finitely generated whenever $k\le g-2$; in particular, $H_2(\I_g;\Z)$ is finitely generated for $g\ge4$.

The case of genus~$3$ has remained completely open. The prevailing view, including that of the author, was that the group~$H_2(\I_3;\Z)$ seemed to be not finitely generated. Nevertheless, attempts to prove this (see e.g.~\cite{Gai21}) led the author to precisely the opposite result. The main result of the present paper is as follows.

\begin{thm}\label{thm_main}
The group~$H_2(\I_3;\Z)$ is finitely generated.
\end{thm}

Note that, by results of Johnson--Millson (unpublished; see~\cite{Mes92}) and Hain (see~\cite{Hai02}), respectively, neither~$H_3(\I_3;\Z)$ nor~$H_4(\I_3;\Z)$ is finitely generated. The group~$H_4(\I_3;\Z)$ was computed explicitly by Spiridonov~\cite{Spi23}. By a theorem of Bestvina, Bux, and Margalit (see~\cite{BBM10}), the Torelli group~$\I_3$ has cohomological dimension~$4$; consequently, $H_k(\I_3;\Z)=0$ for $k\ge5$. A survey of further results concerning which of the groups~$H_k(\I_g;\Z)$ are finitely generated, and which are not, can be found in~\cite{Gai26}.

Along with the Torelli groups of closed surfaces, one can also consider the Torelli subgroups~$\I_g^1$ of the mapping class groups~$\Mod(\Sigma^1_g)$ for compact surfaces of genus~$g$ with one boundary component. Here $\Mod(\Sigma^1_g)$ is the group of isotopy classes of homeomorphisms of~$\Sigma_g^1$ that fix the boundary pointwise. Theorem~\ref{thm_main} can be in a standard way extended to this case.

\begin{corintro}
 The group~$H_2(\I_3^1;\Z)$ is finitely generated.
\end{corintro}

\begin{proof}
 The following version of the Birman short exact sequence is due to Johnson~\cite{Joh83}:
 $$
 1\longrightarrow \pi_1(U\Sigma_3)\longrightarrow \I_3^1\longrightarrow \I_3\longrightarrow 1,
 $$
 where $U\Sigma_3$ is the unit tangent bundle of~$\Sigma_3$. Let
 \begin{equation}\label{eq_LHS}
  E^2_{p,q}=H_p(\I_3;H_q(\pi_1(U\Sigma_3);\Z)) \ \Longrightarrow \ H_{p+q}(\I_3^1;\Z)
 \end{equation}
 be the Lyndon--Hochschild--Serre spectral sequence for this short exact sequence. The long exact sequence of homotopy groups for the fibration $S^1\to U\Sigma_3\to\Sigma_3$ gives that $U\Sigma_3$ is aspherical, and therefore $H_q(\pi_1(U\Sigma_3);\Z)=H_q(U\Sigma_3;\Z)$; note that these coefficient groups in~\eqref{eq_LHS} are twisted. Since $U\Sigma_3$ is compact, it follows that the abelian groups $H_q(U\Sigma_3;\Z)$ are finitely generated. Since $\I_3$ is finitely generated, we obtain that the groups~$E^2_{0,2}=H_0(\I_3;H_2(U\Sigma_3;\Z))$ and~$E^2_{1,1}=H_1(\I_3;H_1(U\Sigma_3;\Z))$ are also finitely generated. Finally, the group $E^2_{2,0}=H_2(\I_3;\Z)$ is finitely generated by Theorem~\ref{thm_main}. Thus $H_2(\I_3^1;\Z)$ is finitely generated as well.
\end{proof}

Throughout the rest of this paper, we consider homology with integer coefficients and omit the coefficient group~$\Z$ from the notation.

\subsection{Finite width and $(P,\CT)$-torsion property}

The short exact sequence~\eqref{eq_ses_I} endows each homology group~$H_k(\I_g)$ with the structure of a $\Z[\Sp(H)]$-module.

For $x\in H\cong\Z^{2g}$, the element $t_x\in\Sp(H)\cong\Sp_{2g}(\Z)$ defined by
$$
t_x(v)=v+(v\cdot x)x
$$
is called a \textit{symplectic transvection}. Here $v\cdot x$ denotes the algebraic intersection number of~$v$ and~$x$.
Our sign convention in the definition of the transvection~$t_x$ differs from the standard convention. However, this choice of sign will be convenient for us, since we want to work with left Dehn twists~$T_{\gamma}$. Under this convention, $T_{\gamma}$ acts on~$H$ by means of the transvection~$t_{[\gamma]}$.

The proof of the fact that $H_k(\I_g)$ is finitely generated for $k\le g-2$ given in~\cite{Gai26} was based on the following \textit{transvection unipotency property}: if $k\le g-2$, then for all $x\in H$, we have
$$
(t_x-1)^{k+1}H_k(\I_g)=0.
$$
To prove Theorem~\ref{thm_main}, we need the following generalization of this property.

\begin{defin}\label{defin_PTtorsion}
 Let $G$ be a group, let $\CT\subseteq G$ be a conjugation-invariant subset, and let $M$ be a left $\Z[G]$-module.
 Suppose that $P\in \Z[t]$. We say that an element $a\in M$ is a \textit{$(P,\CT)$-torsion element} if $P(\tau)a=0$ for all $\tau\in \CT$. A $\Z[G]$-module is called a \textit{$(P,\CT)$-torsion module} if it consists of $(P,\CT)$-torsion elements.

 In the case $G=\Sp(H)\cong\Sp_{2g}(\Z)$, we take for $\CT$ the set of all symplectic transvections~$t_x$. Then a $(P,\CT)$-torsion $\Sp(H)$-module will be called a \textit{transvection $P$-torsion} $\Sp(H)$-module.
\end{defin}

The main result on which Theorem~\ref{thm_main} relies is as follows.

\begin{thm}\label{thm_torsion}
 The group $H_2(\I_3)$ is a transvection $P$-torsion $\Z[\Sp(H)]$-module, where
 $$
 P=\bigl(t^{120}-1\bigr)^{57}\bigl(t^{24}-1\bigr)(t-1)^6.
 $$
 In other words, for each $x\in H$, we have
 $$
 \bigl(t^{120}_x-1\bigr)^{57}\bigl(t_x^{24}-1\bigr)(t_x-1)^6H_2(\I_3)=0.
 $$
 \end{thm}

As in~\cite{Gai26}, the key result that allows us to derive Theorem~\ref{thm_main} from Theorem~\ref{thm_torsion} is the Tavgen theorem on bounded elementary generation of~$\Sp_{2g}(\Z)$ for $g\ge 2$; see~\cite{Tav91}. More precisely, a group $G$ is said to have \textit{finite width} with respect to a subset $\CT\subseteq G$ if there exist elements $\tau_1,\ldots,\tau_m\in\CT$ such that
$$
G=\langle \tau_1\rangle\cdots \langle \tau_m\rangle,
$$
i.e., every element of~$G$ is equal to $\tau_1^{k_1}\cdots\tau_m^{k_m}$ for some $k_1,\ldots,k_m\in\Z$. Then we need the following consequence of the result by Tavgen; see~\cite[Section~5]{Gai26} for the explanation of how to extract this proposition from the results of~\cite{Tav91}.

\begin{propos}[Tavgen~\cite{Tav91}]\label{propos_Tavgen}
If $g\ge 2$, then the group~$\Sp_{2g}(\Z)$ has finite width with respect to the set of all symplectic transvections, i.e., $\Sp_{2g}(\Z)=\langle t_{x_1}\rangle \cdots \langle t_{x_m}\rangle$ for certain sequence of elements $x_1,\ldots,x_m\in H$.
\end{propos}

A key observation is that finite width property for a group~$G$ together with a suitable $(P,\CT)$-torsion property for a finitely generated $\Z[G]$-module~$M$ imply that $M$ is finitely generated as an abelian group. This is a direct generalization of the same observation for $\CT$-unipotent modules made in~\cite{Gai26}. More precisely, we have  the following proposition.

\begin{propos}\label{propos_general_suffcond}
 Let $G$ be a group, and $\CT\subseteq G$ a conjugation-invariant subset. Assume that $G$ has finite width with respect to~$\CT$. Let $M$ be a finitely generated $\Z[G]$-module. Suppose that there exists a monic polynomial $P\in \Z[t]$ with constant term $\pm1$,
 $$
 P=t^n+c_{n-1}t^{n-1}+\cdots+c_1t\pm 1,
 $$
  such that $M$ is a $(P,\CT)$-torsion module.
 Then $M$ is finitely generated as an abelian group.
\end{propos}

\begin{cor}\label{cor_general_suffcond}
 Let $M$ be a finitely generated $\Z[\Sp_{2g}(\Z)]$-module, where $g\ge 2$. Suppose that there exists a monic polynomial $P\in \Z[t]$ with constant term $\pm1$  such that $M$ is a transvection $P$-torsion module.
 Then $M$ is finitely generated as an abelian group.
\end{cor}

\begin{remark}
 Andrew Putman communicated to me that Nate Harman has independently proved the assertion of Corollary~\ref{cor_general_suffcond}  (without any application to the homology of Torelli groups); however, his preprint has not yet appeared.
\end{remark}

Kassabov and Putman~\cite{KaPu20} proved that $H_2(\I_3)$ is finitely generated as a $\Z[\Sp_6(\Z)]$-module. Thus Theorem~\ref{thm_main} follows immediately from Theorem~\ref{thm_torsion} and Corollary~\ref{cor_general_suffcond}.

\subsection{Abelian cycles and relations among them}\label{subsection_ac}
A central role in the proof of Theorem~\ref{thm_torsion} will be played by special homology classes called abelian cycles. Let $g_1$ and~$g_2$ be two commuting elements of a group~$G$. Consider the homomorphism $\varphi\colon \Z\times\Z\to G$ that sends the generators of the first and second factors to~$g_1$ and~$g_2$, respectively. Then the \textit{abelian cycle} $\CA(g_1,g_2)\in H_2(G)$ is defined by $\CA(g_1,g_2)=\varphi_*(\mu)$, where $\mu$ is the standard generator of~$H_2(\Z\times\Z)\cong\Z$.

In $H_2(\I_3)$ there are the following two types of abelian cycles.

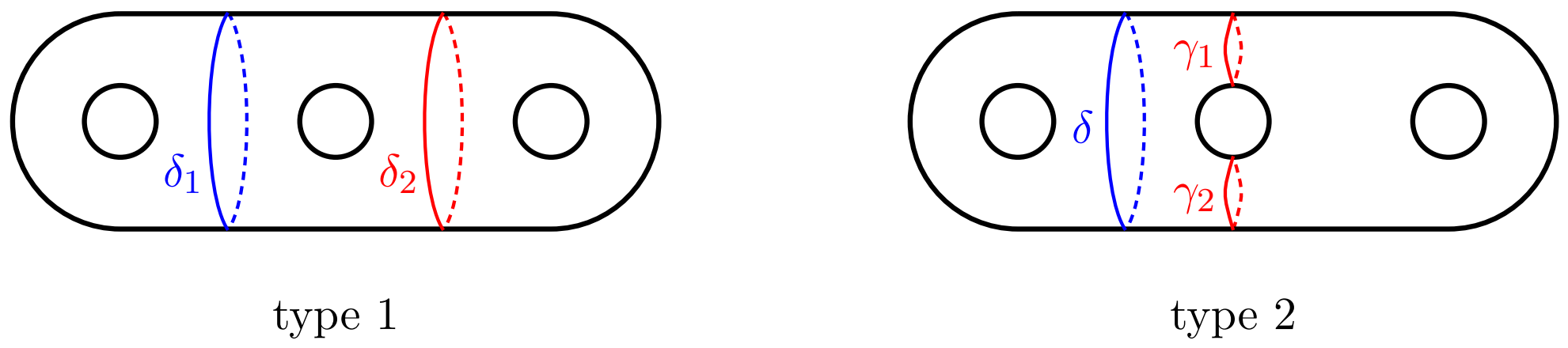
\begin{figure}[!ht]
\begin{tikzpicture}[scale=.3,
        mid arrow/.style={
            postaction={
                decorate,
                decoration={
                    markings,
                    mark=at position 0.55 with {
                        \arrow[#1]{stealth}
                    }
                }
            }
        },
        mid arrow/.default=black,
    ]

\small
% \definecolor{blue}{rgb}{0, 0, 0.7}
% \definecolor{mygreen}{rgb}{0, 0.4, 0}
\tikzset{every path/.append style={line width=.2mm}}
\tikzset{my dash/.style={dash pattern=on 2pt off 1.5pt}}

% Border lines

\begin{scope}
% \draw[color=red, thick, my dash] (.8,3.03) .. controls (1.1,2) .. (.8,.97); % e_0+e_1
%
% \draw[color=red, thick, my dash] (.8,-3.03) .. controls (1.1,-2) .. (.8,-.97); % e_0+e_1

\draw[color=blue, thick, my dash] (-3,3.03) .. controls (-2.3,2) and (-2.3,-2) .. (-3,-3.03);
\draw[color=red, thick, my dash] (3,3.03) .. controls (3.7,2) and (3.7,-2) .. (3,-3.03);

\draw [very thick] (-6,0) circle (1);
\draw [very thick] (0,0) circle (1);
\draw [very thick] (6,0) circle (1);
\draw [very thick]  (-6,3)--(6,3) arc (90:-90:3) -- (-6,-3) arc (270:90:3);

% Visible lines

\draw[color=blue, thick] (-3,3.03) .. controls (-3.7,2) and (-3.7,-2) .. (-3,-3.03) node [pos=0.7,left=-2pt] {$\delta_1$};
\draw[color=red, thick] (3,3.03) .. controls (2.3,2) and (2.3,-2) .. (3,-3.03) node [pos=0.7,left=-2pt] {$\delta_2$};

% \draw[color=red, thick] (.8,3.03) .. controls (.5,2) .. (.8,.97) node [pos=0.6,left=-1pt] {$\gamma_1$}; % e_0+e_1
%
% \draw[color=red, thick] (.8,-3.03) .. controls (.5,-2) .. (.8,-.97) node [pos=0.4,left=-1pt] {$\gamma_2$}; % e_0+e_1

\draw (0, -5.5) node {type 1};
\end{scope}

%%%%%%%
%%%%%%%

\begin{scope}[xshift=25cm]
 \draw[color=red, thick, my dash] (0,3.03) .. controls (.3,2) .. (0,.97); % e_0+e_1

 \draw[color=red, thick, my dash] (0,-3.03) .. controls (.3,-2) .. (0,-.97); % e_0+e_1

\draw[color=blue, thick, my dash] (-3,3.03) .. controls (-2.3,2) and (-2.3,-2) .. (-3,-3.03);

\draw [very thick] (-6,0) circle (1);
\draw [very thick] (0,0) circle (1);
\draw [very thick] (6,0) circle (1);
\draw [very thick]  (-6,3)--(6,3) arc (90:-90:3) -- (-6,-3) arc (270:90:3);

% Visible lines

\draw[color=blue, thick] (-3,3.03) .. controls (-3.7,2) and (-3.7,-2) .. (-3,-3.03) node [pos=0.52,left=-1pt] {$\delta$};

 \draw[color=red, thick] (0,3.03) .. controls (-.3,2) .. (0,.97) node [pos=0.6,left=-2pt] {$\gamma_1$}; % e_0+e_1

 \draw[color=red, thick] (0,-3.03) .. controls (-.3,-2) .. (0,-.97) node [pos=0.4,left=-2pt] {$\gamma_2$}; % e_0+e_1

\draw (0, -5.5) node {type 2};
\end{scope}
\end{tikzpicture}
 \caption{Standard abelian cycles}\label{fig_AC}
\end{figure}

\begin{enumerate}
 \item Let $\delta_1$ and~$\delta_2$ be disjoint non-isotopic separating simple closed curves on~$\Sigma_3$, see the left part of Fig.~\ref{fig_AC}. Then the Dehn twists~$T_{\delta_1}$ and~$T_{\delta_2}$ belong to~$\I_3$ and commute, so we have the abelian cycle $\CA(T_{\delta_1},T_{\delta_2})\in H_2(\I_3)$. Such abelian cycles will be called \textit{type~1 abelian cycles}.
 \item Let $\delta$, $\gamma_1$, and~$\gamma_2$ be pairwise disjoint simple closed curves on~$\Sigma_3$ such that $\delta$ is separating and $\{\gamma_1,\gamma_2\}$ is a bounding pair, see the right part of Fig.~\ref{fig_AC}. Then the mapping classes~$T_{\delta}$ and~$T_{\gamma_1}T_{\gamma_2}^{-1}$ belong to~$\I_3$ and commute, so we have the abelian cycle $\CA\bigl(T_{\delta},T_{\gamma_1}T_{\gamma_2}^{-1}\bigr)\in H_2(\I_3)$. Such abelian cycles will be called \textit{type~2 abelian cycles}.
\end{enumerate}

Abelian cycles of these two types will be referred to as \textit{standard abelian cycles}.

We denote by~$\bA_1$ the subgroup of~$H_2(\I_3)$ generated by all type~1 abelian cycles and by by~$\bA_2$ the subgroup of~$H_2(\I_3)$ generated by all type~2 abelian cycles. Let $\bA=\bA_1+\bA_2$ be the subgroup generated by all standard abelian cycles. Obviously, the subgroups~$\bA_1$, $\bA_2$, and~$\bA$ are $\Z[\Sp(H)]$-submodules of~$H_2(\I_3)$. Vladimirov~\cite{Vla25} showed that all type~1 abelian cycles in~$H_2(\I_3)$ have order~$2$, and therefore $\bA_1$ is a vector space over~$\F_2=\Z/2\Z$.

The proof of Theorem~\ref{thm_torsion} consists of three steps: we will prove the following three assertions in sequence:
\begin{itemize}
 \item $\bA_2$ is a transvection $\bigl(t^{120}-1\bigr)^{57}$-torsion $\Z[\Sp(H)]$-module,
 \item the quotient~$\bA/\bA_2$ is a transvection  $(t^{24}-1)$-torsion $\Z[\Sp(H)]$-module,
 \item the quotient~$H_2(\I_3)/\bA$ is a transvection  $(t-1)^6$-torsion $\Z[\Sp(H)]$-module.
\end{itemize}

To prove the first two assertions we need to study relations among standard abelian cycles.
First, the abelian cycles~$\CA(T_{\delta_1},T_{\delta_2})$ and~$\CA\bigl(T_{\delta},T_{\gamma_1}T_{\gamma_2}^{-1}\bigr)$ depend not on the particular simple closed curves~$\delta_1$, $\delta_2$, $\delta$, $\gamma_1$, and~$\gamma_2$, but only on certain homological data associated with them. Specifically, a type~1 abelian cycle~$\CA(T_{\delta_1},T_{\delta_2})$ is determined solely by the two subgroups $V_1=H_1(S_1)$ and $V_2=H_1(S_2)$ of~$H$, where $S_i\subseteq\Sigma$ is the genus~1 subsurface bounded by~$\delta_i$. We set $\CA(V_1,V_2)=\CA(T_{\delta_1},T_{\delta_2})$. The class~$\CA(V_1,V_2)$ is well-defined whenever $V_1,V_2\subseteq H$ are rank~2 subgroups of~$H$ orthogonal with respect to the intersection form and such that the restrictions of the intersection form to~$V_1$ and~$V_2$ are unimodular.

Similarly, a type~2 abelian cycle~$\CA\bigl(T_{\delta},T_{\gamma_1}T_{\gamma_2}^{-1}\bigr)$ is determined solely by the subgroup $V=H_1(R)\subseteq H$ and the homology class $c=[\gamma_1]=[\gamma_2]\in H$, where $R\subseteq\Sigma$ is the genus~1 subsurface bounded by~$\delta$, and the curves~$\gamma_1$ and~$\gamma_2$ are oriented so that the component of~$\Sigma_3\setminus(\gamma_1\cup\gamma_2)$ containing~$\delta$ lies to the left of~$\gamma_1$ and to the right of~$\gamma_2$. We set $\CA(V,c)=\CA\bigl(T_{\delta},T_{\gamma_1}T_{\gamma_2}^{-1}\bigr)$. The class~$\CA(V,c)$ is well-defined whenever $V\subseteq H$ is a rank~2 subgroup with unimodular restriction of the intersection form to it, and $c\in H$ is a primitive homology class belonging to~$V^{\perp}$. (Recall that an element of an abelian group is called  \textit{primitive} if it is not divisible by any integer greater than~$1$.)

Second, as we have already mentioned, Vladimirov~\cite{Vla25} showed that $2\CA(V_1,V_2)=0$ for all~$V_1$ and~$V_2$. Hence $\CA(V_1,V_2)=\CA(V_2,V_1)$. Vladimirov also proved that if $H=V_1\oplus V_2\oplus V_3$ is a splitting into three rank~2 summands that are mutually orthogonal with respect to the intersection form, then
$$
\CA(V_1,V_2)=\CA(V_2,V_3)=\CA(V_3,V_1);
$$
see Corollary~7.4 in~\cite{Vla25}.

Third, we will show in Section~\ref{section_ac} that the classes~$\CA(V,c)$ satisfy the following \textit{linearity relations}. We extend the function~$\CA(V,c)$ to non-primitive $c=ka$ by setting $\CA(V,ka)=k\CA(V,a)$. Then, for each $V$, the function~$\CA(V,-)$ becomes a linear function on~$V^{\perp}$, see Proposition~\ref{propos_linear}.

Of course, the simple relations listed above are not sufficient to prove the transvection $P$-torsion property for~$\bA$. The following theorem provides far more complicated and deeper relations among standard abelian cycles.

\begin{thm}\label{thm_long_rel}
 Let $a_1,a_2,a_3,b_1,b_2,b_3$ be a symplectic basis of~$H$. Set $U_i=\langle a_i,b_i\rangle$, $t_i=t_{a_j+a_k}$, where $i,j,k$ is a permutation of~$1,2,3$, and
 $$
 \Lambda=\Z\bigl[t_1,t_1^{-1},t_2,t_2^{-1},t_3,t_3^{-1}\bigr].
 $$
 Then

 \textnormal{(1)} the following relation holds in~$H_2(\I_3)$:
\begin{align}
  &(1-t_2)(1-t_3)(t_2-t_3)(1-t_2t_3)\CA(U_1,a_2-a_3)\notag \\
  &{}\quad +(1-t_3)(1-t_1)(t_3-t_1)(1-t_3t_1)\CA(U_2,a_3-a_1)\notag\\
  &{}\quad +(1-t_1)(1-t_2)(t_1-t_2)(1-t_1t_2)\CA(U_3,a_1-a_2)=0;\label{eq_relA2}
 \end{align}

 \textnormal{(2)} the following two homology classes~$\sigma_1$ and~$\sigma_2$ lie in the $\Lambda$-submodule of~$H_2(\I_3)$ generated by the classes~$\CA(U_i,a_j)$ with $i\ne j$:
 \begin{gather*}
  \sigma_1=(1+t_1)(1+t_2)(1+t_3)\CA(U_1,U_2),\\
  \sigma_2=\bigl(t_1+t_1^{-1}+t_2+t_2^{-1}+t_3+t_3^{-1}\bigr)\CA(U_1,U_2);
 \end{gather*}
 in particular, they lie in~$\bA_2$.
\end{thm}

\begin{remark}
 The symplectic transvections~$t_1$, $t_2$, $t_3$ pairwise commute because the corresponding homology classes~$a_2+a_3$, $a_3+a_1$, $a_1+a_2$  have pairwise zero algebraic intersection numbers. For a cyclic permutation of~$1,2,3$, we have
 $$
 t_j^mt_k^n\CA(U_i,a_j-a_k)=\CA(\langle a_i,b_i-na_j-ma_k\rangle,a_j-a_k).
 $$
 After expanding the brackets, relation~\eqref{eq_relA2}  becomes the statement that a certain linear combination of $24$ abelian cycles of this form equals zero. Similarly,
 $$
 t_1^lt_2^mt_3^n\CA(U_1,U_2)=\CA(\langle a_1,b_1-na_2-ma_3\rangle,\langle a_2,b_2-na_1-la_3\rangle).
 $$
 Such classes are naturally indexed by points~$(l,m,n)$ of the lattice~$\Z^3$.
 Then $\sigma_1$ and~$\sigma_2$ are the sums of eight and six abelian cycles of this form corresponding to the vertices of the cube and the regular octahedron, repectively.
\end{remark}

\begin{corintro}\label{cor_short_rel}
 With the notation of Theorem~\ref{thm_long_rel}, the following relations hold in~$H_2(\I_3)$:
 \begin{gather}
 (1-t_1)^2(1-t_2)^2(t_1-t_2)(1-t_1t_2)\CA(U_3,a_1-a_2)=0,\label{eq_rel2_short}\\
  \bigl(1+t_1^2\bigr)\bigl(1+t_2^2\bigr)\bigl(1+t_3^2\bigr)\CA(U_1,U_2)=0,\label{eq_rel1_short}\\
  (1+t_1)(1+t_2)(1+t_3)\bigl(t_1+t_1^{-1}+t_2+t_2^{-1}+t_3+t_3^{-1}\bigr)\CA(U_1,U_2)=0.\label{eq_rel1_short_dop}
 \end{gather}
\end{corintro}

\begin{proof}
 Let $i,j,k$ be a permutation of~$1,2,3$.
 The symplectic transvection~$t_i$ stabilizes~$U_i$ and fixes every homology class~$c$ in~$\langle a_j,a_k\rangle$. Hence $t_i\CA(U_i,c)=\CA(U_i,c)$ for any such~$c$. Therefore, multiplying~\eqref{eq_relA2} on the left  by~$(1-t_1)(1-t_2)$, we obtain relation~\eqref{eq_rel2_short}. Moreover, the element  $(1-t_1)(1-t_2)(1-t_3)$ annihilates the $\Lambda$-submodule of~$H_2(\I_3)$ generated by the classes~$\CA(U_i,a_j)$ with $i\ne j$. Therefore $(1-t_1)(1-t_2)(1-t_3)\sigma_l=0$ for $l=1,2$. Since $2\CA(U_1,U_2)=0$, we obtain relations~\eqref{eq_rel1_short} and~\eqref{eq_rel1_short_dop}.
\end{proof}

\subsection{Notation and conventions}

Throughout the paper, we consider only a closed surface of genus~$3$. Henceforth we omit the subscript~$3$ from the notation. Thus, we denote by~$\Sigma$ a closed  surface of genus~$3$, and by~$\I$ its Torelli group.

We always consider homology with integer coefficients, and omit the coefficient group~$\Z$ from the notation. If $G$ is a group, then $H_1(G)$ is the abelianization of~$G$. For an element $g\in G$, we denote by~$[g]$ the class in~$H_1(G)$ represented by~$g$.

We always set $H=H_1(\Sigma)$. This group is endowed with the intersection form, which we denote by~$x\cdot y$. A \textit{symplectic basis} of~$H$ is a basis $a_1,a_2,a_3,b_1,b_2,b_3$ such that $a_i\cdot a_j=b_i\cdot b_j=0$ and $a_i\cdot b_j=\delta_{ij}$. We denote by~$\Sp(H)$ the group of automorphisms of~$H$ preserving the intersection form. Any symplectic basis of~$H$ gives an isomorphism $\Sp(H)\cong\Sp_6(\Z)$. For a subgroup $V\subseteq H$, we denote by~$V^{\perp}$ the orthogonal complement of~$V$ with respect to the intersection form. A subgroup~$V$ is said to be a \textit{symplectic submodule} if the restriction of the intersection form to~$V$ is unimodular, equivalently, if $H=V\oplus V^{\perp}$. A subgroup~$V$ is said to be \textit{isotropic} if the restriction of the intersection form to~$V$ vanishes, i.e., if $V\subseteq V^{\perp}$. A \textit{Lagrangian submodule} is a maximal  isotropic subgroup with respect to inclusion. An element $x\in H$ is called \textit{primitive} if it is not divisible by any integer greater than~$1$. For any elements $c_1,\ldots,c_k\in H$, we denote by $\langle c_1,\ldots,c_k\rangle$ the subgroup of~$H$ generated by them.

For each $x\in H$, the corresponding \textit{symplectic transvection} is the element $t_x\in\Sp(H)$ defined by $t_x(v)=v+(v\cdot x)x$. Throughout the paper, we will repeatedly use, without further comment, the fact that the symplectic transvections~$t_{x}$ and~$t_{y}$ commute whenever $x\cdot y=0$.

When speaking  of simple closed curves on a surface, we adopt the usual convention in the theory of mapping class groups: we do not distinguish notationally between a simple closed curve and its isotopy class. If $\gamma$ is an oriented simple closed curve on~$\Sigma$,  we denote by~$[\gamma]$ its homology class. We denote by~$T_{\gamma}$ the \textit{left} Dehn twist about~$\gamma$. Then $T_{\gamma}$ acts on~$H$ via the symplectic transvection~$t_{[\gamma]}$.  A \textit{bounding pair} is a pair~$\{\gamma_1,\gamma_2\}$ such that $\gamma_1$ and~$\gamma_2$ are disjoint non-isotopic non-separating simple closed curves on~$\Sigma$, but their union~$\gamma_1\cup\gamma_2$ separates~$\Sigma$. The mapping class~$T_{\gamma_1}T_{\gamma_2}^{-1}$ is called a \textit{bounding pair map}, or a \textit{BP-map}. It is a standard fact that Dehn twists about separating simple closed curves and BP-maps belong to the Torelli group.

By a \textit{multicurve} on a surface, we always mean a union of finitely many pairwise disjoint, pairwise non-isotopic, homotopically nontrivial simple closed curves. A multicurve is said to be \textit{oriented} if all of its components are endowed with orientations. For oriented multicurves, the notation $N\subseteq M$ will always mean that the orientations of the components of~$N$ coincide with their orientations in~$M$.

If a group~$G$ acts on a set~$X$, then we denote by~$G_x$ the stabilizer of an element~$x\in X$.

\subsection{Outline of the paper}

The proofs of our main theorems, Theorem~\ref{thm_main} and Theorem~\ref{thm_torsion}, are based on two hard technical results. The first of them is Theorem~\ref{thm_long_rel}, and the  second is the following theorem, which will allow us to reduce the study of the whole group~$H_2(\I_3)$ to the study of its subgroup~$\bA$ generated by standard abelian cycles.

\begin{thm}\label{thm_toAC}
For each $x\in H$, we have
  $(t_x-1)^6H_2(\I_3)\subseteq \bA$.
 \end{thm}

Thus, the paper naturally divides into three parts.
\begin{description}
 \item[Part 1] Derivation of Theorems~\ref{thm_main}, \ref{thm_torsion} from Theorems~\ref{thm_long_rel} and~\ref{thm_toAC}.
 \item[Part 2] Proof of Theorem~\ref{thm_long_rel}.
 \item[Part 3] Proof of Theorem~\ref{thm_toAC}.
\end{description}
These three parts of the paper are almost independent of each other. Let us discuss them separately.

\subsubsection*{Part 1: Sections~\ref{section_ac}--\ref{section_SpecSeq}}

Section~\ref{section_ac} is preliminary in nature. In it we prove that abelian cycles of the first and second types are uniquely determined by the homological data described in Subsection~\ref{subsection_ac}, and we also prove the linearity relation for abelian cycles of the second type. In Sections~\ref{section_A2} and~\ref{section_A/A2} we derive the transvection $P$-torsion properties for~$\bA_2$ and~$\bA/\bA_2$, respectively, from parts~(1) and~(2) of Theorem~\ref{thm_long_rel}. We note that these derivations, although not very simple, are completely formal: in essence they reduce to certain computations in the group~$\Sp_6(\Z)$ and do not use any additional information about Torelli groups. In Section~\ref{section_proof_main} we prove Proposition~\ref{propos_general_suffcond} and complete the proofs of Theorems~\ref{thm_torsion} and~\ref{thm_main} modulo Theorems~\ref{thm_long_rel} and~\ref{thm_toAC}. Finally, we have included in Part~1 Section~\ref{section_SpecSeq}, which has no logical connection with that part. This section contains preliminaries on the spectral sequence for the action of a group on a CW complex. The reason Section~\ref{section_SpecSeq} is included in Part~1 is that the author wanted Parts~2 and~3 to be as independent of each other as possible, while the contents of Section~\ref{section_SpecSeq} are constantly used in both parts.

\subsubsection*{Part 2: Sections~\ref{section_IM}--\ref{section_end_thm3}}

The aim of this part is to prove Theorem~\ref{thm_long_rel}. In fact, all relations in Theorem~\ref{thm_long_rel} come from certain relations in~$H_2(\I_M)$, where $M\subseteq \Sigma$ is a three-component multicurve that does not separate~$\Sigma$, and $\I_M$ is the stabilizer of it in the Torelli group. An analogue of Theorem~\ref{thm_long_rel} for the stabilzer~$\I_M$ (Theorem~\ref{theorem_3}) is formulated in Section~\ref{section_IM}. The proof of Theorem~\ref{theorem_3} occupies the whole Part~2. It is based on a detailed study of the spectral sequence for the action of~$\I_M$ on the complex of relative cycles~$\CB_{0,6}$ introduced by Spiridonov~\cite{Spi23},~\cite{Spi24}, including hard computations of the differentials~$d^2$ and~$d^3$. Along the way, while studying this spectral sequence, we  also prove that the group~$H_2(\I_M)$ is generated by standard abelian cycles; this result (part~(1) of Theorem~\ref{theorem_3}) will then be used in Part~3 to prove Theorem~\ref{thm_toAC}. This is the only dependency between Parts~2 and~3.

\subsubsection*{Part 3: Sections~\ref{section_modA_first}--\ref{section_modA_end}}

In this part we prove Theorem~\ref{thm_toAC}. The proof relies on the study of the spectral sequence for the action of~$\I$ on the Bestvina--Bux--Margalit complex of cycles~$\CB$, as well as the spectral sequence for the action of the stabilizer~$\I_N$ of a two-component multicurve~$N$ on an auxiliary complex~$\tCB$. This part of the proof is partly conceptually related to the author's proof of the finite generation of the groups~$H_k(\I_g)$ for $k\le g-2$ (see~\cite{Gai26}). This part also makes active use of the author's results from~\cite{Gai21}. In that paper, the author proposed an approach toward a possible proof that the group~$H_2(\I_3)$ is \textit{not} finitely generated. Naturally, this plan cannot be carried out, since in the present paper we prove that $H_2(\I_3)$ \textit{is} finitely generated. Nevertheless, a number of technical results from~\cite{Gai21} prove useful in the study of the spectral sequence for the action of~$\I$ on~$\CB$.
\smallskip

We do not give detailed outlines of Parts~2 and~3 here. Since these parts are almost completely independent of each other, it is more convenient to give an outline of each part at the beginning of that part; these outlines are given at the end of Sections~\ref{section_IM} and~\ref{section_modA_first}, respectively.

\subsection{Some open questions}

Certainly, the main problem is to compute explicitly the group~$H_2(\I_3)$, or at least~$H_2(\I_3;\Q)$. A problem whose solution can provide serious progress towards such a computation is as follows.

\begin{problem}\label{problem_min}
(1) Describe the ideal of~$\Z\bigl[t]$ consisting of all polynomials~$P$ such that $H_2(\I_3)$ is a transvection $P$-torsion $\Z\bigl[\Sp_6(\Z)\bigr]$-module.

(2) Find the minimal polynomial~$P_{\min}\in\Q[t]$ such that $H_2(\I_3;\Q)$ is a transvection $P_{\min}$-torsion $\Q\bigl[\Sp_6(\Z)\bigr]$-module.
\end{problem}

It is also an important question whether the representation of~$\Sp_6(\Z)$ on~$H_2(\I_3;\Q)$ is algebraic or not. Recall that by the author's result~\cite{Gai26} the representation of~$\Sp_{2g}(\Z)$ on~$H_k(\I_g)$ is algebraic whenever $k\le g-2$. Note that from Proposition~10.1 in~\cite{Gai26} it follows that the representation of~$\Sp_6(\Z)$ on~$H_2(\I_3;\Q)$ is algebraic if and only if the minimal polynomial~$P_{\min}$ from part~(2) of Problem~\ref{problem_min} has the form $(t-1)^n$.

\subsection*{Acknowledgements}

The author is grateful to Igor Patsankov, Vasilii Rozhdestvenskii, and Andrei Vladimirov for fruitful discussions.
\smallskip

No AI tools were used in obtaining the results of this paper. AI (DeepSeek, free web version) was used solely for grammar checking.

\part{Derivation of main results from Theorems~\ref{thm_long_rel} and~\ref{thm_toAC}}

\section{Standard abelian cycles}\label{section_ac}

The fact that the abelian cycles~$\CA(V_1,V_2)$ and~$\CA(V,c)$ introduced in Subsection~\ref{subsection_ac} are well-defined is standard (see e.g. Lemma~7.3 in~\cite{Vla25} for type~1). However, for the convenience of the reader, we provide a proof of this fact.

\begin{propos}
 \textnormal{(1)} If\/ $V_1,V_2\subseteq H$ are rank~$2$ symplectic submodules orthogonal with respect to the intersection form, then the type~$1$ abelian cycle~$\CA(V_1,V_2)$ is well-defined.

 \textnormal{(2)} If\/ $V\subseteq H$ is a rank~$2$ symplectic submodule, and $c$ is a primitive element  of~$V^{\perp}$, then the type~$2$ abelian cycle~$\CA(V,c)$ is well-defined.
\end{propos}

\begin{proof}
 For any $V_1$ and~$V_2$ as in part~(1) of the proposition, there exist disjoint separating simple closed curves~$\delta_1$ and~$\delta_2$ such that $\delta_i$ bounds a genus~$1$ subsurface $S_i\subseteq\Sigma$ with $H_1(S_i)=V_i$ for $i=1,2$. We need to show that the abelian cycle $\CA(V_1,V_2)=\CA(T_{\delta_1},T_{\delta_2})$ is independent of the choice of a pair of curves~$(\delta_1,\delta_2)$ with this property. Indeed, let $(\delta_1',\delta_2')$ be another pair of disjoint separating simple closed curves inducing the same splitting in homology. Then a standard application of the Alexander method (see~\cite[Section~2.3]{FaMa12}) yields that there exists a mapping class $h\in\I$ that sends~$\delta_1$ and~$\delta_2$ to~$\delta_1'$ and~$\delta_2'$, respectively. Since any group acts trivially on its homology groups, we obtain that $\CA(\delta_1,\delta_2)=\CA(\delta_1',\delta_2')$.

 Similarly, if $V\subseteq H$ is as in part~(2) of the proposition, then there exists a separating simple closed curve~$\delta$ on~$\Sigma$ giving the splitting $H=V\oplus V^{\perp}$ in homology. Let $S$ be the genus~$2$ component of~$\Sigma\setminus\delta$. Since $c$ is a primitive class in~$V^{\perp}$, there exists an oriented simple closed curve $\gamma_1$ in~$S$ with $[\gamma_1]=c$. Moreover, $\gamma_1$ does not separate~$S$. It follows that there exists a simple closed curve~$\gamma_2\subseteq S\setminus \gamma_1$ such that the curves~$\delta$, $\gamma_1$, and~$\gamma_2$ bound a pair of pants, and moreover, this pair of pants lies to the left of~$\gamma_1$. Again an application of the Alexander method yields that the triple~$(\delta,\gamma_1,\gamma_2)$ is defined uniquely up to the action of~$\I$, and therefore, the class $\CA(V,c)=\CA\bigl(T_{\delta},T_{\gamma_1}T_{\gamma_2}^{-1}\bigr)$ is well-defined.
\end{proof}

\begin{propos}\label{propos_linear}
Let $V\subseteq H$ be a symplectic submodule of rank~$2$. Extend the function~$\CA(V,c)$ to non-primitive $c=ka$ by setting $\CA(V,ka)=k\CA(V,a)$. Then the function~$\CA(V,-)$ becomes a linear function on~$V^{\perp}$.
\end{propos}

\begin{proof}
Let $\delta$ be a separating simple closed curve giving the splitting $H=V\oplus V^{\perp}$, and let $S$ be the  connected components of~$\Sigma\setminus\delta$ of genus~$2$. Then $S$ is a surface of genus~$2$ with one puncture, which we denote by~$p$. Let $\oS$ be the closure of~$S$ in~$\Sigma$, and let~$\hS$ be the closed surface of genus~$2$ obtained by filling in the puncture of~$S$. Let $\I(S)$,  $\I\bigl(\oS\bigr)$, and~$\I\bigl(\hS\,\bigr)$ be the Torelli subgroups of the mapping class groups~$\Mod(S)$,  $\Mod\bigl(\oS\bigr)$, and~$\Mod\bigl(\hS\,\bigr)$, respectively.

Extending homeomorphisms by the identity gives a homomorphism
$$i\colon \I\bigl(\oS\bigr)\to \I(\Sigma)=\I.$$
It is a standard fact that $\CA(g,g_1g_2)=\CA(g,g_1)+\CA(g,g_2)$ whenever $g$, $g_1$, $g_2$ are elements of a group~$G$ such that $g$ commutes with both~$g_1$ and~$g_2$. Therefore,
since $T_{\delta}$ commutes with all elements in the image of~$i$, we obtain a well-defined homomorphism $\bar{\theta}\colon \I\bigl(\oS\bigr)\to H_2(\I)$ given by
$$
\bar{\theta}(h)=\CA(T_{\delta},i(h)).
$$
Moreover, we have a capping short exact sequence (cf.~\cite[Section~3.6]{FaMa12}):
$$
1\longrightarrow \langle T_{\delta}\rangle \longrightarrow \I\bigl(\oS\bigr)\longrightarrow \I(S)\longrightarrow 1.
$$
Since $\bar{\theta}(T_{\delta})=\CA(T_{\delta},T_{\delta})=0$, the homomorphism~$\bar{\theta}$ gives a well-defined homomorphism $\theta\colon \I(S)\to H_2(\I)$.

Further, we have the Birman short exact sequence:
$$
1\longrightarrow \pi_1\bigl(\hS,p\bigr)\stackrel{\varphi}{\longrightarrow}\I(S)\longrightarrow\I\bigl(\hS\,\bigr)\longrightarrow 1.
$$
More precisely, this is the restriction of the classical Birman exact sequence (see~\cite{Bir69}, cf.~\cite[Section~4.2]{FaMa12}) to the Torelli groups. Recall that $\varphi(\gamma)$ is the mapping class obtained by pushing the puncture along the loop~$\gamma^{-1}$. Since the group~$H_2(\I)$ is abelian, the homomorphism $\theta\circ \varphi\colon \pi_1\bigl(\hS,p\bigr)\to H_2(\I)$ factors through a homomorphism
$$
\psi\colon V^{\perp}=H_1\bigl(\hS\,\bigr)\to H_2(\I).
$$

Now, let $c$ be a primitive element of~$V^{\perp}$. Then $c$ can be represented by a simple (i.e., non-self-intersecting) loop~$\gamma$ in~$\hS$ based at~$p$. Let $\gamma_1$ and~$\gamma_2$ be simple closed curves in~$S$ obtained by pushing~$\gamma$ slightly to the right and to the left, respectively. Then $\varphi(\gamma)=T_{\gamma_1}T_{\gamma_2}^{-1}$. Moreover, $[\gamma_1]=[\gamma_2]=c$, and the pair of pants enclosed by~$\delta$, $\gamma_1$, and~$\gamma_2$ lies to the left of~$\gamma_1$. Therefore
$$
\psi(c)=(\theta\circ\varphi)(\gamma)=\theta\bigl(T_{\gamma_1}T_{\gamma_2}^{-1}\bigr)=
 \CA\bigl(T_{\delta},T_{\gamma_1}T_{\gamma_2}^{-1}\bigr)=\CA(V,c).
$$
Thus the function~$\CA(V,-)$ defined on the set of primitive elements of~$V^{\perp}$ extends to the homomorphism~$\psi\colon V^{\perp}\to H_2(\I)$. The proposition follows.
\end{proof}

\section{Transvection $P$-torsion property for~$\bA_2$}
\label{section_A2}

The goal of this section is to derive the following proposition from relation~\eqref{eq_rel2_short} in Corollary~\ref{cor_short_rel}, which follows from part~1 of Theorem~\ref{thm_long_rel}.

\begin{propos}\label{propos_A2}
The group~$\bA_2$ is a transvection $\bigl(t^{120}-1\bigr)^{57}$-torsion $\Z[\Sp(H)]$-module.
\end{propos}

The following properties of symplectic transvections will be repeatedly used not only in this section but throughout the paper:
\begin{itemize}
 \item if $x\cdot y=0$, then the transvections~$t_x$ and~$t_y$ commute,
 \item if $x$ lies in either~$V$ or~$V^{\perp}$, and $x\cdot c=0$, then $t_x$ stabilizes both~$V$ and~$c$, and therefore $t_x\CA(V,c)=\CA(V,c)$,
 \item if $x$ lies in~$V_1$, or~$V_2$, or~$(V_1\oplus V_2)^{\perp}$, then $t_x\CA(V_1,V_2)=\CA(V_1,V_2)$.
\end{itemize}

The following lemma is proved by a direct computation. It will also be repeatedly used  throughout the paper.

\begin{lem}\label{lem_transvec3} Suppose that $a_1,a_2,a_3$ form a basis of a Lagrangian submodule of~$H$, and $x=n_1a_1+n_2a_2+n_3a_3$. Then
 \begin{equation*}
 t_x=t_{a_2+a_3}^{n_2n_3}t_{a_1+a_3}^{n_1n_3}t_{a_1+a_2}^{n_1n_2}t_{a_1}^{n_1(n_1-n_2-n_3)}t_{a_2}^{n_2(n_2-n_1-n_3)}t_{a_3}^{n_3(n_3-n_1-n_2)}.
\end{equation*}
\end{lem}

\begin{cor}\label{cor_transvec3}
 Let $V$ be a rank~$2$ symplectic submodule of~$H$, and let $c$ be an element of\/~$V^{\bot}$. Let $a_1,a_2,a_3$ be a basis of a Lagrangian submodule~$L$ of~$H$ such that $c\in L$, $a_1,a_2\in V^{\perp}$, and $a_3\in V$. Suppose that $x=n_1a_1+n_2a_2+n_3a_3$. Then
 \begin{equation*}
 t_x\CA(V,c)=t_{a_2+a_3}^{n_2n_3}t_{a_1+a_3}^{n_1n_3}\CA(V,c).
 \end{equation*}
\end{cor}

The proof of Proposition~\ref{propos_A2} will be divided into several propositions.

\begin{propos}\label{propos_annih_y+c}
 Let $V$ be a rank~$2$ symplectic submodule of~$H$, let $y$ be a primitive element of\/~$V$, and let~$c$ and~$z$ be primitive elements of\/~$V^{\perp}$ such that $\langle c,z\rangle$ is a Lagrangian submodule of\/~$V^{\perp}$. Then the following relation holds in~$H_2(\I)$:
 \begin{equation}\label{eq_annih_y+c}
  (t_{y+c}-1)f(t_{y+z})\CA(V,c)=0,
 \end{equation}
 where
 $$
 f(t)=(t^8-1)(t^5-1)^2(t^3-1)^4(t^2-1)^4(t-1)^6.
 $$
\end{propos}

To prove this proposition we need the following auxiliary lemma.

\begin{lem}\label{lem_ideal1}
 The ideal of the Laurent polynomial ring~$\Z\bigl[t_1,t_1^{-1},t_2,t_2^{-1}\bigr]$ generated by the two polynomials
 \begin{align*}
  p_1&=(t_1-1)^2(t_2-1)^2(t_1t_2-1),\\
  p_2&=\bigl(t_1^{-1}t_2^2-1\bigr)^2\bigl(t_1^{-2}t_2^{3}-1\bigr)^2\bigl(t_1^{-3}t_2^{5}-1\bigr)
 \end{align*}
 contains the polynomial
 $$
 f(t_1)=(t_1^8-1)(t_1^5-1)^2(t_1^3-1)^4(t_1^2-1)^4(t_1-1)^6.
 $$
\end{lem}

\begin{proof}
Let $I$ be the ideal generated by~$p_1$ and~$p_2$. It is convenient for us to write $p_i(t_1,t_2)$ rather than simply~$p_i$, since we are going to substitute various expressions into these polynomials. We have
\begin{align*}
 p_2(t_1,1)&=p_2(t_1,t_2)+(t_2-1)q_1(t_1,t_2),\\
 p_2(t_1,t_1^{-1})&=p_2(t_1,t_2)+(t_1t_2-1)q_2(t_1,t_2)
\end{align*}
for certain Laurent polynomials $q_1,q_2\in\Z\bigl[t_1,t_1^{-1},t_2,t_2^{-1}\bigr]$. Therefore
$$
p_2(t_1,1)^2p_2(t_1,t_1^{-1})=p_2(t_1,t_2)A(t_1,t_2)+(t_2-1)^2(t_1t_2-1)q_1(t_1,t_2)^2q_2(t_1,t_2)
$$
for a certain Laurent polynomial $A\in\Z\bigl[t_1,t_1^{-1},t_2,t_2^{-1}\bigr]$. Hence the Laurent polynomial
$$
(t_1-1)^2p_2(t_1,1)^2p_2(t_1,t_1^{-1})
$$
belongs to~$I$. The lemma follows, since this Laurent polynomial coincides with~$f(t_1)$ up to multiplication by a monomial.
\end{proof}

\begin{proof}[Proof of Proposition~\ref{propos_annih_y+c}]
Set $a_1=c+z$, $a_2=z$, $a_3=y$, extend $a_1,a_2$ to a symplectic basis $a_1,a_2,b_1,b_2$ of~$V^{\perp}$, and extend~$a_3$ to a symplectic basis $a_3,b_3$ of~$V$. As in Theorem~\ref{thm_long_rel}, we set $U_i=\langle a_i,b_i\rangle$, and $t_i=t_{a_j+a_k}$ whenever $i,j,k$ is a permutation of~$1,2,3$. Then $\CA(V,c)=\CA(U_3,a_1-a_2)$, $t_{y+z}=t_1$, and $t_{y+c}=t_{a_1-a_2+a_3}$. The transvections~$t_{y+z}$, $t_{y+c}$, $t_1$, $t_2$, and~$t_3$ mutually commute. Moreover, by Corollary~\ref{cor_transvec3}, we have
$$
t_{y+c}\CA(V,c)=t_1^{-1}t_2\CA(V,c).
$$
Hence relation~\eqref{eq_rel2_short} from Corollary~\ref{cor_short_rel} reads as
 \begin{equation}\label{eq_rel_t1t2}
  (t_1-1)^2(t_2-1)^2(t_1t_2-1)(t_{y+c}-1)\CA(V,c)=0.
 \end{equation}
 Consider a new symplectic basis $a_1',b_1',a_2',b_2'$ of~$V^{\perp}$, where
 \begin{align*}
  a_1'&=3a_1-2a_2,& b_1'&=-b_1-2b_2\\
  a_2'&=2a_1-a_2,&b_2'&=2b_1+3b_2.
 \end{align*}
 Since $a_1'-a_2'=a_1-a_2=c$, we obtain that relation~\eqref{eq_rel_t1t2} will remain true if we replace the transvections $t_1=t_{a_2+a_3}$ and~$t_2=t_{a_1+a_3}$ by the transvections $t_1'=t_{a_2'+a_3}=t_{2a_1-a_2+a_3}$ and $t_2'=t_{a_1'+a_3}=t_{3a_1-2a_2+a_3}$, respectively. By Corollary~\ref{cor_transvec3}, we have
 \begin{equation*}
  t_1'\CA(V,c)=t_1^{-1}t_2^2\CA(V,c),\qquad
  t_2'\CA(V,c)=t_1^{-2}t_2^{3}\CA(V,c).
 \end{equation*}
 Since the transvections~$t_1$, $t_2$, $t_1'$, $t_2'$, and~$t_{y+c}$ pairwise commute, the relation obtained from~\eqref{eq_rel_t1t2} by replacing~$t_1$ and~$t_2$ by~$t_1'$ and~$t_2'$, respectively, reads as
 \begin{equation}\label{eq_rel_t1t2_2}
  \bigl(t_1^{-1}t_2^2-1\bigr)^2\bigl(t_1^{-2}t_2^{3}-1\bigr)^2\bigl(t_1^{-3}t_2^{5}-1\bigr)(t_{y+c}-1)\CA(V,c)=0.
 \end{equation}
 The required relation~\eqref{eq_annih_y+c} follows immediately from~\eqref{eq_rel_t1t2}, \eqref{eq_rel_t1t2_2}, and Lemma~\ref{lem_ideal1}.
\end{proof}

\begin{propos}\label{propos_annih_wo}
 Let $V$, $c$, $y$, and~$z$ be as in Proposition~\ref{propos_annih_y+c}. Then we have
 \begin{equation}\label{eq_annih_wo}
  g(t_{y+z})\CA(V,c)=0,
 \end{equation}
 where
 $$
 g(t)=(t-1)(t^2-1)f(t)=(t^8-1)(t^5-1)^2(t^3-1)^4(t^2-1)^5(t-1)^7.
 $$
\end{propos}

\begin{proof}
Consider relation~\eqref{eq_annih_y+c} for the abelian cycles~$\CA(V,c)$, $\CA(V,c+z)$, and~$\CA(V,c+2z)$:
\begin{align}
 (t_{y+c}-1)f(t_{y+z})\CA(V,c)&=0,\label{eq_rel_CAV1}\\
 (t_{y+c+z}-1)f(t_{y+z})\CA(V,c+z)&=0,\label{eq_rel_CAV2}\\
 (t_{y+c+2z}-1)f(t_{y+z})\CA(V,c+2z)&=0.\label{eq_rel_CAV3}
\end{align}
Set $t=t_{y+z}$ and $u=t_{y+c}$. Using Proposition~\ref{propos_linear}, Corollary~\ref{cor_transvec3}, and the fact that all transvections under consideration pairwise commute, we can write relations~\eqref{eq_rel_CAV1}--\eqref{eq_rel_CAV3} in the form
\begin{align}
 (u-1)f(t)\CA(V,c)&=0,\label{eq_rel_CAV1'}\\
 (tu-1)f(t)\bigl(\CA(V,c)+\CA(V,z)\bigr)&=0,\label{eq_rel_CAV2'}\\
 (t^2u-1)f(t)\bigl(\CA(V,c)+2\CA(V,z)\bigr)&=0.\label{eq_rel_CAV3'}
\end{align}
Multiply~\eqref{eq_rel_CAV2'} by~$2(t^2u-1)$ and \eqref{eq_rel_CAV3'} by~$tu-1$, and subtract the second resulting relation from the first. We get
$$
(tu-1)(t^2u-1)f(t)\CA(V,c)=0.
$$
This relation and relation~\eqref{eq_rel_CAV1'} immediately imply the required relation~\eqref{eq_annih_wo}.
\end{proof}

We now find it convenient to weaken the statement of Proposition~\ref{propos_annih_wo} as follows.

\begin{cor}\label{cor_annih_120_18}
 Let $V$, $c$, $y$, and~$z$ be as in Proposition~\ref{propos_annih_y+c}. Then we have
 \begin{equation*}
  \bigl(t_{y+z}^{120}-1\bigr)^{19}\CA(V,c)=0,
 \end{equation*}
\end{cor}

\begin{propos}\label{propos_annih_120_53}
 Let $V$ be a rank~$2$ symplectic submodule of~$H$, $y$ a primitive element of\/~$V$, and~$z$ a primitive element of~$V^{\perp}$. Then
 \begin{equation}\label{eq_annih_120_53}
  \bigl(t_{y+z}^{120}-1\bigr)^{56}\CA(V,c)=0
 \end{equation}
 for all $c\in V^{\perp}$.
\end{propos}

\begin{proof}
 Set $a_2=z$, $a_3=y$, extend~$a_2$ to a symplectic basis $a_1,a_2,b_1,b_2$ of~$V^{\perp}$, and extend~$a_3$ to a symplectic basis $a_3,b_3$ of~$V$. As above, set $t_1=t_{a_2+a_3}=t_{y+z}$, $t_2=t_{a_3+a_1}$, and $t_3=t_{a_1+a_2}$. If $c$ is one of the three homology classes $a_1$, $a_1+a_2$, and~$b_1$, then $c,z$ form a basis of a Lagrangian submodule of~$V^{\perp}$, and therefore the required relation~\eqref{eq_annih_120_53} holds by Corollary~\ref{cor_annih_120_18}.
 Let us prove relation~\eqref{eq_annih_120_53} for $c=b_2$. Then, for all $c\in V^{\perp}$, the relation  will follow from the linearity of~$\CA(V,c)$.

 By Corollary~\ref{cor_annih_120_18}, we have:
 \begin{align}
 \bigl(t_1^{120}-1\bigr)^{19}\CA(V,a_1)&=0,\label{eq_18_0}\\
 \bigl(t_1^{120}-1\bigr)^{19}\CA(V,b_1)&=0,\label{eq_18_1}\\
 \bigl(t_2^{120}-1\bigr)^{19}\CA(V,b_2)&=0,\label{eq_18_2}\\
 \bigl(t_{a_1+a_2+a_3}^{120}-1\bigr)^{19}\CA(V,b_2-b_1)&=0.\label{eq_18_3}
 \end{align}
 By Proposition~\ref{propos_linear}, we have $\CA(V,b_2-b_1)=\CA(V,b_2)-\CA(V,b_1)$. Since the symplectic transvections~$t_1$ and~$t_{a_1+a_2+a_3}$ commute, relations~\eqref{eq_18_1} and~\eqref{eq_18_3} imply that
 \begin{equation}\label{eq_120120}
  \bigl(t_1^{120}-1\bigr)^{19}\bigl(t_{a_1+a_2+a_3}^{120}-1\bigr)^{19}\CA(V,b_2)=0.
 \end{equation}
 By Lemma~\ref{lem_transvec3}, we have $t_{a_1+a_2+a_3}=t_1t_2t_3t_{a_1}^{-1}t_{a_2}^{-1}t_{a_3}^{-1}$. Note that now we cannot use Corollary~\ref{cor_transvec3}, since $b_2\notin\langle a_1,a_2,a_3\rangle$. The transvections~$t_{a_1}$ and~$t_{a_3}$ fix~$\CA(V,b_2)$. The transvections $t_3=t_{a_1+a_2}$ and~$t_{a_2}$ stabilize the subgroup~$V$, and $t_3t_{a_2}^{-1}(b_2)=b_2-a_1$. Hence, again using the linearity of~$\CA(V,-)$, we obtain that
 $$
 t_3t_{a_2}^{-1}\CA(V,b_2)= \CA(V,b_2)-\CA(V,a_1).
 $$
Consequently,
  $$
 t_{a_1+a_2+a_3}\CA(V,b_2)= t_1t_2\CA(V,b_2)-t_1t_2\CA(V,a_1).
 $$
Moreover, each of the transvections~$t_3$, $t_{a_1}$, $t_{a_2}$, and~$t_{a_3}$ fixes~$\CA(V,a_1)$. It follows that
$$
\bigl(t_{a_1+a_2+a_3}^{120}-1\bigr)^{19}\CA(V,b_2)=\bigl(t_1^{120}t_2^{120}-1\bigr)^{19}\CA(V,b_2)+p(t_1,t_2)\CA(V,a_1)
$$
for a certain polynomial $p$. Now, relations~\eqref{eq_18_0} and~\eqref{eq_120120} imply that
\begin{equation}\label{eq_120120120}
\bigl(t_1^{120}-1\bigr)^{19}\bigl(t_1^{120}t_2^{120}-1\bigr)^{19}\CA(V,b_2)=0.
\end{equation}
Since
$$
t_1^{120}-1=t_1^{120}t_2^{120}-1-t_1^{120}\bigl(t_2^{120}-1\bigr),
$$
it follows that the polynomial $\bigl(t_1^{120}-1\bigr)^{37}$ belongs to the ideal of~$\Z\bigl[t_1,t_1^{-1},t_2,t_2^{-1}\bigr]$ generated by~$\bigl(t_2^{120}-1\bigr)^{19}$ and $\bigl(t_1^{120}t_2^{120}-1\bigr)^{19}$. Thus the required relation~\eqref{eq_annih_120_53} for $c=b_2$ follows from~\eqref{eq_18_2} and~\eqref{eq_120120120}.
\end{proof}

\begin{propos}\label{propos_annih_120_54}
 Let $\CA(V,c)$ be a type~2 abelian cycle. Then
 \begin{equation*}
  \bigl(t_x^{120}-1\bigr)^{57}\CA(V,c)=0
 \end{equation*}
 for all $x\in H$.
\end{propos}

\begin{proof}
 We have $x=my+nz$, where $y$ is a primitive element of~$V$, $z$ is a primitive element of~$V^{\perp}$, and~$m,n\in\Z$. By Lemma~\ref{lem_transvec3}, we have $t_x=t_{y+z}^{mn}t_y^{m(m-n)}t_z^{n(n-m)}$. Set
 $$
 \tau_1=t_z^{120n(n-m)},\qquad \tau_2=t_{y+z}^{120mn}.
 $$
 The transvections~$t_y$, $t_z$, and~$t_{y+z}$ pairwise commute, and $t_y$ fixes~$\CA(V,c)$. Consequently,
 \begin{equation*}
  \bigl(t_x^{120}-1\bigr)^{57}\CA(V,c)=(\tau_1\tau_2-1)^{57}\CA(V,c).
 \end{equation*}
 Thus we need to prove that
 \begin{equation}\label{eq_annih_57}
  (\tau_1\tau_2-1)^{57}\CA(V,c)=0.
 \end{equation}

 The transvection~$t_z$ stabilizes~$V$, and sends $c$ to $c+(c\cdot z)z$. Using the linearity of~$\CA(V,-)$, we obtain that
$$
t_z^k\CA(V,c)=\CA(V,c)+k(c\cdot z)\CA(V,z).
$$
Therefore
 \begin{equation*}
  (\tau_1-1)^2\CA(V,c)=0.
 \end{equation*}
 On the other hand, from Proposition~\ref{propos_annih_120_53} it follows that
 \begin{equation*}
  (\tau_2-1)^{56}\CA(V,c)=0.
 \end{equation*}
 Since the polynomial $(\tau_1\tau_2-1)^{57}$ belongs to the ideal of~$\Z\bigl[\tau_1,\tau_2^{-1},\tau_2,\tau_2^{-1}\bigr]$ generated by~$(\tau_1-1)^2$ and~$(\tau_2-1)^{56}$, we obtain the required relation~\eqref{eq_annih_57}.
\end{proof}

Proposition~\ref{propos_A2} follows immediately from Proposition~\ref{propos_annih_120_54}.

\section{Transvection $P$-torsion property for~$\bA/\bA_2$}\label{section_A/A2}

The goal of this section is to derive from part~2 of Theorem~\ref{thm_long_rel} the following proposition.

\begin{propos}\label{propos_A/A2}
The quotient group~$\bA/\bA_2$ is a transvection $(t^{24}-1)$-torsion $\Z[\Sp(H)]$-module.
\end{propos}

Let $a_1,a_2,a_3,b_1,b_2,b_3$ be a symplectic basis of~$H$. As in Theorem~\ref{thm_long_rel}, set
$U_i=\langle a_i,b_i\rangle$ and $t_i=t_{a_j+a_k}$, where $i,j,k$ is a pertmutation of~$1,2,3$. Also set
$$
\CA=\CA(U_1,U_2)=\CA(U_2,U_3)=\CA(U_3,U_1),
$$
and denote by~$\oCA$ the image of~$\CA$ in the quotient group~$\bA/\bA_2$. By a result of Vladimirov~\cite{Vla25}, we have $2\CA=0$. Hence $2\oCA=0$. Throughout the rest of this section, we will keep this relation  in mind, and therefore perform all computations  modulo~$2$, that is, regard~$\bA/\bA_2$ as an  $\F_2\bigl[\Sp(H)\bigr]$-module rather than a $\Z\bigl[\Sp(H)\bigr]$-module.

All type~1 abelian cycles in~$H_2(\I)$ lie in the $\Sp(H)$-orbit of~$\CA$. Therefore $\bA/\bA_2$ is an $\F_2\bigl[\Sp(H)\bigr]$-module with a single generator~$\oCA$. Hence to prove Proposition~\ref{propos_A/A2} it suffices to show that
$
t_x^{24}\oCA=\oCA
$
for all $x\in H$.

 By part~2 of Theorem~\ref{thm_long_rel}, we have the following relations in~$\bA/\bA_2$:
\begin{gather}
   (t_1+1)(t_2+1)(t_3+1)\oCA=0,\label{eq_rel2_bar'}\\
   \left(t_1+t_1^{-1}+t_2+t_2^{-1}+t_3+t_3^{-1}\right)\oCA=0.\label{eq_rel1_bar'}
 \end{gather}
Since each of the transvections~$t_{a_1}$, $t_{a_2}$, and~$t_{a_3}$ fixes~$\oCA$, from Lemma~\ref{lem_transvec3} it follows that
\begin{equation}\label{eq_3form_oCA}
 t_{n_1a_1+n_2a_2+n_3a_3}\oCA=t_1^{n_2n_3}t_2^{n_1n_3}t_3^{n_1n_2}\oCA.
\end{equation}

\begin{propos}
 For $i\ne j$, we have
 \begin{equation}\label{eq_rel_aux_bar1}
  (t_i+1)(t_j+1)(t_it_j+1)\bigl(t_it_j^{-1}+1\bigr)\oCA=0.
 \end{equation}
\end{propos}

\begin{proof}
 We may assume that $(i,j)=(1,2)$. Multiplying \eqref{eq_rel1_bar'} by $(t_1+1)(t_2+1)$ and using~\eqref{eq_rel2_bar'}, we obtain that
 $$
 (t_1+1)(t_2+1)\left(t_1+t_1^{-1}+t_2+t_2^{-1}\right)\oCA=0.
 $$
 Since
 $$
 t_1+t_1^{-1}+t_2+t_2^{-1}=t_1^{-1}(t_1t_2+1)\bigl(t_1t_2^{-1}+1\bigr),
 $$
 this immediately yields the required relation.
\end{proof}

\begin{propos}
 For $i\ne j$, we have
 \begin{equation}\label{eq_rel_aux_bar2}
  (t_i+1)^2(t_j+1)(t_it_j+1)\bigl(t_i^2t_j+1\bigr)\oCA=0.
 \end{equation}
\end{propos}

\begin{proof}
We may assume that $(i,j)=(1,2)$.
Theorem~\ref{thm_long_rel} holds for any symplectic basis of~$H$. So relation~\eqref{eq_rel2_bar'} will remain true if we replace the symplectic transvections~$t_1$, $t_2$, and~$t_3$ by the transvections $t_{c_2+c_3}$, $t_{c_3+c_1}$, and~$t_{c_1+c_2}$, respectively, where $c_i\in U_i$ are arbitrary primitive elements. In particular, for the triple $(c_1,c_2,c_3)=(b_1,a_2,a_3)$, we have
$$
(t_{b_1+a_2}+1)(t_{b_1+a_3}+1)(t_{a_2+a_3}+1)\oCA=0.
$$
The transvection~$t_{b_1}$ commutes with $t_{b_1+a_2}$, $t_{b_1+a_3}$, and~$t_{a_2+a_3}$, and it fixes~$\oCA$. Hence the last relation takes the form:
\begin{equation}\label{eq_g_1g_2g_3}
(g_1+g_2+g_3+1)(t_1+1)\oCA=0,
\end{equation}
where
\begin{align*}
 g_1&=t_{b_1+a_2}t_{b_1+a_3}t_{b_1}^{-2},&
 g_2&=t_{b_1+a_2}t_{b_1}^{-1},&
 g_3&=t_{b_1+a_3}t_{b_1}^{-1}.
\end{align*}
Multiply this relation on the left by the following element $q\in\F_2\bigl[\Sp_6(\Z)\bigr]$:
$$
q=\left(t_{a_1+a_3}t_{a_1}^{-1}+1\right)
\left(t_{a_2+a_3}t_{a_1+a_3}t_{a_1}^{-1}+1\right)
\left(t_{a_2+a_3}^2t_{a_1+a_3}t_{a_1}^{-1}+1\right)
\left(t_{a_2+a_3}t_{a_1+a_2}t_{a_1}^{-1}+1\right).
$$
We have $g^{-1}t_yg=t_{g^{-1}(y)}$ for all $g\in\Sp(H)$ and $y\in H$. Each of the elements~$g_1^{-1}$, $g_2^{-1}$, and~$g_3^{-1}$ fixes~$a_2$ and~$a_3$; moreover,
\begin{align*}
 g_1^{-1}(a_1)&=a_1-a_2-a_3,&g_2^{-1}(a_1)&=a_1-a_2,& g_3^{-1}(a_1)&=a_1-a_3.
\end{align*}
A direct computation using these formulas yields:
\begin{align*}
g_1^{-1}qg_1&=\left(t_{a_1-a_2}t_{a_1-a_2-a_3}^{-1}+1\right)
\left(t_{a_2+a_3}t_{a_1-a_2}t_{a_1-a_2-a_3}^{-1}+1\right)\\
&{}\hspace{4.5cm}
\left(t_{a_2+a_3}^2t_{a_1-a_2}t_{a_1-a_2-a_3}^{-1}+1\right)\left(t_{a_2+a_3}t_{a_1-a_3}t_{a_1-a_2-a_3}^{-1}+1\right),\\
g_2^{-1}qg_2&=\left(t_{a_1-a_2+a_3}t_{a_1-a_2}^{-1}+1\right)
\left(t_{a_2+a_3}t_{a_1-a_2+a_3}t_{a_1-a_2}^{-1}+1\right)\\
&{}\hspace{4.5cm}
\left(t_{a_2+a_3}^2t_{a_1-a_2+a_3}t_{a_1-a_2}^{-1}+1\right)\left(t_{a_2+a_3}t_{a_1}t_{a_1-a_2}^{-1}+1\right),\\
g_3^{-1}qg_3&=\left(t_{a_1}t_{a_1-a_3}^{-1}+1\right)
\left(t_{a_2+a_3}t_{a_1}t_{a_1-a_3}^{-1}+1\right)\\
&{}\hspace{4.5cm}
\left(t_{a_2+a_3}^2t_{a_1}t_{a_1-a_3}^{-1}+1\right)
\left(t_{a_2+a_3}t_{a_1+a_2-a_3}t_{a_1-a_3}^{-1}+1\right).
\end{align*}
Further, applying~\eqref{eq_3form_oCA} to all transvections in these three formulas, we obtain:
\begin{align*}
 qg_1(t_1+1)\oCA&=g_1\left(t_1^{-1}t_2+1\right)(t_2+1)(t_1t_2+1)(t_3+1)(t_1+1)\oCA=0,\\
 qg_2(t_1+1)\oCA&=g_2\left(t_1^{-1}t_2+1\right)(t_2+1)(t_1t_2+1)(t_1t_3+1)(t_1+1)\oCA=0,\\
 qg_3(t_1+1)\oCA&=g_3(t_2+1)(t_1t_2+1)\bigl(t_1^2t_2+1\bigr)(t_3+1)(t_1+1)\oCA=0.
\end{align*}
Here the vanishing of the entries in the first and third rows follows from~\eqref{eq_rel2_bar'}, while the vanishing of the entry in the second row from~\eqref{eq_rel_aux_bar1}. Therefore relation~\eqref{eq_g_1g_2g_3} multiplied on the left by~$q$ reads as
$$
q(t_1+1)\oCA=0.
$$
Using~\eqref{eq_3form_oCA} once more, we obtain:
\begin{equation*}
(t_2+1)(t_1t_2+1)\bigl(t_1^2t_2+1\bigr)(t_1t_3+1)(t_1+1)\oCA=0.
\end{equation*}
Thus the element~$t_1t_3+1$ annihilates the element
$$
u=(t_1+1)(t_2+1)(t_1t_2+1)\bigl(t_1^2t_2+1\bigr)\oCA.
$$
On the other hand, from~\eqref{eq_rel2_bar'} it follows that the element $t_3+1$ annihilates~$u$. Hence, the element $t_1+1$ also annihilates~$u$, which yields the required relation~\eqref{eq_rel_aux_bar2}.
\end{proof}

\begin{propos}\label{propos_rel_aux_bar3}
 For $i\ne j$, we have
 \begin{equation*}
  (t_i+1)^4(t_i^3+1)(t_j+1)\oCA=0.
 \end{equation*}
\end{propos}

\begin{proof}
 Consider the element
 $$
 v=(t_i+1)^2(t_j+1)(t_it_j+1)\oCA.
 $$
 By~\eqref{eq_rel_aux_bar2}, the element~$t_i^2t_j+1$ annihilates~$v$. On the other hand, from~\eqref{eq_rel_aux_bar1} it follows that $t_i^{-1}t_j+1$ annihilates~$v$. Hence $t_i^3+1$ also annihilates~$v$, that is,
 \begin{equation*}
 (t_i+1)^2(t_i^3+1)(t_j+1)(t_it_j+1)\oCA=0.
 \end{equation*}
 Now, note that this relation will remain true if we replace the homology classes $a_1$, $a_2$, and~$a_3$ with any other triple of primitive homology classes from~$U_1$, $U_2$, and~$U_3$, respectively, and replace~$t_i$ and~$t_j$ accordingly. In particular, changing the sign of~$a_i$ results in the replacement of~$t_j$ with~$t_j^{-1}$. Therefore,
 \begin{equation*}
 (t_i+1)^2\bigl(t_i^3+1\bigr)\bigl(t_j^{-1}+1\bigr)\bigl(t_it_j^{-1}+1\bigr)\oCA=0.
 \end{equation*}
  It follows that both $t_it_j+1$ and~$t_it_j^{-1}+1$ annihilate the element
  $$w=(t_i+1)^2(t_i^3+1)(t_j+1)\oCA.$$
  Thus the element $t_i^2+1=(t_i+1)^2$ also annihilates~$w$. The proposition follows.
\end{proof}

\begin{propos}\label{propos_rel_aux_bar5}
 For $i=1,2,3$, we have
 \begin{equation}\label{eq_rel_aux_bar5}
  (t_i+1)^6(t_i^3+1)\oCA=0.
 \end{equation}
\end{propos}

\begin{proof}
 Multiply relation~\eqref{eq_rel1_bar'} by $(t_i+1)^4(t_i^3+1)$. From Proposition~\ref{propos_rel_aux_bar3} it follows that the element $(t_i+1)^4(t_i^3+1)$ annihilates both the element $\left(t_j+t_j^{-1}\right)\oCA$ and the element $\left(t_k+t_k^{-1}\right)\oCA$, where $i,j,k$ is a permutation of~$1,2,3$. Therefore we obtain the required relation~\eqref{eq_rel_aux_bar5}.
\end{proof}

\begin{propos}\label{propos_tx24}
 We have $t_x^{24}\oCA=\oCA$ for all $x\in H$.
\end{propos}

\begin{proof}
 The polynomial~$(t+1)^6(t^3+1)$ divides $t^{24}+1$ in $\F_2[t]$. Therefore Proposition~\ref{propos_rel_aux_bar5} implies that $t_{a_j+a_k}^{24}\oCA=\oCA$ for $j\ne k$. Moreover, these relations hold for any choice of the primitive homology classes~$a_i\in U_i$, $i=1,2,3$.

 Consider an arbitrary homology class $x\in H$. We have $x=n_1a_1+n_2a_2+n_3a_3$ for some primitive homology classes $a_i\in U_i$. Then $t_{a_j+a_k}^{24}\oCA=\oCA$ for $j\ne k$. The proposition now follows from~\eqref{eq_3form_oCA}, since the symplectic transvections~$t_{a_j+a_k}$ commute with each other and with~$t_x$.
\end{proof}

Proposition~\ref{propos_A/A2} follows from Proposition~\ref{propos_tx24}.

\section{Proofs of Theorems~\ref{thm_torsion} and~\ref{thm_main}}\label{section_proof_main}

In this section, we complete the  derivation of Theorems~\ref{thm_torsion} and~\ref{thm_main} from Theorems~\ref{thm_long_rel} and~\ref{thm_toAC}, which will be proved in Parts~2 and~3 of this paper, respectively.

Theorem~\ref{thm_torsion} follows immediately from Propositions~\ref{propos_A2} and~\ref{propos_A/A2}, and Theorem~\ref{thm_toAC}.

To derive Theorem~\ref{thm_main} from Theorem~\ref{thm_torsion}, we need to prove the following proposition.

\begin{propos}[$=$ Proposition~\ref{propos_general_suffcond}]\label{propos_general_suffcond'}
 Let $G$ be a group, and $\CT\subseteq G$ a conjugation-invariant subset. Assume that $G$ has finite width with respect to~$\CT$. Let $M$ be a finitely generated $\Z[G]$-module. Suppose that there exists a monic polynomial $P\in \Z[t]$ with constant term $\pm1$,
 $$
 P=t^n+c_{n-1}t^{n-1}+\cdots+c_1t\pm 1,
 $$
  such that $M$ is a $(P,\CT)$-torsion module.
 Then $M$ is finitely generated as an abelian group.
\end{propos}

\begin{proof}
 The proof mimics the proof of Proposition~4.3 in~\cite{Gai26}.

 For each $\tau\in\CT$ and each element $v\in M$, we have $P(\tau)v=0$. Since $P$ is monic and has constant term~$\pm1$, it follows that any element $\tau^kv$ with $k\in\Z$ can be expressed as a $\Z$-linear combination of the elements $v,\tau v,\ldots,\tau^{n-1}v$.

The group~$G$ has finite width with respect to~$\CT$. Therefore
$$
G=\langle \tau_1\rangle\cdots\langle\tau_m\rangle
$$
for some $\tau_1,\ldots,\tau_m\in\CT$.

Let $u_1,\ldots,u_q$ form a finite set of generators for~$M$ as a $\Z[G]$-module.

Each element of~$G$ can be written as $g=\tau_1^{k_1}\cdots\tau_m^{k_m}$.
Taking $\tau=\tau_1$ and $v=\tau_2^{k_2}\cdots\tau_m^{k_m}u_s$, we can write $gu_s$ as a $\Z$-linear combination of elements of the form $\tau_1^{l_1}\tau_2^{k_2}\cdots\tau_m^{k_m}u_s$ with $0\le l_1<n$. Then, using
$$
\tau_1^{l_1}\tau_2^{k_2}\cdots\tau_m^{k_m}u_s=\left(\tau_1^{l_1}\tau_2\tau_1^{-l_1}\right)^{k_2}\tau_1^{l_1}\tau_3^{k_3}\cdots\tau_m^{k_m}u_s
$$
and setting $\tau=\tau_1^{l_1}\tau_2\tau_1^{-l_1}$ and $v=\tau_1^{l_1}\tau_3^{k_3}\cdots\tau_m^{k_m}u_s$, we can further write $gu_s$ as a $\Z$-linear combination of elements of the form $\tau_1^{l_1}\tau_2^{l_2}\tau_3^{k_3}\cdots\tau_m^{k_m}u_s$ with $0\le l_1,l_2<n$. Proceeding inductively, we eventually express $gu_s$ as a $\Z$-linear combination of elements $\tau_1^{l_1}\cdots\tau_m^{l_m}u_s$ where $0\le l_i<n$ for each~$i$. Thus $M$ is generated as an abelian group by these $n^mq$ elements.
\end{proof}

\begin{proof}[Proof of Theorem~\ref{thm_main}]
Let $G=\Sp(H)\cong\Sp_6(\Z)$, and let $\CT$ be the set of all symplectic transvections~$t_x$, where $x\in H$. Then the set~$\CT$ is conjugation-invariant. Moreover, by a result of Tavgen~\cite{Tav91} (see Proposition~\ref{propos_Tavgen} in the Introduction), the group~$G$ has finite width with respect to~$\CT$. Theorem~\ref{thm_main} now follows from Theorem~\ref{thm_torsion} and Proposition~\ref{propos_general_suffcond'}.
\end{proof}

\section{Spectral sequence for the action of a group on a CW complex}\label{section_SpecSeq}

Let $X$ be a CW complex. Suppose that a group~$G$ acts cellularly and without rotations on~$X$. The latter condition means that an element~$g\in G$ fixes a cell~$Q$ of~$X$ pointwise whenever $g$ fixes~$Q$ setwise. We denote by~$G_Q$ the stabilizer of~$Q$ in~$G$.

In this paper, we will make constant use of the spectral sequence
\begin{equation}\label{eq_SpecSeq}
 E^1_{p,q}=\bigoplus_{Q\in\fX_p}H_q(G_Q)\quad\Longrightarrow\quad H_{p+q}^G(X)
\end{equation}
associated with the action of~$G$ on~$X$. Here $\fX_p$ is any set of representatives of the $G$-orbits of $p$-cells of~$X$, and $H_*^G(X)$ is the equivariant homology of~$X$. The coefficient group is always taken to be the (untwisted) group~$\Z$; we omit it from the notation. The necessary background on spectral sequence~\eqref{eq_SpecSeq} can be found in the book~\cite{Bro82}, where it is called the \textit{second spectral sequence} (see also~\cite{Eve91}). For the applications of this spectral sequence to the study of the Torelli groups and related groups, see~\cite{BBM10}, \cite{Gai21}, \cite{Gai22}, \cite{Gai24}, \cite{Gai26}, \cite{Spi23}, \cite{Spi24}. In this section we recall the known facts about spectral sequence~\eqref{eq_SpecSeq} that we need for the subsequent computations; all these facts can be found in~\cite[Chapter~VII]{Bro82} and~\cite[Chapter~7]{Eve91}. After that, we derive some consequences of these facts, which also appear to be standard, but the author is not aware of precise references in the literature where they can be found in exactly this form.

We find it convenient to work throughout with left $\Z[G]$-modules. Recall that the tensor product over~$G$ of two left $\Z[G]$-modules~$A$ and~$B$ is, by definition, the abelian group
$$
A\otimes_GB=A\otimes_{\Z}B/\langle(ga)\otimes (gb)-a\otimes b\mid a\in A, b\in B, g\in G\rangle.
$$
We will use repeatedly  the relation
$$
(g^{-1}a)\otimes b=a\otimes(gb)
$$
in $A\otimes_GB$.

 Let $(\CR_*,d)$ be the \textit{normalized bar resolution} of~$\Z$ (with trivial $G$-action) over~$\Z[G]$. Recall that $\CR_q$ is a free left $\Z[G]$-module with basis consisting of all symbols $[g_1|\ldots|g_q]$, where $g_i\in G\setminus\{1\}$, and the differential $d\colon\CR_q\to\CR_{q-1}$ is given by
\begin{multline*}
d[g_1|\ldots|g_q]=g_1[g_2|\ldots|g_q]+\sum_{i=1}^{q-1}(-1)^i[g_1|\ldots|g_{i-1}|g_ig_{i+1}|g_{i+2}|\ldots |g_q]\\{}+(-1)^q[g_1|\ldots|g_{q-1}].
\end{multline*}
(By definition, $[g_1|\ldots|g_q]=0$ whenever $g_i=1$ for some~$i$.) In particular, $\CR_0$ is a free $\Z[G]$-module with a single generator~$\ebr$. The chain complex~$\CR_*$ is a resolution of~$\Z$, that is, the following sequence is exact:
$$
 \ldots\longrightarrow\CR_2\stackrel{d}{\longrightarrow}\CR_1\stackrel{d}{\longrightarrow}\CR_0
 \longrightarrow\Z\longrightarrow 0.
 $$
In most constructions below, the normalized bar resolution can be replaced by any projective resolution of~$\Z$ over~$\Z[G]$. Nevertheless, for specific computations, it is convenient for us to work specifically with the normalized bar resolution. Recall that $H_*(G)$ is the homology of the chain complex $\Z\otimes_G\CR_*$.

Now, let $C_*(X)$ be the cellular chain complex of~$X$, and let~$\partial$ denote its differential. We consider the bicomplex
\begin{equation*}
B_{p,q}=C_p(X)\otimes_G\CR_q
\end{equation*}
with the differentials
$$
\partial'\colon B_{p,q}\to B_{p-1,q},\qquad \partial''\colon B_{p,q}\to B_{p,q-1}
$$
given by
\begin{align*}
\partial'(Q\otimes\xi)&=(\partial Q)\otimes \xi,\\ \partial''(Q\otimes\xi)&=(-1)^p(Q\otimes d\xi).
\end{align*}
Our sign convention is to include the sign~$(-1)^p$ in the latter formula, so that $\partial'$ and~$\partial''$ anticommute.

By definition, the equivariant homology $H_*^G(X)$ is the homology of the total complex
$$
\Tot_s(B_{*,*})=\bigoplus_{p+q=s}B_{p,q}
$$
with differential
$$
\partial_{\tot}=\partial'+\partial''.
$$
If $X$ is contractible, then we have $H_*^G(X)\cong H_*(G)$. Moreover, this isomorphism is canonical, since it is induced by the map $X\to\mathrm{pt}$. Hence we may identify these  homology groups.

Spectral sequence~\eqref{eq_SpecSeq} is the spectral sequence associated  with the filtration by columns in the bicomplex~$B_{*,*}$. This means that we can start the spectral sequence from $E^0_{p,q}=B_{p,q}$ with $d^0=\partial''$. The following are standard facts about this spectral sequence.

\begin{fact}\label{fact_E1}
 If\/ $\fX_p$ is a set of $p$-cells of~$X$ containing exactly one cell in every $G$-orbit, then
 \begin{equation}\label{eq_SpecSeqE1}
 E^1_{p,q}= \bigoplus_{Q\in\fX_p}H_q(G_Q).
\end{equation}
The embedding
$$
\iota_Q\colon H_q(G_Q)\hookrightarrow E^1_{p,q}
$$
is induced by sending each homology class $x\in H_q(G_Q)$ to $Q\otimes \xi\in B_{p,q}$, where $\xi$ is the cycle representing~$x$ in $\Z\otimes_{G_Q}\CR_q(G_Q)$. (Here $\CR_*(G_Q)$ is the normalized bar resolution of\/~$\Z$ over~$\Z[G_Q]$; it is a subcomplex of the normalized bar resolution~$\CR_*$ of\/~$\Z$ over~$\Z[G]$.) The homomorphism~$\iota_Q$ changes sign whenever we change the orientation of~$Q$.

Moreover, if we replace~$Q$ by another representative $gQ$ of the same $G$-orbit, then the corresponding inclusions~$\iota_Q$ and~$\iota_{gQ}$ are related by the commutative diagram
\begin{equation}\label{eq_cd_iota}
\begin{tikzcd}[row sep=small]
 H_q(G_Q) \arrow[dd,"{g_*}","{\cong}"'] \arrow[rd,hookrightarrow,"{\iota_Q}"] & \\
 & E^1_{p,q}\\
 H_q(G_{gQ})  \arrow[ru,hookrightarrow,"{\iota_{gQ}}"'] &
\end{tikzcd}
\end{equation}
where $g_*$ is the isomorphism induced by the conjugation by~$g$. In particular, the images of~$\iota_Q$ and~$\iota_{gQ}$ coincide, and therefore the summands of the direct sum decomposition~\eqref{eq_SpecSeqE1} do not depend on the choice of the set of representatives~$\fX_p$.
\end{fact}

\begin{fact}\label{fact_dr}
 The spectral sequence~\eqref{eq_SpecSeq} have differentials~$d^r$, $r=1,2,\ldots$, of bidegrees~$(-r,r-1)$, respectively. Suppose that $z\in E^1_{p,q}$. Let $r\ge 1$ and let
 $$
 \theta_{p,q}\in B_{p,q},\ \theta_{p-1,q+1}\in B_{p-1,q+1},\ \ldots,\ \theta_{p-r+1,q+r-1}\in B_{p-r+1,q+r-1}
 $$
 be elements satisfying the following conditions:
 \begin{itemize}
  \item $\partial''\theta_{p,q}=0$ and $\theta_{p,q}$ represents~$z$ in~$E^1_{p,q}$
  \item $\partial'\theta_{p-i,q+i}+\partial''\theta_{p-i-1,q+i+1}=0$ for $i=0,\ldots,r-2$.
 \end{itemize}
 Then $d^sz=0$ for $s=1,\ldots,r-1$, and the element $\partial'\theta_{p-r+1,q+r-1}$ represents in~$E^r_{p-r,q+r-1}$ the class~$d^rz$.

 In particular, the differential~$d^1$ is induced by~$\partial'$, that is, if $\partial Q=\sum\limits_P [Q:P]P$, then
 \begin{equation*}
 d^1\iota_Q(y)=\sum_P [Q:P]\,\iota_P\left(i_{Q,P}y\right),
 \end{equation*}
 where $i_{Q,P}\colon H_q(G_Q)\to H_q(G_P)$ is the homomorphism induced by the inclusion $G_Q\le G_P$.
\end{fact}

\begin{remark}
 In practice, to compute~$d^rz$  one typically needs to construct not only elements $\theta_{p,q},\ldots,\theta_{p-r+1,q+r-1}$ satisfying the conditions of Fact~\ref{fact_dr}, but also an additional element $\theta_{p-r,q+r}\in B_{p-r,q+r}$. Then the class~$d^rz$ is represented by the element $\partial'\theta_{p-r+1,q+r-1}+\partial''\theta_{p-r,q+r}$. At first glance, introducing such an element~$\theta_{p-r,q+r}$ changes nothing, since adding an element from the image of~$\partial''$ does not change the class in~$E^r_{p-r,q+r-1}$. Nevertheless, introducing the element~$\theta_{p-r,q+r}$ is useful. The point is that from the element $\partial'\theta_{p-r+1,q+r-1}$ it is often very difficult to tell which class in~$E^r_{p-r,q+r-1}$ it represents. Adding a suitable~$\partial''\theta_{p-r,q+r}$ is sometimes necessary in order to put~$\partial'\theta_{p-r+1,q+r-1}$ into a form whose homology class is easy to identify.
\end{remark}

\begin{fact}\label{fact_Einfty}
 The $E^{\infty}$-page of the spectral sequence~\eqref{eq_SpecSeq} is associated with a filtration in~$H_*^G(X)$. Specifically, for each~$n$, there is a filtration
 \begin{equation}\label{eq_filtration}
 0=\CF_{-1,n+1}\subseteq\CF_{0,n}\subseteq\CF_{1,n-1}\subseteq\cdots\subseteq \CF_{n,0}=H_n^G(X)
 \end{equation}
 with (canonical) isomorphisms
 \begin{equation}\label{eq_Einfty_isomorphism}
 \CF_{p,q}/\CF_{p-1,q+1}\cong E^{\infty}_{p,q}.
 \end{equation}
 Moreover, this is precisely the filtration induced by the filtration by columns in the bicomplex~$B_{*,*}$. This means that a class $z\in H_n^G(X)$ belongs to~$\CF_{p,n-p}$ if and only if $z$ can be represented by a cycle lying in the subgroup
 $$
 \bigoplus_{i\le p}B_{i,n-i}\subseteq \Tot_n(B_{*,*}).
 $$
\end{fact}

Recall that, if $X$ is contractible, then $H_*^G(X)\cong H_*(G)$, hence the spectral sequence converges to~$H_*(G)$.

Note that, since $\CR_0$ is a free $\Z[G]$-module with a single generator~$\ebr$, we have a canonical isomorphism $C_p(X)\otimes_G\CR_0\cong C_p(X)$. We always identify these groups along this isomorphism. Then the differential $\partial'$ in the $0$th row of~$B_{*,*}$ is identified with the usual boundary operator~$\partial$ in~$C_*(X)$.

\begin{fact}\label{fact_E0row}
 There is a canonical isomorphism~$E^1_{p,0}\cong C_p(X/G)$ such that the projection $B_{p,0}\twoheadrightarrow E^1_{p,0}$ coincides with the standard projection $C_p(X)\twoheadrightarrow C_p(X/G)$. Moreover, this isomorphism identifies the differential~$d^1$ in~$E^1_{*,0}$ with the boundary operator~$\partial$ in~$C_*(X/G)$. Thus there is a canonical isomorphism $E^2_{p,0}\cong H_p(X/G)$.
\end{fact}

Now, assume that the action of~$G$ on~$X$ extends to a cellular action of a larger group $\CG>G$ such that $G$ is a normal subgroup of~$\CG$. (We do not require that $\CG$ acts without rotations.) We set $\Gamma=\CG/G$. Then the group~$\Gamma$ acts on~$H_*^G(X)$.

\begin{fact}\label{fact_Gamma}
 Consider the action of~$\CG$ on the bicomplex
 $
 B_{p,q}=C_p(X)\otimes_G\CR_q
 $
 defined by
 $$
 f\bullet (Q\otimes g_0[g_1|\ldots|g_q])=(fQ)\otimes \bigl(fg_0f^{-1}\bigr)\bigl[fg_1f^{-1}\big|\ldots\big|fg_qf^{-1}\bigr].
 $$
 Then the differentials~$\partial'$ and~$\partial''$ become homomorphisms of~$\Z[\G]$-modules. Hence this action induces the action of\/~$\CG$ on the spectral sequence~$E^r_{p,q}$. Starting at the $E^1$-page, the action of the subgroup $G<\CG$ is trivial, so $E^r_{p,q}$ becomes a spectral sequence of (left) $\Z[\Gamma]$-modules. Moreover, the filtration~\eqref{eq_filtration} is $\Gamma$-equivariant, and~\eqref{eq_Einfty_isomorphism} is an isomorphism of $\Z[\Gamma]$-modules. Furthermore, the $\Gamma$-action on the $0$th row~$E^1_{*,0}$ coincides with the standard $\Gamma$-action on~$C_*(X/G)$.
\end{fact}

In this paper, we will explicitly compute differentials~$d^1$, $d^2$, and~$d^3$ using Fact~\ref{fact_dr}. Such computations are rather cumbersome. To shorten them, we would like to be able to perform computations in the bicomplex~$B_{*,*}$ modulo some $G$-equivariant subcomplex $Y\subseteq X$. The following proposition provides this possibility. Let
$$
B_{p,q}(Y)=C_p(Y)\otimes_G\CR_q\subseteq C_p(X)\otimes_G\CR_q=B_{p,q}
$$
be the subbicomplex corresponding to~$Y$. We will use the symbol~$\Yequivs$ to denote equality in~$B_{p,q}$ modulo~$B_{p,q}(Y)$.

\begin{lem}\label{lem_dr}
 Let $Y\subseteq X$ be a $G$-equivariant subcomplex such that $G$ acts freely on~$Y$. Suppose that $z\in E^1_{p,q}$, where $p,q\ge 0$. When $q=0$, we additionally assume that $d^1z=0$.  Let $r>1$, and let
 \begin{gather*}
 \theta_{p,q}\in B_{p,q},\ \theta_{p-1,q+1}\in B_{p-1,q+1},\ \ldots,\ \theta_{p-r,q+r}\in B_{p-r,q+r},\\
 \psi\in B_{p-r,q+r-1}
 \end{gather*}
 be elements satisfying the following conditions:
 \begin{itemize}
  \item $\partial''\theta_{p,q}=0$ and $\theta_{p,q}$ represents~$z$ in~$E^1_{p,q}$,
  \item $\partial'\theta_{p-i,q+i}+\partial''\theta_{p-i-1,q+i+1}\Yequiv 0$ for $i=0,\ldots,r-2$,
  \item $\partial'\theta_{p-r+1,q+r-1}+\partial''\theta_{p-r,q+r}\Yequiv \psi$,
  \item $\partial''\psi=0$.
 \end{itemize}
 Then $d^sz=0$ for $s=1,\ldots,r-1$, and  $\psi$ represents in~$E^r_{p-r,q+r-1}$ the class~$d^rz$.
\end{lem}

The proof of this lemma will be based on the following standard assertion.

\begin{lem}\label{lem_acyclic}
 If $G$ acts freely on~$Y$, then, for each~$p$, the chain complex $\bigl(B_{p,*}(Y),\partial''\bigr)$ is acyclic in positive degrees.
\end{lem}

\begin{proof}
Since  $G$ acts freely on~$Y$, we have that $C_p(Y)$ is a free $\Z[G]$-module. Since the resolution~$\CR_*$ is acyclic in positive degrees, it now follows that the chain complex $B_{p,*}(Y)=C_p(Y)\otimes_G\CR_*$ is also acyclic in positive degrees.
\end{proof}

\begin{proof}[Proof of Lemma~\ref{lem_dr}]
Let us show that we can add to the elements~$\theta_{p-i,q+i}$ certain elements $\varphi_{p-i,q+i}\in B_{p-i,q+i}(Y)$ so that the resulting elements $\tilde{\theta}_{p-i,q+i}=\theta_{p-i,q+i}+\varphi_{p-i,q+i}$ satisfy exact equalities
\begin{align}
\partial'\tilde\theta_{p-i,q+i}+\partial''\tilde\theta_{p-i-1,q+i+1}&=0,& i&=0,\ldots,r-2,\label{eq_p'p''1} \\
\partial'\tilde\theta_{p-r+1,q+r-1}+\partial''\tilde\theta_{p-r,q+r}&=\psi.&&\label{eq_p'p''2}
\end{align}
Then the assertion of the lemma will follow from Fact~\ref{fact_dr}.

We begin by setting $\varphi_{p,q}=0$, hence $\tilde\theta_{p,q}=\theta_{p,q}$, and then inductively construct $\varphi_{p-i,q+i}$ for $i=1,\ldots,r$. Assume that we have already constructed $\varphi_{p,q},\ldots, \varphi_{p-s+1,q+s-1}$ so that \eqref{eq_p'p''1} holds for $i=0,\ldots,s-2$. We need to find an element~$\varphi_{p-s,q+s}$ such that the next equality, namely~\eqref{eq_p'p''1} for $i=s-1$  if $s<r$ or~\eqref{eq_p'p''2} if $s=r$, is also satisfied. The equation for~$\varphi_{p-s,q+s}$ is as follows:
\begin{equation*}
\partial''\varphi_{p-s,q+s}=-\partial''\theta_{p-s,q+s}-\partial'\theta_{p-s+1,q+s-1}-\partial'\varphi_{p-s+1,q+s-1}+\left\{
\begin{aligned}
 &0&&\text{if}\ s<r,\\
 &\psi&&\text{if}\ s=r.
\end{aligned}
\right.
\end{equation*}
From the conditions of Lemma~\ref{lem_dr} and equality~\eqref{eq_p'p''1} for $i=s-2$, it follows that the right-hand side of this equation lies in~$B_{p-s,q+s-1}(Y)$ and in the kernel of~$\partial''$. If $q+s>1$, then the existence of the required element~$\varphi_{p-s,q+s}\in B_{p-s,q+s}(Y)$ is guaranteed by Lemma~\ref{lem_acyclic}.

The only case in which the above argument fails is the case $q=0$ and~$s=1$. This is why  Lemma~\ref{lem_dr} includes the additional assumption $d^1z=0$ when $q=0$. In this case, the equation for~$\varphi_{p-1,1}$ is as follows:
$$
\partial''\varphi_{p-1,1}=-\partial''\theta_{p-1,1}-\partial'\theta_{p,0}.
$$
The right-hand side of this equation lies both in~$B_{p-1,0}(Y)$ and in the image of~$\partial''$, since $\theta_{p,0}$ represents~$z$ and~$d^1z=0$. The existence of the required element~$\varphi_{p-1,1}$ now follows from the equality $B_{p-1,0}(Y)\cap\partial''B_{p-1,1}=\partial''B_{p-1,1}(Y)$.
\end{proof}

In this paper we will have to work extensively with abelian cycles. Suppose that the CW complex~$X$ on which the group~$G$ acts is contractible, so that the spectral sequence converges to~$H_*(G)$. Let $\CA(g_1,g_2)$ be an abelian cycle in~$H_2(G)$. Then $\CA(g_1,g_2)$ is represented by the cycle
$$
[g_1|g_2]-[g_2|g_1]\in\Z\otimes_G\CR_2.
$$
First, assume that $X$ contains a vertex~$v$ such that both~$g_1$ and~$g_2$ belong to~$G_v$. Then the abelian cycle~$\CA(g_1,g_2)$ is represented by the element $\iota_v(\CA_v(g_1,g_2))\in E^1_{0,2}$, where $\CA_v(g_1,g_2)$ is the abelian cycle in~$H_2(G_v)$ corresponding to the pair~$(g_1,g_2)$. Therefore $\CA(g_1,g_2)$ lies in the subgroup $\CF_{0,2}\subseteq H_2(G)$. Moreover, $\CA(g_1,g_2)$ can be represented by the cycle
$$
v\otimes([g_1|g_2]-[g_2|g_1])\in C_0(X)\otimes_G\CR_2\subseteq \Tot_2(B_{*,*}).
$$

However, we will also have to work with abelian cycles for which no vertex~$v$ with the described property exists. We will need the following lemma, which allows one, under certain conditions, to work with abelian cycles lying in~$\CF_{1,1}$; this lemma is related to Lemma 5.1 in~\cite{Gai24}.

\begin{lem}\label{lem_ab_F11}
Suppose  that $X$ is contractible. Let $g$ and~$h$ be commuting elements of~$G$. Suppose that $X$ contains an oriented $1$-cell~$P$ with the following properties:
 \begin{itemize}
  \item the element~$g$ fixes~$P$ pointwise,
  \item the element~$h$ sends the initial vertex~$v$ of~$P$ to its terminal vertex.
 \end{itemize}
 Then the abelian cycle~$\CA(g,h)$ belongs to the subgroup $\CF_{1,1}\subseteq H_2(G)$, and the projection $\CF_{1,1}\twoheadrightarrow E^{\infty}_{1,1}$ sends~$\CA(g,h)$ to the class represented by the element $-\iota_P([g])\in E^1_{1,1}$. Moreover, $\CA(g,h)$ is represented by the  element
 $$
 \fA=v\otimes\left(\bigl[h^{-1}\big|g\bigr]-\bigl[g\big|h^{-1}\bigr]\right)-P\otimes [g]\in \Tot_2(B_{*,*}).
 $$
\end{lem}

\begin{proof}
 We  have $\partial P=hv-v$, $gP=P$, and~$gv=v$. A direct computation shows that $\partial_{\tot}\fA=0$. Since $X$ is contractible, the augmentation $\eta\colon C_*(X)\to\Z$ induces a chain homotopy equivalence
$$
\Tot_*(B_{*,*})=\Tot_*\bigl(C_*(X)\otimes_G\CR_*\bigr)\stackrel{\eta_*}{\longrightarrow} \Z\otimes_G\CR_*.
$$
This chain homotopy equivalence sends~$\fA$ to the element
$$
\bigl[h^{-1}\big|g\bigr]-\bigl[g\big|h^{-1}\bigr]\in \Z\otimes_G\CR_2,
$$
which represents the abelian cycle~$\CA\bigl(h^{-1},g\bigr)=\CA(g,h)$. Hence $\fA$ also represents~$\CA(g,h)$. But $\fA$ lies in $B_{0,2}\oplus B_{1,1}$. Consequently, by Fact~\ref{fact_Einfty}, we have $\CA(g,h)\in\CF_{1,1}$. Moreover, by Fact~\ref{fact_E1}, the element $P\otimes [g]$ represents the class $\iota_P([g])$ in~$E^1_{1,1}$. Thus $\CA(g,h)$ projects to the class in~$E^{\infty}_{1,1}$ represented by~$-\iota_P([g])$.
\end{proof}

\part{Stabilizer of a three-component multicurve}\label{section_stab3}

\section{Main theorem}\label{section_IM}

The main goal of this part of the paper (Sections~\ref{section_IM}--\ref{section_end_thm3}) is  to prove Theorem~\ref{thm_long_rel}. In fact, all relations in Theorem~\ref{thm_long_rel} come from certain relations in~$H_2(\I_M)$, where $M\subseteq \Sigma$ is a three-component multicurve that does not separate~$\Sigma$, and $\I_M$ is the stabilizer of it in the Torelli group. So in this part we study in details the group~$H_2(\I_M)$. The main tool is the spectral sequence~\eqref{eq_SpecSeq} for the action of~$\I_M$ on the \textit{complex of relative cycles}~$\CB_{0,6}$ introduced by Spiridonov~\cite{Spi23},~\cite{Spi24} (see Section~\ref{section_RCC}). Along the way, while studying this spectral sequence, we will also prove that the group~$H_2(\I_M)$ is generated by standard abelian cycles; this result will then be used in Part~3 to prove Theorem~\ref{thm_toAC}. Note, however, that the study of the homology of~$\I_M$ is not a goal in itself for us. In this part we intentionally include only the computations in the spectral sequence and the results on the homology of~$\I_M$ that are required for the proof of Theorem~\ref{thm_long_rel}. In fact, the same methods (but with considerably more computation) can yield far more refined results on the homology of~$\I_M$; the author plans to address these in a separate paper.

To formulate an analogue of Theorem~\ref{thm_long_rel} for the stabilizer~$\I_M$ (from which Theorem~\ref{thm_long_rel} itself will follow easily), we need to introduce some notation for specific simple closed curves on~$\Sigma$.

Let $M=\alpha_1\cup\alpha_2\cup\alpha_3$ be a three-component multicurve on~$\Sigma=\Sigma_3$ that does not separate~$\Sigma$. Endow the components of~$M$ with some orientations, and let $a_1,a_2,a_3\in H$ be the homology classes of $\alpha_1,\alpha_2,\alpha_3$, respectively. Then $a_1,a_2,a_3$ form a basis of a Lagrangian submodule of~$H$. Extend it to a symplectic basis $a_1,a_2,a_3,b_1,b_2,b_3$ of~$H$, and identify $\Sp(H)=\Sp_6(\Z)$ using this basis. Let $\CX_M$ be the subgroup of~$\Mod(\Sigma)$ consisting of all mapping classes that fix every component of~$M$, and preserve the orientation of each component. A restriction of the short exact sequence
$$
1\longrightarrow \I\longrightarrow\Mod(\Sigma)\longrightarrow\Sp_6(\Z)\longrightarrow 1
$$
is the short exact sequence
\begin{equation}\label{eq_ses_IXW}
1\longrightarrow \I_M\longrightarrow\CX_M\longrightarrow W\longrightarrow 1,
\end{equation}
where $W$ is the subgroup of~$\Sp_6(\Z)$ consisting of all matrices
$$
\begin{pmatrix}
 I & S\\
 0 & I
\end{pmatrix}
$$
such that $I$ is the  $3\times 3$ identity matrix, and $S$ is a symmetric $3\times 3$ matrix. The group~$W$ is free abelian of rank~$6$ with basis  $t_{a_1},t_{a_2},t_{a_3},t_{a_2+a_3},t_{a_3+a_1},t_{a_1+a_2}$. This group acts on~$H_*(\I_M)$. However, the symplectic transvections~$t_{a_1},t_{a_2},t_{a_3}$ are induced by the Dehn twists~$T_{\alpha_1},T_{\alpha_2},T_{\alpha_3}$, respectively, which commute with all elements of~$\I_M$. Therefore these transvections act trivially on~$H_*(\I_M)$. So only the action of the group $V=\langle t_{a_2+a_3},t_{a_3+a_1},t_{a_1+a_2}\rangle\cong \Z^3$ is of interest.
Set $t_i=t_{a_j+a_k}$, where $i,j,k$ is a permutation of~$1,2,3$, and
$$
\Lambda=\Z[V]=\Z\left[t_1,t_1^{-1},t_2,t_2^{-1},t_3,t_3^{-1}\right].
$$
Then $H_*(\I_M)$ is a $\Lambda$-module and the natural map $H_*(\I_M)\to H_*(\I)$ is a homomorphism of $\Lambda$-modules.

It is convenient to arrange the surface~$\Sigma$ as follows. Consider the standard unit sphere~$\bS\subseteq\R^3$ centered at the origin, and a regular octahedron inscribed in it. It is convenient for us to denote the vertices of this regular octahedron by~$1,2,3,\bar{1},\bar{2},\bar{3}$ in such a way that $i$ and~$\bar{i}$ are antipodes for every~$i$. Our convention on orientations is that the vertices~$1,2,3$ appear in this order as one goes around the face~$123$ clockwise. Now make a small round hole around each vertex of the octahedron on the sphere~$\bS$, and glue the boundary circles of the holes corresponding to each pair of opposite vertices~$i$ and~$\bar{i}$ by parallel transport. Identify the genus~$3$ surface obtained as a result of this gluing with the surface~$\Sigma$; moreover, identify the curve formed by gluing the holes around the vertices~$i$ and~$\bar{i}$ with~$\alpha_i$. The orientation of~$\alpha_i$ is chosen so that it goes clockwise around~$i$ and counterclockwise around~$\bar{i}$. Connect every pair of vertices of the octahedron, except for the three pairs~$\bigl\{i,\bar{i}\bigr\}$, by the shortest geodesic segment. Then the obtained configuration of geodesic segments is the radial projection onto the sphere of the $1$-skeleton of the octahedron~$123\bar{1}\bar{2}\bar{3}$. Such octahedral representation of the genus-three surface was used by the author in~\cite{Gai25} to study the genus-three Johnson kernel~$\K_3$.

The described ``octahedral'' realization of~$\Sigma$ is convenient because symmetries of the octahedron will provide an additional group of automorphisms of the spectral sequence that we are going to study. However, it is rather inconvenient to draw configurations of curves on the unit sphere, primarily because half of the sphere is always invisible. Therefore, in the figures we will always unfold the sphere~$\bS$ onto the plane compactified by a point~$\infty$, corresponding to the center of the triangle~$\bar{1}\bar{2}\bar{3}$; see Fig.~\ref{fig_sphere}, where the geodesic segments between punctures in this representation are shown. It is worth noting that $\infty$ is not a puncture, so ``moving an arc through infinity'' is an admissible isotopy.

\begin{figure}[!ht]
 \begin{center}
\begin{tikzpicture}[scale=1.25,
         mid arrow/.style={
            postaction={
                decorate,
                decoration={
                    markings,
                    mark=at position 0.55 with {
                        \arrow[#1]{stealth}
                    }
                }
            }
        },
        mid arrow/.default=black,
    ]

 %\definecolor{myblue}{rgb}{0, 0, 0.7}
 \unitlength=1cm
 \fill [black!4] (0,0) circle (3.2);
 \thicklines

 \begin{scope}[white]
 \fill  (-.6,-.346) circle (.27);
\fill  (0,.693) circle (.27);
\fill  (.6,-.346) circle (.27);
\fill  (2,1.155) circle (.27);
\fill  (0,-2.309) circle (.27);
\fill  (-2,1.155) circle (.27);
 \end{scope}

 \begin{scope}
  \draw (-.6,-.346) -- (0,.693) -- (.6,-.346) -- (-.6,-.346) -- (-2,1.155) -- (0,.693) -- (2,1.155) -- (.6,-.346) -- (0,-2.309) -- (-.6,-.346);
 \begin{scope}
 \draw (-2,1.155) .. controls (-2.3,2.4) and (-.8,2.7)  .. (0,2.7) .. controls (.8,2.7) and (2.3,2.4) .. (2,1.155);
 \end{scope}
 \begin{scope}[rotate=120]
 \draw (-2,1.155) .. controls (-2.3,2.4) and (-.8,2.7)  .. (0,2.7) .. controls (.8,2.7) and (2.3,2.4) .. (2,1.155);
 \end{scope}
 \begin{scope}[rotate=-120]
 \draw (-2,1.155) .. controls (-2.3,2.4) and (-.8,2.7)  .. (0,2.7) .. controls (.8,2.7) and (2.3,2.4) .. (2,1.155);
 \end{scope}
 \end{scope}

\begin{scope}[white, opacity=.6]
 \fill  (-.6,-.346) circle (.27);
\fill  (0,.693) circle (.27);
\fill  (.6,-.346) circle (.27);
\fill  (2,1.155) circle (.27);
\fill  (0,-2.309) circle (.27);
\fill  (-2,1.155) circle (.27);
 \end{scope}

\draw [very thick] (-.6,-.346) circle (.27);
\draw [very thick] (0,.693) circle (.27);
\draw [very thick] (.6,-.346) circle (.27);
\draw [very thick] (2,1.155) circle (.27);
\draw [very thick] (0,-2.309) circle (.27);
\draw [very thick] (-2,1.155) circle (.27);

 \fill[black!30] (0,-2.309) circle (1pt) node [black,below=-1.5pt,yshift=.5pt] {\tiny $\bar{2}$} node [black,below=3.5mm,xshift=1mm] {$\alpha_2$};
 \fill[black!30] (-2,1.155) circle (1pt) node [black,above left=-3pt,yshift=-.5pt] {\tiny $\bar{3}$} node [black,below=3.5mm] {$\alpha_3$};
 \fill[black!30] (2,1.155) circle (1pt) node [black,above right=-3pt,yshift=-.5pt] {\tiny $\bar{1}$} node [black,below=3.5mm,xshift=2mm] {$\alpha_1$};
 \fill[black!30] (0,.693) circle (1pt) node [black,above=-1.5pt] {\tiny $2$} node [black,above=3.5mm] {$\alpha_2$};
 \fill[black!30] (-.6,-.346) circle (1pt) node [black,below left=-3pt] {\tiny $1$} node [black,below left=1.6mm] {$\alpha_1$};
 \fill[black!30] (.6,-.346) circle (1pt)node [black,below right=-3pt] {\tiny $3$} node [black,below right=1.6mm] {$\alpha_3$};

 %\draw (0,.693) arc(90:270:1.501);
 %\draw (-.6,-.346) arc(-150:30:1.501);
 %\draw (.6,-.346) arc(-30:150:1.501);

%  \begin{scope}[blue]
%  \draw [thick] (0,1.1) arc (90:-90:.407)-- +(-.5,0) arc (90:150:.794) arc(150:210:1.394) arc (210:270:.794) -- +(.5,0) arc (90:-90:.407)  arc (270:90:1.908) node [pos=.5, left=-.5pt] {$\delta_2$};
%  \end{scope}
%  \begin{scope}[blue,rotate=120]
%  \draw [thick] (0,1.1) arc (90:-90:.407)-- +(-.5,0) arc (90:150:.794) arc(150:210:1.394)  arc (210:270:.794) -- +(.5,0) arc (90:-90:.407)  arc (270:90:1.908) node [pos=.5, below] {$\delta_1$};
%  \end{scope}
%  \begin{scope}[blue,rotate=-120]
%  \draw [thick] (0,1.1) arc (90:-90:.407)-- +(-.5,0) arc (90:150:.794) arc(150:210:1.394) arc (210:270:.794) -- +(.5,0) arc (90:-90:.407)  arc (270:90:1.908) node [pos=.5, above right=-2pt] {$\delta_3$};
%  \end{scope}

 %1
 \begin{scope}[xshift=-.6cm, yshift=-.346cm]
 \begin{scope}[rotate=-68]
 \draw [thick](0,.27)-- (-.1,.31);
\draw [thick](0,.27)-- (-.08,.2);
\end{scope}
\end{scope}

%2
\begin{scope}[yshift=.693cm]
 \begin{scope}[rotate=164]
 \draw [thick](0,.27)-- (-.1,.31);
\draw [thick](0,.27)-- (-.08,.2);
\end{scope}
\begin{scope}[rotate=180]
 \draw [thick](0,.27)-- (-.1,.31);
\draw [thick](0,.27)-- (-.08,.2);
\end{scope}
\end{scope}

%3
\begin{scope}[xshift=.6cm, yshift=-.346cm]
 \begin{scope}[rotate=66]
 \draw [thick](0,.27)-- (-.1,.31);
\draw [thick](0,.27)-- (-.08,.2);
\end{scope}
 \begin{scope}[rotate=52]
 \draw [thick](0,.27)-- (-.1,.31);
\draw [thick](0,.27)-- (-.08,.2);
\end{scope}
\begin{scope}[rotate=38]
 \draw [thick](0,.27)-- (-.1,.31);
\draw [thick](0,.27)-- (-.08,.2);
\end{scope}
\end{scope}

%bar1
 \begin{scope}[xshift=2cm, yshift=1.155cm]
 \begin{scope}[rotate=128]
 \draw [thick](0,.27)-- (.1,.31);
\draw [thick](0,.27)-- (.08,.2);
\end{scope}
\end{scope}

%bar2
\begin{scope}[yshift=-2.309cm]
 \begin{scope}[rotate=16]
 \draw [thick](0,.27)-- (.1,.31);
\draw [thick](0,.27)-- (.08,.2);
\end{scope}
\begin{scope}
 \draw [thick](0,.27)-- (.1,.31);
\draw [thick](0,.27)-- (.08,.2);
\end{scope}
\end{scope}

%bar3
\begin{scope}[xshift=-2cm, yshift=1.155cm]
 \begin{scope}[rotate=-98]
 \draw [thick](0,.27)-- (.1,.31);
\draw [thick](0,.27)-- (.08,.2);
\end{scope}
\begin{scope}[rotate=-112]
 \draw [thick](0,.27)-- (.1,.31);
\draw [thick](0,.27)-- (.08,.2);
\end{scope}
\begin{scope}[rotate=-126]
 \draw [thick](0,.27)-- (.1,.31);
\draw [thick](0,.27)-- (.08,.2);
\end{scope}
\end{scope}
\end{tikzpicture}
\end{center}
 \caption{The surface~$\Sigma$}\label{fig_sphere}
 \bigskip
 \bigskip

  \begin{center}
\begin{tikzpicture}[scale=1.1,
         mid arrow/.style={
            postaction={
                decorate,
                decoration={
                    markings,
                    mark=at position 0.55 with {
                        \arrow[#1]{stealth}
                    }
                }
            }
        },
        mid arrow/.default=black,
    ]

 %\definecolor{myblue}{rgb}{0, 0, 0.7}
 \unitlength=1cm

 \begin{scope}
 \fill [black!4] (0,0) circle (3.2);
 \thicklines

 \begin{scope}[white]
 \fill  (-.6,-.346) circle (.27);
\fill  (0,.693) circle (.27);
\fill  (.6,-.346) circle (.27);
\fill  (2,1.155) circle (.27);
\fill  (0,-2.309) circle (.27);
\fill  (-2,1.155) circle (.27);
 \end{scope}

\draw [thick] (-.6,-.346) circle (.27);
\draw [thick] (0,.693) circle (.27);
\draw [thick] (.6,-.346) circle (.27);
\draw [thick] (2,1.155) circle (.27);
\draw [thick] (0,-2.309) circle (.27);
\draw [thick] (-2,1.155) circle (.27);

 \draw (0,-2.309) node {\footnotesize $\bar{2}$};
 \draw (-2,1.155) node {\footnotesize $\bar{3}$};
 \draw (2,1.155) node  {\footnotesize $\bar{1}$};
 \draw (0,.693) node {\footnotesize $2$};
 \draw (-.6,-.346) node  {\footnotesize $1$};
 \draw (.6,-.346)node {\footnotesize $3$};

 %\draw (0,.693) arc(90:270:1.501);
 %\draw (-.6,-.346) arc(-150:30:1.501);
 %\draw (.6,-.346) arc(-30:150:1.501);

 \begin{scope}[blue]
 \draw [thick] (0,1.1) arc (90:-90:.407)-- +(-.5,0) arc (90:150:.794) arc(150:210:1.394) arc (210:270:.794) -- +(.5,0) arc (90:-90:.407)  arc (270:90:1.908) node [pos=.5, left=-.5pt] {$\delta_2$};
 \end{scope}
 \begin{scope}[red,rotate=120]
 \draw [thick] (0,1.1) arc (90:-90:.407)-- +(-.5,0) arc (90:150:.794) arc(150:210:1.394)  arc (210:270:.794) -- +(.5,0) arc (90:-90:.407)  arc (270:90:1.908) node [pos=.5, below] {$\delta_1$};
 \end{scope}
 \begin{scope}[orange,rotate=-120]
 \draw [thick] (0,1.1) arc (90:-90:.407)-- +(-.5,0) arc (90:150:.794) arc(150:210:1.394) arc (210:270:.794) -- +(.5,0) arc (90:-90:.407)  arc (270:90:1.908) node [pos=.5, above right=-2pt] {$\delta_3$};
 \end{scope}

 %1
 \begin{scope}[xshift=-.6cm, yshift=-.346cm]
 \begin{scope}[rotate=-68]
 \draw (0,.27)-- (-.1,.31);
\draw (0,.27)-- (-.08,.2);
\end{scope}
\end{scope}

%2
\begin{scope}[yshift=.693cm]
 \begin{scope}[rotate=164]
 \draw (0,.27)-- (-.1,.31);
\draw (0,.27)-- (-.08,.2);
\end{scope}
\begin{scope}[rotate=180]
 \draw (0,.27)-- (-.1,.31);
\draw (0,.27)-- (-.08,.2);
\end{scope}
\end{scope}

%3
\begin{scope}[xshift=.6cm, yshift=-.346cm]
 \begin{scope}[rotate=66]
 \draw (0,.27)-- (-.1,.31);
\draw (0,.27)-- (-.08,.2);
\end{scope}
 \begin{scope}[rotate=52]
 \draw(0,.27)-- (-.1,.31);
\draw (0,.27)-- (-.08,.2);
\end{scope}
\begin{scope}[rotate=38]
 \draw (0,.27)-- (-.1,.31);
\draw (0,.27)-- (-.08,.2);
\end{scope}
\end{scope}

%bar1
 \begin{scope}[xshift=2cm, yshift=1.155cm]
 \begin{scope}[rotate=128]
 \draw (0,.27)-- (.1,.31);
\draw (0,.27)-- (.08,.2);
\end{scope}
\end{scope}

%bar2
\begin{scope}[yshift=-2.309cm]
 \begin{scope}[rotate=16]
 \draw (0,.27)-- (.1,.31);
\draw (0,.27)-- (.08,.2);
\end{scope}
\begin{scope}
 \draw (0,.27)-- (.1,.31);
\draw (0,.27)-- (.08,.2);
\end{scope}
\end{scope}

%bar3
\begin{scope}[xshift=-2cm, yshift=1.155cm]
 \begin{scope}[rotate=-98]
 \draw (0,.27)-- (.1,.31);
\draw (0,.27)-- (.08,.2);
\end{scope}
\begin{scope}[rotate=-112]
 \draw (0,.27)-- (.1,.31);
\draw (0,.27)-- (.08,.2);
\end{scope}
\begin{scope}[rotate=-126]
 \draw (0,.27)-- (.1,.31);
\draw (0,.27)-- (.08,.2);
\end{scope}
\end{scope}

\end{scope}

%%%%%%%%%%%
%%%%%%%%%%%

 \begin{scope}[xshift=7.5cm]
 \fill [black!4] (0,0) circle (3.2);
 \thicklines

 \begin{scope}[white]
 \fill  (-.6,-.346) circle (.27);
\fill  (0,.693) circle (.27);
\fill  (.6,-.346) circle (.27);
\fill  (2,1.155) circle (.27);
\fill  (0,-2.309) circle (.27);
\fill  (-2,1.155) circle (.27);
 \end{scope}

\draw [thick] (-.6,-.346) circle (.27);
\draw [thick] (0,.693) circle (.27);
\draw [thick] (.6,-.346) circle (.27);
\draw [thick] (2,1.155) circle (.27);
\draw [thick] (0,-2.309) circle (.27);
\draw [thick] (-2,1.155) circle (.27);

 \draw (0,-2.309) node {\footnotesize $\bar{2}$};
 \draw (-2,1.155) node {\footnotesize $\bar{3}$};
 \draw (2,1.155) node  {\footnotesize $\bar{1}$};
 \draw (0,.693) node {\footnotesize $2$};
 \draw (-.6,-.346) node  {\footnotesize $1$};
 \draw (.6,-.346)node {\footnotesize $3$};

 \begin{scope}[violet,rotate=-120]
  \draw[thick] (-.6,.061) arc (90:270:.407) node[pos=.65,above] {$\gamma_{12}$} --+(1.2,0) arc(-90:90:.407) -- +(-1.2,0);
 \end{scope}
 \begin{scope}[orange]
 \draw[thick] (.302,-.068) -- +(1.4,1.501) arc (137:-43:.407) -- +(-1.4,-1.501) node[pos=.5,below right=-2pt] {$\gamma_{3\bar{1}}$} arc (-43:-223:.407);
  \end{scope}
 \begin{scope}[red]
 \draw[thick] (0,-2.9) arc (270:150:2.9) arc (150:-30:.55) arc (150:270:1.8) node[pos=.5,below left=-2pt] {$\gamma_{\bar{2}\bar{3}}$} arc (90:-90:.55);
 \end{scope}
 \begin{scope}[xshift=.6cm, yshift=-.346cm]
 \begin{scope}[rotate=-154]
 \begin{scope}[xshift=-.6cm, yshift=.346cm, blue]
 \draw[thick] (.302,-.068) -- +(1.4,1.501) node[pos=.5,right] {$\gamma_{\bar{2}3}$} arc (137:-43:.407) -- +(-1.4,-1.501)  arc (-43:-223:.407);
  \end{scope}
  \end{scope}
  \end{scope}

 %\draw (0,.693) arc(90:270:1.501);
 %\draw (-.6,-.346) arc(-150:30:1.501);
 %\draw (.6,-.346) arc(-30:150:1.501);

%  \begin{scope}[blue,rotate=120]
%  \draw [thick] (0,1.1) arc (90:-90:.407)-- +(-.5,0) arc (90:150:.794) arc(150:210:1.394)  arc (210:270:.794) -- +(.5,0) arc (90:-90:.407)  arc (270:90:1.908) node [pos=.5, above] {$\delta_1$};
%  \end{scope}
%
%  \draw[thick,red] (2,1.705) arc (90:30:.55) arc (30:-150:2.051) node [pos=.55, below] {$\talpha_3$}  arc (210:75:.55) arc(-105:-15:1.57) arc (165:80:.55);

 %1
 \begin{scope}[xshift=-.6cm, yshift=-.346cm]
 \begin{scope}[rotate=-68]
 \draw (0,.27)-- (-.1,.31);
\draw (0,.27)-- (-.08,.2);
\end{scope}
\end{scope}

%2
\begin{scope}[yshift=.693cm]
 \begin{scope}[rotate=164]
 \draw (0,.27)-- (-.1,.31);
\draw (0,.27)-- (-.08,.2);
\end{scope}
\begin{scope}[rotate=180]
 \draw (0,.27)-- (-.1,.31);
\draw (0,.27)-- (-.08,.2);
\end{scope}
\end{scope}

%3
\begin{scope}[xshift=.6cm, yshift=-.346cm]
 \begin{scope}[rotate=66]
 \draw (0,.27)-- (-.1,.31);
\draw (0,.27)-- (-.08,.2);
\end{scope}
 \begin{scope}[rotate=52]
 \draw(0,.27)-- (-.1,.31);
\draw (0,.27)-- (-.08,.2);
\end{scope}
\begin{scope}[rotate=38]
 \draw (0,.27)-- (-.1,.31);
\draw (0,.27)-- (-.08,.2);
\end{scope}
\end{scope}

%bar1
 \begin{scope}[xshift=2cm, yshift=1.155cm]
 \begin{scope}[rotate=128]
 \draw (0,.27)-- (.1,.31);
\draw (0,.27)-- (.08,.2);
\end{scope}
\end{scope}

%bar2
\begin{scope}[yshift=-2.309cm]
 \begin{scope}[rotate=16]
 \draw (0,.27)-- (.1,.31);
\draw (0,.27)-- (.08,.2);
\end{scope}
\begin{scope}
 \draw (0,.27)-- (.1,.31);
\draw (0,.27)-- (.08,.2);
\end{scope}
\end{scope}

%bar3
\begin{scope}[xshift=-2cm, yshift=1.155cm]
 \begin{scope}[rotate=-98]
 \draw (0,.27)-- (.1,.31);
\draw (0,.27)-- (.08,.2);
\end{scope}
\begin{scope}[rotate=-112]
 \draw (0,.27)-- (.1,.31);
\draw (0,.27)-- (.08,.2);
\end{scope}
\begin{scope}[rotate=-126]
 \draw (0,.27)-- (.1,.31);
\draw (0,.27)-- (.08,.2);
\end{scope}
\end{scope}

\end{scope}

\end{tikzpicture}
\end{center}
 \caption{Simple closed curves ~$\delta_1$, $\delta_2$, $\delta_3$, $\gamma_{12}$, $\gamma_{\bar{2}\bar{3}}$, $\gamma_{\bar{2}3}$, and $\gamma_{3\bar{1}}$}\label{fig_deltagamma}
 \end{figure}

Let us now define several useful simple closed curves on~$\Sigma$. Let~$\delta_1$, $\delta_2$, and~$\delta_3$ be the simple closed curves shown in the left part of Fig.~\ref{fig_deltagamma}. The curves~$\delta_i$ are separating and disjoint from $M$. Hence $T_{\delta_i}\in \I_M$. Therefore, we have the abelian cycles  $\CA(T_{\delta_i},T_{\delta_j})\in H_2(\I_M)$.  We shall assume that the basis $a_1,a_2,a_3,b_1,b_2,b_3$ of~$H$ is chosen so that the first homology group of the genus~$1$ subsurface bounded by $\delta_i$ is~$\langle a_i,b_i\rangle$.

Let $p$ and~$q$ be distinct elements of the set~$\bigl\{1,2,3,\bar{1},\bar{2},\bar{3}\bigr\}$ such that $\{p,q\}\ne \bigl\{i,\bar{i}\bigr\}$. Then we denote by $\gamma_{pq}$ the circuit around the geodesic segment~$pq$. Certainly, we have $\gamma_{qp}=\gamma_{pq}$. For instance, the curves~$\gamma_{12}$, $\gamma_{\bar{2}\bar{3}}$, $\gamma_{\bar{2}3}$, and~$\gamma_{3\bar{1}}$ are shown in the right part of Fig.~\ref{fig_deltagamma}. All other curves~$\gamma_{pq}$ are obtained from these four by the clockwise and counterclockwise rotations by~$2\pi/3$ around the center~$o$ of the triangle~$123$.

Let $i,j,k$ be a cyclic permutation of~$1,2,3$. Then $\bigl\{\gamma_{j\bar{k}},\gamma_{\bar{j}k}\bigr\}$ is a bounding pair. However, this bounding pair is not disjoint from~$\delta_i$. To obtain a bounding pair disjoint from~$\delta_i$ we need to replace~$\gamma_{\bar{j}k}$ with another simple closed curve enclosing the same vertices~$\bar{j}$ and~$k$ of the octahedron. Specifically, consider the simple closed curve~$\tgamma_{\bar{2}3}$ depicted in the left part of Fig.~\ref{fig_AC2}, together with the simple closed curves~$\gamma_{2\bar{3}}$ and~$\delta_1$. Then $\bigl\{\gamma_{2\bar{3}},\tgamma_{\bar{2}3}\bigr\}$ is a bounding pair disjoint from~$\delta_1$ and from~$M$, hence we obtain the abelian cycle $\CA\bigl(T_{\delta_1},T_{\gamma_{2\bar{3}}}T_{\tgamma_{\bar{2}3}}^{-1}\bigr)\in H_2(\I_M)$.
Similarly, we have the abelian cycles~$\CA\bigl(T_{\delta_2},T_{\gamma_{3\bar{1}}}T_{\tgamma_{\bar{3}1}}^{-1}\bigr)$ and~$\CA\bigl(T_{\delta_3},T_{\gamma_{1\bar{2}}}T_{\tgamma_{\bar{1}2}}^{-1}\bigr)$, where the simple closed curves $\tgamma_{\bar{3}1}$ and~$\tgamma_{\bar{1}2}$ are obtained from~$\tgamma_{\bar{2}3}$ by the clockwise and counterclockwise rotations by~$2\pi/3$, respectively, around the center~$o$ of the triangle~$123$.

\begin{figure}[!ht]
 \begin{center}
\begin{tikzpicture}[scale=1.1,
         mid arrow/.style={
            postaction={
                decorate,
                decoration={
                    markings,
                    mark=at position 0.55 with {
                        \arrow[#1]{stealth}
                    }
                }
            }
        },
        mid arrow/.default=black,
    ]

 %\definecolor{myblue}{rgb}{0, 0, 0.7}
 \unitlength=1cm

 \begin{scope}
 \fill [black!4] (0,0) circle (3.2);
 \thicklines

 \begin{scope}[white]
 \fill  (-.6,-.346) circle (.27);
\fill  (0,.693) circle (.27);
\fill  (.6,-.346) circle (.27);
\fill  (2,1.155) circle (.27);
\fill  (0,-2.309) circle (.27);
\fill  (-2,1.155) circle (.27);
 \end{scope}

\draw [thick] (-.6,-.346) circle (.27);
\draw [thick] (0,.693) circle (.27);
\draw [thick] (.6,-.346) circle (.27);
\draw [thick] (2,1.155) circle (.27);
\draw [thick] (0,-2.309) circle (.27);
\draw [thick] (-2,1.155) circle (.27);

 \draw (0,-2.309) node {\footnotesize $\bar{2}$};
 \draw (-2,1.155) node {\footnotesize $\bar{3}$};
 \draw (2,1.155) node  {\footnotesize $\bar{1}$};
 \draw (0,.693) node {\footnotesize $2$};
 \draw (-.6,-.346) node  {\footnotesize $1$};
 \draw (.6,-.346)node {\footnotesize $3$};

 %\draw (0,.693) arc(90:270:1.501);
 %\draw (-.6,-.346) arc(-150:30:1.501);
 %\draw (.6,-.346) arc(-30:150:1.501);

 \begin{scope}[blue,rotate=120]
 \draw [thick] (0,1.1) arc (90:-90:.407)-- +(-.5,0) arc (90:150:.794) arc(150:210:1.394)  arc (210:270:.794) -- +(.5,0) arc (90:-90:.407)  arc (270:90:1.908) node [pos=.5, below] {$\delta_1$};
 \end{scope}

 \begin{scope}[red,rotate=120]
 \draw[thick] (.302,-.068) -- +(1.4,1.501) arc (137:-43:.407) -- +(-1.4,-1.501) node[pos=.5,above] {$\gamma_{2\bar{3}}$} arc (-43:-223:.407);
  \end{scope}

 \draw [red,thick] (0,-2.716) arc (-90:56:.407) node [pos=.4, below right=-2pt] {$\tgamma_{\bar{2}3}$} arc (-124:-150:4.15) arc (210:30:.52) arc(-150:-135:2.9) arc (-130:86:.44) arc (-94:-105:5.5) arc (75:195.8:.65) arc (-164.2:-149.2:8.235) arc (-149.2:-85:.407);

 %1
 \begin{scope}[xshift=-.6cm, yshift=-.346cm]
 \begin{scope}[rotate=-68]
 \draw (0,.27)-- (-.1,.31);
\draw (0,.27)-- (-.08,.2);
\end{scope}
\end{scope}

%2
\begin{scope}[yshift=.693cm]
 \begin{scope}[rotate=164]
 \draw (0,.27)-- (-.1,.31);
\draw (0,.27)-- (-.08,.2);
\end{scope}
\begin{scope}[rotate=180]
 \draw (0,.27)-- (-.1,.31);
\draw (0,.27)-- (-.08,.2);
\end{scope}
\end{scope}

%3
\begin{scope}[xshift=.6cm, yshift=-.346cm]
 \begin{scope}[rotate=66]
 \draw (0,.27)-- (-.1,.31);
\draw (0,.27)-- (-.08,.2);
\end{scope}
 \begin{scope}[rotate=52]
 \draw(0,.27)-- (-.1,.31);
\draw (0,.27)-- (-.08,.2);
\end{scope}
\begin{scope}[rotate=38]
 \draw (0,.27)-- (-.1,.31);
\draw (0,.27)-- (-.08,.2);
\end{scope}
\end{scope}

%bar1
 \begin{scope}[xshift=2cm, yshift=1.155cm]
 \begin{scope}[rotate=128]
 \draw (0,.27)-- (.1,.31);
\draw (0,.27)-- (.08,.2);
\end{scope}
\end{scope}

%bar2
\begin{scope}[yshift=-2.309cm]
 \begin{scope}[rotate=16]
 \draw (0,.27)-- (.1,.31);
\draw (0,.27)-- (.08,.2);
\end{scope}
\begin{scope}
 \draw (0,.27)-- (.1,.31);
\draw (0,.27)-- (.08,.2);
\end{scope}
\end{scope}

%bar3
\begin{scope}[xshift=-2cm, yshift=1.155cm]
 \begin{scope}[rotate=-98]
 \draw (0,.27)-- (.1,.31);
\draw (0,.27)-- (.08,.2);
\end{scope}
\begin{scope}[rotate=-112]
 \draw (0,.27)-- (.1,.31);
\draw (0,.27)-- (.08,.2);
\end{scope}
\begin{scope}[rotate=-126]
 \draw (0,.27)-- (.1,.31);
\draw (0,.27)-- (.08,.2);
\end{scope}
\end{scope}

\end{scope}

%%%%%%%%%%%
%%%%%%%%%%%

 \begin{scope}[xshift=7.5cm]
 \fill [black!4] (0,0) circle (3.2);
 \thicklines

 \begin{scope}[white]
 \fill  (-.6,-.346) circle (.27);
\fill  (0,.693) circle (.27);
\fill  (.6,-.346) circle (.27);
\fill  (2,1.155) circle (.27);
\fill  (0,-2.309) circle (.27);
\fill  (-2,1.155) circle (.27);
 \end{scope}

\draw [thick] (-.6,-.346) circle (.27);
\draw [thick] (0,.693) circle (.27);
\draw [thick,red] (.6,-.346) circle (.27) node [right=2mm] {$\alpha_3$};
\draw [thick] (2,1.155) circle (.27);
\draw [thick] (0,-2.309) circle (.27);
\draw [thick,red] (-2,1.155) circle (.27) node [below=2.5mm,xshift=1mm] {$\alpha_3$};

 \draw (0,-2.309) node {\footnotesize $\bar{2}$};
 \draw (-2,1.155) node {\footnotesize $\bar{3}$};
 \draw (2,1.155) node  {\footnotesize $\bar{1}$};
 \draw (0,.693) node {\footnotesize $2$};
 \draw (-.6,-.346) node  {\footnotesize $1$};
 \draw (.6,-.346)node {\footnotesize $3$};

 %\draw (0,.693) arc(90:270:1.501);
 %\draw (-.6,-.346) arc(-150:30:1.501);
 %\draw (.6,-.346) arc(-30:150:1.501);

 \begin{scope}[blue,rotate=120]
 \draw [thick] (0,1.1) arc (90:-90:.407)-- +(-.5,0) arc (90:150:.794) arc(150:210:1.394)  arc (210:270:.794) -- +(.5,0) arc (90:-90:.407)  arc (270:90:1.908) node [pos=.5, above] {$\delta_1$};
 \end{scope}

 \draw[thick,red] (2,1.705) arc (90:30:.55) arc (30:-150:2.051) node [pos=.55, below] {$\talpha_3$}  arc (210:75:.55) arc(-105:-15:1.57) arc (165:80:.55);

 %1
 \begin{scope}[xshift=-.6cm, yshift=-.346cm]
 \begin{scope}[rotate=-68]
 \draw (0,.27)-- (-.1,.31);
\draw (0,.27)-- (-.08,.2);
\end{scope}
\end{scope}

%2
\begin{scope}[yshift=.693cm]
 \begin{scope}[rotate=164]
 \draw (0,.27)-- (-.1,.31);
\draw (0,.27)-- (-.08,.2);
\end{scope}
\begin{scope}[rotate=180]
 \draw (0,.27)-- (-.1,.31);
\draw (0,.27)-- (-.08,.2);
\end{scope}
\end{scope}

%3
\begin{scope}[red,xshift=.6cm, yshift=-.346cm]
 \begin{scope}[rotate=66]
 \draw (0,.27)-- (-.1,.31);
\draw (0,.27)-- (-.08,.2);
\end{scope}
 \begin{scope}[rotate=52]
 \draw(0,.27)-- (-.1,.31);
\draw (0,.27)-- (-.08,.2);
\end{scope}
\begin{scope}[rotate=38]
 \draw (0,.27)-- (-.1,.31);
\draw (0,.27)-- (-.08,.2);
\end{scope}
\end{scope}

%bar1
 \begin{scope}[xshift=2cm, yshift=1.155cm]
 \begin{scope}[rotate=128]
 \draw (0,.27)-- (.1,.31);
\draw (0,.27)-- (.08,.2);
\end{scope}
\end{scope}

%bar2
\begin{scope}[yshift=-2.309cm]
 \begin{scope}[rotate=16]
 \draw (0,.27)-- (.1,.31);
\draw (0,.27)-- (.08,.2);
\end{scope}
\begin{scope}
 \draw (0,.27)-- (.1,.31);
\draw (0,.27)-- (.08,.2);
\end{scope}
\end{scope}

%bar3
\begin{scope}[red,xshift=-2cm, yshift=1.155cm]
 \begin{scope}[rotate=-98]
 \draw (0,.27)-- (.1,.31);
\draw (0,.27)-- (.08,.2);
\end{scope}
\begin{scope}[rotate=-112]
 \draw (0,.27)-- (.1,.31);
\draw (0,.27)-- (.08,.2);
\end{scope}
\begin{scope}[rotate=-126]
 \draw (0,.27)-- (.1,.31);
\draw (0,.27)-- (.08,.2);
\end{scope}
\end{scope}

\end{scope}

\end{tikzpicture}
\end{center}
 \caption{The abelian cycles $\CA\bigl(T_{\delta_1},T_{\gamma_{2\bar{3}}}T_{\tgamma_{\bar{2}3}}^{-1}\bigr)$ and $\CA\bigl(T_{\delta_1},T_{\talpha_{3}}T_{\alpha_3}^{-1}\bigr)$}\label{fig_AC2}
\end{figure}

Further, let $\talpha_1$ be the circuit around the  union of the geodesic segments~$12$ and~$1\bar{2}$, let $\talpha_2$ be the circuit around the  union of the geodesic segments~$23$ and~$2\bar{3}$, and let $\talpha_3$ be the circuit around the  union of the geodesic segments~$31$ and~$3\bar{1}$. The simple closed curve~$\talpha_3$ is depicted in the right part of Fig.~\ref{fig_AC2}, and the curves~$\talpha_1$ and~$\talpha_2$ are obtained from it by the clockwise and counterclockwise rotations by~$2\pi/3$, respectively, around the center~$o$ of the triangle~$123$. If $k=i-1\pmod 3$, then $\bigl\{\talpha_k,\alpha_k\bigr\}$ is a bounding pair disjoint from~$\delta_i$, and we obtain the abelian cycle $\CA\bigl(T_{\delta_i},T_{\talpha_k}T_{\alpha_k}^{-1}\bigr)\in H_2(\I_M)$.

We are now ready to formulate the main result of this part of the paper.

\begin{theorem}\label{theorem_3}
 \textnormal{(1)} The group $H_2(\I_M)$ is generated by standard abelian cycles. Moreover, $H_2(\I_M)$ is generated as a $\Lambda$-module by the $9$ abelian cycles
 \begin{align*}
  &\CA(T_{\delta_i},T_{\delta_j}),&& \CA\bigl(T_{\delta_i},T_{\gamma_{j\bar{k}}}T_{\tgamma_{\bar{j}k}}^{-1}\bigr),&& \CA\bigl(T_{\delta_i},T_{\talpha_{k}}T_{\alpha_k}^{-1}\bigr),
 \end{align*}
 where $i,j,k$ is a cyclic permutation of~$1,2,3$.
 \smallskip

 \textnormal{(2)} The following relation holds in $H_2(\I_M)$:
 \begin{align}
  &(1-t_2)(1-t_3)(t_2-t_3)(1-t_2t_3)\CA\bigl(T_{\delta_1},T_{\gamma_{2\bar{3}}}T_{\tgamma_{\bar{2}3}}^{-1}\bigr)\notag \\
  &{}\quad +(1-t_3)(1-t_1)(t_3-t_1)(1-t_3t_1)\CA\bigl(T_{\delta_2},T_{\gamma_{3\bar{1}}}T_{\tgamma_{\bar{3}1}}^{-1}\bigr)\notag\\
  &{}\quad +(1-t_1)(1-t_2)(t_1-t_2)(1-t_1t_2)\CA\bigl(T_{\delta_3},T_{\gamma_{1\bar{2}}}T_{\tgamma_{\bar{1}2}}^{-1}\bigr)=0.\label{eq_thm3_relA2}
 \end{align}

 \textnormal{(3)} The following two homology classes~$\sigma^M_1$ and~$\sigma^M_2$ lie in the $\Lambda$-submodule of~$H_2(\I_M)$ generated by the $6$ abelian cycles $\CA\bigl(T_{\delta_i},T_{\gamma_{j\bar{k}}}T_{\tgamma_{\bar{j}k}}^{-1}\bigr)$ and $\CA\bigl(T_{\delta_i},T_{\talpha_{k}}T_{\alpha_k}^{-1}\bigr)$,
 where $i,j,k$ is a cyclic permutation of~$1,2,3$:
 \begin{align*}
  \sigma^M_1&=(1+t_1+t_2+t_3-t_1t_2-t_1t_3-t_2t_3-t_1t_2t_3)\\
  &\hspace{5cm}\bigl(\CA(T_{\delta_1},T_{\delta_2})+\CA(T_{\delta_2},T_{\delta_3})+\CA(T_{\delta_3},T_{\delta_1})\bigr),\\
  \sigma^M_2&=\bigl(t_1-t_1^{-1}+t_2-t_2^{-1}+t_3-t_3^{-1}\bigr)\bigl(\CA(T_{\delta_1},T_{\delta_2})+\CA(T_{\delta_2},T_{\delta_3})+\CA(T_{\delta_3},T_{\delta_1})\bigr).
 \end{align*}
\end{theorem}

\begin{remark}
 In contrast to the group~$H_2(\I)$, in the group~$H_2(\I_M)$ it is generally not true that the abelian cycles~$\CA(T_{\delta_i},T_{\delta_j})$ are elements of order~$2$, nor that these abelian cycles are equal to one another. Hence the formulas for the classes~$\sigma_l^M$ look more complicated than those for the classes~$\sigma_l$ in Theorem~\ref{thm_long_rel}, and the signs in them matter.
\end{remark}

Let us now derive Theorem~\ref{thm_long_rel} from Theorem~\ref{theorem_3}.

\begin{proof}[Proof of Theorem~\ref{thm_long_rel}]
Set $U_i=\langle a_i,b_i\rangle$ for $i=1,2,3$. Let $i,j,k$ be a cyclic permutation of~$1,2,3$. Then the natural map $H_2(\I_M)\to H_2(\I)$ sends the abelian cycles~$\CA(T_{\delta_i},T_{\delta_j})$, $\CA\bigl(T_{\delta_i},T_{\gamma_{j\bar{k}}}T_{\tgamma_{\bar{j}k}}^{-1}\bigr)$, and~$\CA\bigl(T_{\delta_i},T_{\talpha_{k}}T_{\alpha_k}^{-1}\bigr)$ to $\CA(U_i,U_j)$, $\CA(U_i,a_j-a_k)$, and~$-\CA(U_i,a_k)$. (The minus sign in the last case arises because, if we orient~$\talpha_{k}$ so that its homology class becomes~$a_k$, then the pair of pants bounded by~$\delta_i$, $\talpha_{k}$, and~$\alpha_k$ lies to the right of~$\talpha_{k}$ rather than to its left.) Therefore, assertion~(1) of Theorem~\ref{thm_long_rel} follows immediately from assertion~(2) of Theorem~\ref{theorem_3}. Moreover, since $2\CA(U_i,U_j)=0$ and $\CA(U_1,U_2)=\CA(U_2,U_3)=\CA(U_3,U_1)$ in~$H_2(\I)$, we see that the map $H_2(\I_M)\to H_2(\I)$ sends~$\sigma_l^M$ to~$\sigma_l$ for $l=1,2$. Thus assertion~(2) of Theorem~\ref{thm_long_rel} follows immediately from assertion~(3) of Theorem~\ref{theorem_3}.
\end{proof}

The complement $\Sigma_{0,6}=\Sigma\setminus M$ is a $6$-punctured sphere. It is convenient to shrink the holes in~$\bS$ around the vertices of the octahedron~$123\bar{1}\bar{2}\bar{3}$ to zero and to regard~$\Sigma_{0,6}$ as the standard sphere~$\bS$ punctured at the vertices of this octahedron.

Since $M$ contains no bounding pairs, we have an injection
$$
\I_{M}\hookrightarrow \PMod(\Sigma\setminus M)=\PMod(\Sigma_{0,6})
$$
(cf.~\cite[Section~6.2]{BBM10}).
Let $\CS<\PMod(\Sigma_{0,6})$ be the image of this injection. It is convenient for us to identify $\I_M$ with~$\CS$ and in fact to work further with~$\CS$ rather than with~$\I_M$. The only difference, however, is that the BP-maps $T_{\talpha_{k}}T_{\alpha_k}^{-1}\in \I_M$ turn into the usual Dehn twists $T_{\talpha_{k}}\in \CS$.

This part of the paper (Sections~\ref{section_RCC}--\ref{section_end_thm3}) is entirely devoted to the proof of Theorem~\ref{theorem_3}. It is organized as follows. In Section~\ref{section_RCC} we recall the construction of the complex of relative cycles~$\CB_{0,2g}$ due to Spiridonov~\cite{Spi23,Spi24}. In Section~\ref{section_CG}, we introduce the spectral sequence
\begin{equation}\label{eq_SS_brief}
E^1_{p,q}=\bigoplus_{Q\in\mathfrak{Q}_p}H_q(\CS_Q) \ \Longrightarrow \ H_{p+q}(\CS)
\end{equation}
for the action of~$\CS$ on~$\CB_{0,6}$ and the corresponding filtration
$$
0\subseteq\CF_{0,2}\subseteq\CF_{1,1}\subseteq\CF_{2,0}=H_2(\CS).
$$
Moreover, we extend the action of~$\CS$ on~$\CB_{0,6}$ to a larger group~$\CG\rhd\CS$, so that \eqref{eq_SS_brief} becomes a spectral sequence of $\Z[\Gamma]$-modules, where $\Gamma=\CG/\CS$. In Section~\ref{section_curves}, we define several useful simple closed curves on~$\Sigma_{0,6}$, and esteblish a lot of lantern relations among the corresponding Dehn twists.

Section~\ref{section_representatives} is central to this part of the paper. In it, we explicitly find representatives of all $\CG$-orbits of cells of the complex~$\CB_{0,6}$ and compute their stabilizers in the groups~$\CS$ and~$\CG$. This result allows us to write down explicitly generators and defining relations for the $\Z[\Gamma]$-modules~$E^1_{p,q}$. This is done for the module~$E^1_{0,2}$ in Section~\ref{section_E102}, for the modules~$E^1_{p,1}$ in Section~\ref{section_first_row}, and for~$E^1_{p,0}$ in Section~\ref{section_zero_row}. In Sections~\ref{section_first_row} and~\ref{section_zero_row}, we also compute the differential~$d^1$ in the first and zeroth rows of the spectral sequence, respectively. These results are already sufficient to prove that the subgroup $\CF_{1,1}\subseteq H_2(\CS)=H_2(\I_M)$ is generated as a $\Lambda$-module by the $9$ abelian cycles listed in part~(1) of  Theorem~\ref{theorem_3}. The next step is the proof that the group~$E^3_{2,0}=E^{\infty}_{2,0}$ is trivial; this implies that $\CF_{1,1}$ is the whole group~$H_2(\I_M)$, which completes the proof of part~(1) of Theorem~\ref{theorem_3}. However, proving the triviality of the group~$E^3_{2,0}$ is technically very difficult. It occupies Sections~\mbox{\ref{section_E220}--\ref{section_triv_E320}}. Namely, in Section~\ref{section_E220} we compute the group~$E^2_{2,0}$, in Section~\ref{section_d2rho} we compute the differential $d^2\colon E^2_{2,0}\to E^2_{0,1}$, and finally in Section~\ref{section_triv_E320} we show that this differential is injective and thus~$E^3_{2,0}=0$.

In Sections~\ref{section_E230}--\ref{section_end_thm3} we prove parts~(2) and~(3) of Theorem~\ref{theorem_3}. These proofs are based on a partial computation of the differentials $d^2\colon E^2_{3,0}\to E^2_{1,1}$ and $d^3\colon E^3_{3,0}\to E^3_{0,2}$. In Section~\ref{section_E230} we present three elements~$\lambda_1$, $\lambda_2$, and~$\lambda_3$ lying in the kernel of the differential $d^1\colon E^1_{3,0}\to E^1_{2,0}$, i.e., in the group~$E^2_{3,0}$. In Section~\ref{section_d3} we show that $d^2\lambda_3=0$ and compute the element~$d^3\lambda_3$; this element, which thus lies in the kernel of the map $E^3_{0,2}\to E^{\infty}_{0,2}=\CF_{0,2}$, gives exactly relation~\eqref{eq_thm3_relA2} from part~(2) of Theorem~\ref{theorem_3}. The differential~$d^3$ is computed directly in the bicomplex $B_{*,*}=C_*(\CB_{0,6})\otimes_{\CS}\CR_*$ using Lemma~\ref{lem_dr}.

In Section~\ref{section_end_thm3} we similarly extract relations in~$H_2(\CS)=H_2(\I_M)$ from the elements~$\lambda_1$ and~$\lambda_2$. The difficulty, however, is that these elements do not lie in the kernel of~$d^2$. Nevertheless, a computation in the bicomplex $B_{*,*}$ similar to the computation of~$d^3\lambda_3$ carries information not only about the differentials of the spectral sequence~$E^*_{*,*}$ but also about how the associated factors~$E^{\infty}_{0,2}$ and~$E^{\infty}_{1,1}$ are assembled into the whole group~$H_2(\I_M)$. Therefore such a computation can actually be used to obtain precise relations in~$H_2(\I_M)$. We will carry out the corresponding computation not in full, but only modulo a certain special subcomplex $L\subseteq\CB_{0,6}$. This will not give us explicit relations in~$H_2(\I_M)$ but will yield exactly part~(3) of Theorem~\ref{theorem_3}.

\section{Complex of relative cycles~$\CB_{0,2g}$}\label{section_RCC}

Following Spiridonov~\cite{Spi23,Spi24}, we consider the following \textit{complex of relative cycles}~$\CB_{0,2g}$. Let $\Sigma_{0,2g}$ be a sphere with $2g$ punctures. We divide the punctures into two parts, $\mathbf{p}=\{p_1,\ldots,p_g\}$ and $\mathbf{q}=\{q_1,\ldots,q_g\}$, each consisting of $g$ punctures. By an \textit{arc} we mean an embedded oriented curve with endpoints at punctures. A \textit{multiarc} is a finite union of arcs that are not isotopic to each other and have pairwise disjoint interiors (common endpoints are allowed). We will usually not distinguish between an arc (or a multiarc) and its isotopy class.

Let $\fD$ be the set of isotopy classes of arcs on~$\Sigma_{0,2g}$ that start at some point from~$\mathbf{p}$ and finish at some point from~$\mathbf{q}$. The \textit{complex of relative cycles}~$\CB_{0,2g}$ is the subset of~$\R_{\ge0}^{\fD}$ consisting of all points (finite formal sums) $\bbeta=\sum_{i=1}^kn_i\beta_i$, where $\beta_i\in\fD$ and $n_i\in\R_{\ge0}$, that satisfy the following conditions:
\begin{enumerate}
 \item $\partial \left(\sum_{i=1}^k n_i \beta_i\right)=\sum_{i=1}^g(q_i-p_i)$, and
 \item the isotopy classes $\beta_1,\ldots,\beta_k$ contain pairwise disjoint representatives.
\end{enumerate}
The set~$\CB_{0,2g}$ has the following structure of a regular CW complex. By definition, the \textit{support} of a relative cycle $\bbeta\in \CB_{0,2g}$ is the multiarc obtained as the union of all arcs~$\beta_i$ that enter~$\bbeta$ with non-zero coefficients. Let $\CL$ be the set of all multiarcs that occur as the supports of relative cycles $\bbeta\in \CB_{0,2g}$. Then, for each $L\in\CL$, the set of all relative cycles $\bbeta\in \CB_{0,2g}$ with supports contained in~$L$ is a (finite-dimensional) convex polytope $Q_L\subseteq\R_{\ge0}^{\fD}$. Moreover, these convex polytopes meet face-to-face and form a regular CW-decomposition of~$\CB_{0,2g}$, see~\cite[Section~5.1]{Spi24}. We will need the following propositions.

\begin{propos}[\protect{\cite[Lemma~5.6]{Spi24}}]\label{propos_dim_rel}
We have
$$
\dim(Q_L)=F-1,
$$
where $F$ is the number of connected components of~$\Sigma_{0,2g}\setminus L$.
\end{propos}

\begin{propos}[\protect{\cite[Proposition~5.4]{Spi24}}]\label{propos_rel_wor}
 The group $\PMod(\Sigma_{0,2g})$ acts on~$\CB_{0,2g}$ without rotations.
\end{propos}

The vertices of~$\CB_{0,2g}$ are exactly \textit{relative basic cycles}, i.e., formal sums $\bbeta=\beta_1+\cdots+\beta_g$ such that $\beta_i\in\fD$, $\partial\bbeta=\sum_{i=1}^g(q_i-p_i)$, and $\beta_1,\ldots,\beta_g$ contain pairwise disjoint representatives. Moreover, the set~$\CL$ indexing the cells of~$\CB_{0,2g}$ admits the following characterization.

\begin{propos}[\protect{\cite[Section~5.1]{Spi24}}]\label{propos_rel_L_char}
 Let~$L$ be a multiarc with all components starting at a point from~$\mathbf{p}$ and finishing at a point from~$\mathbf{q}$. Then $L$ belongs to~$\CL$ if and only if, for each component~$\beta$ of~$L$, there exists a relative basic cycle whose support is contained in~$L$ and contains~$\beta$.
\end{propos}

The main property of the complex of relative cycles is the following theorem due to Spiridonov.

\begin{theorem}[\protect{\cite[Theorem~5.3]{Spi24}}]
 Let $g\ge 1$. Then $\CB_{0,2g}$ is contractible.
\end{theorem}

Finally, note that the maximal cells of~$\CB_{0,2g}$ correspond to bipartite graphs on the sphere all of whose faces are quadrangles. A standard computation of the Euler characteristic yields that such a graph has $2g-2$ faces. Therefore $\dim(\CB_{0,2g})=2g-3$. In particular, for $g=3$, which is the  only case of interest to us in this paper, we have $\dim(\CB_{0,6})=3$.

\section{Spectral sequence and its symmetries}\label{section_CG}

Consider the surface $\Sigma_{0,6}=\Sigma\setminus M$, which is identified with the sphere~$\bS$ punctured at the vertices of the regular octahedron~$123\bar{1}\bar{2}\bar{3}$ as in Section~\ref{section_IM}. Set $\mathbf{p}=\{1,2,3\}$ and $\mathbf{q}=\{\bar{1},\bar{2},\bar{3}\}$, and consider the corresponding complex  of relative cycles~$\CB_{0,6}$. By Proposition~\ref{propos_rel_wor}, the group~$\CS$ acts on~$\CB_{0,6}$ without rotations.
Throughout this part of the paper, we study the spectral sequence~\eqref{eq_SpecSeq} for the action of~$\CS$ on~$\CB_{0,6}$, i.e., the spectral sequence
\begin{equation}\label{eq_SSS}
E^1_{p,q}=\bigoplus_{Q\in\mathfrak{Q}_p}H_q(\CS_Q) \ \Longrightarrow \ H_{p+q}(\CS),
\end{equation}
where $\mathfrak{Q}_p$ is a set of representatives of $\CS$-orbits of $p$-cells of~$\CB_{0,6}$. We write
$$
\iota_Q\colon  H_q(\CS_Q)\hookrightarrow E^1_{p,q}
$$
for the inclusion of the summand. Let
$$
0\subseteq\CF_{0,2}\subseteq\CF_{1,1}\subseteq\CF_{2,0}=H_2(\CS)
$$
be the filtration associated with the spectral sequence~\eqref{eq_SSS}, i.e., such that $$\CF_{p,q}/\CF_{p-1,q+1}=E^{\infty}_{p,q}.$$

A standard capping short exact sequence (cf.~\cite[Section~3.6]{FaMa12}) implies that
the group~$\CX_M$ defined in Section~\ref{section_IM} surjects onto~$\PMod(\Sigma_{0,6})$, and the kernel of this surjection is the free abelian group with basis $T_{\alpha_1}, T_{\alpha_2},T_{\alpha_3}$. Therefore, the short exact sequence~\eqref{eq_ses_IXW} turns into the short exact sequence
\begin{equation}\label{eq_semidir_PMod}
   1\longrightarrow \CS \longrightarrow \PMod(\Sigma_{0,6}) \stackrel{\kappa}{\longrightarrow} V \longrightarrow 1,
 \end{equation}
where $V$ is the free abelian group of rank~$3$ with basis~$t_1=t_{a_2+a_3}$, $t_2=t_{a_3+a_1}$, $t_3=t_{a_1+a_2}$. Moreover, Lemma~\ref{lem_transvec3} gives the following formula for~$\kappa$: if $\gamma$ is a simple closed curve on~$\Sigma_{0,6}=\Sigma\setminus M$ with $[\gamma]=n_1a_1+n_2a_2+n_3a_3$, then
 \begin{equation}\label{eq_kappa}
  \kappa(T_{\gamma})=t_1^{n_2n_3}t_2^{n_3n_1}t_3^{n_1n_2}.
 \end{equation}

To perform computations in the spectral sequence~\eqref{eq_SSS}, it is convenient to find a larger group $\CG> \CS$ such that $\CS$ is normal in~$\CG$ and the action of~$\CS$ on~$\CB_{0,6}$ extends to an action of~$\CG$. Then, by Fact~\ref{fact_Gamma}, \eqref{eq_SSS} would become a spectral sequence of $\Z[\CG/\CS]$-modules. We could take for~$\CG$ the group~$\PMod(\Sigma_{0,6})$. However, it is more convenient to take even larger group~$\CG$.

We define~$\CG$ as follows. First, let $\Mod^{\pm}(\Sigma_{0,6})$ be the \textit{extended mapping class group} of~$\Sigma_{0,6}$, that is, the group consisting  of the isotopy classes of both orientation-preserving and orientation-reversing homeomorphisms of~$\Sigma_{0,6}$. Second, let $\tG<\Mod^{\pm}(\Sigma_{0,6})$ be the subgroup consisting of all mapping classes~$f$ that satisfy the following condition:
\begin{enumerate}
 \item[($*$)] there is a permutation $\sigma\in S_3$ such that either $f(i)=\sigma(i)$ and $f\bigl(\bar{i}\bigr)=\overline{\sigma(i)}$ for all~$i$ or $f(i)=\overline{\sigma(i)}$ and $f\bigl(\bar{i}\bigr)=\sigma(i)$ for all~$i$.
\end{enumerate}
We have three homomorphisms $e_1,e_2,e_3\colon \tG\to\{\pm 1\}$ given by
\begin{itemize}
 \item $e_1(f)=1$ if $f$ preserves the orientation, and $e_1(f)=-1$ if $f$ reverses the orientation;
 \item $e_2(f)=1$ if $f$ stabilizes each of the two sets $\mathbf{p}=\{1,2,3\}$ and $\mathbf{q}=\{\bar{1},\bar{2},\bar{3}\}$, and  $e_2(f)=-1$ if $f$ swaps the sets~$\mathbf{p}$ and~$\mathbf{q}$;
 \item $e_3(f)=\sgn(\sigma)$, where $\sigma$ is the permutation from condition~($*$).
\end{itemize}
We define $\G$ to be the kernel of the product homomorphism
$$e_1e_2e_3\colon\tG\to\{\pm1\}.$$

 Let $G$ be the group of all isometries of the octahedron~$123\bar{1}\bar{2}\bar{3}$ that either preserve each of the sets~$\mathbf{p}$ and~$\mathbf{q}$ or swap~$\mathbf{p}$ and~$\mathbf{q}$; then $G$ is isomorphic to~$S_3\times C_2$, where $S_3$ is the symmetric group of degree~$3$ and $C_2$ is the cyclic group of order~$2$. Hence $|G|=12$. We consider the following elements of~$G$:
\begin{itemize}
 \item let $\nu$ be the antipodal involution on the sphere~$\bS$; then $\nu(i)=\bar{i}$ and $\nu\bigl(\bar{i}\bigr)=i$ for all~$i$;
 \item let $\omega$ be the rotation by~$2\pi/3$ around the axis through the centers of the spherical triangles~$123$ and~$\bar{1}\bar{2}\bar{3}$, where the direction of rotation is chosen so that $\omega(i)=i+1$ and~$\omega\bigl(\bar{i}\bigr)=\overline{i+1}$ for all~$i$ (the sums of indices are taken modulo~$3$);
 \item for each permutation $i,j,k$ of~$1,2,3$, let $r_i$ be the rotation by~$\pi$ around the axis~$\ell_i$ through the midpoints of the segments~$j\bar{k}$ and~$\bar{j}k$, and let $s_i$ be the orthogonal reflection in the plane~$\ell_i^{\bot}$; then $r_i$ swaps the punctures by $i\leftrightarrow \bar{i}$, $j\leftrightarrow \bar{k}$, and~$k\leftrightarrow \bar{j}$, and $s_i$ fixes each of the punctures~$i$ and~$\bar{i}$ and swaps $j\leftrightarrow k$ and $\bar{j}\leftrightarrow \bar{k}$.
\end{itemize}
These elements of~$G$ satisfy the following relations:
\begin{gather*}
 \nu^2=1,\qquad \omega^3=1,\qquad r_i^2=1,\qquad s_i^2=1,\\
 s_i=r_i\nu=\nu r_i,\qquad \omega\nu=\nu\omega,\\
 \omega r_i=r_{i-1},\qquad r_i\omega=r_{i+1},\qquad r_ir_{i+1}=\omega,\qquad r_{i+1}r_i=\omega^2.
 \end{gather*}
Hereafter, the sums of indices are taken modulo~$3$. Moreover, $1$, $\omega$, $\omega^2$, $r_1$, $r_2$, $r_3$, $\nu$, $\omega\nu$, $\omega^2\nu$, $r_1\nu$, $r_2\nu$, and~$r_3\nu$ are precisely the $12$ elements of~$G$. Throughout this part of the paper, we will repeatedly use these relations without further comment.

As already mentioned in Section~\ref{section_IM}, in figures we always depict the sphere~$\bS$ as the plane compactified by a single point~$\infty$, corresponding to the center of the triangle~$\bar{1}\bar{2}\bar{3}$, as in Fig.~\ref{fig_sphere}. Under this representation, the antipodal involution~$\nu$ becomes the composition of the central symmetry about the center~$o$ of the triangle~$123$ and the inversion in a certain circle centered at~$o$ and enclosing the punctures~$1$, $2$, and~$3$, the map~$\omega$ becomes the clockwise rotation by~$2\pi/3$ about~$o$, and $s_i$ becomes the reflection in the line~$i\bar{i}$.

Combining the definition of~$\G$ and the short exact sequence~\eqref{eq_semidir_PMod}, we obtain the following proposition.

\begin{propos}\label{propose_semidir_G}
 We have a semidirect product decomposition:
 \begin{equation}\label{eq_semi_dir}
 \G=\PMod(\Sigma_{0,6})\rtimes G.
 \end{equation}
 Moreover, we have a short exact sequence
 \begin{equation*}
   1\longrightarrow  \CS \longrightarrow \G \stackrel{\kappa}{\longrightarrow} V\rtimes G \longrightarrow 1,
 \end{equation*}
where the action of\/~$G$ on~$V=\langle t_1,t_2,t_3\rangle$ is given by
$$
\nu t_i\nu=t_i^{-1},\qquad \omega t_i\omega^{-1}=t_{i+1}, \qquad r_i t_i r_i=t_i,\qquad r_i t_{i+1} r_i=t_{i-1},\qquad r_i t_{i-1} r_i=t_{i+1},
$$
the restriction of~$\kappa$ to~$\PMod(\Sigma_{0,6})$ is given by~\eqref{eq_kappa}, and $\kappa(g)=g$ for all $g\in G$. In particular, $\CS$ is a normal subgroup of~$\G$.
\end{propos}

\begin{remark}
 The element~$\nu t_i\nu$ is equal to~$t_i^{-1}$ rather than~$t_i$ because $\nu$ reverses the sign of the intersection form.
\end{remark}

The action of~$\tG$ on~$\CB_{0,6}$ is defined as follows. Consider a mapping class $f\in\tG$. If $f(\mathbf{p})=\mathbf{p}$ and~$f(\mathbf{q})=\mathbf{q}$, then $f$ acts on~$\CB_{0,6}$ in a naive way, that is, by sending each relative cycle $\sum n_i\beta_i$ to $\sum n_if(\beta_i)$. If  $f(\mathbf{p})=\mathbf{q}$ and~$f(\mathbf{q})=\mathbf{p}$, then we need in addition to reverse the orientations of all arcs, that is, $f$ sends $\sum n_i\beta_i$ to $\sum n_i\widehat{f(\beta_i)}$, where $\widehat{f(\beta_i)}$ is the arc~$f(\beta_i)$ with the orientation reversed.

It is more convenient to consider the subgroup~$\G$ rather than the whole group~$\tG$ because of the existence of the semidirect product decomposition~\eqref{eq_semi_dir}.
We set $$\Gamma=V\rtimes G\cong \G/\CS.$$ Then~\eqref{eq_SSS} is a spectral sequence of $\Z[\Gamma]$-modules.

\section{More curves on~$\Sigma_{0,6}$ and lantern relations}\label{section_curves}

In Section~\ref{section_IM} we defined the simple closed curves~$\delta_i$ and~$\talpha_i$ for $i=1,2,3$, $\gamma_{pq}$ for distinct $p,q\in \bigl\{1,2,3,\bar{1},\bar{2},\bar{3}\bigr\}$ with $\{p,q\}\ne \bigl\{i,\bar{i}\bigr\}$, and $\tgamma_{{}\,\bar{i}j}$ for $j=i+1\pmod 3$; see Figs.~\ref{fig_deltagamma}--\ref{fig_AC2}. To simplify notation, we set
\begin{equation*}
 T_{pq}=T_{\gamma_{pq}}.
\end{equation*}
The curves~$\delta_i$ are separating on~$\Sigma$, and the curves~$\talpha_i$ are homologous to components of~$M$. Therefore the twists~$T_{\delta_i}$ and $T_{\talpha_i}$ belong to~$\CS$, and  hence to the kernel of~$\kappa$. If $i,j,k$ is a cyclic permutation of~$1,2,3$, then the curves~$\gamma_{ij}$ and~$\gamma_{{}\,\bar{i}\bar{j}}$ (with appropriate orientations) have homology class~$a_i+a_j$, and the curves~$\gamma_{i\bar{j}}$, $\gamma_{{}\,\bar{i}j}$, and~$\tgamma_{{}\,\bar{i}j}$ have homology class~$a_i-a_j$. Therefore, by~\eqref{eq_kappa} we  have:
$$
\kappa(T_{ij})=\kappa(T_{\bar{i}\bar{j}})=t_k,\qquad \kappa(T_{i\bar{j}})=\kappa(T_{\bar{i}j})=\kappa\bigl(T_{\tgamma_{{}\,\bar{i}j}}\bigr)=t_k^{-1}.
$$

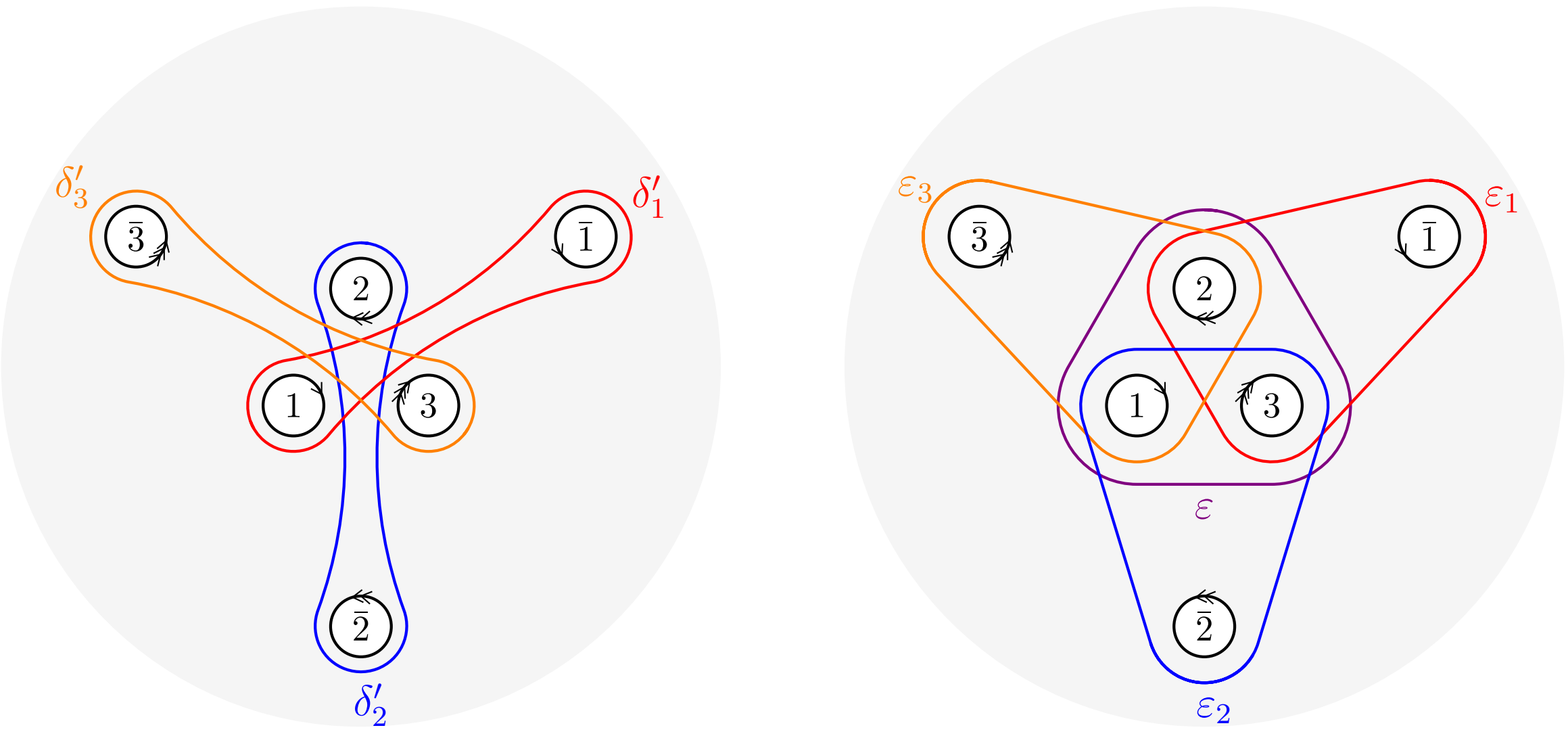
\begin{figure}[!ht]
 \begin{center}
\begin{tikzpicture}[scale=1.1,
         mid arrow/.style={
            postaction={
                decorate,
                decoration={
                    markings,
                    mark=at position 0.55 with {
                        \arrow[#1]{stealth}
                    }
                }
            }
        },
        mid arrow/.default=black,
    ]

 %\definecolor{myblue}{rgb}{0, 0, 0.7}
 \unitlength=1cm

 \begin{scope}
 \fill [black!4] (0,0) circle (3.2);
 \thicklines

 \begin{scope}[white]
 \fill  (-.6,-.346) circle (.27);
\fill  (0,.693) circle (.27);
\fill  (.6,-.346) circle (.27);
\fill  (2,1.155) circle (.27);
\fill  (0,-2.309) circle (.27);
\fill  (-2,1.155) circle (.27);
 \end{scope}

\draw [thick] (-.6,-.346) circle (.27);
\draw [thick] (0,.693) circle (.27);
\draw [thick] (.6,-.346) circle (.27);
\draw [thick] (2,1.155) circle (.27);
\draw [thick] (0,-2.309) circle (.27);
\draw [thick] (-2,1.155) circle (.27);

 \draw (0,-2.309) node {\footnotesize $\bar{2}$};
 \draw (-2,1.155) node {\footnotesize $\bar{3}$};
 \draw (2,1.155) node  {\footnotesize $\bar{1}$};
 \draw (0,.693) node {\footnotesize $2$};
 \draw (-.6,-.346) node  {\footnotesize $1$};
 \draw (.6,-.346)node {\footnotesize $3$};

 %\draw (0,.693) arc(90:270:1.501);
 %\draw (-.6,-.346) arc(-150:30:1.501);
 %\draw (.6,-.346) arc(-30:150:1.501);

 \begin{scope}[blue]
 \draw [thick] (0,1.1) arc (90:-20:.407) arc(160:200:3.976) arc (20:-200:.407) node [pos=.5,below=-1pt,xshift=1mm] {$\delta_2'$} arc(-20:20:3.976) arc (200:80:.407);
 \end{scope}

\begin{scope}[red,rotate=120]
 \draw [thick] (0,1.1) arc (90:-20:.407) arc(160:200:3.976) arc (20:-200:.407) node [pos=.5,above right=-2.5pt,yshift=-1mm] {$\delta_1'$} arc(-20:20:3.976) arc (200:80:.407);
 \end{scope}
 \begin{scope}[orange,rotate=-120]
 \draw [thick] (0,1.1) arc (90:-20:.407) arc(160:200:3.976) arc (20:-200:.407) node [pos=.5,above left=-2.5pt] {$\delta_3'$} arc(-20:20:3.976) arc (200:80:.407);
 \end{scope}

 %1
 \begin{scope}[xshift=-.6cm, yshift=-.346cm]
 \begin{scope}[rotate=-68]
 \draw (0,.27)-- (-.1,.31);
\draw (0,.27)-- (-.08,.2);
\end{scope}
\end{scope}

%2
\begin{scope}[yshift=.693cm]
 \begin{scope}[rotate=164]
 \draw (0,.27)-- (-.1,.31);
\draw (0,.27)-- (-.08,.2);
\end{scope}
\begin{scope}[rotate=180]
 \draw (0,.27)-- (-.1,.31);
\draw (0,.27)-- (-.08,.2);
\end{scope}
\end{scope}

%3
\begin{scope}[xshift=.6cm, yshift=-.346cm]
 \begin{scope}[rotate=66]
 \draw (0,.27)-- (-.1,.31);
\draw (0,.27)-- (-.08,.2);
\end{scope}
 \begin{scope}[rotate=52]
 \draw(0,.27)-- (-.1,.31);
\draw (0,.27)-- (-.08,.2);
\end{scope}
\begin{scope}[rotate=38]
 \draw (0,.27)-- (-.1,.31);
\draw (0,.27)-- (-.08,.2);
\end{scope}
\end{scope}

%bar1
 \begin{scope}[xshift=2cm, yshift=1.155cm]
 \begin{scope}[rotate=128]
 \draw (0,.27)-- (.1,.31);
\draw (0,.27)-- (.08,.2);
\end{scope}
\end{scope}

%bar2
\begin{scope}[yshift=-2.309cm]
 \begin{scope}[rotate=16]
 \draw (0,.27)-- (.1,.31);
\draw (0,.27)-- (.08,.2);
\end{scope}
\begin{scope}
 \draw (0,.27)-- (.1,.31);
\draw (0,.27)-- (.08,.2);
\end{scope}
\end{scope}

%bar3
\begin{scope}[xshift=-2cm, yshift=1.155cm]
 \begin{scope}[rotate=-98]
 \draw (0,.27)-- (.1,.31);
\draw (0,.27)-- (.08,.2);
\end{scope}
\begin{scope}[rotate=-112]
 \draw (0,.27)-- (.1,.31);
\draw (0,.27)-- (.08,.2);
\end{scope}
\begin{scope}[rotate=-126]
 \draw (0,.27)-- (.1,.31);
\draw (0,.27)-- (.08,.2);
\end{scope}
\end{scope}

\end{scope}

%%%%%%%%%%%
%%%%%%%%%%%

 \begin{scope}[xshift=7.5cm]
 \fill [black!4] (0,0) circle (3.2);
 \thicklines

 \begin{scope}[white]
 \fill  (-.6,-.346) circle (.27);
\fill  (0,.693) circle (.27);
\fill  (.6,-.346) circle (.27);
\fill  (2,1.155) circle (.27);
\fill  (0,-2.309) circle (.27);
\fill  (-2,1.155) circle (.27);
 \end{scope}

\draw [thick] (-.6,-.346) circle (.27);
\draw [thick] (0,.693) circle (.27);
\draw [thick] (.6,-.346) circle (.27);
\draw [thick] (2,1.155) circle (.27);
\draw [thick] (0,-2.309) circle (.27);
\draw [thick] (-2,1.155) circle (.27);

 \draw (0,-2.309) node {\footnotesize $\bar{2}$};
 \draw (-2,1.155) node {\footnotesize $\bar{3}$};
 \draw (2,1.155) node  {\footnotesize $\bar{1}$};
 \draw (0,.693) node {\footnotesize $2$};
 \draw (-.6,-.346) node  {\footnotesize $1$};
 \draw (.6,-.346)node {\footnotesize $3$};

 \begin{scope}
 \draw[thick,violet] (0,1.393) arc (90:30:.7) -- +(.6,-1.039) arc (30:-90:.7) -- +(-1.2,0) node [pos=.5,below] {$\varepsilon$} arc (-90:-210:.7) -- +(.6,1.039) arc (150:80:.7);
 \end{scope}

 \begin{scope}[red]
 \draw[thick] (2,1.655) arc (90:-43:.5) -- +(-1.4,-1.501) arc (-43:-150:.5) -- +(-.6,1.039) arc (210:103:.5) -- +(2,.462) arc (103:-43:.5) node [pos=.5,above right=-2.5pt,yshift=-1mm] {$\varepsilon_1$} ;
\end{scope}

 \begin{scope}[orange,rotate=120]
 \draw[thick] (2,1.655) arc (90:-43:.5) -- +(-1.4,-1.501) arc (-43:-150:.5) -- +(-.6,1.039) arc (210:103:.5) -- +(2,.462) arc (103:-43:.5) node [pos=.5,above left=-2.5pt,xshift=1mm] {$\varepsilon_3$} ;
\end{scope}

\begin{scope}[blue,rotate=-120]
 \draw[thick] (2,1.655) arc (90:-43:.5) -- +(-1.4,-1.501) arc (-43:-150:.5) -- +(-.6,1.039) arc (210:103:.5) -- +(2,.462) arc (103:-43:.5) node [pos=.5,below,xshift=1mm] {$\varepsilon_2$} ;
\end{scope}

 %\draw (0,.693) arc(90:270:1.501);
 %\draw (-.6,-.346) arc(-150:30:1.501);
 %\draw (.6,-.346) arc(-30:150:1.501);

%  \begin{scope}[blue,rotate=120]
%  \draw [thick] (0,1.1) arc (90:-90:.407)-- +(-.5,0) arc (90:150:.794) arc(150:210:1.394)  arc (210:270:.794) -- +(.5,0) arc (90:-90:.407)  arc (270:90:1.908) node [pos=.5, above] {$\delta_1$};
%  \end{scope}
%
%  \draw[thick,red] (2,1.705) arc (90:30:.55) arc (30:-150:2.051) node [pos=.55, below] {$\talpha_3$}  arc (210:75:.55) arc(-105:-15:1.57) arc (165:80:.55);

 %1
 \begin{scope}[xshift=-.6cm, yshift=-.346cm]
 \begin{scope}[rotate=-68]
 \draw (0,.27)-- (-.1,.31);
\draw (0,.27)-- (-.08,.2);
\end{scope}
\end{scope}

%2
\begin{scope}[yshift=.693cm]
 \begin{scope}[rotate=164]
 \draw (0,.27)-- (-.1,.31);
\draw (0,.27)-- (-.08,.2);
\end{scope}
\begin{scope}[rotate=180]
 \draw (0,.27)-- (-.1,.31);
\draw (0,.27)-- (-.08,.2);
\end{scope}
\end{scope}

%3
\begin{scope}[xshift=.6cm, yshift=-.346cm]
 \begin{scope}[rotate=66]
 \draw (0,.27)-- (-.1,.31);
\draw (0,.27)-- (-.08,.2);
\end{scope}
 \begin{scope}[rotate=52]
 \draw(0,.27)-- (-.1,.31);
\draw (0,.27)-- (-.08,.2);
\end{scope}
\begin{scope}[rotate=38]
 \draw (0,.27)-- (-.1,.31);
\draw (0,.27)-- (-.08,.2);
\end{scope}
\end{scope}

%bar1
 \begin{scope}[xshift=2cm, yshift=1.155cm]
 \begin{scope}[rotate=128]
 \draw (0,.27)-- (.1,.31);
\draw (0,.27)-- (.08,.2);
\end{scope}
\end{scope}

%bar2
\begin{scope}[yshift=-2.309cm]
 \begin{scope}[rotate=16]
 \draw (0,.27)-- (.1,.31);
\draw (0,.27)-- (.08,.2);
\end{scope}
\begin{scope}
 \draw (0,.27)-- (.1,.31);
\draw (0,.27)-- (.08,.2);
\end{scope}
\end{scope}

%bar3
\begin{scope}[xshift=-2cm, yshift=1.155cm]
 \begin{scope}[rotate=-98]
 \draw (0,.27)-- (.1,.31);
\draw (0,.27)-- (.08,.2);
\end{scope}
\begin{scope}[rotate=-112]
 \draw (0,.27)-- (.1,.31);
\draw (0,.27)-- (.08,.2);
\end{scope}
\begin{scope}[rotate=-126]
 \draw (0,.27)-- (.1,.31);
\draw (0,.27)-- (.08,.2);
\end{scope}
\end{scope}

\end{scope}

\end{tikzpicture}
\end{center}
 \caption{Simple closed curves~$\delta_1'$, $\delta_2'$, $\delta_3'$, $\varepsilon$, $\varepsilon_1$, $\varepsilon_2$, and~$\varepsilon_3$}\label{fig_deltaeps}
\end{figure}

Let us define several other useful simple closed curves on~$\Sigma_{0,6}=\Sigma\setminus M$.
First, consider the simple closed curves~$\delta_1'$, $\delta_2'$, and~$\delta_3'$ depicted in the left part of Fig.~\ref{fig_deltaeps}. Since $\delta_i'$ is separating on~$\Sigma$, we have $T_{\delta_i'}\in\CS$.

Now, let $\varepsilon$ be the simple closed curve that separates from each other the faces~$123$ and~$\bar{1}\bar{2}\bar{3}$ of the octahedron~$123\bar{1}\bar{2}\bar{3}$. Also, for a permutation~$i,j,k$ of~$1,2,3$, let~$\varepsilon_i$ be the simple closed curve that separates from each other the faces~$i\bar{j}\bar{k}$ and~$\bar{i}jk$ of the octahedron (see the right part of Fig.~\ref{fig_deltaeps}). The homology classes of the curves~$\varepsilon$ and~$\varepsilon_i$ (with suitable orientations) are~$a_1+a_2+a_3$ and~$a_i-a_j-a_k$, respectively. Hence, by~\eqref{eq_kappa} we have
$$
\kappa(T_{\varepsilon})=t_1t_2t_3,\qquad \kappa(T_{\varepsilon_i})=t_it_j^{-1}t_k^{-1}.
$$

The Dehn twists about the introduced simple closed curves satisfy a lot of relations. Let us describe some of them. We will repeatedly use that Dehn twists about disjoint simple closed curves commute. In particular, $T_{p_1q_1}T_{p_2q_2}=T_{p_2q_2}T_{p_1q_1}$ whenever $p_1$, $p_2$, $q_1$, and~$q_2$ are pairwise distinct.

\begin{figure}[!ht]
 \begin{center}
\unitlength=1mm
\tikzset{->-/.style={decoration={
     markings,
     mark=at position #1 with {\arrow{>}}},postaction={decorate}}}

  %  \draw[->-=.5] (0,0) to [bend left] (2,4);
    %\draw[->-=.8] (0,0) to [bend right] (2,4);

\definecolor{mygreen}{rgb}{0, 0.4, 0}
\definecolor{myblue}{rgb}{0, 0, 0.7}
\definecolor{mypink}{rgb}{1, 0.3, 0.8}
\definecolor{myyellow}{rgb}{1,1, 0}

\colorlet{cola}{violet}
\colorlet{cold}{red}
\colorlet{colg}{myblue}
\colorlet{colkf}{mypink}
\colorlet{colks}{violet}
\colorlet{coll}{orange}
\colorlet{colx}{mygreen}
\colorlet{coly}{orange}
\colorlet{colm}{mygreen}
\colorlet{colz}{brown}

\begin{tikzpicture}
\small
\tikzset{every path/.append style={line width=.3mm}}

\filldraw [fill=black!4, very thick] (0,0) circle (1.9); 
%\fill [color=] (0,0) circle (2);
\filldraw [fill=white, very thick] (0,1) circle (.3);
\filldraw [fill=white, very thick] (0.866,-0.5) circle (.3);
\filldraw [fill=white, very thick] (-0.866,-0.5) circle (.3);

\draw [thick] (0.1,-.5) ellipse (1.4 and .7);
\draw [thick, rotate=120] (0.1,-.5) ellipse (1.4 and .7);
\draw [thick, rotate=240] (0.1,-.5) ellipse (1.4 and .7);
\path (0.03,1) node {$b_2$};
\path (-0.836,-0.5) node {$b_1$};
\path (0.896,-0.5) node {$b_3$};
\path (0.13,-1.45) node {$c_2$};
\path (1.27,0.7) node {$c_1$};
\path (-1.27,0.7) node {$c_3$};
\path (1.56,-1.56) node {$b_4$};

\end{tikzpicture}
\end{center}
\caption{Lantern relation}\label{fig_lantern}
\end{figure}
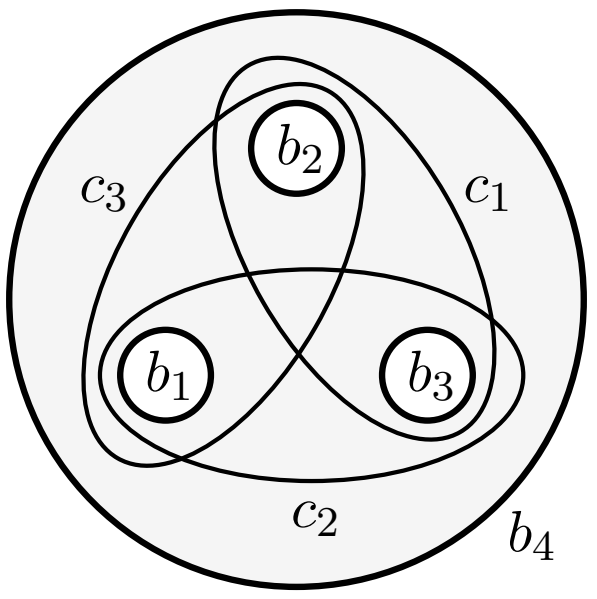

Recall the \textit{lantern relation} in the mapping class group of a genus~$0$ surface with~$4$ boundary components:
$$
T_{c_1}T_{c_2}T_{c_3}=T_{b_1}T_{b_2}T_{b_3}T_{b_4}
$$
for the curves $c_i$ and~$b_j$ shown in Fig.~\ref{fig_lantern}; see~\cite[Section~5.1]{FaMa12}. Any $7$-tuple of curves as in this figure will be referred to as a \textit{lantern configuration}.

There are several lantern configurations among  the introduced simple closed curves. First, for each face of the octahedron~$123\bar{1}\bar{2}\bar{3}$, the circuit around this face, the three circuits around the edges of this face, and the three circuits around the vertices  of this face form a lantern configuration. Since the Dehn twists about circuits around vertices  are trivial in~$\PMod(\Sigma_{0,6})$, the corresponding lantern relations are as follows:
\begin{equation}\label{eq_lantern_faces}
 \begin{aligned}
  T_{\varepsilon}&=T_{12}T_{23}T_{31}=T_{23}T_{31}T_{12}=T_{31}T_{12}T_{23}
  =T_{\bar{1}\bar{2}}T_{\bar{3}\bar{1}}T_{\bar{2}\bar{3}}=T_{\bar{3}\bar{1}}T_{\bar{2}\bar{3}}T_{\bar{1}\bar{2}}=T_{\bar{2}\bar{3}}T_{\bar{1}\bar{2}}T_{\bar{3}\bar{1}},\\
  T_{\varepsilon_1}&=T_{1\bar{2}}T_{\bar{2}\bar{3}}T_{\bar{3}1}=T_{\bar{2}\bar{3}}T_{\bar{3}1}T_{1\bar{2}}=T_{\bar{3}1}T_{1\bar{2}}T_{\bar{2}\bar{3}}
  =T_{\bar{1}2}T_{3\bar{1}}T_{23}=T_{3\bar{1}}T_{23}T_{\bar{1}2}=T_{23}T_{\bar{1}2}T_{3\bar{1}},\\
  T_{\varepsilon_2}&=T_{\bar{1}2}T_{2\bar{3}}T_{\bar{3}\bar{1}}=T_{2\bar{3}}T_{\bar{3}\bar{1}}T_{\bar{1}2}=T_{\bar{3}\bar{1}}T_{\bar{1}2}T_{2\bar{3}}
  =T_{1\bar{2}}T_{31}T_{\bar{2}3}=T_{31}T_{\bar{2}3}T_{1\bar{2}}=T_{\bar{2}3}T_{1\bar{2}}T_{31},\\
  T_{\varepsilon_3}&=T_{\bar{1}\bar{2}}T_{\bar{2}3}T_{3\bar{1}}=T_{\bar{2}3}T_{3\bar{1}}T_{\bar{1}\bar{2}}=T_{3\bar{1}}T_{\bar{1}\bar{2}}T_{\bar{2}3}
  =T_{12}T_{\bar{3}1}T_{2\bar{3}}=T_{\bar{3}1}T_{2\bar{3}}T_{12}=T_{2\bar{3}}T_{12}T_{\bar{3}1}.
 \end{aligned}
\end{equation}

Further relations are provided by the following lemma.

\begin{lem}
 Let $i,j,k$ be a cyclic permutation of~$1,2,3$. Then we have the following relations in~$\PMod(\Sigma_{0,6})$:
 \begin{gather}
  \label{eq_lantern_eps}
 T_{\delta_i}T_{\varepsilon}T_{\varepsilon_i}=T_{\varepsilon}T_{\varepsilon_i}T_{\delta_i}=T_{\varepsilon_i}T_{\delta_i}T_{\varepsilon}=T_{jk}T_{\bar{j}\bar{k}},\\
 \label{eq_lantern_eps_jk}
 T_{\delta_i'}T_{\varepsilon_k}T_{\varepsilon_j}=T_{\varepsilon_k}T_{\varepsilon_j}T_{\delta_i'}=T_{\varepsilon_j}T_{\delta_i'}T_{\varepsilon_k}=T_{j\bar{k}}T_{\bar{j}k},\\
 T_{\delta_i'}T_{k\bar{i}}T_{ki}=T_{k\bar{i}}T_{ki}T_{\delta_i'}=T_{ki}T_{\delta_i'}T_{k\bar{i}}=T_{\talpha_k},\label{eq_131}
 \\
 \label{eq_Tdeltai}
 T_{\delta_i}=T_{i\bar{j}}^{-1}T_{\bar{k}i}^{-1}T_{ij}^{-1}T_{ki}^{-1}=T_{k\bar{i}}^{-1}T_{\bar{i}j}^{-1}T_{\bar{k}\bar{i}}^{-1}T_{\bar{i}\bar{j}}^{-1}.
 \end{gather}
\end{lem}

\begin{proof}
 The curves~$\delta_i$, $\varepsilon$, $\varepsilon_i$, $\gamma_{jk}$, $\gamma_{\bar{j}\bar{k}}$, and the circuites around the punctures~$i$ and~$\bar{i}$ form a lantern configuration; the corresponding lantern relation is~\eqref{eq_lantern_eps}. The curves~$\delta_i'$, $\varepsilon_j$, $\varepsilon_k$, $\gamma_{j\bar{k}}$, $\gamma_{\bar{j}k}$, and the circuits around the punctures~$i$ and~$\bar{i}$ also form a lantern configuration; the corresponding lantern relation is~\eqref{eq_lantern_eps_jk}. One more lantern configuration is formed by the curves~$\delta_i'$, $\gamma_{k\bar{i}}$, $\gamma_{ki}$, $\talpha_k$, and the circuits around the punctures~$i$, $\bar{i}$, and~$k$; the corresponding relation is~\eqref{eq_131}.

 Finally, let us prove~\eqref{eq_Tdeltai}. It suffices to prove the first equality in~\eqref{eq_Tdeltai}, since the conjugation by~$r_i$ fixes~$T_{\delta_i}$ and sends $T_{i\bar{j}}^{-1}T_{\bar{k}i}^{-1}T_{ij}^{-1}T_{ki}^{-1}$ to~$T_{k\bar{i}}^{-1}T_{\bar{i}j}^{-1}T_{\bar{k}\bar{i}}^{-1}T_{\bar{i}\bar{j}}^{-1}$. Combining~\eqref{eq_lantern_eps} and relations~\eqref{eq_lantern_faces}, we have
 $$
 T_{\delta_i}=T_{jk}T_{\bar{j}\bar{k}}T_{\varepsilon_i}^{-1}T_{\varepsilon}^{-1}=T_{jk}T_{\bar{j}\bar{k}}\left(T_{\bar{j}\bar{k}}^{-1}T_{i\bar{j}}^{-1}T_{\bar{k}i}^{-1}\right)\left(T_{jk}^{-1}T_{ij}^{-1}T_{ki}^{-1}\right)=T_{i\bar{j}}^{-1}T_{\bar{k}i}^{-1}T_{ij}^{-1}T_{ki}^{-1},
 $$
 as required.
\end{proof}

We will also need several relations among the elements of~$G$ and the Dehn twists about the introduced simple closed curves. If $g\in G$ and $\gamma$ is a simple closed curve on~$\Sigma_{0,6}$, then $gT_{\gamma}g^{-1}=T_{g(\gamma)}$ if $g$ preserves the orientation of~$\Sigma_{0,6}$, and  $gT_{\gamma}g^{-1}=T_{g(\gamma)}^{-1}$ if $g$ reverses the orientation of~$\Sigma_{0,6}$. This simple observation immediately implies the following lemma.

\begin{lem}
 Let $g\in G$, and let $\sigma\in S_3$ be the corresponding permutation from condition~$(*)$. Then the following  relations hold in~$\CG$:
 \begin{align*}
  gT_{pq}g^{-1}&=T_{g(p)g(q)}^{{}\,e_1(g)},&
  gT_{\varepsilon}g^{-1}&=T_{\varepsilon}^{e_1(g)},&
  gT_{\varepsilon_i}g^{-1}&=T_{\varepsilon_{\sigma(i)}}^{e_1(g)}.
 \end{align*}
 If $e_1(g)=1$, i.e., $g\in\{1,\omega,\omega^2,r_1,r_2,r_3\}$, then we additionally have
 \begin{align*}
 gT_{\delta_i}g^{-1}&=T_{\delta_{\sigma(i)}},&
 gT_{\talpha_i}g^{-1}&=T_{\talpha_{\sigma(i)}}.
 \end{align*}
\end{lem}

Throughout this part of the paper,  we will repeatedly use these relations without furhter comment.

\section{Representatives of $\G$-orbits of cells}\label{section_representatives}

In this section, we describe an explicit set of representatives of~$\G$-orbits of the cells of~$\CB_{0,6}$. Also, for each such representative~$Q$, we compute its algebraic boundary~$\partial Q$ in~$C_*(\CB_{0,6})$, and its stabilizers~$\CS_Q$ and~$\G_Q$ in the groups~$\CS$ and~$\G$, respectively. The result is as follows.

\begin{propos}\label{propos_tables}
\textnormal{(1)} The complex $\CB_{0,6}$ contains exactly three $\G$-orbits of vertices, nine $\G$-orbits of\/ $1$-cells, eleven $\G$-orbits of\/ $2$-cells, and seven $\G$-orbits of\/ $3$-cells, representatives are listed in Tables~\ref{table_0cells}--\ref{table_3cells}, respectively.

\textnormal{(2)} For each of the three pairs of cells~$Z_3^{\pm}$, $Z_4^{\pm}$, and~$W_1^{\pm}$, the corresponding two $\G$-orbits form a single $\tG$-orbit. For any other cell in Tables~\ref{table_0cells}--\ref{table_3cells}, its $\tG$-orbit coincides with its $\G$-orbit.

\textnormal{(3)} All cells of~$\CB_{0,6}$ are simplices. For each cell~$Q$ in Tables~\ref{table_0cells}--\ref{table_3cells}, its vertex set~$\CV(Q)$ and its algebraic boundary~$\partial Q$ are as depicted in the tables. (We choose the orientation of~$Q$ given by the specified order of vertices.)

\textnormal{(4)} For every cell~$Q$ in Tables~\ref{table_0cells}--\ref{table_2cells}, the stabilizer~$\CS_Q$ is a free abelian group whose basis is depicted in the tables. The stabilizers~$\CS_Q$ are trivial for all $3$-cells of~$\CB_{0,6}$, so Table~\ref{table_3cells} does not contain the corresponding column.

\textnormal{(5)} For each cell~$Q$ in Tables~\ref{table_0cells}--\ref{table_3cells}, the last column contains elements $g_1,\ldots,g_s\in\G_Q$ that generate $\G_Q$ modulo~$\CS_Q$ (i.e., together with $\CS_Q$ they generate~$\G_Q$). We use the following convention: we write $\boxed{g}$, if $g$ stabilizes~$Q$ but reverses its orientation; elements preserving the orientation of~$Q$ are not boxed.
\end{propos}

%\clearpage

\afterpage{
\clearpage

\begin{longtable}[p!]{|c|m{5cm}|c|c|}
\caption{Representatives of $\G$-orbits of $0$-cells}\label{table_0cells}\\
\hline
\multirow{2}{*}{$Q$} & \multicolumn{1}{c|}{\multirow{2}{*}{Multiarc}}  & $\CS_Q$ & $\G_Q\!\!\!\mod \CS_Q$\\
& & (basis) & (generators)
\\
\hline
\endhead
$X_1$ &
\begin{tikzpicture}
 \thicklines
 \fill (0,-2.309) circle (1.5pt) node [below=-.5pt] {\tiny $\bar{2}$};;
 \fill (-2,1.155) circle (1.5pt) node [above left=-2pt] {\tiny $\bar{3}$};;
 \fill (2,1.155) circle (1.5pt) node [above right=-2pt] {\tiny $\bar{1}$};
 \fill (0,.693) circle (1.5pt) node [above=-.5pt] {\tiny $2$};
 \fill (-.6,-.346) circle (1.5pt) node [below left=-2pt] {\tiny $1$};
 \fill (.6,-.346) circle (1.5pt)node [below right=-2pt] {\tiny $3$};

 \draw (0,.693) arc(90:270:1.501);
 \draw (-.6,-.346) arc(-150:30:1.501);
 \draw (.6,-.346) arc(-30:150:1.501);
\end{tikzpicture}
 & $T_{\delta_1}$, $T_{\delta_2}$, $T_{\delta_3}$ & $\omega$, $r_1$, $T_{\varepsilon}\nu$ \\
 \hline
$X_2$ &
\begin{tikzpicture}
 \thicklines
 \fill (0,-2.309) circle (1.5pt) node [below=-.5pt] {\tiny $\bar{2}$};;
 \fill (-2,1.155) circle (1.5pt) node [above left=-2pt] {\tiny $\bar{3}$};;
 \fill (2,1.155) circle (1.5pt) node [above right=-2pt] {\tiny $\bar{1}$};
 \fill (0,.693) circle (1.5pt) node [above=-.5pt] {\tiny $2$};
 \fill (-.6,-.346) circle (1.5pt) node [below left=-2pt] {\tiny $1$};
 \fill (.6,-.346) circle (1.5pt)node [below right=-2pt] {\tiny $3$};

 \draw (0,.693)--(-2,1.155);
 \draw (-.6,-.346)--(2,1.155);
 \draw (.6,-.346)--(0,-2.309);
\end{tikzpicture}
 & $T_{\delta'_1}$, $T_{2\bar{3}}T_{\bar{2}3}^{-1}$ & $T_{2\bar{3}}$, $s_1$, $T_{\varepsilon_2}^{-1}r_1$ \\
 \hline
$X_3$ &
\begin{tikzpicture}
 \thicklines
 \fill (0,-2.309) circle (1.5pt) node [below=-.5pt] {\tiny $\bar{2}$};;
 \fill (-2,1.155) circle (1.5pt) node [above left=-2pt] {\tiny $\bar{3}$};;
 \fill (2,1.155) circle (1.5pt) node [above right=-2pt] {\tiny $\bar{1}$};
 \fill (0,.693) circle (1.5pt) node [above=-.5pt] {\tiny $2$};
 \fill (-.6,-.346) circle (1.5pt) node [below left=-2pt] {\tiny $1$};
 \fill (.6,-.346) circle (1.5pt)node [below right=-2pt] {\tiny $3$};

 \draw (0,.693)--(-2,1.155);
 \draw (-.6,-.346)--(0,-2.309);
 \draw (.6,-.346)--(2,1.155);
\end{tikzpicture}
 & $\varnothing$ & $T_{1\bar{2}}$, $\omega$, $r_1$ \\
 \hline
\end{longtable}

% \newpage
\begin{longtable}[p!]{|c|m{4.5cm}|c|c|c|c|}
\caption{Representatives of $\G$-orbits of $1$-cells}\label{table_1cells}\\
\hline
\multirow{2}{*}{$Q$} & \multicolumn{1}{c|}{\multirow{2}{*}{Multiarc}} & \multirow{2}{*}{$\CV(Q)$} & \multirow{2}{*}{$\partial Q$} & $\CS_Q$ &$\G_Q\!\!\!\mod \CS_Q$ \\
&&&& (basis) & (generators)
\\
\hline
\endhead
$Y_1$ &
\begin{tikzpicture}[scale=.9]
 \thicklines
 \fill (0,-2.309) circle (1.5pt) node [below=-.5pt] {\tiny $\bar{2}$};;
 \fill (-2,1.155) circle (1.5pt) node [above left=-2pt] {\tiny $\bar{3}$};;
 \fill (2,1.155) circle (1.5pt) node [above right=-2pt] {\tiny $\bar{1}$};
 \fill (0,.693) circle (1.5pt) node [above=-.5pt] {\tiny $2$};
 \fill (-.6,-.346) circle (1.5pt) node [below left=-2pt] {\tiny $1$};
 \fill (.6,-.346) circle (1.5pt)node [below right=-2pt] {\tiny $3$};

 \draw (0,.693) arc(90:270:1.501);
 \draw (-.6,-.346) arc(-150:30:1.501);
 \draw (.6,-.346) arc(-30:150:1.501);
 \draw (0,.693)--(-2,1.155);
 \draw (0,-2.309) .. controls (-2.4,.4) and (0,.4) .. (.6,-.346);
\end{tikzpicture}
& \multicolumn{1}{m{1.2cm}|}{\centering $X_1$, $T_{31}X_2$}  & $T_{31}X_2-X_1$ & $T_{\delta_1}$ & $T_{\delta_2}T_{\varepsilon}\nu$\\
\hline
$Y_2$ &
\begin{tikzpicture}[scale=.9]
 \thicklines
 \fill (0,-2.309) circle (1.5pt) node [below=-.5pt] {\tiny $\bar{2}$};;
 \fill (-2,1.155) circle (1.5pt) node [above left=-2pt] {\tiny $\bar{3}$};;
 \fill (2,1.155) circle (1.5pt) node [above right=-2pt] {\tiny $\bar{1}$};
 \fill (0,.693) circle (1.5pt) node [above=-.5pt] {\tiny $2$};
 \fill (-.6,-.346) circle (1.5pt) node [below left=-2pt] {\tiny $1$};
 \fill (.6,-.346) circle (1.5pt)node [below right=-2pt] {\tiny $3$};

 \draw (0,.693)--(-2,1.155);
 \draw (-.6,-.346)--(2,1.155);
 \draw (.6,-.346)--(0,-2.309);
 \draw (-.6,-.346)--(0,-2.309);
 \draw (.6,-.346) -- (2,1.155);
\end{tikzpicture}
& \multicolumn{1}{m{1.2cm}|}{\centering $X_2$, $X_3$}  & $X_3-X_2$ & $\varnothing$ & $T_{2\bar{3}}$\\
\hline
$Y_3$ &
\begin{tikzpicture}[scale=.9]
 \thicklines
 \fill (0,-2.309) circle (1.5pt) node [below=-.5pt] {\tiny $\bar{2}$};;
 \fill (-2,1.155) circle (1.5pt) node [above left=-2pt] {\tiny $\bar{3}$};;
 \fill (2,1.155) circle (1.5pt) node [above right=-2pt] {\tiny $\bar{1}$};
 \fill (0,.693) circle (1.5pt) node [above=-.5pt] {\tiny $2$};
 \fill (-.6,-.346) circle (1.5pt) node [below left=-2pt] {\tiny $1$};
 \fill (.6,-.346) circle (1.5pt)node [below right=-2pt] {\tiny $3$};

 \draw (0,.693) arc(90:270:1.501);
 \draw (-.6,-.346) arc(-150:30:1.501);
 \draw (.6,-.346) arc(-30:150:1.501);
 \draw (0,.693)--(-2,1.155);
 \draw (-.6,-.346)--(0,-2.309);
 \draw (.6,-.346)--(2,1.155);
\end{tikzpicture}
& \multicolumn{1}{m{1.2cm}|}{\centering $X_1$, $X_3$}  & $X_3-X_1$ & $\varnothing$ & $\omega$, $r_1$\\
\hline
$Y_4$ &
\begin{tikzpicture}[scale=.9]
 \thicklines
 \fill (0,-2.309) circle (1.5pt) node [below=-.5pt] {\tiny $\bar{2}$};;
 \fill (-2,1.155) circle (1.5pt) node [above left=-2pt] {\tiny $\bar{3}$};;
 \fill (2,1.155) circle (1.5pt) node [above right=-2pt] {\tiny $\bar{1}$};
 \fill (0,.693) circle (1.5pt) node [above=-.5pt] {\tiny $2$};
 \fill (-.6,-.346) circle (1.5pt) node [below left=-2pt] {\tiny $1$};
 \fill (.6,-.346) circle (1.5pt)node [below right=-2pt] {\tiny $3$};

 \draw (0,.693)--(2,1.155) ;
 \draw (-.6,-.346)--(-2,1.155);
 \draw (.6,-.346)--(0,-2.309);
 \draw (0,.693)--(-2,1.155);
 \draw (-.6,-.346)--(0,-2.309);
 \draw (.6,-.346)--(2,1.155);
\end{tikzpicture}
& \multicolumn{1}{m{1.2cm}|}{\centering $X_3$, $\nu X_3$}  & $(\nu-1)X_3$ & $\varnothing$ & $\omega$, $r_1$, $\boxed{\nu}$\\
\hline
$Y_5$ &
\begin{tikzpicture}[scale=.9]
 \path[use as bounding box] (-2.5,-2.8) rectangle (2.5,1.8);
 \thicklines
 \fill (0,-2.309) circle (1.5pt) node [below=-.5pt] {\tiny $\bar{2}$};;
 \fill (-2,1.155) circle (1.5pt) node [above left=-2pt] {\tiny $\bar{3}$};;
 \fill (2,1.155) circle (1.5pt) node [above right=-2pt] {\tiny $\bar{1}$};
 \fill (0,.693) circle (1.5pt) node [above=-.5pt] {\tiny $2$};
 \fill (-.6,-.346) circle (1.5pt) node [below left=-2pt] {\tiny $1$};
 \fill (.6,-.346) circle (1.5pt)node [below right=-2pt] {\tiny $3$};

 \draw (0,.693)--(2,1.155) ;
 \draw (-.6,-.346)--(-2,1.155);
 \draw (.6,-.346)--(2,1.155);
 \draw (.6,-.346)--(-2,1.155);
 \draw (-.6,-.346)--(0,-2.309);
 \draw (0,-2.309) .. controls (-3.5,1.2) and (-3,2.7) .. (0,.693);

\end{tikzpicture}
& \multicolumn{1}{m{1.2cm}|}{\centering $\omega^2X_2$, $r_3 X_2$}  & $(r_3-\omega^2)X_2$ & $\varnothing$ & $T_{\varepsilon_1}^{-1}\nu$, $\boxed{r_1}$\\
\hline
$Y_6$ &
\begin{tikzpicture}[scale=.9]
 \path[use as bounding box] (-2.5,-2.8) rectangle (2.5,2.4);
 \thicklines
 \fill (0,-2.309) circle (1.5pt) node [below=-.5pt] {\tiny $\bar{2}$};;
 \fill (-2,1.155) circle (1.5pt) node [above left=-2pt] {\tiny $\bar{3}$};;
 \fill (2,1.155) circle (1.5pt) node [above right=-2pt] {\tiny $\bar{1}$};
 \fill (0,.693) circle (1.5pt) node [above=-.5pt] {\tiny $2$};
 \fill (-.6,-.346) circle (1.5pt) node [below left=-2pt] {\tiny $1$};
 \fill (.6,-.346) circle (1.5pt)node [below right=-2pt] {\tiny $3$};

 \draw (0,.693) arc(90:270:1.501);
 \draw (-.6,-.346) arc(-150:30:1.501);
 \draw (.6,-.346) arc(-30:150:1.501);
 \draw (-.6,-.346) .. controls (0,0) and (.33,.3) .. (.33,.693) .. controls (.33,1.9) and (-2.4,-.1) .. (-2.4,1.3) .. controls (-2.4,2.4) and (1,2.4) .. (2,1.155);

\end{tikzpicture}
& \multicolumn{1}{m{1.2cm}|}{\centering $X_1$, $T_{\talpha_1} X_1$}  & $\bigl(T_{\talpha_1}-1\bigr)X_1$ & $T_{\delta_2}$, $T_{\delta_3}$ & $r_1$, $T_{\varepsilon}\nu$ \\
\hline
$Y_7$ &
\begin{tikzpicture}[scale=.9]
 \path[use as bounding box] (-2.5,-2.8) rectangle (2.5,1.8);
 \thicklines
 \fill (0,-2.309) circle (1.5pt) node [below=-.5pt] {\tiny $\bar{2}$};;
 \fill (-2,1.155) circle (1.5pt) node [above left=-2pt] {\tiny $\bar{3}$};;
 \fill (2,1.155) circle (1.5pt) node [above right=-2pt] {\tiny $\bar{1}$};
 \fill (0,.693) circle (1.5pt) node [above=-.5pt] {\tiny $2$};
 \fill (-.6,-.346) circle (1.5pt) node [below left=-2pt] {\tiny $1$};
 \fill (.6,-.346) circle (1.5pt)node [below right=-2pt] {\tiny $3$};

 \draw (0,.693)--(-2,1.155);
 \draw (-.6,-.346)--(2,1.155);
 \draw (.6,-.346)--(0,-2.309);
 \draw (2,1.155) .. controls (.711,-3.631) and (-.752,-3.898) .. (-.6,-.346);

\end{tikzpicture}
& \multicolumn{1}{m{1.2cm}|}{\centering $X_2$, $r_1 X_2$}  & $(r_1-1)X_2$ & $T_{2\bar{3}}T_{\bar{2}3}^{-1}$ & $T_{2\bar{3}}$, $s_1$, $\boxed{r_1}$\\
\hline
$Y_8$ &
\begin{tikzpicture}[scale=.9]
\path[use as bounding box] (-2.5,-2.8) rectangle (2.5,1.8);
 \thicklines
 \fill (0,-2.309) circle (1.5pt) node [below=-.5pt] {\tiny $\bar{2}$};;
 \fill (-2,1.155) circle (1.5pt) node [above left=-2pt] {\tiny $\bar{3}$};;
 \fill (2,1.155) circle (1.5pt) node [above right=-2pt] {\tiny $\bar{1}$};
 \fill (0,.693) circle (1.5pt) node [above=-.5pt] {\tiny $2$};
 \fill (-.6,-.346) circle (1.5pt) node [below left=-2pt] {\tiny $1$};
 \fill (.6,-.346) circle (1.5pt)node [below right=-2pt] {\tiny $3$};

 %\draw (0,.693) arc(90:270:1.501);
 \draw (-.6,-.346) arc(-150:30:1.501);
 %\draw (.6,-.346) arc(-30:150:1.501);
 \draw (0,.693)--(-2,1.155);
 \draw (0,-2.309) .. controls (-2.4,.4) and (0,.4) .. (.6,-.346);
 \draw (0,-2.309) .. controls (4.3,-.5) and (2.3,4.3) .. (.6,-.346);

\end{tikzpicture}
& \multicolumn{1}{m{1.2cm}|}{\centering $T_{31}X_2$, $T_{3\bar{1}}^{-1}X_2$}  & $\bigl(T_{3\bar{1}}^{-1}-T_{31}\bigr)X_2$ & $T_{\delta_1}$ & $T_{2\bar{3}}$, $\boxed{r_1}$\\
\hline
$Y_9$ &
\begin{tikzpicture}[scale=.9]
 \path[use as bounding box] (-2.5,-2.8) rectangle (2.5,1.8);
 \thicklines
 \fill (0,-2.309) circle (1.5pt) node [below=-.5pt] {\tiny $\bar{2}$};;
 \fill (-2,1.155) circle (1.5pt) node [above left=-2pt] {\tiny $\bar{3}$};;
 \fill (2,1.155) circle (1.5pt) node [above right=-2pt] {\tiny $\bar{1}$};
 \fill (0,.693) circle (1.5pt) node [above=-.5pt] {\tiny $2$};
 \fill (-.6,-.346) circle (1.5pt) node [below left=-2pt] {\tiny $1$};
 \fill (.6,-.346) circle (1.5pt)node [below right=-2pt] {\tiny $3$};

 \draw (0,.693)--(-2,1.155);
 \draw (-.6,-.346)--(0,-2.309);
 \draw (.6,-.346)--(2,1.155);
 \draw (-2,1.155) .. controls (-.711,-3.631) and (1,-4.4) .. (0,.693);

\end{tikzpicture}
& \multicolumn{1}{m{1.2cm}|}{\centering $X_3$, $T_{\talpha_1}^{-1} X_3$}  & $\bigl(T_{\talpha_1}^{-1}-1\bigr)X_3$ & $\varnothing$ & $T_{1\bar{2}}$, $T_{3\bar{1}}$, $r_1$\\
\hline
\end{longtable}

% \newpage
\begin{longtable}[p!]{|c|m{4.5cm}|c|c|c|c|}
\caption{Representatives of $\G$-orbits of $2$-cells}\label{table_2cells}\\
\hline
\multirow{2}{*}{$Q$} & \multicolumn{1}{c|}{\multirow{2}{*}{Multiarc}} & \multirow{2}{*}{$\CV(Q)$} & \multirow{2}{*}{$\partial Q$} & $\CS_Q$ &$\G_Q\!\!\!\mod \CS_Q$ \\
&&&& (basis) & (generators)
\\
\hline
\endhead
$Z_1$ &
\begin{tikzpicture}[scale=.9]
 \thicklines
 \fill (0,-2.309) circle (1.5pt) node [below=-.5pt] {\tiny $\bar{2}$};;
 \fill (-2,1.155) circle (1.5pt) node [above left=-2pt] {\tiny $\bar{3}$};;
 \fill (2,1.155) circle (1.5pt) node [above right=-2pt] {\tiny $\bar{1}$};
 \fill (0,.693) circle (1.5pt) node [above=-.5pt] {\tiny $2$};
 \fill (-.6,-.346) circle (1.5pt) node [below left=-2pt] {\tiny $1$};
 \fill (.6,-.346) circle (1.5pt)node [below right=-2pt] {\tiny $3$};

 \draw (0,.693) arc(90:270:1.501);
 \draw (-.6,-.346) arc(-150:30:1.501);
 \draw (.6,-.346) arc(-30:150:1.501);
 \draw (0,.693)--(-2,1.155);
 \draw (-.6,-.346)--(0,-2.309);
 \draw (.6,-.346)--(2,1.155);
 \draw (0,-2.309) .. controls (-2.4,.4) and (0,.4) .. (.6,-.346);
\end{tikzpicture}
& \multicolumn{1}{m{1.5cm}|}{\centering $X_1$, $T_{31}X_2$, $X_3$}  &
% \multicolumn{1}{p{2.2cm}|}{$Y_1+T_{1\bar{2}}^{-1}r_1Y_2$ ${}\ \ \,{}-Y_3\strut$}
$Y_1+T_{1\bar{2}}^{-1}r_1Y_2-Y_3$
& $\varnothing$ & $\varnothing$\\
\hline
$Z_2$ &
\begin{tikzpicture}[scale=.9]
 \thicklines
 \fill (0,-2.309) circle (1.5pt) node [below=-.5pt] {\tiny $\bar{2}$};;
 \fill (-2,1.155) circle (1.5pt) node [above left=-2pt] {\tiny $\bar{3}$};;
 \fill (2,1.155) circle (1.5pt) node [above right=-2pt] {\tiny $\bar{1}$};
 \fill (0,.693) circle (1.5pt) node [above=-.5pt] {\tiny $2$};
 \fill (-.6,-.346) circle (1.5pt) node [below left=-2pt] {\tiny $1$};
 \fill (.6,-.346) circle (1.5pt)node [below right=-2pt] {\tiny $3$};

 \draw (0,.693)--(2,1.155) ;
 \draw (-.6,-.346)--(-2,1.155);
 \draw (.6,-.346)--(0,-2.309);
 \draw (0,.693)--(-2,1.155);
 \draw (-.6,-.346)--(0,-2.309);
 \draw (.6,-.346)--(2,1.155);
 \draw (-.6,-.346)--(2,1.155);
\end{tikzpicture}
& \multicolumn{1}{m{1.5cm}|}{\centering \mbox{{}\ $X_2$,\ {}} \mbox{{}\ $X_3$,\ {}} $\nu X_3$}  & $(1-s_1)Y_2+Y_4$ & $\varnothing$ &  $\boxed{s_1}$\\
\hline
$Z_3^+$ &
\begin{tikzpicture}[scale=.9]
 \path[use as bounding box] (-2.5,-2.8) rectangle (2.5,1.8);
 \thicklines
 \fill (0,-2.309) circle (1.5pt) node [below=-.5pt] {\tiny $\bar{2}$};;
 \fill (-2,1.155) circle (1.5pt) node [above left=-2pt] {\tiny $\bar{3}$};;
 \fill (2,1.155) circle (1.5pt) node [above right=-2pt] {\tiny $\bar{1}$};
 \fill (0,.693) circle (1.5pt) node [above=-.5pt] {\tiny $2$};
 \fill (-.6,-.346) circle (1.5pt) node [below left=-2pt] {\tiny $1$};
 \fill (.6,-.346) circle (1.5pt)node [below right=-2pt] {\tiny $3$};

 \draw (0,.693)--(2,1.155) ;
 \draw (-.6,-.346)--(-2,1.155);
 \draw (.6,-.346)--(2,1.155);
 \draw (.6,-.346)--(-2,1.155);
 \draw (-.6,-.346)--(0,-2.309);
 \draw (-.6,-.346) arc(-150:30:1.501);
 \draw (0,-2.309) .. controls (-3.5,1.2) and (-3,2.7) .. (0,.693);

\end{tikzpicture}
& \multicolumn{1}{m{1.5cm}|}{\centering $T_{2\bar{3}}X_1$, $\omega^2X_2$, $r_3 X_2$}  & $(1-r_1)T_{2\bar{3}}r_2Y_1+Y_5$ & $\varnothing$ & $T_{\varepsilon_1}^{-1}\nu$, $\boxed{r_1}$\\
\hline
$Z_3^-$ &
\begin{tikzpicture}[scale=.9]
 \path[use as bounding box] (-2.5,-2.8) rectangle (2.5,1.8);
 \thicklines
 \fill (0,-2.309) circle (1.5pt) node [below=-.5pt] {\tiny $\bar{2}$};;
 \fill (-2,1.155) circle (1.5pt) node [above left=-2pt] {\tiny $\bar{3}$};;
 \fill (2,1.155) circle (1.5pt) node [above right=-2pt] {\tiny $\bar{1}$};
 \fill (0,.693) circle (1.5pt) node [above=-.5pt] {\tiny $2$};
 \fill (-.6,-.346) circle (1.5pt) node [below left=-2pt] {\tiny $1$};
 \fill (.6,-.346) circle (1.5pt)node [below right=-2pt] {\tiny $3$};

 \draw (0,.693)--(2,1.155) ;
 \draw (-.6,-.346)--(-2,1.155);
 \draw (.6,-.346)--(2,1.155);
 \draw (.6,-.346)--(-2,1.155);
 \draw (-.6,-.346)--(0,-2.309);
 \draw (-.6,-.346) .. controls (-3,.9) and (-3.05,2.8) .. (-.05,.4) .. controls (.45,0) and  (1,.5) .. (2,1.155);
 \draw (0,-2.309) .. controls (-3.6,1.2) and (-3,2.9) .. (0,.693);

\end{tikzpicture}
& \multicolumn{1}{m{1.5cm}|}{\centering $hX_1$, $\omega^2X_2$, $r_3 X_2$}  & \multicolumn{1}{m{3.8cm}|}{\centering \phantom{$(r_1-1)hr_3Y_1+Y_5$} \phantom{$(r_1-1)hr_3Y_1+Y_5$} $(1-r_1)h\omega^2Y_1+Y_5$ \phantom{$(r_1-1)hr_3Y_1+Y_5$} \mbox{where $h=T_{\bar{3}1}^{-1}T_{\bar{1}2}^{-1}\nu$}} & $\varnothing$ & $T_{\varepsilon_1}^{-1}\nu$, $\boxed{r_1}$\\
\hline
$Z_4^+$ &
\begin{tikzpicture}[scale=.9]
 \path[use as bounding box] (-2.5,-2.8) rectangle (2.5,1.8);
 \thicklines
 \fill (0,-2.309) circle (1.5pt) node [below=-.5pt] {\tiny $\bar{2}$};;
 \fill (-2,1.155) circle (1.5pt) node [above left=-2pt] {\tiny $\bar{3}$};;
 \fill (2,1.155) circle (1.5pt) node [above right=-2pt] {\tiny $\bar{1}$};
 \fill (0,.693) circle (1.5pt) node [above=-.5pt] {\tiny $2$};
 \fill (-.6,-.346) circle (1.5pt) node [below left=-2pt] {\tiny $1$};
 \fill (.6,-.346) circle (1.5pt)node [below right=-2pt] {\tiny $3$};

 \draw (0,.693)--(2,1.155) ;
 \draw (-.6,-.346)--(-2,1.155);
 \draw (.6,-.346)--(2,1.155);
 \draw (.6,-.346)--(-2,1.155);
 \draw (-.6,-.346)--(0,-2.309);
 \draw (0,.693) -- (-2,1.155);
 \draw (0,-2.309) .. controls (-3.5,1.2) and (-3,2.7) .. (0,.693);

\end{tikzpicture}
& \multicolumn{1}{m{1.5cm}|}{\centering $\omega^2X_2$, $r_3 X_2$, $X_3$}  & $(r_3-\omega^2)Y_2+Y_5$ & $\varnothing$ & $\boxed{r_1}$\\
\hline
$Z_4^-$ &
\begin{tikzpicture}[scale=.9]
 \path[use as bounding box] (-2.5,-2.8) rectangle (2.5,1.8);
 \thicklines
 \fill (0,-2.309) circle (1.5pt) node [below=-.5pt] {\tiny $\bar{2}$};;
 \fill (-2,1.155) circle (1.5pt) node [above left=-2pt] {\tiny $\bar{3}$};;
 \fill (2,1.155) circle (1.5pt) node [above right=-2pt] {\tiny $\bar{1}$};
 \fill (0,.693) circle (1.5pt) node [above=-.5pt] {\tiny $2$};
 \fill (-.6,-.346) circle (1.5pt) node [below left=-2pt] {\tiny $1$};
 \fill (.6,-.346) circle (1.5pt)node [below right=-2pt] {\tiny $3$};

 \draw (0,.693)--(2,1.155) ;
 \draw (-.6,-.346)--(-2,1.155);
 \draw (.6,-.346)--(2,1.155);
 \draw (.6,-.346)--(-2,1.155);
 \draw (-.6,-.346)--(0,-2.309);
 \draw (0,-2.309)--(.6,-.346);
 \draw (0,-2.309) .. controls (-3.5,1.2) and (-3,2.7) .. (0,.693);

\end{tikzpicture}
& \multicolumn{1}{m{1.5cm}|}{\centering $\omega^2X_2$, $r_3 X_2$, $\nu X_3$}  & $(\omega-r_2)\nu Y_2+Y_5$ & $\varnothing$ & $\boxed{r_1}$\\
\hline
$Z_5$ &
\begin{tikzpicture}[scale=.9]
 \path[use as bounding box] (-2.5,-2.8) rectangle (2.5,2.4);
 \thicklines
 \fill (0,-2.309) circle (1.5pt) node [below=-.5pt] {\tiny $\bar{2}$};;
 \fill (-2,1.155) circle (1.5pt) node [above left=-2pt] {\tiny $\bar{3}$};;
 \fill (2,1.155) circle (1.5pt) node [above right=-2pt] {\tiny $\bar{1}$};
 \fill (0,.693) circle (1.5pt) node [above=-.5pt] {\tiny $2$};
 \fill (-.6,-.346) circle (1.5pt) node [below left=-2pt] {\tiny $1$};
 \fill (.6,-.346) circle (1.5pt)node [below right=-2pt] {\tiny $3$};

 \draw (0,.693) arc(90:270:1.501);
 \draw (-.6,-.346) arc(-150:30:1.501);
 \draw (.6,-.346) arc(-30:150:1.501);
 \draw (-.6,-.346)--(0,-2.309);
 \draw (2,1.155) .. controls (1.3,2.6) and (-2.55,2.6) .. (-2.55, 1.3) .. controls (-2.55,0) and (-1,.8) .. (0,.693);
 \draw (-.6,-.346) .. controls (0,0) and (.33,.3) .. (.33,.693) .. controls (.33,1.9) and (-2.4,-.1) .. (-2.4,1.3) .. controls (-2.4,2.4) and (1,2.4) .. (2,1.155);

\end{tikzpicture}
& \multicolumn{1}{m{1.5cm}|}{\centering $X_1$, $T_{\talpha_1} X_1$, ${}\!T_{2\bar{3}}^{-1}\omega^2X_2$\!{}}  & $\bigl(T_{\talpha_1}\,\omega^2-r_2\bigr)Y_1+Y_6$ & $T_{\delta_3}$ &  $T_{\delta_2}T_{\varepsilon}\nu$ \\
\hline
$Z_6$ &
\begin{tikzpicture}[scale=.9]
 \path[use as bounding box] (-2.5,-2.8) rectangle (2.5,1.8);
 \thicklines
 \fill (0,-2.309) circle (1.5pt) node [below=-.5pt] {\tiny $\bar{2}$};;
 \fill (-2,1.155) circle (1.5pt) node [above left=-2pt] {\tiny $\bar{3}$};;
 \fill (2,1.155) circle (1.5pt) node [above right=-2pt] {\tiny $\bar{1}$};
 \fill (0,.693) circle (1.5pt) node [above=-.5pt] {\tiny $2$};
 \fill (-.6,-.346) circle (1.5pt) node [below left=-2pt] {\tiny $1$};
 \fill (.6,-.346) circle (1.5pt)node [below right=-2pt] {\tiny $3$};

 \draw (0,.693)--(-2,1.155);
 \draw (-.6,-.346)--(2,1.155);
 \draw (.6,-.346)--(0,-2.309);
 \draw (-.6,-.346)--(0,-2.309);
 \draw (.6,-.346)--(2,1.155);
 \draw (2,1.155) .. controls (.711,-3.631) and (-.752,-3.898) .. (-.6,-.346);

\end{tikzpicture}
& \multicolumn{1}{m{1.5cm}|}{\centering $X_2$, $r_1 X_2$, $X_3$}  & $(r_1-1)Y_2+Y_7$ & $\varnothing$ & $T_{2\bar{3}}$,  $\boxed{r_1}$\\
\hline
$Z_7$ &
\begin{tikzpicture}[scale=.9]
\path[use as bounding box] (-2.5,-2.8) rectangle (2.5,2.3);
 \thicklines
 \fill (0,-2.309) circle (1.5pt) node [below=-.5pt] {\tiny $\bar{2}$};;
 \fill (-2,1.155) circle (1.5pt) node [above left=-2pt] {\tiny $\bar{3}$};;
 \fill (2,1.155) circle (1.5pt) node [above right=-2pt] {\tiny $\bar{1}$};
 \fill (0,.693) circle (1.5pt) node [above=-.5pt] {\tiny $2$};
 \fill (-.6,-.346) circle (1.5pt) node [below left=-2pt] {\tiny $1$};
 \fill (.6,-.346) circle (1.5pt)node [below right=-2pt] {\tiny $3$};

 \draw (0,.693) arc(90:270:1.501);
 \draw (-.6,-.346) arc(-150:30:1.501);
 \draw (.6,-.346) arc(-30:150:1.501);
 \draw (0,.693)--(-2,1.155);
%  \draw (2,1.155)--(.6,-.346); %lishnee
 \draw (0,-2.309) .. controls (-2.4,.4) and (0,.4) .. (.6,-.346);
 \draw (0,-2.309) .. controls (4,-.5) and (2.4,3.15) .. (1.36,.95) .. controls (1.1,.4)   .. (.6,-.346);

\end{tikzpicture}
& \multicolumn{1}{m{1.5cm}|}{\centering $X_1$, $T_{31}X_2$, $T_{3\bar{1}}^{-1}X_2$}  & $(1-r_1)Y_1+Y_8$ & $T_{\delta_1}$ &  $\boxed{r_1}$\\
\hline
$Z_8$ &
\begin{tikzpicture}[scale=.9]
\path[use as bounding box] (-2.5,-2.8) rectangle (2.5,1.8);
 \thicklines
 \fill (0,-2.309) circle (1.5pt) node [below=-.5pt] {\tiny $\bar{2}$};;
 \fill (-2,1.155) circle (1.5pt) node [above left=-2pt] {\tiny $\bar{3}$};;
 \fill (2,1.155) circle (1.5pt) node [above right=-2pt] {\tiny $\bar{1}$};
 \fill (0,.693) circle (1.5pt) node [above=-.5pt] {\tiny $2$};
 \fill (-.6,-.346) circle (1.5pt) node [below left=-2pt] {\tiny $1$};
 \fill (.6,-.346) circle (1.5pt)node [below right=-2pt] {\tiny $3$};

 %\draw (0,.693) arc(90:270:1.501);
 \draw (-.6,-.346) arc(-150:30:1.501);
 \draw (2,1.155)--(.6,-.346);
 \draw (0,.693)--(-2,1.155);
 \draw (-.6,-.346)--(0,-2.309);
 \draw (0,-2.309) .. controls (-2.4,.4) and (0,.4) .. (.6,-.346);
 \draw (0,-2.309) .. controls (4,-.5) and (2.4,3.15) .. (1.36,.95) .. controls (1.1,.4)   .. (.6,-.346);
\end{tikzpicture}
& \multicolumn{1}{m{1.5cm}|}{\centering $T_{31}X_2$, $T_{3\bar{1}}^{-1}X_2$, $X_3$}  & $(1-r_1)T_{3\bar{1}}^{-1}Y_2+Y_8$ & $\varnothing$ & $T_{2\bar{3}}$, $\boxed{r_1}$\\
\hline
$Z_9$ &
\begin{tikzpicture}[scale=.9]
 \path[use as bounding box] (-2.5,-2.8) rectangle (2.5,1.8);
 \thicklines
 \fill (0,-2.309) circle (1.5pt) node [below=-.5pt] {\tiny $\bar{2}$};;
 \fill (-2,1.155) circle (1.5pt) node [above left=-2pt] {\tiny $\bar{3}$};;
 \fill (2,1.155) circle (1.5pt) node [above right=-2pt] {\tiny $\bar{1}$};
 \fill (0,.693) circle (1.5pt) node [above=-.5pt] {\tiny $2$};
 \fill (-.6,-.346) circle (1.5pt) node [below left=-2pt] {\tiny $1$};
 \fill (.6,-.346) circle (1.5pt)node [below right=-2pt] {\tiny $3$};

 \draw (0,.693)--(-2,1.155);
 \draw (-.6,-.346)--(0,-2.309);
 \draw (.6,-.346)--(2,1.155);
 \draw (-.6,-.346)--(-2,1.155);
 \draw (0,.693)--(0,-2.309);
 \draw (-2,1.155) .. controls (-.711,-3.631) and (1,-4.4) .. (0,.693);

\end{tikzpicture}
& \multicolumn{1}{m{1.5cm}|}{\centering $\omega X_2$, $X_3$, $T_{\talpha_1}^{-1} X_3$}  & $(\omega-T_{\varepsilon_1}r_3)Y_2+Y_9$ & $\varnothing$ &  $T_{3\bar{1}}$\\
\hline
\end{longtable}

\begin{longtable}[p!]{|c|m{4.5cm}|c|c|c|c|}
\caption{Representatives of $\G$-orbits of $3$-cells}\label{table_3cells}\\
\hline
\multirow{2}{*}{$Q$} & \multicolumn{1}{c|}{\multirow{2}{*}{Multiarc}} & \multirow{2}{*}{$\CV(Q)$} & \multirow{2}{*}{$\partial Q$}  &$\G_Q$ \\
&&&&  (generators)
\\
\hline
\endhead
$W_1^+$ &
\begin{tikzpicture}[scale=.9]
 \path[use as bounding box] (-2.5,-2.8) rectangle (2.5,1.8);
 \thicklines
 \fill (0,-2.309) circle (1.5pt) node [below=-.5pt] {\tiny $\bar{2}$};;
 \fill (-2,1.155) circle (1.5pt) node [above left=-2pt] {\tiny $\bar{3}$};;
 \fill (2,1.155) circle (1.5pt) node [above right=-2pt] {\tiny $\bar{1}$};
 \fill (0,.693) circle (1.5pt) node [above=-.5pt] {\tiny $2$};
 \fill (-.6,-.346) circle (1.5pt) node [below left=-2pt] {\tiny $1$};
 \fill (.6,-.346) circle (1.5pt)node [below right=-2pt] {\tiny $3$};

 \draw (0,.693)--(2,1.155) ;
 \draw (-.6,-.346)--(-2,1.155);
 \draw (.6,-.346)--(2,1.155);
 \draw (.6,-.346)--(-2,1.155);
 \draw (-.6,-.346)--(0,-2.309);
 \draw (0,.693)--(-2,1.155);
 \draw (-.6,-.346) arc(-150:30:1.501);
 \draw (0,-2.309) .. controls (-3.5,1.2) and (-3,2.7) .. (0,.693);

\end{tikzpicture}
& \multicolumn{1}{m{2cm}|}{\centering $T_{2\bar{3}}X_1$, $\omega^2X_2$, \mbox{{}\ \ $r_3 X_2$,\ \ {}} $X_3$}  & $T_{2\bar{3}}(r_2-\omega)Z_1-Z_3^++Z_4^+$ &   $\boxed{r_1}$\\
\hline
$W_1^-$ &
\begin{tikzpicture}[scale=.9]
 \path[use as bounding box] (-2.5,-2.8) rectangle (2.5,1.8);
 \thicklines
 \fill (0,-2.309) circle (1.5pt) node [below=-.5pt] {\tiny $\bar{2}$};;
 \fill (-2,1.155) circle (1.5pt) node [above left=-2pt] {\tiny $\bar{3}$};;
 \fill (2,1.155) circle (1.5pt) node [above right=-2pt] {\tiny $\bar{1}$};
 \fill (0,.693) circle (1.5pt) node [above=-.5pt] {\tiny $2$};
 \fill (-.6,-.346) circle (1.5pt) node [below left=-2pt] {\tiny $1$};
 \fill (.6,-.346) circle (1.5pt)node [below right=-2pt] {\tiny $3$};

 \draw (0,.693)--(2,1.155) ;
 \draw (-.6,-.346)--(-2,1.155);
 \draw (.6,-.346)--(2,1.155);
 \draw (.6,-.346)--(-2,1.155);
 \draw (-.6,-.346)--(0,-2.309);
 \draw (.6,-.346)--(0,-2.309);
 \draw (-.6,-.346) .. controls (-3,.9) and (-3.05,2.8) .. (-.05,.4) .. controls (.45,0) and  (1,.5) .. (2,1.155);
 \draw (0,-2.309) .. controls (-3.6,1.2) and (-3,2.9) .. (0,.693);

\end{tikzpicture}
& \multicolumn{1}{m{2cm}|}{\centering ${}\!\!T_{\bar{3}1}^{-1}T_{\bar{1}2}^{-1}\nu X_1$,\!\!\!{} \mbox{{}\ $\omega^2X_2$,\ {}} \mbox{{}\ \ $r_3 X_2$,\ \ {}} $\nu X_3$}  & \multicolumn{1}{m{4.4cm}|}{\centering $T_{\bar{3}1}^{-1}T_{\bar{1}2}^{-1}(\omega^2-r_3)\nu Z_1\qquad {}$ ${}\qquad\qquad\ -Z_3^-+Z_4^-$} &   $\boxed{r_1}$\\
\hline
$W_2$ &
\begin{tikzpicture}[scale=.9]
 \path[use as bounding box] (-2.5,-2.8) rectangle (2.5,1.8);
 \thicklines
 \fill (0,-2.309) circle (1.5pt) node [below=-.5pt] {\tiny $\bar{2}$};;
 \fill (-2,1.155) circle (1.5pt) node [above left=-2pt] {\tiny $\bar{3}$};;
 \fill (2,1.155) circle (1.5pt) node [above right=-2pt] {\tiny $\bar{1}$};
 \fill (0,.693) circle (1.5pt) node [above=-.5pt] {\tiny $2$};
 \fill (-.6,-.346) circle (1.5pt) node [below left=-2pt] {\tiny $1$};
 \fill (.6,-.346) circle (1.5pt)node [below right=-2pt] {\tiny $3$};

 \draw (0,.693)--(2,1.155) ;
 \draw (-.6,-.346)--(-2,1.155);
 \draw (.6,-.346)--(2,1.155);
 \draw (.6,-.346)--(-2,1.155);
 \draw (-.6,-.346)--(0,-2.309);
 \draw (0,.693) -- (-2,1.155);
 \draw (0,-2.309)--(.6,-.346);
 \draw (0,-2.309) .. controls (-3.5,1.2) and (-3,2.7) .. (0,.693);

\end{tikzpicture}
& \multicolumn{1}{m{2cm}|}{\centering $\omega^2X_2$, \mbox{{}\quad $r_3 X_2$,\quad{}} \mbox{{}\quad $X_3$,\quad{}} $\nu X_3$}  & $(-\omega\nu-\omega^2)Z_2-Z_4^++Z_4^-$  & $\boxed{r_1}$\\
\hline
$W_3$ &
\begin{tikzpicture}[scale=.9]
 \path[use as bounding box] (-2.5,-2.8) rectangle (2.5,2.4);
 \thicklines
 \fill (0,-2.309) circle (1.5pt) node [below=-.5pt] {\tiny $\bar{2}$};;
 \fill (-2,1.155) circle (1.5pt) node [above left=-2pt] {\tiny $\bar{3}$};;
 \fill (2,1.155) circle (1.5pt) node [above right=-2pt] {\tiny $\bar{1}$};
 \fill (0,.693) circle (1.5pt) node [above=-.5pt] {\tiny $2$};
 \fill (-.6,-.346) circle (1.5pt) node [below left=-2pt] {\tiny $1$};
 \fill (.6,-.346) circle (1.5pt)node [below right=-2pt] {\tiny $3$};

 \draw (0,.693) arc(90:270:1.501);
 \draw (-.6,-.346) arc(-150:30:1.501);
 \draw (.6,-.346) arc(-30:150:1.501);
 \draw (-.6,-.346)--(0,-2.309);
 \draw (.6,-.346)--(2,1.155);
 \draw (2,1.155) .. controls (1.3,2.6) and (-2.55,2.6) .. (-2.55, 1.3) .. controls (-2.55,0) and (-1,.8) .. (0,.693);
 \draw (-.6,-.346) .. controls (0,0) and (.33,.3) .. (.33,.693) .. controls (.33,1.9) and (-2.4,-.1) .. (-2.4,1.3) .. controls (-2.4,2.4) and (1,2.4) .. (2,1.155);
 \draw (-.6,-.346) .. controls (0,-.2) and (.5,.3) .. (.5,.693) .. controls (.5,1.8) and (-1.2,1.4) .. (-2,1.155);

\end{tikzpicture}
& \multicolumn{1}{m{2cm}|}{\centering \mbox{{}\quad$X_1$,\quad{}} $T_{\talpha_1} X_1$, $T_{2\bar{3}}^{-1}\omega^2X_2$, ${}\!\!T_{2\bar{3}}^{-1}r_3X_2$}  &\multicolumn{1}{m{4.4cm}|}{\centering $T_{2\bar{3}}^{-1}(-Z_3^++Z_3^-)\qquad\qquad{}$ ${}\quad\quad\ +(r_1-1)Z_5$}
&  \multicolumn{1}{m{1.8cm}|}{\centering $T_{\delta_2}T_{\delta_3}T_{\varepsilon}\nu$, $\boxed{r_1}$} \\
\hline
$W_4$ &
\begin{tikzpicture}[scale=.9]
 \path[use as bounding box] (-2.5,-2.8) rectangle (2.5,1.8);
 \thicklines
 \fill (0,-2.309) circle (1.5pt) node [below=-.5pt] {\tiny $\bar{2}$};;
 \fill (-2,1.155) circle (1.5pt) node [above left=-2pt] {\tiny $\bar{3}$};;
 \fill (2,1.155) circle (1.5pt) node [above right=-2pt] {\tiny $\bar{1}$};
 \fill (0,.693) circle (1.5pt) node [above=-.5pt] {\tiny $2$};
 \fill (-.6,-.346) circle (1.5pt) node [below left=-2pt] {\tiny $1$};
 \fill (.6,-.346) circle (1.5pt)node [below right=-2pt] {\tiny $3$};

 \draw (-.6,-.346)-- (-2,1.155);
 \draw (2,1.155)--(0,.693);
 \draw (0,.693)--(-2,1.155);
 \draw (-.6,-.346)--(2,1.155);
 \draw (.6,-.346)--(0,-2.309);
 \draw (-.6,-.346)--(0,-2.309);
 \draw (.6,-.346)--(2,1.155);
 \draw (2,1.155) .. controls (.711,-3.631) and (-.752,-3.898) .. (-.6,-.346);

\end{tikzpicture}
& \multicolumn{1}{m{2cm}|}{\centering \mbox{{}\quad$X_2$,\quad{}} \mbox{{}\quad$r_1 X_2$,\quad{}} \mbox{{}\quad$X_3$,\quad{}} $\nu X_3$}  & $-(\nu+1)(Z_2+Z_6)$ &  $\nu$,  $\boxed{r_1}$\\
\hline
$W_5$ &
\begin{tikzpicture}[scale=.9]
\path[use as bounding box] (-2.5,-2.8) rectangle (2.5,2.3);
 \thicklines
 \fill (0,-2.309) circle (1.5pt) node [below=-.5pt] {\tiny $\bar{2}$};;
 \fill (-2,1.155) circle (1.5pt) node [above left=-2pt] {\tiny $\bar{3}$};;
 \fill (2,1.155) circle (1.5pt) node [above right=-2pt] {\tiny $\bar{1}$};
 \fill (0,.693) circle (1.5pt) node [above=-.5pt] {\tiny $2$};
 \fill (-.6,-.346) circle (1.5pt) node [below left=-2pt] {\tiny $1$};
 \fill (.6,-.346) circle (1.5pt)node [below right=-2pt] {\tiny $3$};

 \draw (0,.693) arc(90:270:1.501);
 \draw (-.6,-.346) arc(-150:30:1.501);
 \draw (.6,-.346) arc(-30:150:1.501);
 \draw (0,.693)--(-2,1.155);
 \draw (2,1.155)--(.6,-.346);
 \draw (-.6,-.346)--(0,-2.309);
 \draw (0,-2.309) .. controls (-2.4,.4) and (0,.4) .. (.6,-.346);
 \draw (0,-2.309) .. controls (4,-.5) and (2.4,3.15) .. (1.36,.95) .. controls (1.1,.4)   .. (.6,-.346);

\end{tikzpicture}
& \multicolumn{1}{m{2cm}|}{\centering \mbox{{}\quad$X_1$,\quad{}} \mbox{{}\ $T_{31}X_2$,\ {}} \mbox{{}\ $T_{3\bar{1}}^{-1}X_2$,\ {}} $X_3$}  & $(1-r_1)Z_1-Z_7+Z_8$  &  $\boxed{r_1}$\\
\hline
$W_6$ &
\begin{tikzpicture}[scale=.9]
 \path[use as bounding box] (-2.5,-3) rectangle (2.5,1.8);
 \thicklines
 \fill (0,-2.309) circle (1.5pt) node [below=-.5pt] {\tiny $\bar{2}$};;
 \fill (-2,1.155) circle (1.5pt) node [above left=-2pt] {\tiny $\bar{3}$};;
 \fill (2,1.155) circle (1.5pt) node [above right=-2pt] {\tiny $\bar{1}$};
 \fill (0,.693) circle (1.5pt) node [above=-.5pt] {\tiny $2$};
 \fill (-.6,-.346) circle (1.5pt) node [below left=-2pt] {\tiny $1$};
 \fill (.6,-.346) circle (1.5pt)node [below right=-2pt] {\tiny $3$};

 \draw (0,.693)--(-2,1.155);
 \draw (-.6,-.346)--(0,-2.309);
 \draw (.6,-.346)--(2,1.155);
 \draw (-.6,-.346)--(-2,1.155);
 \draw (0,.693)--(0,-2.309);
 \draw (2,1.155)--(0,.693);
 \draw (-2,1.155) .. controls (-.711,-3.631) and (1,-4.4) .. (0,.693);
 \draw (-2,1.155) .. controls (-1.1,-3.8) and (.9,-4.2) .. (.6,-.346);

\end{tikzpicture}
& \multicolumn{1}{m{2cm}|}{\centering \mbox{{}\quad$\omega X_2$,\quad{}} \mbox{{}\quad$r_2X_2$,\quad{}} \mbox{{}\quad$X_3$,\quad{}} $T_{\talpha_1}^{-1} X_3$}  &
\multicolumn{1}{m{4.4cm}|}{\centering $-T_{\varepsilon_1}Z_4^++\nu Z_4^-\qquad\qquad{}$ ${}\quad\quad\ +(r_1-1)Z_9$}
&   $\boxed{r_1}$\\
\hline
\end{longtable}
\clearpage
}

\begin{remark}
 All arcs in all pictures in Tables~\ref{table_0cells}--\ref{table_3cells} are oriented from a point in $\mathbf{p}=\{1,2,3\}$ to a point in $\mathbf{q}=\bigl\{\bar{1},\bar{2},\bar{3}\bigr\}$; hence, these orientations are not shown explicitly in the pictures.
\end{remark}

\begin{proof}
Let $Q$ be a cell of~$\CB_{0,6}$. Then the corresponding multiarc~$L\in\CL$ is a bipartite graph on sphere with three vertices in either part (the graph is allowed to have multiple edges but is not allowed to have isotopic edges). Let us prove that this graph satisfies the following conditions:
\begin{itemize}
 \item[(i)] $L$ has no isolated vertices,
 \item[(ii)] every connected component of~$L$ either is a single edge or contains no leaves.
\end{itemize}
(Recall that an \textit{isolated vertex} is a vertex with no edges adjacent to it and a \textit{leaf} is a vertex with exactly one edge adjacent to it.) Since $L\in\CL$, there exists a relative cycle~$\bbeta\in\CB_{0,6}$ with
\begin{equation*}
 \partial\bbeta=\sum_{i=1}^3\bigl(\bar{i}-i\bigr)
\end{equation*}
and support precisely~$L$. (In the formula for~$\partial\bbeta$, $i$ should be considered just  as a symbol denoting a puncture, not as a number.) Condition~(i) follows immediately from the fact that all vertices~$i$ and~$\bar{i}$ enter the boundary~$\partial\bbeta$. Let us prove condition~(ii). Suppose that $i$ is a leaf of~$L$. (The case of a leaf~$\bar{i}$ is completely similar.) Let $\beta$ be the edge outgoing from~$i$, and let $\bar{j}$ be the endpoint of~$\beta$. Since $i$ enters~$\partial\bbeta$ with coefficient~$-1$, we obtain that $\beta$ enters~$\bbeta$ with coefficient~$1$. If $L$ contained another edge adjacent to~$\bar{j}$, then $\bar{j}$ would enter~$\partial\bbeta$ with coefficient strictly greater than~$1$, which is not the case. Therefore $\bar{j}$ is also a leaf, and hence the edge~$\beta$ forms a connected component of~$L$. Condition~(ii) follows.

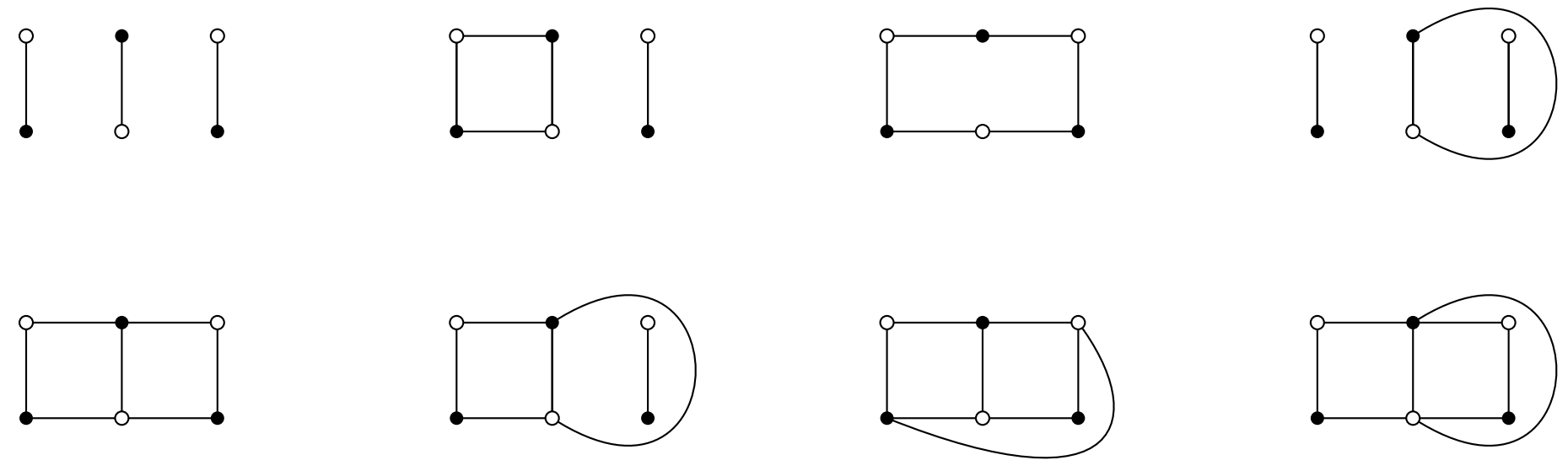
\begin{figure}[!ht]
\begin{tikzpicture}
 \begin{scope}[scale=.75]
  \path[use as bounding box] (0,-.3) rectangle (2.4,1.2);
  \thicklines
  \draw(0,0)--(0,1);
  \draw(1,0)--(1,1);
  \draw(2,0)--(2,1);
  \fill(0,0) circle(2pt);
  \fill(1,1) circle(2pt);
  \fill(2,0) circle(2pt);
  \fill[white](0,1) circle(2pt);
  \fill[white](1,0) circle(2pt);
  \fill[white](2,1) circle(2pt);
  \draw(0,1) circle(2pt);
  \draw(1,0) circle(2pt);
  \draw(2,1) circle(2pt);
 \end{scope}

 \begin{scope}[scale=.75,xshift=4.5cm]
  \path[use as bounding box] (0,-.3) rectangle (2.4,1.2);
  \thicklines
  \draw(0,0)--(0,1);
  \draw(1,0)--(1,1);
  \draw(2,0)--(2,1);
  \draw(0,0)--(1,0);
  \draw(0,1)--(1,1);
  \draw(0,0)--(0,1);
  \draw(1,0)--(1,1);
  \draw(2,0)--(2,1);
  \fill(0,0) circle(2pt);
  \fill(1,1) circle(2pt);
  \fill(2,0) circle(2pt);
  \fill[white](0,1) circle(2pt);
  \fill[white](1,0) circle(2pt);
  \fill[white](2,1) circle(2pt);
  \draw(0,1) circle(2pt);
  \draw(1,0) circle(2pt);
  \draw(2,1) circle(2pt);
 \end{scope}

 \begin{scope}[scale=.75,xshift=9cm]
  \path[use as bounding box] (0,-.3) rectangle (2.4,1.2);
  \thicklines
  \draw(0,0)--(2,0)--(2,1)--(0,1)--(0,0);
  \fill(0,0) circle(2pt);
  \fill(1,1) circle(2pt);
  \fill(2,0) circle(2pt);
  \fill[white](0,1) circle(2pt);
  \fill[white](1,0) circle(2pt);
  \fill[white](2,1) circle(2pt);
  \draw(0,1) circle(2pt);
  \draw(1,0) circle(2pt);
  \draw(2,1) circle(2pt);
 \end{scope}

 \begin{scope}[scale=.75,xshift=13.5cm]
  \path[use as bounding box] (0,-.3) rectangle (2.4,1.2);
  \thicklines
  \draw(0,0)--(0,1);
  \draw(1,0)--(1,1);
  \draw(2,0)--(2,1);
  \draw(1,0) .. controls (3,-1.3) and (3,2.3) .. (1,1);
  \draw(0,0)--(0,1);
  \draw(1,0)--(1,1);
  \draw(2,0)--(2,1);
  \fill(0,0) circle(2pt);
  \fill(1,1) circle(2pt);
  \fill(2,0) circle(2pt);
  \fill[white](0,1) circle(2pt);
  \fill[white](1,0) circle(2pt);
  \fill[white](2,1) circle(2pt);
  \draw(0,1) circle(2pt);
  \draw(1,0) circle(2pt);
  \draw(2,1) circle(2pt);
 \end{scope}

 \begin{scope}[scale=.75,yshift=-3cm]
  \path[use as bounding box] (0,-.3) rectangle (2.4,1.2);
  \thicklines
  \draw(0,0)--(2,0)--(2,1)--(0,1)--(0,0);
  \draw(1,0)--(1,1);
  \fill(0,0) circle(2pt);
  \fill(1,1) circle(2pt);
  \fill(2,0) circle(2pt);
  \fill[white](0,1) circle(2pt);
  \fill[white](1,0) circle(2pt);
  \fill[white](2,1) circle(2pt);
  \draw(0,1) circle(2pt);
  \draw(1,0) circle(2pt);
  \draw(2,1) circle(2pt);
 \end{scope}

 \begin{scope}[scale=.75,xshift=4.5cm,yshift=-3cm]
  \path[use as bounding box] (0,-.3) rectangle (2.4,1.2);
  \draw(0,0)--(0,1);
  \draw(1,0)--(1,1);
  \draw(2,0)--(2,1);
  \draw(0,0)--(1,0);
  \draw(0,1)--(1,1);
  \draw(1,0)--(1,1);
  \draw(1,0) .. controls (3,-1.3) and (3,2.3) .. (1,1);

  \fill(0,0) circle(2pt);
  \fill(1,1) circle(2pt);
  \fill(2,0) circle(2pt);
  \fill[white](0,1) circle(2pt);
  \fill[white](1,0) circle(2pt);
  \fill[white](2,1) circle(2pt);
  \draw(0,1) circle(2pt);
  \draw(1,0) circle(2pt);
  \draw(2,1) circle(2pt);
 \end{scope}

 \begin{scope}[scale=.75,xshift=9cm,yshift=-3cm]
  \path[use as bounding box] (0,-.3) rectangle (2.4,1.2);
  \draw(0,0)--(2,0)--(2,1)--(0,1)--(0,0);
  \draw(1,0)--(1,1);
  \draw(0,0) .. controls (2.5,-1) and (2.75,0) .. (2,1);
  \fill(0,0) circle(2pt);
  \fill(1,1) circle(2pt);
  \fill(2,0) circle(2pt);
  \fill[white](0,1) circle(2pt);
  \fill[white](1,0) circle(2pt);
  \fill[white](2,1) circle(2pt);
  \draw(0,1) circle(2pt);
  \draw(1,0) circle(2pt);
  \draw(2,1) circle(2pt);
 \end{scope}

 \begin{scope}[scale=.75,xshift=13.5cm,yshift=-3cm]
  \path[use as bounding box] (0,-.3) rectangle (2.4,1.2);
  \draw(0,0)--(2,0)--(2,1)--(0,1)--(0,0);
  \draw(1,0)--(1,1);
  \draw(1,0) .. controls (3,-1.3) and (3,2.3) .. (1,1);
  \fill(0,0) circle(2pt);
  \fill(1,1) circle(2pt);
  \fill(2,0) circle(2pt);
  \fill[white](0,1) circle(2pt);
  \fill[white](1,0) circle(2pt);
  \fill[white](2,1) circle(2pt);
  \draw(0,1) circle(2pt);
  \draw(1,0) circle(2pt);
  \draw(2,1) circle(2pt);
 \end{scope}

\end{tikzpicture}

\caption{Bipartite graphs satisfying conditions~(i) and~(ii)}\label{fig_graphs}
\end{figure}

By a direct enumeration we obtain that up to isomorphism there are exactly $8$ bipartite graphs on the sphere that have three vertices in either part, contain no isotopic edges, and satisfy conditions~(i) and~(ii), see Fig.~\ref{fig_graphs}. (By ``isomorphism'' one can mean either an abstract isomorphism of graphs or an isomorphism that is induced by a homeomorphism of the sphere; nevertheless, it is easy to check that for the graphs under consideration these two types of isomorphisms are equivalent.) Consequently, the $\tG$-orbits of the bipartite graphs satisfying conditions~(i) and~(ii) correspond to the possible ways to label the six vertices of each of the graphs on Fig.~\ref{fig_graphs} with labels~$1$, $2$, $3$, $\bar{1}$, $\bar{2}$, $\bar{3}$ so that $\{1,2,3\}$ and~$\{\bar{1},\bar{2},\bar{3}\}$ are the two parts. Such labellings are considered up to graph isomorphisms (induced by homeomorphisms of~$\bS$) and permutations of labels that are induced by the $12$ elements of~$G$. Again, a direct enumeration shows that there are exactly $27$~\ $\tG$-orbits of bipartite graphs satisfying~(i) and~(ii), with representatives listed in Tables~\ref{table_0cells}--\ref{table_3cells}, where graphs in each of the three pairs~$Z_3^{\pm}$, $Z_4^{\pm}$, and~$W_1^{\pm}$ lie in the same $\tG$-orbit.

Note that the above reasoning does not guarantee that all multiarcs listed in Tables~\mbox{\ref{table_0cells}--\ref{table_3cells}} do indeed belong to~$\CL$, i.e., give cells of~$\CB_{0,6}$: we have only proved that all $\tG$-orbits of multiarcs in~$\CL$ are among the $27$ orbits with representatives listed in these tables. However, we obviously have that $X_1$, $X_2$, and~$X_3$ are vertices of~$\CB_{0,6}$, and all vertices of~$\CB_{0,6}$ lie in the union of their $\tG$-orbits. We denote the multiarc corresponding to~$X_i$ by~$L_{X_i}$. One can now check directly that, for each multiarc~$L$ in Tables~\mbox{\ref{table_1cells}--\ref{table_3cells}}, the set of all its submultiarcs belonging to $\tG L_{X_1}\cup\tG L_{X_2}\cup\tG L_{X_3}$ is exactly the set~$\CV(Q)$ specified in the tables (more precisely, the set of the multiarcs corresponding to this set of vertices); moreover, every component of~$L$ is contained in one of these submultiarcs.  From Proposition~\ref{propos_rel_L_char} it follows that $L$ belongs to~$\CL$, that is, indeed defines a cell~$Q$ of~$\CB_{0,6}$, and $\CV(Q)$ is the vertex set of this cell. Moreover, Proposition~\ref{propos_dim_rel} implies that the dimensions of the cells corresponding to multicurves in Tables~\ref{table_1cells}, \ref{table_2cells}, and~\ref{table_3cells} are equal to~$1$, $2$, and~$3$, respectively. Now, we see that every cell~$Q$ from these tables has exactly $\dim Q+1$ vertices, and therefore is a simplex. Hence the facets of~$Q$ correspond to the unions of the multiarcs corresponding to all but one vertices of~$Q$. It can now be checked directly in each of the cases that the boundary~$\partial Q$ is exactly as specified in the tables (including the signs). Thus we obtain assertion~(3).

Let us now prove assertions~(4) and~(5). Let $Q$ be a cell from Tables~\mbox{\ref{table_0cells}--\ref{table_3cells}} and $L$ the corresponding multiarc. We set $\CH=\PMod(\Sigma_{0,6})$; then $\CS<\CH<\G$. Along with the stabilizers~$\CS_Q$ and~$\G_Q$, we consider the stabilizer~$\CH_Q$ of~$Q$ in~$\CH$. By Proposition~\ref{propos_rel_wor}, the group~$\CH$ acts on~$\CB_{0,6}$ without rotations, i.e., each element $h\in \CH_Q$ fixes every component of~$L$. Every connected component of~$\Sigma_{0,6}\setminus L$ is either an open disc, or an open annulus~$\Sigma_{0,2}$, or a pair of pants, i.e., a three-punctured sphere~$\Sigma_{0,3}$. A standard application of the Alexander method (see~\cite[Section~2.3]{FaMa12}) yields the following assertions.
\begin{enumerate}
 \item[(a)] The stabilizer~$\CH_Q$ is free abelian, generated by the Dehn twists about the circuites around the punctures of the connected components of~$\Sigma_{0,6}\setminus L$ homeomorphic to either~$\Sigma_{0,2}$ or~$\Sigma_{0,3}$.
 \item[(b)] The quotient $\tG_Q/\CH_Q$ coincides with the group of all automorphisms~$\varphi$ of the graph~$L$ such that
 \begin{itemize}
 \item $\varphi$ extends to a homeomorphism of the sphere (either preserving or reversing the orientation), which we also denote by~$\varphi$;
 \item $\varphi$ either stabilizes each of the sets~$\{1,2,3\}$ and~$\{\bar{1},\bar{2},\bar{3}\}$ or swaps these two sets;
 \item $\varphi$ takes every set~$\{i,\bar{i}\}$ to some set~$\{j,\bar{j}\}$.
 \end{itemize}
 Moreover, $\G_Q/\CH_Q$ is the kernel of the homomorphism $\tG_Q/\CH_Q\to\{-1,1\}$ induced by the homomorphism~$e_1e_2e_3$ from the definition of the group~$\G$.
\end{enumerate}

Fact~(a) allows us to compute explicitly the group~$\CH_Q$ for each~$Q$ in Tables~\ref{table_0cells}--\ref{table_3cells}. The subgroup $\CS_Q<\CH_Q$ is then easily computed by applying~\eqref{eq_kappa}. The obtained bases of the groups~$\CS_Q$ are given in the tables. Thus we obtain assertion~(4).

To find generators of~$\G_Q$ modulo~$\CS_Q$, we need to write some generators of~$\CH_Q$ modulo~$\CS_Q$ (which are already found as described above) and some generators of~$\G_Q$ modulo~$\CH_Q$. The latter can be found as follows. First, one can check directly in each case that the elements of~$\G$ listed in the last column of Tables~\ref{table_0cells}--\ref{table_3cells} do indeed stabilize the multiarc~$L$ corresponding to~$Q$, and hence belong to~$\G_Q$. Moreover, one can find out for every such element~$g$ whether it preserves or reverses the orientation of~$Q$ by looking on the sign of the permutation of vertices of~$Q$ given by~$g$. Second, the group~$\G_Q/\CH_Q$, i.e., the group of the automorphisms of~$L$ satisfying conditions from~(b), is also easily computed exlicitly in each case. Then one can check that the images of the elements in the last column generate the whole group~$\G_Q/\CH_Q$. Thus we obtain assertion~(5).

To complete the proofs of assertions~(1) and~(2), we need to show that:
\begin{itemize}
 \item for each of the three pairs of cells $Q^{\pm}=Z_3^{\pm}$, or~$Z_4^{\pm}$, or~$W_1^{\pm}$, the cells~$Q^+$ and~$Q^-$ lie in different $\G$-orbits;
 \item for each of the other cells $Q$ in the tables, its $\G$-orbit coincides with its $\tG$-orbit.
\end{itemize}

First, let us prove that~$Z_3^+$ and~$Z_3^-$ lie in different $\CH$-orbits, and $Z_4^+$ and~$\nu Z_4^-$ lie in different $\CH$-orbits. Indeed, in the multiarc~$L_{Z_3^+}$ corresponding to~$Z_3^+$, if we travel from~$1$ to~$\bar{1}$ along the arc with these two endpoints, then the connected component of $\Sigma_{0,6}\setminus L_{Z^+_3}$ lying to the left is adjacent to the punctures~$1$, $\bar{1}$, $3$, and~$\bar{3}$. However, if we travel in the same way in the multiarc~$L_{Z_3^-}$ corresponding to~$Z_3^-$, then the connected component of $\Sigma_{0,6}\setminus L_{Z_3^-}$ lying to the left is adjacent to the punctures~$1$, $\bar{1}$, $2$, and~$\bar{2}$. Consequently, an element from $\CH=\PMod(\Sigma_{0,6})$ cannot take~$L_{Z_3^+}$ to~$L_{Z_3^-}$. Thus $Z_3^+$ and~$Z_3^-$ lie in different $\CH$-orbits. The proof for the pair $\bigl\{Z_4^+,\nu Z_4^-\bigr\}$ is literally the same, but one needs to travel from~$2$ to~$\bar{3}$.

Second, let us show that~$Z_3^+$ and~$Z_3^-$ lie in different $\G$-orbits. Suppose to the contrary that there is an element $g\in\G$ such that $gZ_3^+=Z_3^-$. The punctures~$1$ and~$\bar{1}$ have degree~$3$ in the graphs corresponding to~$Z_3^{\pm}$, while the punctures~$2$, $\bar{2}$, $3$, and~$\bar{3}$ have degree~$2$ in these graphs. Hence $g$ preserves each of the  sets~$\{1,\bar{1}\}$ and~$\{2,\bar{2},3,\bar{3}\}$. We already know that the elements~$T_{\varepsilon_1}^{-1}\nu$ and~$r_1$ lie in the stabilizer~$\G_{Z_3^+}$. Therefore, precomposing~$g$ with an element of~$\G_{Z_3^+}$, we may obtain an element $h\in\G$ that sends~$Z_3^+$ to~$Z_3^-$ and fixes every puncture of~$\Sigma_{0,6}$. Then $e_2(h)=e_3(h)=1$. Since $h\in\G$, we obtain that $h$ preserves the orientation of~$\Sigma_{0,6}$, and hence belongs to $\CH$. We obtain a contradiction, since we have already proved that $Z_3^+$ and~$Z_3^-$ lie in different $\CH$-orbits. Thus $Z_3^+$ and~$Z_3^-$ lie in different $\G$-orbits.

Let us similarly prove that~$Z_4^+$ and~$Z_4^-$ lie in different $\G$-orbits. Suppose to the contrary that they lie in the same $\G$-orbit. Then there is an element $g\in\G$ such that $gZ_4^+=\nu Z_4^-$. (It is more convenient here to take the representative~$\nu Z_4^-$ rather than~$Z_4^-$ itself.) The punctures~$2$ and~$\bar{3}$ have degree~$3$ in the graphs corresponding to~$Z_4^+$ and~$\nu Z_4^-$, while the punctures~$1$, $\bar{1}$, $\bar{2}$, and~$3$ have degree~$2$ in these graphs. Hence $g$ preserves the set~$\{2,\bar{3}\}$. But $r_1\in\G_{Z_4^+}$. Therefore, precomposing~$g$ with~$r_1$ if necessary, we may obtain an element $h\in\G$ that sends~$Z_4^+$ to~$\nu Z_4^-$ and fixes each of the punctures~$2$ and~$\bar{3}$. Then $h$ fixes every puncture of~$\Sigma_{0,6}$, and therefore $e_2(h)=e_3(h)=1$. Since $h\in\G$, we obtain that $h$ preserves the orientation of~$\Sigma_{0,6}$, and hence belongs to $\CH$. We obtain a contradiction, since we have already proved that $Z_4^+$ and~$\nu Z_4^-$ lie in different $\CH$-orbits. Thus $Z_4^+$ and~$Z_4^-$ lie in different $\G$-orbits.

The cell $W_1^+$ contains~$Z_3^+$ and contains no cells from the $\G$-orbit of~$Z_3^-$, while $W_1^-$ contains~$Z_3^-$ and contains no cells from the $\G$-orbit of~$Z_3^+$. Thus the cells $W_1^+$ and~$W_1^-$ lie in different $\G$-orbits as well.

Finally, let us prove that $\tG Q=\G Q$ for all cells~$Q$ in Tables~\mbox{\ref{table_0cells}--\ref{table_3cells}}, except~$Z_3^{\pm}$, $Z_4^{\pm}$, and~$W_1^{\pm}$. First, assume that $Q=W_i$ with $i\in\{2,\ldots,6\}$. Since $\G$ has index~$2$ in~$\tG$, it suffices to find an element $f\in \tG$ such that $f(Q)=Q$ and $f\notin \G$, that is, $e_1(f)e_2(f)e_3(f)=-1$.

Suppose that $i\in\{3,\ldots,6\}$. Then the multiarc~$L_{W_i}$ corresponding to~$W_i$ contains a (unique) pair of arcs~$\{\beta_1,\beta_2\}$ with the same pair of endpoints. We define the required mapping class~$f$ as follows: $f$ fixes every puncture of~$\Sigma_{0,6}$, swaps the two edges~$\beta_1$ and~$\beta_2$, and fixes every other edge of~$L_{W_i}$; then $f$ extends to an orientation-reversing homeomorphism of~$\Sigma_{0,6}$. We have $e_1(f)=-1$ and $e_2(f)=e_3(f)=1$. Hence $f\in\tG_{W_i}$ but $f\notin\G_{W_i}$. Thus $\tG W_i=\G W_i$ for $i=3,\ldots,6$.

The multiarc~$L_{W_2}$ corresponding to~$W_2$ contains no pairs of edges with the same pair of endpoints. Let $f$ fix each of the punctures~$1$ and~$\bar{1}$ and swap the pairs of punctures $2\leftrightarrow 3$ and $\bar{2}\leftrightarrow \bar{3}$. Then $f$ can be extended (uniquely up to isotopy) first to the edges of~$L_{W_2}$ and then to the connected components of~$\Sigma_{0,6}\setminus L_{W_2}$. The obtained mapping class~$f$ is orientation-preserving. Hence, we have $e_1(f)=e_2(f)=1$ and $e_3(f)=-1$, and therefore $f\in\tG_{W_2}$ but $f\notin\G_{W_2}$. Thus $\tG W_2=\G W_2$.

The cell $W_2$ contains the cells $r_3Z_2$, and $\omega^2Z_2$ and cointains no other cells in the $\tG$-orbits of~$Z_2$. Hence equality $\tG Z_2=\G Z_2$ follows immediately from $\tG W_2=\G W_2$. We abbreviate this by writing $W_2\Rightarrow Z_2$. The same argument applies to many other pairs of cells. In particular, we have $W_5\Rightarrow Z_1$, $W_2\Rightarrow Z_2$, $W_3\Rightarrow Z_5$, $W_4\Rightarrow Z_6$, $W_5\Rightarrow Z_7$, $W_5\Rightarrow Z_8$, $W_6\Rightarrow Z_9$, $Z_1\Rightarrow Y_1$, $Z_1\Rightarrow Y_2$, $Z_1\Rightarrow Y_3$, $Z_2\Rightarrow Y_4$, $W_3\Rightarrow Y_5$, $Z_5\Rightarrow Y_6$, $Z_6\Rightarrow Y_7$, $Z_7\Rightarrow Y_8$, $Z_9\Rightarrow Y_9$, $Y_1\Rightarrow X_1$, $Y_1\Rightarrow X_2$, and $Y_2\Rightarrow X_3$. Thus $\tG Q=\G Q$ for all cells~$Q$ in Tables~\mbox{\ref{table_0cells}--\ref{table_3cells}}, except~$Z_3^{\pm}$, $Z_4^{\pm}$, and~$W_1^{\pm}$. This completes the proof of the proposition.
\end{proof}

\section{The group~$E^1_{0,2}$}\label{section_E102}

In this section we start studying the spectral sequence~\eqref{eq_SSS} associated with the action of~$\CS$ on~$\CB_{0,6}$. Recall that this spectral sequence~$E^*_{*,*}$ is a spectral sequence of left $\Z[\Gamma]$-modules.

\begin{propos}\label{propos_E102}
 The group $E^1_{0,2}$ is a left $\Z[\Gamma]$-module with two generators
 \begin{align*}
  A_1&=\iota_{X_1}\bigl(\CA(T_{\delta_2},T_{\delta_3})\bigr),\\
  A_2&=\iota_{T_{31}X_2}\left(\CA\bigl(T_{\delta_1},T_{\gamma_{2\bar{3}}}T_{\tgamma_{\bar{2}3}}^{-1}\bigr)\right),
 \end{align*}
  and a complete set of defining relations
  \begin{gather}
   (r_1+1)A_1=0,\qquad (t_1t_2t_3\nu-1)A_1=0,\label{eq_rel_AE1}\\
   (t_1-1)A_2=0,\qquad (r_1-1)A_2=0,\qquad (t_1t_2t_3\nu+1)A_2=0.\label{eq_rel_BE1}
  \end{gather}
\end{propos}

\begin{remark}
 We have $\delta_1=T_{31}(\delta_1')$, and hence $T_{\delta_1}=T_{31}T_{\delta_1'}T_{31}^{-1}$.
 So instead of~$A_2$ we could take for a generator the element
 $$
 t_2^{-1}A_2=\iota_{X_2}\left(T_{\delta_1'},T_{2\bar{3}}T_{\bar{2}3}^{-1}\right).
 $$
Nevertheless, our choice of the generator is more convenient for further computations.
\end{remark}

\begin{proof}[Proof of Proposition~\ref{propos_E102}]
 First, let us prove that relations~\eqref{eq_rel_AE1} and~\eqref{eq_rel_BE1} hold in~$E^1_{0,2}$. The mapping class~$r_1$ fixes~$X_1$ and swaps~$T_{\delta_2}$ and~$T_{\delta_3}$; hence, $r_1A_1=-A_1$. The mapping class~$T_{\varepsilon}\nu$ fixes~$X_1$ and stabilizes each of the curves~$\delta_2$ and~$\delta_3$ but reverses the orientation of the surface. Hence the conjugation by~$T_{\varepsilon}\nu$ sends~$T_{\delta_2}$ and~$T_{\delta_3}$ to~$T_{\delta_2}^{-1}$ and~$T_{\delta_3}^{-1}$, respectively. Since $\CA\bigl(T_{\delta_2}^{-1},T_{\delta_3}^{-1}\bigr)=\CA(T_{\delta_2},T_{\delta_3})$ and $\kappa(T_{\varepsilon})=t_1t_2t_3$, we obtain that $t_1t_2t_3\nu A_1=A_1$.

 Each of the mapping classes~$T_{2\bar{3}}$ and~$T_{\talpha_3}r_1$  fixes~$T_{31}X_2$ and stabilizes each of the curves~$T_{\delta_1}$, $\gamma_{2\bar{3}}$, and~$\tgamma_{\bar{2}3}$. Since $\kappa(T_{2\bar{3}})=t_1^{-1}$ and $\kappa(T_{\talpha_3})=1$, we obtain that $t_1A_2=A_2$ and $r_1A_2=A_2$. The mapping class~$T_{\delta_2}T_{\varepsilon}\nu$ fixes~$T_{31}X_2$, stabilizes~$\delta_1$, swaps~$\gamma_{2\bar{3}}$ and~$\tgamma_{\bar{2}3}$, and reverses the orientation of the surface. Therefore the conjugation by~$T_{\delta_2}T_{\varepsilon}\nu$ sends~$T_{\delta_1}$ to~$T_{\delta_1}^{-1}$ and fixes~$T_{\gamma_{2\bar{3}}}T_{\tgamma_{\bar{2}3}}^{-1}$.
But $\kappa(T_{\delta_2}T_{\varepsilon})=t_1t_2t_3$. Hence $t_1t_2t_3\nu A_2=-A_2$.

 Let us now prove that $A_1$ and~$A_2$ generate the $\Z[\Gamma]$-module~$E^1_{0,2}$, and~\eqref{eq_rel_AE1} and~\eqref{eq_rel_BE1} form a complete set of relations. Choose some mapping classes~$\tau_1,\tau_2,\tau_3\in\G$ with $\kappa(\tau_i)=t_i$. From Proposition~\ref{propos_tables} it follows that for the representatives of $\CS$-orbits of vertices of~$\CB_{0,6}$ one can take all vertices $\tau_1^{k_1}\tau_2^{k_2}\tau_3^{k_3}X_1$ with $k_1,k_2,k_3\in\Z$, all vertices $\omega^l\tau_2^{k_2}\tau_3^{k_3}T_{31}X_2$ with $l\in\{0,1,2\}$  and~$k_2,k_3\in\Z$, and some vertices of the form~$gX_3$, where $g\in\G$. Moreover, we have
 \begin{itemize}
  \item $H_2(\CS_{X_1})\cong\Z^3$ with basis~$\CA(T_{\delta_2},T_{\delta_3})$, $\CA(T_{\delta_3},T_{\delta_1})$, $\CA(T_{\delta_1},T_{\delta_2})$,
  \item $H_2(\CS_{T_{31}X_2})\cong\Z$ with generator~$\CA\left(T_{\delta_1},T_{\gamma_{2\bar{3}}}T_{\tgamma_{\bar{2}3}}^{-1}\right)$, and
  \item  $H_2(\CS_{X_3})=0$.
 \end{itemize}
   Hence the following elements form a $\Z$-basis of the group~$E_{0,2}^1$: the elements $\omega^lt_1^{k_1}t_2^{k_2}t_3^{k_3}A_1$ with $l\in\{0,1,2\}$ and $k_1,k_2,k_3\in\Z$, and the elements $\omega^lt_2^{k_2}t_3^{k_3}A_2$ with $l\in\{0,1,2\}$  and~$k_2,k_3\in\Z$. The proposition follows.
\end{proof}

The group~$E^{\infty}_{0,2}$, which embeds as a subgroup in~$H_2(\CS)$, is the quotient group of~$E^1_{0,2}$ by the images of the differentials~$d^1$, $d^2$, and~$d^3$. Factoring by the images of these differentials leads to additional relations on the generators~$A_1$ and~$A_2$. However, in what follows we will only need the relations coming from the differential~$d^3$. Therefore we will not explicitly compute the differentials~$d^1$ and~$d^2$ (although this is not difficult to do). It will suffice for us to use the following statement, which is an immediate consequence of Proposition~\ref{propos_E102}. With some abuse of notation, we continue to denote the images of~$A_1$ and~$A_2$ in~$E^r_{0,2}$ for $r\ge 2$ by~$A_1$ and~$A_2$, respectively.

\begin{cor}\label{cor_Er02}
For each $r\ge 2$, the $\Z[\Gamma]$-module $E^r_{0,2}$ is generated by~$A_1$ and~$A_2$. These generators satisfy relations~\eqref{eq_rel_AE1} and~\eqref{eq_rel_BE1} (and possibly some other relations).
\end{cor}

\begin{remark}
 In fact, one can show that the differential $d^1\colon E^1_{1,2}\to E^1_{0,2}$ is trivial and the image of $d^2\colon E^2_{2,1}\to E^2_{0,2}$ is generated by a single element $(t_1-1)\bigl(\omega+\omega^2)A_1$. We will not prove this, since it will not be needed for the proof of Theorem~\ref{theorem_3}.
\end{remark}

Let us now consider the following abelian cycles in~$H_2(\CS)$:
\begin{align}
 \CA_1&=\CA(T_{\delta_2},T_{\delta_3}),\label{eq_CA1}\\
 \CA_2&=\CA\left(T_{\delta_1},T_{\gamma_{2\bar{3}}}T_{\tgamma_{\bar{2}3}}^{-1}\right).\label{eq_CA2}
\end{align}
Recall that
$$
0\subseteq \CF_{0,2}\subseteq \CF_{1,1}\subseteq \CF_{2,0}=H_2(\CS)
$$
is the filtration associated with the action of~$\CS$ on~$\CB_{0,6}$, i.e., such that $\CF_{p,q}/\CF_{p-1,q+1}=E^{\infty}_{p,q}$. The following proposition is straightforward.

\begin{propos}\label{propos_A1A2}
 The classes~$\CA_1$ and~$\CA_2$ lie in~$\CF_{0,2}$. The natural isomorphism $\CF_{0,2}\cong E^{\infty}_{0,2}$ sends them to~$A_1$ and~$A_2$, respectively.
\end{propos}

\begin{cor}\label{cor_F02}
 The group~$\CF_{0,2}$ is generated as a $\Lambda$-module by the $6$ abelian cycles~$\CA_1$, $\omega\CA_1$, $\omega^2\CA_1$, $\CA_2$, $\omega\CA_2$, and~$\omega^2\CA_2$.
\end{cor}

\begin{proof}
 By Corollary~\ref{cor_Er02}, the group $\CF_{0,2}$ is generated as a $\Z[\Gamma]$-module by~$\CA_1$ and~$\CA_2$. Combining this with relations~\eqref{eq_rel_AE1} and~\eqref{eq_rel_BE1}, we get the required assertion.
\end{proof}

\section{First row of the spectral sequence}\label{section_first_row}
The next proposition follows from the data in Tables~\mbox{\ref{table_0cells}--\ref{table_3cells}} in exactly the same way as Proposition~\ref{propos_E102} does.

\begin{propos}\label{propos_E1p1}
 \textnormal{(1)} The group $E^1_{0,1}$ is a left $\Z[\Gamma]$-module with three generators
 \begin{align*}
  \xi_1&=\iota_{X_1}([T_{\delta_1}]),\\
  \xi_2&=\iota_{T_{31}X_2}([T_{\delta_1}]),\\
  \xi_3&=\iota_{X_2}\bigl(\bigl[T_{2\bar{3}}T_{\bar{2}3}^{-1}\bigr]\bigr),
 \end{align*}
  and a complete set of defining relations
  \begin{gather}
   (r_1-1)\xi_1=0,\qquad (t_1t_2t_3\nu+1)\xi_1=0,\label{eq_E1xi_rel1}\\
   (t_1-1)\xi_2=0,\qquad\quad\quad\, (r_1-1)\xi_2=0,\!\qquad \bigl(t_1t_2t_3\nu+1\bigr)\xi_2=0,\label{eq_E1xi_rel2}\\
   (t_1-1)\xi_3=0,\qquad \bigl(t_2^{-1}t_3r_1-1\bigr)\xi_3=0,\qquad \bigl(t_2^{-1}t_3\nu-1\bigr)\xi_3=0.\label{eq_E1xi_rel3}
  \end{gather}

  \textnormal{(2)} The group $E^1_{1,1}$ is a left $\Z[\Gamma]$-module with four generators
 \begin{align*}
  \eta_1&=\iota_{Y_1}([T_{\delta_1}]),\\
  \eta_2&=\iota_{\omega Y_6}([T_{\delta_1}]),\\
  \eta_3&=\iota_{Y_7}\bigl(\bigl[T_{2\bar{3}}T_{\bar{2}3}^{-1}\bigr]\bigr),\\
  \eta_4&=\iota_{Y_8}([T_{\delta_1}]).
 \end{align*}
  and a complete set of defining relations
  \begin{gather}
   (t_1t_2t_3\nu+1)\eta_1=0,\label{eq_E1eta_rel1}\\
   (t_1t_2t_3\nu+1)\eta_2=0,\label{eq_E1eta_rel2}\\
   (t_1-1)\eta_3=0,\qquad (r_1+1)\eta_3=0,\qquad (\nu+1)\eta_3=0,\label{eq_E1eta_rel3}\\
   (t_1-1)\eta_4=0,\qquad  (r_1+1)\eta_4=0.\label{eq_E1eta_rel4}
  \end{gather}

   \textnormal{(3)} The group $E^1_{2,1}$ is a left $\Z[\Gamma]$-module with two generators
 \begin{align*}
  \zeta_1&=\iota_{\omega Z_5}([T_{\delta_1}]),\\
  \zeta_2&=\iota_{Z_7}([T_{\delta_1}]).
 \end{align*}
  and a complete set of defining relations
  \begin{gather*}
   (t_1t_2t_3\nu+1)\zeta_1=0,\\
   (r_1+1)\zeta_2=0.
   \end{gather*}
\end{propos}

\begin{remark}
 Again, as with~$A_2$ in $E^1_{0,2}$, we could take for generators the classes $t_2^{-1}\xi_2=\iota_{X_2}([T_{\delta_1'}])$, $\omega^2\eta_2=\iota_{Y_6}([T_{\delta_3}])$, and $\omega^2\zeta_1=\iota_{Z_5}([T_{\delta_3}])$ instead of~$\xi_2$, $\eta_2$, and~$\zeta_1$, respectively. However, our choice of generators is more convenient for further computations.
\end{remark}

The next proposition follows immediately from Fact~\ref{fact_dr} and the formulas for the boundary maps $\partial\colon C_k(\CB_{0,6})\to C_{k-1}(\CB_{0,6})$ for $k=1,2$ from Tables~\ref{table_1cells} and~\ref{table_2cells}.

\begin{propos}\label{propos_d1_1row}
The differentials $d^1\colon E^1_{1,1}\to E^1_{0,1}$ and $d^1\colon E^1_{2,1}\to E^1_{1,1}$ are given  by
 \begingroup
 \allowdisplaybreaks
 \begin{align}
  d^1\eta_1&=\xi_2-\xi_1,\notag\\
  d^1\eta_2&=0,\notag\\
  d^1\eta_3&=(r_1-1)\xi_3=\bigl(t_2t_3^{-1}-1)\xi_3,\notag\\
  d^1\eta_4&=0,\notag\\
  d^1\zeta_1&=(1-r_1)\eta_1+\eta_2,\notag\\
  d^1\zeta_2&=(1-r_1)\eta_1+\eta_4.\label{eq_d1zeta2}
 \end{align}
 \endgroup
 \end{propos}

 The group~$E^2_{0,1}$ is the quotient group of~$E^1_{0,1}$ by the image of the differential~$d^1$. With some abuse of notation, we denote the image of~$\xi_i$ in~$E^2_{0,1}$ again by~$\xi_i$. Similarly, we denote the classes in~$E^2_{1,1}$ represented by the cycles~$\eta_2$ and~$\eta_4$ again by~$\eta_2$ and~$\eta_4$, respectively.

 \begin{propos}\label{propos_E201}
  The group $E^2_{0,1}$ is a left $\Z[\Gamma]$-module with two generators
 $\xi_1$ and~$\xi_3$, and a complete set of defining relations
  \begin{gather*}
   (t_1-1)\xi_1=0,\qquad (r_1-1)\xi_1=0,\qquad (t_2t_3\nu+1)\xi_1=0,\\
   (t_1-1)\xi_3=0,\qquad (t_2-t_3)\xi_3=0,\qquad (r_1-1)\xi_3=0,\qquad (\nu-1)\xi_3=0.
  \end{gather*}
  Moreover, $\xi_2=\xi_1$ in~$E^2_{0,1}$.
 \end{propos}

 \begin{proof}
  By Proposition~\ref{propos_d1_1row}, $E^2_{0,1}$ is the quotient of the $\Z[\Gamma]$-module~$E^1_{0,1}$ by the submodule generated by $\xi_2-\xi_1$ and  $(t_2-t_3)\xi_3$. Substituting $\xi_2=\xi_1$ to relations~\eqref{eq_E1xi_rel1}--\eqref{eq_E1xi_rel3} and adding the relation $(t_2-t_3)\xi_3=0$, we get precisely the listed set of relations.
 \end{proof}

\begin{propos}\label{propos_E211}
 The group $E^2_{1,1}$ is a left $\Z[\Gamma]$-module with a single generator
 $\eta_4$, and a complete set of defining relations
 \begin{equation}\label{eq_rel_eta2E2}
  (t_1-1)\eta_4=0,\qquad  (r_1+1)\eta_4=0,\qquad (t_2t_3\nu+1)\eta_4=0.
 \end{equation}
 Moreover, $\eta_2=\eta_4$ in~$E^2_{1,1}$.
\end{propos}

\begin{proof}
Since $d^1(\zeta_1-\zeta_2)=\eta_2-\eta_4$, it follows that $\eta_2=\eta_4$ in~$E^2_{1,1}$.

To prove that $E^2_{1,1}$ is generated by~$\eta_4$ as a $\Z[\Gamma]$-module, it suffices to show that any element $\eta=c_1\eta_1+c_3\eta_3$ (with $c_1,c_3\in\Z[\Gamma]$) lying in the kernel of~$d^1$ is $d^1$-homologous to some element $c_4\eta_4$. We have
$$
0=d^1(c_1\eta_1+c_3\eta_3)=-c_1\xi_1+c_1\xi_2+c_3\bigl(t_2t_3^{-1}-1\bigr)\xi_3.
$$
Since $E^1_{0,1}$ is the direct sum of the three $\Z[\Gamma]$-modules generated by~$\xi_1$, $\xi_2$, and~$\xi_3$, respectively, we obtain that $c_1\xi_1=0$, $c_1\xi_2=0$, and $c_3\bigl(t_2t_3^{-1}-1\bigr)\xi_3=0$.

Since relations~\eqref{eq_E1eta_rel1} and~\eqref{eq_E1eta_rel3} hold in~$E^1_{0,2}$, we may assume that
\begin{align}
 c_1&=\sum_{l=0}^2\omega^l\bigl(p_l+\tilde{p}_lr_1\bigr),&&p_l,\tilde{p}_l\in \Lambda,\label{eq_c1}\\
 c_3&=\sum_{l=0}^2\omega^l  q_l,&&q_l\in \Z\bigl[t_2,t_2^{-1},t_3,t_3^{-1}\bigr]
\end{align}
The annihilator of $\xi_1\in E^1_{0,1}$ is the left ideal of~$\Z[\Gamma]$ generated by $r_1-1$ and $t_1t_2t_3\nu+1$. Consequently, an element of the form~\eqref{eq_c1} annihilates~$\xi_1$ if and only if $\tilde{p}_l=-p_l$ for all~$l$, that is,
$$
c_1=\sum_{l=0}^2\omega^lp_l(1-r_1).
$$

The equality $c_3\bigl(t_2t_3^{-1}-1\bigr)\xi_3=0$ reads as
$$
\sum_{l=0}^2\omega^l \bigl(t_2t_3^{-1}-1\bigr) q_l\xi_3=0.
$$
But, by~\eqref{eq_E1xi_rel3}, the annihilator of the element $\xi_3\in E^1_{0,1}$ is the left $\Z[\Gamma]$-module generated by~$t_1-1$, $t_2^{-1}t_3r_1-1$, and~$t_2^{-1}t_3\nu-1$. Consequently, the elements $\omega^lt_2^{k_2}t_3^{k_3}\xi_3$ (where $l\in\{0,1,2\}$, $k_2,k_3\in\Z$) are linearly independent in~$E^1_{0,1}$. Therefore $q_0=q_1=q_2=0$, i.e., $c_3=0$. Hence
$$
\eta=\sum_{l=0}^2\omega^lp_l(1-r_1)\eta_1.
$$
By~\eqref{eq_d1zeta2}, the element~$\eta$  is $d^1$-homologous to $-\sum_{l=0}^2\omega^lp_l\eta_4$. Thus $E^2_{1,1}$ is generated as a $\Z[\Gamma]$-module by~$\eta_4$.

Relations~\eqref{eq_rel_eta2E2} follow immediately from~\eqref{eq_E1eta_rel2} and~\eqref{eq_E1eta_rel4}. Finally, we need to prove that relations~\eqref{eq_rel_eta2E2} form a complete set of defining relations for the $\Z[\Gamma]$-module~$E^2_{1,1}$. Indeed, assume that $c\eta_4=0$ in~$E^2_{1,1}$ for some $c\in\Z[\Gamma]$. Then the following holds in~$E^1_{1,1}$ for some $b_1,b_2\in\Z[\Gamma]$:
$$
c\eta_4=d^1(b_1\zeta_1+b_2\zeta_2)=(b_1+b_2)(1-r_1)\eta_1+b_1\eta_2+b_2\eta_4.
$$
Therefore
\begin{align*}
 (b_1+b_2)(1-r_1)\eta_1&=0,\\
 b_1\eta_2&=0,\\
 (c-b_2)\eta_4&=0.
\end{align*}
Since the annihilator of~$\eta_1$ is generated by $t_1t_2t_3\nu+1$ and the elements~$t_1t_2t_3\nu$ and~$r_1$ have order~$2$ and commute, it follows that the annihilator of~$(1-r_1)\eta_1$ is generated by~$r_1+1$ and $t_1t_2t_3\nu+1$. Hence $b_1+b_2$ lies in the ideal $$\Z[\Gamma](r_1+1)+\Z[\Gamma](t_1t_2t_3\nu+1).$$
Further, the element~$b_1$ lies in the annihilator of~$\eta_2$, i.e., in the ideal
$\Z[\Gamma](t_1t_2t_3\nu+1)$.
Finally,
the element~$c-b_2$ lies in the annihilator of~$\eta_4$, i.e., in the ideal
$$\Z[\Gamma](t_1-1)+\Z[\Gamma](r_1+1).$$
It follows that $c$ lies in the ideal
$$\Z[\Gamma](t_1-1)+\Z[\Gamma](r_1+1)+\Z[\Gamma](t_1t_2t_3\nu+1).$$
Thus \eqref{eq_rel_eta2E2} is a complete set of defining relations for the $\Z[\Gamma]$-module~$E^2_{1,1}$.
\end{proof}

\begin{cor}\label{cor_Einfty11}
The group $E^{\infty}_{1,1}$ is generated as a $\Z[\Gamma]$-module by the element~$\eta_4$, which satisfies relations~\eqref{eq_rel_eta2E2} (and possibly some other relations). Hence, $E^{\infty}_{1,1}$ is generated as a $\Lambda$-module by the three elements~$\eta_4$, $\omega\eta_4$, and~$\omega^2\eta_4$.
\end{cor}

In addition to the abelian cycles~$\CA_1$ and~$\CA_2$ defined in Section~\ref{section_E102}, consider the abelian cycle
\begin{equation}\label{eq_CA3}
 \CA_3=\CA(T_{\delta_1},T_{\talpha_3})\in H_2(\CS).
\end{equation}
Recall that we denote by~$\{\CF_{p,q}\}$ the filtration in~$H_*(\CS)$ associated with the action of~$\CS$ on~$\CB_{0,6}$.

\begin{propos}\label{propos_A3}
 The abelian cycle~$\CA_3$ belongs to the group~$\CF_{1,1}$, and the projection $\CF_{1,1}\twoheadrightarrow E^{\infty}_{1,1}$ sends~$\CA_3$ to~$\eta_4$. Moreover, $\CA_3$ is represented by the following cycle in~$\Tot_2\bigl(C_*(\CB_{0,6})\otimes_{\CS}\CR_*\bigr)$:
  \begin{equation*}
   \fA_3=(T_{31}X_2)\otimes\left(\bigl[T_{\talpha_3}\big|T_{\delta_1}^{-1}\bigr]-\bigl[T_{\delta_1}^{-1}\big|T_{\talpha_3}\bigr]\right)-Y_8\otimes\bigl[T_{\delta_1}^{-1}\bigr].
  \end{equation*}
 \end{propos}

 \begin{proof}
Set $g=T_{\delta_1}^{-1}$ and $h=T_{\talpha_3}^{-1}$. We have $\partial Y_8=\bigl(T_{3\bar{1}}^{-1}-T_{31}\bigr)X_2$. By relation~\eqref{eq_131}, we have $T_{3\bar{1}}^{-1}=T_{\talpha_3}^{-1}T_{31}T_{\delta_1'}$. Since $T_{\delta_1'}$ fixes~$X_2$, we obtain that $\partial Y_8=(h-1)T_{31}X_2$. Therefore, the elements~$g$ and~$h$, the $1$-cell~$P=Y_8$, and the vertex $v=T_{31}X_2$ satisfy all conditions of Lemma~\ref{lem_ab_F11}. The proposition now follows from this lemma.
 \end{proof}

\begin{cor}\label{cor_F11}
 The group $\CF_{1,1}$ is generated as a $\Lambda$-module by the $9$ abelian cycles~$\omega^l\CA_m$, where $l\in\{0,1,2\}$ and $m\in\{1,2,3\}$.
\end{cor}

\begin{proof}
 We have a short exact sequence
 $$
 0\longrightarrow \CF_{0,2}\longrightarrow\CF_{1,1}\longrightarrow E^{\infty}_{1,1}\longrightarrow 0.
 $$
 Therefore the assertion of the corollary follows from Corollaries~\ref{cor_F02} and~\ref{cor_Einfty11}, and Proposition~\ref{propos_A3}.
\end{proof}

\section{Groups $E^1_{p,0}$ and differential~$d^1$}\label{section_zero_row}

Recall that the $0$th row~$E^1_{*,0}$ with differential~$d^1$ is canonically isomorphic to the cellular chain complex~$C_*(\CB_{0,6}/\CS)$ with the usual cellular differential~$\partial$. For the cells~$X_1,\ldots,W_6$, we denote their images under the quotient map $\CB_{0,6}\to \CB_{0,6}/\CS$ by the corresponding lowercase letters~$x_1,\ldots,w_6$. The data in Tables~\ref{table_0cells}--\ref{table_3cells} immediately yield the following two propositions.

\begingroup
 \allowdisplaybreaks
\begin{propos}\label{propos_E1p0}
 \textnormal{(1)} The group $E^1_{1,0}$ is a left\/ $\Z[\Gamma]$-module with $9$ generators $y_1,\ldots,y_9$ and a complete set of defining relations
 \begin{gather*}
  (t_1t_2t_3\nu - 1)y_1=0,\\
  (t_1-1)y_2=0,\\
  (\omega-1)y_3=0,\qquad (r_1-1)y_3=0,\\
  (\omega-1)y_4=0,\qquad (r_1-1)y_4=0,\qquad (\nu+1)y_4=0,\\
  (r_1+1)y_5=0,\qquad \bigl(t_1^{-1}t_2t_3\nu-1\bigr)y_5=0, \\
  (r_1-1)y_6=0,\qquad (t_1t_2t_3\nu-1)y_6=0,\\
  (t_1-1)y_7=0,\qquad (r_1+1)y_7=0,\qquad (\nu+1)y_7=0,\\
  (t_1-1)y_8=0,\qquad (r_1+1)y_8=0,\\
  (t_2-1)y_9=0,\qquad (t_3-1)y_9=0,\qquad (r_1-1)y_9=0.
 \end{gather*}

\textnormal{(2)} The group $E^1_{2,0}$ is a left\/ $\Z[\Gamma]$-module with $11$ generators $z_1,z_2,z_3^{\pm},z_4^{\pm},z_5,\ldots,z_9$ and a complete set of defining relations
 \begin{gather*}
  (s_1+1)z_2=0,\\
  (r_1+1)z_3^{\pm}=0,\qquad\bigl(t_1^{-1}t_2t_3\nu-1\bigr)z_3^{\pm}=0,\\
  (r_1+1)z_4^{\pm}=0,\\
  (t_1t_2t_3\nu-1)z_5=0,\\
  (t_1-1)z_6=0,\qquad (r_1+1)z_6=0,\\
  (r_1+1)z_7=0,\\
  (t_1-1)z_8=0,\qquad (r_1+1)z_8=0,\\
  (t_2-1)z_9=0.
 \end{gather*}

 \textnormal{(3)} The group $E^1_{3,0}$ is a left\/ $\Z[\Gamma]$-module with $7$ generators $w_1^{\pm},w_2,\ldots,w_6$ and a complete set of defining relations
 \begin{gather*}
  (r_1+1)w_1^{\pm}=0,\\
  (r_1+1)w_i=0,\qquad i=2,\ldots,6,\\
  (t_1t_2t_3\nu-1)w_3=0,\\
  (\nu-1)w_4=0.
 \end{gather*}
\end{propos}
\endgroup

\begin{propos}\label{propos_d1_0row}
 The differentials $d^1\colon E^1_{2,0}\to E^1_{1,0}$ and $d^1\colon E^1_{3,0}\to E^1_{2,0}$ are given by
 \begingroup
 \allowdisplaybreaks
 \begin{align*}
  d^1z_1&=y_1+t_3r_1y_2-y_3,\\
  d^1z_2&=(1-s_1)y_2+y_4,\\
  d^1z_3^+&=t_1^{-1}(r_2-\omega)y_1+y_5,\\
  d^1z_3^-&=t_2t_3\bigl(\omega^2-r_3\bigr)\nu y_1+y_5=t_1^{-1}\bigl(\omega^2-r_3\bigr)y_1+y_5,\\
  d^1z_4^+&=\bigl(r_3-\omega^2\bigr)y_2+y_5,\\
  d^1z_4^-&=(\omega-r_2)\nu y_2+y_5,\\
  d^1z_5&=\bigl(\omega^2-r_2\bigr)y_1+y_6,\\
  d^1z_6&=(r_1-1)y_2+y_7,\\
  d^1z_7&=(1-r_1)y_1+y_8,\\
  d^1z_8&=(1-r_1)t_2y_2+y_8,\\
  d^1z_9&=\bigl(\omega-t_1t_3^{-1}r_3\bigr)y_2+y_9,\\
  d^1w_1^+&= t_1^{-1}(r_2-\omega)z_1-z_3^+ +z_4^+,\\
  d^1w_1^-&= t_2t_3\bigl(\omega^2-r_3)\nu z_1-z_3^- +z_4^-,\\
  d^1w_2&=\bigl(-\omega\nu-\omega^2\bigr)z_2-z_4^++z_4^-,\\
  d^1w_3&=-t_1z_3^+ +t_1z_3^- + (r_1-1)z_5,\\
  d^1w_4&=-(\nu+1)(z_2+z_6),\\
  d^1w_5&=(1-r_1)z_1-z_7+z_8,\\
  d^1w_6&=-t_1t_2^{-1}t_3^{-1}z_4^++\nu z_4^-+(r_1-1)z_9.
 \end{align*}
 \endgroup
\end{propos}

\begin{remark}
 We could similarly write formulas for the differential $d^1\colon E^1_{1,0}\to E^1_{0,0}$, but we do not need them.
\end{remark}

\section{The group~$E^2_{2,0}$}\label{section_E220}

Let $\Delta<\Gamma$ be the subgroup generated by~$t_1$, $t_2$, $t_3$, and~$\omega$. The subgroup~$\Delta$ is normal of index~$4$ and $\Gamma/\Delta\cong C_2\times C_2$, where $C_2$ denotes the cyclic group of order~$2$. Moreover, $1$, $r_1$, $\nu$, and $s_1=r_1\nu$ is a set of representatives of cosets of~$\Delta$ in~$\Gamma$. By some technical reasons it is sometimes more convenient to study~$E^2_{2,0}$ as a $\Z[\Delta]$-module rather than as a $\Z[\Gamma]$-module.   The goal of this section is to find generators of~$E^2_{2,0}$ as a $\Z[\Delta]$-module, and also as an abelian group.

Let $I$ be the image of the differential $d^1\colon E^1_{3,0}\to E^1_{2,0}$. We write~$\Iequiv$ to denote equality modulo~$I$. For an element~$a$ of an arbitrary left $\Z[\Gamma]$-module, we denote by~$\Ann(a)$ the annihilator of~$a$. Let us start with two auxiliary lemmas.

\begin{lem}\label{lem_z2z9}
 Let $a\in E^1_{2,0}$ be an element satisfying $d^1a=0$. Then $$a\Iequiv c_2z_2+c_9z_9$$ for some $c_2,c_9\in\Z[\Gamma]$.
\end{lem}

\begin{proof}
 We write
 $$
 a=\sum_{i\in\{1,2,5,6,7,8,9\}}c_iz_i+c_3^+z_3^++c_3^-z_3^-+c_4^+z_4^++c_4^-z_4^-,
 $$
 where the coefficients belong to~$\Z[\Gamma]$. Our aim is to show that, by adding to~$a$ elements from~$I$, we can kill all coefficients except~$c_2$ and~$c_9$.

 First, adding to~$a$ the element~$d^1c_7w_5$, we kill the coefficient~$c_7$. When we have already achieved that $c_7=0$, we obtain that $y_8$ enters~$d^1a$ with coefficient~$c_8$. Since $d^1a=0$, we have $c_8y_8=0$, that is, $c_8\in \Ann(y_8)$. But by Proposition~\ref{propos_E1p0}, $\Ann(z_8)=\Ann(y_8)$. Consequently, $c_8z_8=0$. Thus we may assume that $c_8=0$.

 Second, $y_6$ enters~$d^1a$ with coefficient~$c_5$. Therefore $c_5\in\Ann(y_6)$. By Proposition~\ref{propos_E1p0}, $\Ann(y_6)$ is the left ideal of~$\Z[\Gamma]$ generated by $r_1-1$ and $t_1t_2t_3\nu-1$. Since $t_1t_2t_3\nu-1$ annihilates~$z_5$, we may assume that $c_5=c_5'(r_1-1)$ for some $c_5'\in\Z[\Gamma]$. Subtracting from~$a$ the element~$d^1c_5'w_3$, we kill the coefficient~$c_5$. Similarly, considering the coefficient of~$y_7$ in~$d^1a$, we conclude that $c_6z_6=c_6'(1+\nu)z_6$ for some $c_6'\in\Z[\Gamma]$. Then the coefficient~$c_6$ is killed by adding to~$a$ the element~$d^1c_6'w_4$.

 Third, subtracting from~$a$ the element $c_4^+w_1^++c_4^-w_1^-$, we can kill the coefficients~$c_4^+$ and~$c_4^-$. When we have already achieved that $c_4^{\pm}=0$, we obtain that $y_5$ enters~$d^1a$ with coefficient~$c_3^++c_3^-$, and therefore $c_3^++c_3^-\in\Ann(y_5)$. But by Proposition~\ref{propos_E1p0}, $\Ann\bigl(z_3^+\bigr)=\Ann\bigl(z_3^-\bigr)=\Ann(y_5)$. Hence $c_3^+z_3^++c_3^-z_3^-=c_3\bigl(z_3^+-z_3^-)$ for some $c_3\in\Z[\Gamma]$. We have
 $$
 d^1\bigl(w_1^+-w_1^-+w_2\bigr)=-z_3^++z_3^-+b_1z_1+b_2z_2
 $$
 with certain $b_1,b_2\in\Z[\Gamma]$, which we do not need to compute explicitly. Thus, adding to~$a$ the element $c_3\bigl(w_1^+-w_1^-+w_2\bigr)$, we kill the coefficients~$c_3^{\pm}$.

 Finally, we may assume that
 $$
 a=c_1z_1+c_2z_2+c_9z_9,
 $$
 and we need to kill the coefficient~$c_1$. The coefficients of~$y_1$ and~$y_3$ in~$d^1a$ are equal to~$c_1$ and~$-c_1$, respectively. Hence $c_1\in\Ann(y_1)\cap\Ann(y_3)$. By Proposition~\ref{propos_E1p0}, $\Ann(y_1)$ is the left ideal of~$\Z[\Gamma]$ generated by~$t_1t_2t_3\nu-1$. Note that $(t_1t_2t_3\nu)^2=1$. It follows that $c_1=c_1'(t_1t_2t_3\nu-1)$ for some $c_1'\in\Z[\Gamma_1]$, where $\Gamma_1<\Gamma$ is the index~$2$ subgroup generated by~$t_1$, $t_2$, $t_3$, $\omega$, and~$r_1$. We write
 $$
 c_1'=p_1+p_2\omega+p_3\omega^2+q_1r_1+q_2\omega r_1+q_3\omega^2r_1,\qquad p_i,q_i\in\Lambda.
 $$
 Since $\Ann(y_3)$ is generated by~$\omega-1$ and~$r_1-1$, the condition $c_1\in\Ann(y_3)$ reads as
 \begin{equation}\label{eq_sum_pq}
 p_1+p_2+p_3+q_1+q_2+q_3=0.
 \end{equation}

 Consider the following elements of~$E^1_{3,0}$:
 \begin{align}\label{eq_u1}
 u_1&=(t_1-t_2t_3\nu)\bigl(w_1^+-w_1^-+w_2\bigr),\\
 u_2&=(t_1-t_2t_3\nu)w_1^--t_1w_2+t_2t_3w_6.\nonumber
 \end{align}
 A direct computation shows that
 \begin{align}\label{eq_d1u1}
 d^1u_1&=\bigl(\omega-\omega^2\bigr)(1+r_1)(t_1t_2t_3\nu-1)z_1-(t_1-t_2t_3\nu)\bigl(\omega\nu+\omega^2\bigr)z_2,\\
 d^1u_2&=\bigl(\omega^2-\omega r_1)(t_1t_2t_3\nu-1)z_1+m_2z_2+m_9z_9\nonumber
 \end{align}
 for some $m_2,m_9\in\Z[\Gamma]$, which we do not need to compute explicitly. (We do not need the coefficient of~$z_2$ in~\eqref{eq_d1u1} yet either. We have included it, since it will be needed later.) Therefore, adding $d^1u_1$ to~$a$ simultaneously adds~$1$ to both~$p_2$ and~$q_2$, and subtracts~$1$ from both~$p_3$ and~$q_3$. In other words, adding~$d^1u_1$ to~$a$ amounts to adding the vector $e_1=(0,1,-1,0,1,-1)$ to $(p_1,p_2,p_3,q_1,q_2,q_3)$. Consequently, adding~$d^1\omega u_1$ and $d^1\omega^2u_1$ to~$a$ amounts to adding the vectors $e_2=(-1,0,1,-1,0,1)$ and $e_3=(1,-1,0,1,-1,0)$, respectively, to  $(p_1,p_2,p_3,q_1,q_2,q_3)$. Similarly, adding~$d^1u_2$, $d^1\omega u_2$, and~$d^1\omega^2u_2$ to~$a$ amounts to adding  the vectors $e_4=(0,0,1,0,-1,0)$, $e_5=(1,0,0,0,0,-1)$, and $e_6=(0,1,0,-1,0,0)$, respectively, to $(p_1,p_2,p_3,q_1,q_2,q_3)$. Since the vectors $e_1,\ldots,e_6$ generate over~$\Lambda$ the subgroup of~$\Lambda^6$ consisting of all vectors satisfying~\eqref{eq_sum_pq}, we see that adding to~$a$ an appropriate linear combination of $u_1$, $\omega u_1$, $\omega^2u_1$, $u_2$, $\omega u_2$, and $\omega^2u_2$ with coefficients in~$\Lambda$, we can kill the coefficient~$c_1$.
\end{proof}

\begin{lem}\label{lem_t1-1nu}
 We have
 \begin{equation}\label{eq_t1-1nu}
  (t_1-1)\nu z_2\Iequiv -(t_1-1)z_2.
 \end{equation}
\end{lem}
\begin{proof}
 Since $d^1w_4=-(\nu+1)(z_2+z_6)$ and $t_1-1$ annihilates~$z_6$, we obtain that
 $$
 d^1(t_1-1)w_4=-(t_1-1)(\nu+1)z_2.
 $$
 The lemma follows.
\end{proof}

The following lemma is straightforward.

\begin{lem}\label{lem_FH}
 Let $F_1,F_2,F_3\in \Z[\tau_1,\tau_2,\tau_3]$ be polynomials in independent variables~$\tau_1$, $\tau_2$, and~$\tau_3$ such that
 $$
 \tau_1F_1+\tau_2F_2+\tau_3F_3=0.
 $$
 Then there exists polynomials $H_1,H_2,H_3\in \Z[\tau_1,\tau_2,\tau_3]$ such that $F_1=\tau_2H_3-\tau_3H_2$, $F_2=\tau_3H_1-\tau_1H_3$, and $F_3=\tau_1H_2-\tau_2H_1$.
\end{lem}

\begin{cor}\label{cor_FH}
 Let $F_1,F_2,F_3\in \Lambda$ be Laurent polynomials satisfying
 $$
 (t_1-1)F_1+(t_2-1)F_2+(t_3-1)F_3=0.
 $$
 Then there exists Laurent polynomials $H_1,H_2,H_3\in \Lambda$ such that
 \begin{align*}
  F_1&=(t_2-1)H_3-(t_3-1)H_2,\\
  F_2&=(t_3-1)H_1-(t_1-1)H_3,\\
  F_3&=(t_1-1)H_2-(t_2-1)H_1.
 \end{align*}
\end{cor}

\begin{proof}
 By multiplying~$F_1$, $F_2$, $F_3$ by $(t_1t_2t_3)^n$ for a sufficiently large~$n$, we may reduce the assertion of the corollary to the case where $F_1, F_2,F_3\in\Z[t_1,t_2,t_3]$. The  assertion then follows from Lemma~\ref{lem_FH} by the change of variables $\tau_i=t_i-1$.
\end{proof}

\begin{propos}\label{propos_E220_gen_mod}
 The group~$E^2_{2,0}$ is generated as a $\Z[\Delta]$-module by the two classes
 \begin{align*}
  \rho_1&=\bigl[(t_1-1)(t_2-1)(1-\omega)z_2\bigr],\\
  \rho_2&=\bigl[(t_1-1)(1-\omega)(1+r_2)t_1^{-1}z_9\bigr],
 \end{align*}
 where square brackets denote the homology class with respect to~$d^1$.
 \end{propos}

 \begin{proof}
 We set
 \begin{align*}
  R_1&=(t_1-1)(t_2-1)(1-\omega)z_2,\\
  R_2&=(t_1-1)(1-\omega)(1+r_2)t_1^{-1}z_9.
 \end{align*}
 A direct computation shows that $d^1R_1=0$ and~$d^1R_2=0$. Therefore, the classes~$\rho_1$ and~$\rho_2$ are well defined.

Recall that we  denote by~$I$ the image of the differential $d^1\colon E^1_{3,0}\to E^1_{2,0}$. We set
 $$
 J=I+\Z[\Delta]R_1+\Z[\Delta]R_2\subseteq E^1_{2,0}.
 $$
Our goal is to show that $J$ coincides with the kernel of the differential $d^1\colon E^1_{2,0}\to E^1_{1,0}$. We already know that $J\subseteq\ker d^1$. To prove the inverse inclusion, we take an arbitrary element $a\in E^1_{2,0}$ with $d^1a=0$, and modify it by adding elements of~$J$. If this process can be used to reduce $a$ to~$0$, then  the proposition follows. The reduction of $a$ to~$0$ will proceed in several steps.
\smallskip

\textsl{Step 1.}  By Lemma~\ref{lem_z2z9}, any element in the kernel of $d^1\colon E^1_{2,0}\to E^1_{1,0}$ can be reduced by adding an element of~$I$ to an element $a=c_2z_2+c_9z_9$ with $c_2,c_9\in\Z[\Gamma]$ and $d^1a=0$. Moreover, since $s_1+1$ annihilates~$z_2$ and $t_2-1$ annihilates~$z_9$, we may assume that
  \begin{align*}
   c_2&=\sum_{l=0}^2\omega^l\bigl(p_l+\tilde{p}_l\nu\bigr),&&p_l,\tilde{p}_l\in\Lambda,\\
   c_9&=\sum_{l=0}^2\omega^l\bigl(q_l+\tilde{q}_l\nu+\hat{q}_lr_2+\bar{q}_lr_2\nu\bigr)
   &&q_l,\tilde{q}_l,\hat{q}_l,\bar{q}_l\in\Z\bigl[t_1,t_1^{-1},t_3,t_3^{-1}\bigr].
  \end{align*}
\smallskip

\textsl{Step 2.}  Consider the element
  $$
  v=(t_1+t_2t_3\nu)(w_2-\nu w_6).
  $$
  We have
  \begin{align}
  d^1v&=-(t_1+t_2t_3\nu)\bigl(\omega\nu+\omega^2\bigr)z_2+(t_1\nu+t_2t_3)(1-r_1)z_9\nonumber\\
  {}&=-\bigl(t_1\omega^2+t_2t_3\omega+t_1\omega\nu+t_2t_3\omega^2\nu\bigr)z_2+(t_3+t_1\nu-\omega t_1r_2-\omega t_3r_2\nu)z_9.\label{eq_d1v}
  \end{align}
  Therefore, adding to~$a$ an appropriate $\Z[\Delta]$-multiple of~$d^1v$, we can kill the coefficients~$\bar{q}_l$ for $l=0,1,2$. Hence below we assume that $\bar{q}_0=\bar{q}_1=\bar{q}_2=0$.
\smallskip

\textsl{Step 3.} By Lemma~\ref{lem_t1-1nu}, the element $(t_1-1)(\nu+1)z_2$ belongs to~$I$. Adding to~$a$  an appropriate $\Z[\Delta]$-multiple of this element, we can add simultaneously to both~$p_l$ and~$\tilde{p}_l$ any Laurent polynomial divisible by~$t_1-1$. In this way, we can arrange that $\tilde{p}_l$ is independent of~$t_1$, i.e., belongs to~$\Z\bigl[t_2,t_2^{-1},t_3,t_3^{-1}\bigr]$ for every~$l$.
\smallskip

\textsl{Step 4.} Multiplying~\eqref{eq_d1v} by~$t_1^{-1}(1-t_3)$ and using that $(t_2-1)z_9=0$, we obtain:
  $$
  d^1t_1^{-1}(1-t_3)v=(t_3-1)\bigl(\omega^2+t_1^{-1}t_2t_3\omega+\omega\nu+t_1^{-1}t_2t_3\omega^2\nu\bigr)z_2-(t_3-1)\bigl(\nu+t_1^{-1}t_2t_3\bigr)z_9.
  $$
  The term involving~$z_2$ in this equality reads as
  \begin{multline*}
   \bigl(\omega^2(t_1-1)+\omega t_1t_2t_3^{-1}(t_2-1)+\omega(t_2-1)\nu+\omega^2t_1t_2^{-1}t_3(t_1-1)\nu\bigr)z_2\\
   {}\Iequiv
   \bigl(\omega^2\bigl(1-t_1t_2^{-1}t_3\bigr)(t_1-1)+\omega t_1t_2t_3^{-1}(t_2-1)+\omega(t_2-1)\nu\bigr)z_2,
  \end{multline*}
  where we have used~\eqref{eq_t1-1nu}. Therefore,
  \begin{multline*}
  0\Iequiv d^1\omega^2t_1^{-1}(1-t_3)v\Iequiv \bigl(\omega\bigl(1-t_1t_2^{-1}t_3\bigr)(t_1-1)+ t_1t_2t_3^{-1}(t_2-1)+(t_2-1)\nu\bigr)z_2\\
  {}-\omega^2(t_3-1)\bigl(\nu+t_1^{-1}t_2t_3\bigr)z_9.
  \end{multline*}
  Adding to~$a$ the latter element multiplied on the left by a certain coefficient from $\omega^l\Z\bigl[t_2,t_2^{-1},t_3,t_3^{-1}\bigr]$, we can add to~$\tilde{p}_l$ any Laurent polynomial in~$t_2$ and~$t_3$ divisible by~$t_2-1$. In this way, we can arrange that $\tilde{p}_l$ is independent of~$t_2$, i.e., belongs to~$\Z\bigl[t_3,t_3^{-1}\bigr]$, for every~$l$. Note that the described procedure does not violate the condition $\bar{q}_l=0$.
\smallskip

\textsl{Step 5.} For a Laurent polynomial $F\in\Lambda$, we set $F^{\omega}=\omega^2F\omega$ and $F^{\omega^2}=\omega F\omega^2$. In other words, the Laurent polynomials~$F^{\omega}$ and~$F^{\omega^2}$ are obtained from~$F$ by cyclic permutations of the variables~$t_1$, $t_2$, and~$t_3$.

Take a Laurent polynomial~$\varphi\in \Z\bigl[t_1,t_1^{-1},t_3,t_3^{-1}\bigr]$, and consider the element
\begin{align}
\varphi R_2&=\varphi(t_1-1)(1-\omega)(1+r_2)t_1^{-1}z_9\notag\\
{}&=\Bigl(\bigl(1-t_1^{-1}\bigr)\varphi -\omega t_1^{-1}(t_3-1)\varphi^{\omega}
+t_3^{-1}(t_1-1)\varphi r_2-\omega\bigl(1-t_3^{-1}\bigr)\varphi^{\omega}r_2\Bigr)z_9\label{eq_R2_subtract_new}
\end{align}
Note that $\varphi^{\omega}$ is a Laurent polynomial in~$t_2$ and~$t_3$. Since $t_2z_9=z_9$, the right-hand side of~\eqref{eq_R2_subtract_new} does not change if we replace all occurrences of~$t_2$ in~$\varphi^{\omega}$ by~$1$. Adding to~$a$ a suitable element of the form~$\varphi R_2$, we can change $\hat{q}_0$ by any Laurent polynomial in~$t_1$ and~$t_3$ divisible by~$t_1-1$. Thus we can arrange that $\hat{q}_0$ is a Laurent polynomial in~$t_3$ alone. Moreover, this operation does not change~$\hat{q}_2$, and changes~$\hat{q}_1$  only by a Laurent polynomial in~$t_3$, specifically, by the polynomial $(t_3^{-1}-1)\varphi^{\omega}$ with all occurrences of~$t_2$ replaced by~$1$. This observation allows us to repeat the same operation for~$\hat{q}_1$ and~$\hat{q}_2$. Namely, we first arrange that $\hat{q}_1$ is independent of~$t_1$ by adding to~$a$ a suitable element of the form~$\omega\varphi R_2$, and then arrange that $\hat{q}_2$ is also independent of~$t_1$ by adding to~$a$ a suitable element of the form~$\omega^2\varphi R_2$. The last operation changes~$\hat{q}_0$ only by a Laurent polynomial in~$t_3$. Therefore we obtain that all three polynomials~$\hat{q}_0$, $\hat{q}_1$, and~$\hat{q}_2$ belong to~$\Z\bigl[t_3,t_3^{-1}\bigr]$.
\smallskip

\textsl{Step 6.} At the previous steps, we have arrived at
\begin{multline}\label{eq_a_middle}
 a=\sum_{l=0}^2\omega^l\Bigl(\bigl(p_l+\tilde{p}_l\nu\bigr)z_2+\bigl(q_l+\tilde{q}_l\nu+\hat{q}_lr_2\bigr)z_9\Bigr),\\
   p_l\in\Lambda,\quad q_l,\tilde{q}_l\in\Z\bigl[t_1,t_1^{-1},t_3,t_3^{-1}\bigr],\quad \tilde{p}_l,\hat{q}_l\in\Z\bigl[t_3,t_3^{-1}\bigr].
 \end{multline}
Let us now study equations on the coefficients $p_l,\tilde{p}_l,q_l,\tilde{q}_l,\hat{q}_l$ arising from the condition~$d^1a=0$.

Since $r_1$ fixes~$y_9$ and $r_2r_1=\omega^2$, the term of~$d^1a$ containing~$y_9$ is equal to
$$
\sum_{l=0}^2\omega^l\bigl(q_l+\tilde{q}_l\nu+\hat{q}_l\omega^2\bigr)y_9=\sum_{l=0}^2\omega^l\bigl(q_l+\hat{q}_{l+1}^{\,\omega^2}+\tilde{q}_l\nu\bigr)y_9.
$$
(Hereafter, sums of indices are always taken modulo~$3$.)
The annihilator of~$y_9$ is generated by~$t_2-1$, $t_3-1$, and~$r_1-1$. Since $q_l+\hat{q}_{l+1}^{\,\omega^2}$ and~$\tilde{q}_l$ are  Laurent polynomials in~$t_1$ and~$t_3$, we obtain that both $q_l+\hat{q}_{l+1}^{\,\omega^2}$ and~$\tilde{q}_l$ are divisible by $t_3-1$.

The coefficient of~$y_2$ in~$d^1a$ is equal to
\begin{multline*}
\sum_{l=0}^2\omega^l\Bigl(\bigl(p_l+\tilde{p}_l\nu\bigr)(1-s_1)+\bigl(q_l+\tilde{q}_l\nu+\hat{q}_lr_2\bigr)\bigl(\omega-t_1t_3^{-1}r_3\bigr)\Bigr)\\
{}=\sum_{l=0}^2\omega^l\Bigl(\bigl(p_l+q_{l-1}^{\omega}-t_2t_3^{-1}\hat{q}_{l-1}^{\,\omega} \bigr)+\bigl(-\tilde{p}_l-t_2^{-1}t_3 q_{l-1}^{\omega}+\hat{q}_{l-1}^{\,\omega}\bigr)r_1+\bigl(\tilde{p}_l+\tilde{q}_{l-1}^{\,\omega}\bigr)\nu\\{}+\bigl(-p_l-t_2t_3^{-1}\tilde{q}_{l-1}^{\,\omega}\bigr)s_1\Bigr).
\end{multline*}
 Since the annihilator of~$y_2$ is generated by~$t_1-1$, we obtain that each of the following Laurent polynomials must be divisible by~$t_1-1$:
\begin{gather}
 p_l+q_{l-1}^{\omega}-t_2t_3^{-1}\hat{q}_{l-1}^{\,\omega},\label{eq_pl+nul}\\
 -\tilde{p}_l-t_2^{-1}t_3q_{l-1}^{\omega}+\hat{q}_{l-1}^{\,\omega},\nonumber\\
 \tilde{p}_l+\tilde{q}_{l-1}^{\,\omega},\nonumber\\
 -p_l-t_2t_3^{-1}\tilde{q}_{l-1}^{\,\omega}.\nonumber
\end{gather}
But the Laurent polynomials $\tilde{p}_l$, $q_{l-1}^{\omega}$, $\hat{q}_{l-1}^{\,\omega}$ and~$\tilde{q}_{l-1}^{\,\omega}$ do not involve~$t_1$. Therefore
\begin{gather}
-\tilde{p}_l-t_2^{-1}t_3q_{l-1}^{\omega}+\hat{q}_{l-1}^{\,\omega}=0,\label{eq_chi_l20}\\
\tilde{p}_l+\tilde{q}_{l-1}^{\,\omega}=0.\nonumber
\end{gather}
We have already proved that $\tilde{q}_{l-1}$ is divisible by~$t_3-1$. Hence $\tilde{q}_{l-1}^{\,\omega}$ is divisible by~$t_2-1$. On the other hand, $\tilde{p}_l$ is a Laurent polynomial in~$t_3$ alone.
Consequently,
$\tilde{p}_l=0$ and~$\tilde{q}_{l-1}^{\,\omega}=0$ for all~$l$, and hence $\tilde{q}_l=0$ for all~$l$.

Now, \eqref{eq_chi_l20} implies that $q_l=t_1^{-1}t_3\hat{q}_l$. But we have already proved that $q_l+\hat{q}_{l+1}^{\,\omega^2}$ is divisible  by~$t_3-1$. Therefore $t_3\hat{q}_l+\bigl(t_3\hat{q}_{l+1}\bigr)^{\omega^2}$ is divisible by~$t_3-1$. Since $t_3\hat{q}_l\in\Z\bigl[t_3,t_3^{-1}\bigr]$ and $\bigl(t_3\hat{q}_{l+1}\bigr)^{\omega^2}\in\Z\bigl[t_1,t_1^{-1}\bigr]$, it follows that $m_{l+1}=t_3\hat{q}_{l+1}$ is a constant, i.e., lies in~$\Z$, for every~$l$. Now, the condition that $t_3\hat{q}_l+\bigl(t_3\hat{q}_{l+1}\bigr)^{\omega^2}$ is divisible by~$t_3-1$ means that $m_l+m_{l+1}=0$. Consequently, $m_0=m_1=m_2=0$. Thus $q_l=0$ and $\hat{q}_l=0$ for all~$l$.

The condition that the Laurent polynomial~\eqref{eq_pl+nul} is divisible by~$t_1-1$ now implies that $p_l$ is divisible by~$t_1-1$ for every~$l$.
\smallskip

\textsl{Step 7.} At the previous steps, we have arranged that all coefficients in~\eqref{eq_a_middle}, except for~$p_l$, vanish, and $p_l$ is divisible by~$t_1-1$ for every~$l$. It is now convenient for us to rewrite~$a$ in the form
$$
a=\bigl((t_1-1)F_1+(t_2-1)F_2\omega+(t_3-1)F_3\omega^2\bigr)z_2,\qquad  F_1,F_2,F_3\in\Lambda.
$$
Since $\omega$ fixes~$y_4$, the coefficient of~$y_4$ in~$d^1a$ is equal to $\sum_{i=1}^3(t_i-1)F_i$. But the elements~$t_1^{k_1}t_2^{k_2}t_3^{k_3}y_4$ with $k_1,k_2,k_3\in\Z$ are linearly independent in~$E^1_{1,0}$. Hence
$$
(t_1-1)F_1+(t_2-1)F_2+(t_3-1)F_3=0.
$$
Let $H_1$, $H_2$, and~$H_3$ be the Laurent polynomials from Corollary~\ref{cor_FH}. Then
$$
a=\bigl(H_3+H_1\omega+H_2\omega^2\bigr)R_1,
$$
which completes the proof of the proposition.
\end{proof}

 \begin{propos}\label{propos_E220_rel}
 The classes~$\rho_1$ and~$\rho_2$ satisfy the relations:
\begin{gather}
(t_1t_2-t_3)\rho_1=0,\label{eq_rel_rho1}\\
(t_1-1)(t_2-1)\rho_1=0,\label{eq_rel_rho2}\\
 \bigl(1+\omega+\omega^2\bigr)\rho_1=0,\label{eq_rel_rho3}\\
 \bigl(1+\omega+\omega^2\bigr)t_3\rho_1=0,\label{eq_rel_rho4}\\
 (t_2-1)(t_3-1)\rho_2=0,\label{eq_rel_rho5}\\
 \bigl(t_2-1+(t_1-1)\omega\bigr)\rho_2=0.\label{eq_rel_rho6}
\end{gather}
\end{propos}

\begin{proof}
 Since $t_2$ fixes~$z_9$ and~$t_3$ fixes~$r_1z_9$, we have
 \begin{multline*}
 d^1(t_2-1)(t_3-1)(t_1+t_2t_3\nu)(w_2-\nu w_6)=-(t_2-1)(t_3-1)(t_1+t_2t_3\nu)\bigl(\omega\nu+\omega^2\bigr)z_2\\
 {}=-\omega(t_1-1)(t_2-1)(t_1t_2+t_3\nu) z_2-\omega^2(t_3-1)(t_1-1)(t_2+t_1t_3\nu)z_2.
 \end{multline*}
 Multiplying this by~$-\omega^2$ and using~\eqref{eq_t1-1nu}, we obtain:
 \begin{align*}
 0&\Iequiv(t_1-1)(t_2-1)(t_1t_2+t_3\nu)z_2+\omega (t_3-1)(t_1-1)(t_2+t_1t_3\nu)z_2\\
 {}&\Iequiv(t_1-1)(t_2-1)(t_1t_2-t_3)z_2+\omega (t_3-1)(t_1-1)(t_2-t_1t_3)z_2\\
 {}&=(t_1t_2-t_3)(t_1-1)(t_2-1)(1-\omega)z_2,
 \end{align*}
 which implies~\eqref{eq_rel_rho1}.

Consider the element~$u_1$ from the proof of Lemma~\ref{lem_z2z9}; see~\eqref{eq_u1}. Multiplying \eqref{eq_d1u1}  by $1+\omega+\omega^2$, we obtain:
\begin{align*}
 0&\Iequiv -\bigl(1+\omega+\omega^2\bigr)(t_1-t_2t_3\nu)\bigl(\omega\nu+\omega^2\bigr)z_2&&\\
 {}&=-\bigl(1+\omega+\omega^2\bigr)\bigl(\omega(t_3\nu-t_1t_2)+\omega^2(t_2-t_1t_3\nu)\bigr)z_2&&\\
 {}&=\bigl(1+\omega+\omega^2\bigr)\bigl(-t_2+t_1t_2-t_3\nu+t_1t_3\nu\bigr)z_2&&\\
 {}&\Iequiv \bigl(1+\omega+\omega^2\bigr)\bigl(-t_2+t_1t_2+t_3-t_1t_3\bigr)z_2&&\text{(by~\eqref{eq_t1-1nu})}\\
 {}&= \bigl(1+\omega+\omega^2\bigr)\bigl((t_1-1)(t_2-1)-(t_1-1)(t_3-1)\bigr)z_2\\
 {}&= \bigl(1+\omega+\omega^2\bigr)(t_1-1)(t_2-1)(1-\omega)z_2.
\end{align*}
This gives~\eqref{eq_rel_rho3}.

Since $1+\omega+\omega^2$ commutes with $(t_1-1)(t_2-1)(t_3-1)$, we have
$$
\bigl(1+\omega+\omega^2\bigr)(t_3-1)\rho_1=0.
$$
Combining this with~\eqref{eq_rel_rho3}, we obtain~\eqref{eq_rel_rho4}.

Now, a direct computation gives that
\begin{multline*}
(t_1-1)(t_2-1)=-(t_1+t_2)\bigl(1+\omega+\omega^2\bigr)+\bigl(1+t_3^{-1}\bigr)\bigl(1+\omega+\omega^2\bigr)t_3\\
{}+ \bigl(1+\omega t_2^{-1}+\omega^2t_1^{-1}\bigr)(t_1t_2-t_3).
\end{multline*}
Therefore~\eqref{eq_rel_rho2} follows from~\eqref{eq_rel_rho1}, ~\eqref{eq_rel_rho3}, and~\eqref{eq_rel_rho4}.

Since the product $(t_1-1)(t_2-1)(t_3-1)$ commutes with~$\omega$ and~$r_2$, and annihilates~$z_9$, we obtain~\eqref{eq_rel_rho5}.

Finally, we have
\begin{align*}
 &\bigl(t_2-1+(t_1-1)\omega\bigr)(t_1-1)(1-\omega)(1+r_2)t_1^{-1}z_9\\
 &\hspace{2cm}=(t_1-1)(t_2-1)\bigl(1-\omega^2\bigr)(1+r_2)t_1^{-1}z_9=0,
\end{align*}
since $t_2-1$ annihilates~$z_9$ and~$r_2z_9$, while $t_1-1$ annihilates~$\omega^2z_9$ and~$\omega^2r_2z_9$.
Thus we obtain~\eqref{eq_rel_rho6}.
\end{proof}

\begin{propos}
\label{propos_E220_gen_group}
 As an abelian group, $E^2_{2,0}$ is generated by the following elements:
 \begin{itemize}
  \item $t_1^k\rho_1$, where $k\in\Z$,
  \item $t_2^k\rho_1$, where $k\in\Z\setminus\{0\}$,
  \item $\omega t_1^k\rho_1$, where $k\in\Z$,
  \item $\omega^lt_1^{k_1}t_2^{k_2}\rho_2$, where $l\in\{0,1,2\}$ and $k_1,k_2\in\Z$.
 \end{itemize}
\end{propos}

\begin{proof}
 By Proposition~\ref{propos_E220_gen_mod}, $E^2_{2,0}$ is generated as an abelian group by the elements $\omega^lt_1^{k_1}t_2^{k_2}t_3^{k_3}\rho_m$, where $l\in\{0,1,2\}$, $k_1,k_2,k_3\in\Z$, and $m\in\{1,2\}$. Our aim is to show that, using relations~\eqref{eq_rel_rho1}--\eqref{eq_rel_rho6}, all these elements can be written as $\Z$-linear combinations of the elements listed in Proposition~\ref{propos_E220_gen_group}.

 First, relations~\eqref{eq_rel_rho3} and~\eqref{eq_rel_rho1} allow us to write all elements $\omega^lt_1^{k_1}t_2^{k_2}t_3^{k_3}\rho_1$ as $\Z$-linear combinations of elements of the form $t_1^{k_1}t_2^{k_2}\rho_1$ and $\omega t_1^{k_1}t_2^{k_2}\rho_1$.

 Second, multiply~\eqref{eq_rel_rho3} by~$t_2$ and then subtract from it~\eqref{eq_rel_rho4}. We obtain:
 $$
 (t_2-t_3)\rho_1+\omega(t_1-t_3)\rho_1=0.
 $$
 Multiplying this by~$t_2^pt_3^q$ and using~\eqref{eq_rel_rho1}, we obtain
 $$
 t_1^qt_2^{p+q+1}(1-t_1)\rho_1+\omega t_1^{p+1}t_2^q(1-t_2)\rho_1=0.
 $$
 These relations allow us to write all elements of the form~$\omega  t_1^{k_1}t_2^{k_2}\rho_1$ as $\Z$-linear combinations of elements of the form~$t_1^{k_1}t_2^{k_2}\rho_1$ and elements of the form~$\omega t_1^{k}\rho_1$.

 Third, relation~\eqref{eq_rel_rho2} allows us to write all elements~$t_1^{k_1}t_2^{k_2}\rho_1$ as $\Z$-linear combinations  of elements of the form~$t_1^k\rho_1$ and~$t_2^k\rho_1$.

 Now we set $\sigma=(t_2-1)\rho_2$. By~\eqref{eq_rel_rho5} we have $t_3\sigma=\sigma$. Multiplying~\eqref{eq_rel_rho6} by~$\omega^2$, we obtain
 $$
 (t_3-1)\rho_2=-\omega^2\sigma.
 $$
 This relation allows us to write elements of the form~$\omega^l  t_1^{k_1}t_2^{k_2}t_3^{k_3}\rho_2$ as $\Z$-linear combinations of elements of the form~$\omega^l  t_1^{k_1}t_2^{k_2}\rho_2$ and elements of the form~$\omega^l  t_1^{k_1}t_2^{k_2}\sigma$. The latter elements are also $\Z$-linear combinations of elements of the form~$\omega^l  t_1^{k_1}t_2^{k_2}\rho_2$.
\end{proof}

\begin{remark}
 In fact, in Section~\ref{section_triv_E320} we will see that the group~$E^2_{2,0}$ is free abelian with basis consisting of the elements listed in Proposition~\ref{propos_E220_gen_group}; see Remark~\ref{remark_basis_E220}.
\end{remark}

\section{Differential $d^2\colon E^2_{2,0}\to E^2_{0,1}$}\label{section_d2rho}

In this section, we compute the differential $d^2\colon E^2_{2,0}\to E^2_{0,1}$. Let $\CR_*$ be the normalized bar resolution of~$\Z$ over~$\Z[\CS]$. We will perform the computation in the bicomplex $C_*(\CB_{0,6})\otimes_{\CS}\CR_*$ modulo $C_*(K)\otimes_{\CS}\CR_*$ for certain $\G$-equivariant subcomplex $K\subseteq\CB_{0,6}$ on which $\CS$ acts freely, and then use Lemma~\ref{lem_dr}. Specifically, let $K$ be the one-dimensional subcomplex of~$\CB_{0,6}$ consisting of all vertices in the $\G$-orbit~$\G X_3$ and all edges in the $\G$-orbits~$\G Y_4$ and~$\G Y_9$. (The reason for considering this subcomplex is that it is precisely the subcomplex consisting of all simplices having all vertices in~$\G X_3$.) By the data in Tables~\ref{table_0cells} and~\ref{table_1cells}, the stabilizers~$\CS_{X_3}$, $\CS_{Y_4}$, and~$\CS_{Y_9}$ are trivial. This yields the following proposition.

\begin{propos}
 The group~$\CS$ acts freely on~$K$.
\end{propos}

This proposition allows us to apply Lemma~\ref{lem_dr}. We will use the symbol~$\Kequiv$ to denote equality modulo~$C_*(K)$ in~$C_*(\CB_{0,6})$, as well as equality modulo $C_*(K)\otimes_{\CS}\CR_*$ in $C_*(\CB_{0,6})\otimes_{\CS}\CR_*$. We always identify $C_*(\CB_{0,6})\otimes_{\CS}\CR_0$ with $C_*(\CB_{0,6})$. Recall that, under this identification, the differential~$\partial'$ coincides with the usual boundary operator~$\partial$.
In the computations below, we will repeatedly use the following lemma.

\begin{lem}\label{lem_d'd''}
 Let $Q$ be a cell of~$\CB_{0,6}$ and let $h\in\CS$. Then
 \begin{equation*}
   \bigl(h^{-1}-1\bigr)Q=(-1)^{\dim Q}\,\partial''(Q\otimes[h]).
 \end{equation*}
\end{lem}

\begin{proof}
 We have
 \begin{multline*}
\bigl(h^{-1}-1\bigr)Q= \bigl(h^{-1}Q\bigr)\otimes\ebr-Q\otimes\ebr=Q\otimes (h\ebr-\ebr)=Q\otimes d[h]\\{}=(-1)^{\dim Q}\,\partial''(Q\otimes[h]).
 \end{multline*}
 The lemma follows.
\end{proof}

\begin{propos}\label{propos_d2rho1}
 We have
 \begin{equation*}%\label{eq_d2rho1}
 d^2\rho_1=(1-\omega)(t_2-1)\xi_3.
 \end{equation*}
\end{propos}

\begin{proof}
 By Lemma~\ref{lem_dr}, to compute~$d^2\rho_1$ it suffices to find elements
 $$
 \theta_{2,0}\in C_2(\CB_{0,6})\otimes_{\CS}\CR_0=C_2(\CB_{0,6}),\qquad\theta_{1,1}\in C_1(\CB_{0,6})\otimes_{\CS}\CR_1,\qquad \psi\in C_0(\CB_{0,6})\otimes_{\CS}\CR_1
 $$
 such that $\theta_{2,0}$ represents~$\rho_1$, $\partial'\theta_{2,0}+\partial''\theta_{1,1}\Kequiv0$, $\partial'\theta_{1,1}\Kequiv\psi$, and~$\partial''\psi=0$. Then the element~$\psi$ will represent~$d^2\rho_1$. (In the  statement of Lemma~\ref{lem_dr} for the case of the differential $d^2\colon E^2_{2,0}\to E^2_{0,1}$, a class $\theta_{0,2}\in C_0(\CB_{0,6})\otimes_{\CS}\CR_2$ also appears; however, this class will be zero in our computation.)

 We set
 $$
 \theta_{2,0}=\bigl(T_{\bar{2}3}^{-1}-1\bigr)\bigl(T_{\bar{3}1}^{-1}-1\bigr)(1-\omega)Z_2.
 $$
 Since $\kappa(T_{\bar{2}3})=t_1^{-1}$ and $\kappa(T_{\bar{3}1})=t_2^{-1}$, we see that $\theta_{2,0}$ represents~$\rho_1$.

 Since $\partial Z_2=(1-s_1)Y_2+Y_4$ and $\omega$ fixes~$Y_4$, we have
 $$
 \partial'\theta_{2,0}=\bigl(T_{\bar{2}3}^{-1}-1\bigr)\bigl(T_{\bar{3}1}^{-1}-1\bigr)(1-\omega)(1-s_1)Y_2.
 $$
 The conjugation by~$s_1$ takes~$\omega$ to~$\omega^2$, $T_{\bar{2}3}$ to~$T_{2\bar{3}}^{-1}$, and~$T_{\bar{3}1}$ to~$T_{1\bar{2}}^{-1}$. Hence
 \begin{align*}
 \partial'\theta_{2,0}&=\bigl(T_{\bar{2}3}^{-1}-1\bigr)\bigl(T_{\bar{3}1}^{-1}-1\bigr)(1-\omega)Y_2-s_1(T_{2\bar{3}}-1)(T_{1\bar{2}}-1)\bigl(1-\omega^2\bigr)Y_2\\
 {}&=\bigl(T_{\bar{2}3}^{-1}-1\bigr)\bigl(T_{\bar{3}1}^{-1}-1\bigr)(1-\omega)Y_2,
 \end{align*}
 where the latter equality holds because $T_{2\bar{3}}$ fixes~$Y_2$, $T_{1\bar{2}}$ fixes~$\omega^2Y_2$, and $T_{2\bar{3}}$ and~$T_{1\bar{2}}$ commute. Further we have
 \begin{align*}
 \partial'\theta_{2,0}&=\Bigl(\bigl(T_{\bar{2}3}^{-1}-1\bigr)\bigl(T_{\bar{3}1}^{-1}-1\bigr)-\omega\bigl(T_{\bar{1}2}^{-1}-1\bigr)\bigl(T_{\bar{2}3}^{-1}-1\bigr)\Bigr)Y_2\\
 &=\Bigl(T_{\bar{3}1}^{-1}-1-\omega\bigl(T_{\bar{1}2}^{-1}-1\bigr)\Bigr)\bigl(T_{\bar{2}3}^{-1}T_{2\bar{3}}-1\bigr)Y_2.
 \end{align*}
 Since $T_{2\bar{3}}^{-1}T_{\bar{2}3}\in\CS$, by Lemma~\ref{lem_d'd''}, we have
 $
 \partial'\theta_{2,0}=-\partial''\theta_{1,1}
 $,
 where
 $$
 \theta_{1,1}=\Bigl(T_{\bar{3}1}^{-1}-1-\omega\bigl(T_{\bar{1}2}^{-1}-1\bigr)\Bigr)\bullet\bigl(Y_2\otimes\bigl[T_{2\bar{3}}^{-1}T_{\bar{2}3}\bigr]\bigl).
 $$
 Hereafter, $\bullet$ denotes the action of~$\G$ on the bicomplex $C_*(\CB_{0,6})\otimes_{\CS}\CR_*$ from Fact~\ref{fact_Gamma}. Recall that the differentials~$\partial'$ and~$\partial''$ are $\G$-equivariant.

 Let us now compute $\partial'\theta_{1,1}$ up to elements of~$C_0(K)\otimes_{\CS}\CR_1$. Since $\partial Y_2\Kequiv -X_2$, we have $\partial'\theta_{1,1}\Kequiv\psi$, where
\begin{gather*}
     \psi=-\Bigl(T_{\bar{3}1}^{-1}-1-\omega\bigl(T_{\bar{1}2}^{-1}-1\bigr)\Bigr)\bullet\bigl(X_2\otimes\bigl[T_{2\bar{3}}^{-1}T_{\bar{2}3}\bigr]\bigl).
\end{gather*}
The element~$X_2\otimes\bigl[T_{2\bar{3}}^{-1}T_{\bar{2}3}\bigr]$ lies in the kernel of~$\partial''$ and represents~$-\xi_3$ in~$E^1_{0,1}$. Therefore, $\partial''\psi=0$ and~$\psi$ represents in~$E^2_{0,1}$ the element
$$
\bigl(t_2-1-\omega(t_3-1)\bigr)\xi_3=(1-\omega)(t_2-1)\xi_3,
$$
where the latter equality holds, since, by Proposition~\ref{propos_E201}, $t_3\xi_3=t_2\xi_3$ in~$E^2_{0,1}$. The proposition follows.
\end{proof}

\begin{propos}\label{propos_d2rho2}
 We have
 \begin{equation*}%\label{eq_d2rho2}
 d^2\rho_2=(t_1-1)\bigl(\omega-\omega^2\bigr)t_3^{-1}\bigl(t_2^{-1}\xi_1+\xi_3\bigr).
 \end{equation*}
\end{propos}

\begin{proof}
 We proceed in the same way as in the proof of Proposition~\ref{propos_d2rho1}. The only difference is that now the element $\theta_{0,2}\in C_0(\CB_{0,6})\otimes\CR_2$ will also be nonzero. We set
 $$
 \theta_{2,0}=\bigl(T_{2\bar{3}}^{-1}-1\bigr)(1-\omega)(1+r_2)T_{2\bar{3}}Z_9.
 $$
 Then $\theta_{2,0}$ represents~$\rho_2$. Since $\partial Z_9\Kequiv(\omega-T_{\varepsilon_1}r_3)Y_2$, we have
 \begin{align*}
  \partial'\theta_{2,0}&\Kequiv \bigl(T_{2\bar{3}}^{-1}-1\bigr)(1-\omega)(1+r_2)T_{2\bar{3}}(\omega-T_{\varepsilon_1}r_3)Y_2\\
  &= \bigl(T_{2\bar{3}}^{-1}-1\bigr)(1-\omega)(1+r_2)(\omega T_{1\bar{2}}-r_3T_{3\bar{1}}T_{\varepsilon_2})Y_2.
 \end{align*}
But
$$
(1-\omega)(1+r_2)\omega=(1-\omega)(1+r_2)r_3=\bigl(\omega-\omega^2\bigr)(1+r_1).
$$
Moreover, $T_{2\bar{3}}$ fixes~$Y_2$ and $T_{1\bar{2}}T_{2\bar{3}}T_{\varepsilon_2}^{-1}T_{3\bar{1}}^{-1}\in\CS$. Therefore
$$
\partial'\theta_{2,0}\Kequiv \bigl(T_{2\bar{3}}^{-1}-1\bigr)\bigl(\omega-\omega^2\bigr)(1+r_1)\bigl(1-T_{3\bar{1}}T_{\varepsilon_2}T_{2\bar{3}}^{-1}T_{1\bar{2}}^{-1}\bigr)T_{1\bar{2}}Y_2.
$$
By Lemma~\ref{lem_d'd''}, we have
$$
\partial'\theta_{2,0}\Kequiv-\partial''\theta_{1,1},
$$
where
$$
\theta_{1,1}=-\bigl(T_{2\bar{3}}^{-1}-1\bigr)\bigl(\omega-\omega^2\bigr)(1+r_1)
\bullet\bigl((T_{1\bar{2}}Y_2)\otimes
\bigl[T_{1\bar{2}}T_{2\bar{3}}T_{\varepsilon_2}^{-1}T_{3\bar{1}}^{-1}\bigr]\bigl).
$$
Since $\partial Y_2\Kequiv -X_2$, we have
$$
\partial'\theta_{1,1}\Kequiv\bigl(T_{2\bar{3}}^{-1}-1\bigr)\bigl(\omega-\omega^2\bigr)(1+r_1)
\bullet\bigl((T_{1\bar{2}}X_2)\otimes
\bigl[T_{1\bar{2}}T_{2\bar{3}}T_{\varepsilon_2}^{-1}T_{3\bar{1}}^{-1}\bigr]\bigl).
$$
We set
\begin{align*}
Q&=T_{1\bar{2}}X_2,&h_1&=T_{1\bar{2}}T_{2\bar{3}}T_{\varepsilon_2}^{-1}T_{3\bar{1}}^{-1},&
h_2&=r_1h_1r_1=T_{3\bar{1}}T_{2\bar{3}}T_{\varepsilon_3}^{-1}T_{1\bar{2}}^{-1}.
\end{align*}
Since the elements~$T_{\varepsilon_2}^{-1}r_1$ and~$T_{2\bar{3}}$ fix~$X_2$, we obtain that
$$
r_1Q=r_1T_{1\bar{2}}X_2=T_{3\bar{1}}r_1X_2=T_{3\bar{1}}T_{\varepsilon_2}X_2=h_1^{-1}Q.
$$
Therefore
\begin{align*}
(1+r_1)\bullet\bigl(Q\otimes [h_1]\bigr)&=Q\otimes[h_1]+\bigl(h_1^{-1}Q\bigr)\otimes [h_2]=Q\otimes([h_1]+h_1[h_2])\\
{}&=\partial''\bigl(Q\otimes[h_1|h_2]\bigr)+Q\otimes[h_1h_2].
\end{align*}
Moreover, using relation~\eqref{eq_lantern_eps_jk} and the fact that twists about disjoint simple closed curves commute, we obtain that
$$
h_1h_2=T_{1\bar{2}}T_{2\bar{3}}^2T_{\varepsilon_2}^{-1}T_{\varepsilon_3}^{-1}T_{1\bar{2}}^{-1}=T_{1\bar{2}}T_{2\bar{3}}T_{\bar{2}3}^{-1}T_{\delta_1'}T_{1\bar{2}}^{-1}.
$$
Now, set
$$
\theta_{0,2}=-\bigl(T_{2\bar{3}}^{-1}-1\bigr)\bigl(\omega-\omega^2\bigr)
\bullet(Q\otimes[h_1|h_2]).
$$
Then we have $\partial'\theta_{1,1}+\partial''\theta_{0,2}\Kequiv \psi$, where
\begin{align*}
 \psi&=\bigl(T_{2\bar{3}}^{-1}-1\bigr)\bigl(\omega-\omega^2\bigr)
\bullet(Q\otimes[h_1h_2])\\
{}&=\bigl(T_{2\bar{3}}^{-1}-1\bigr)\bigl(\omega-\omega^2\bigr)T_{1\bar{2}}\bullet \bigl(X_2\otimes\bigl[T_{2\bar{3}}T_{\bar{2}3}^{-1}T_{\delta_1'}\bigr]\bigr)
\end{align*}
The mapping class~$T_{2\bar{3}}T_{\bar{2}3}^{-1}T_{\delta_1'}$ lies in the stabilizer~$\CS_{X_2}$. Consequently, the element $X_2\otimes\bigl[T_{2\bar{3}}T_{\bar{2}3}^{-1}T_{\delta_1'}\bigr]$ lies in the kernel of~$\partial''$ and represents in~$E^1_{0,1}$ the class
$$
\iota_{X_2}\bigl([T_{\delta_1'}]+\bigr[T_{2\bar{3}}T_{\bar{2}3}^{-1}\bigr]\bigr)=t_2^{-1}\xi_2+\xi_3.
$$
Here we have  used that $\delta_1'=T_{31}^{-1}(\delta_1)$, and hence $T_{\delta_1'}=T_{31}^{-1}T_{\delta_1}T_{31}$.

Thus $\partial''\psi=0$ and~$\psi$ represents in~$E^1_{0,1}$ the class
$$
(t_1-1)\bigl(\omega-\omega^2\bigr)t_3^{-1}\bigl(t_2^{-1}\xi_2+\xi_3\bigr).
$$
The proposition follows, since, by Proposition~\ref{propos_E201}, we have $\xi_1=\xi_2$ in~$E^2_{0,1}$.
\end{proof}

\section{Triviality of the group $E^{3}_{2,0}$}\label{section_triv_E320}

\begin{propos}\label{propos_E320}
 We have $E^3_{2,0}=0$.
\end{propos}

\begin{proof}
 To prove the proposition, we need to show that the differential $d^2\colon E^2_{2,0}\to E^2_{0,1}$ is injective. Suppose that $a$ is an element lying in the kernel of the differential $d^2\colon E^2_{2,0}\to E^2_{0,1}$. By Proposition~\ref{propos_E220_gen_group}, we have
$$
a=\bigl(q_1(t_1)+q_2(t_2)+\omega q_3(t_1)\bigr)\rho_1+\sum_{l=0}^2\omega^lp_l(t_1,t_2)\rho_2
$$
for some Laurent polynomials $p_l\in\Z\bigl[x^{\pm1},y^{\pm1}\bigr]$, $q_l\in\Z\bigl[x^{\pm1}\bigr]$, where $q_2$ has zero constant term. Here it is convenient to modify our usual notation and specify explicitly the variables that Laurent polynomials depend on, for instance, write $p_l(t_1,t_2)$ rather than~$p_l$. This is because polynomials with certain expressions substituted instead of variables will occur in our formulas.

From Proposition~\ref{propos_E201} it follows that there is an isomorphism of left~$\Z[\Delta]$-modules
 \begin{equation*}
  E^2_{0,1}\cong \frac{\Z[\Delta]}{(t_1-1)\Z[\Delta]}\oplus\frac{\Z[\Delta]}{(t_1-1)\Z[\Delta]+(t_2-t_3)\Z[\Delta]}
 \end{equation*}
where the generators of the summands are~$\xi_1$ and~$\xi_3$, respectively. We denote the summands by~$\Xi_1$ and~$\Xi_3$, respectively, and write
$$
d^2a=b_1+b_3,\qquad b_i\in\Xi_i.
$$
Since  $d^2a=0$, we have $b_1=0$ and~$b_3=0$.

By Propositions~\ref{propos_d2rho1} and~\ref{propos_d2rho2}, we have
\begin{align*}
b_1&=\sum_{l=0}^2\omega^lp_l(t_1,t_2)(t_1-1)\bigl(\omega-\omega^2\bigr)t_2^{-1}t_3^{-1}\xi_1\\
{}&=\sum_{l=0}^2\omega^l\bigl(
(t_3-1)p_{l-1}(t_3,t_1)-(t_2-1)p_{l+1}(t_2,t_3)\bigr)t_2^{-1}t_3^{-1}\xi_1,
\end{align*}
where sums of indices are taken modulo~$3$. Since $t_1\xi_1=\xi_1$, this formula reads  as
\begin{equation}\label{eq_b1}
b_1=\sum_{l=0}^2\omega^l\bigl(
(t_3-1)p_{l-1}(t_3,1)-(t_2-1)p_{l+1}(t_2,t_3)\bigr)t_2^{-1}t_3^{-1}\xi_1.
\end{equation}
Moreover, by Proposition~\ref{propos_E201}, the elements $\omega_lt_2^{k_2}t_3^{k_3}\xi_1$, where $l\in\{0,1,2\}$ and $k_2,k_3\in\Z$, are linearly independent in~$E^2_{0,1}$. Therefore all coefficients in~\eqref{eq_b1} vanish, that is,
\begin{equation*}
(t_2-1)p_{l+1}(t_2,t_3)=(t_3-1)p_{l-1}(t_3,1),\qquad l=0,1,2.
\end{equation*}
The left-hand side of  this equality is divisible by~$t_2-1$, and  the right-hand side is a Laurent polynomial of~$t_3$ alone. It follows that both sides of the equality vanish, and therefore $p_l=0$ for $l=0,1,2$. Thus
$$
a=\bigl(q_1(t_1)+q_2(t_2)+\omega q_3(t_1)\bigr)\rho_1.
$$

Now, by Proposition~\ref{propos_d2rho1}, we have
\begin{align*}
b_3&=\bigl(q_1(t_1)+q_2(t_2)+\omega q_3(t_1)\bigr)(1-\omega)(t_2-1)\xi_3\\
{}&=\bigl(q_1(t_1)+q_2(t_2)+\omega(-q_1(t_3)-q_2(t_1)+q_3(t_1))-\omega^2q_3(t_3)\bigr)(t_2-1)\xi_3\\
{}&=\bigl(q_1(1)+q_2(t_2)+\omega(-q_1(t_2)-q_2(1)+q_3(1))-\omega^2q_3(t_2)\bigr)(t_2-1)\xi_3,
\end{align*}
where the last equality holds because $t_1\xi_3=\xi_3$ and~$t_3\xi_3=t_2\xi_3$.
By Proposition~\ref{propos_E201}, the elements $\omega_lt_2^{k}\xi_3$, where $l\in\{0,1,2\}$ and $k\in\Z$, are linearly independent in~$E^2_{0,1}$. Therefore
\begin{align*}
 q_1(1)+q_2(t_2)&=0,\\
 q_1(t_2)+q_2(1)-q_3(1)&=0,\\
 q_3(t_2)&=0.
\end{align*}
Since $q_2$ has zero constant term, these equations imply that $q_i=0$ for $i=1,2,3$, and hence $a=0$.
\end{proof}

\begin{remark}\label{remark_basis_E220}
This prove, in fact, shows that the images under~$d^2$ of the generators listed in Proposition~\ref{propos_E220_gen_group} are linearly independent in~$E^2_{0,1}$. Therefore the group~$E^2_{2,0}$ is free abelian, and these elements form a basis of it.
\end{remark}

\begin{proof}[Proof of part~(1) of Theorem~\ref{theorem_3}]
 Recall that
 $$
 0\subseteq \CF_{0,2}\subseteq\CF_{1,1}\subseteq \CF_{2,0}=H_2(\CS)
 $$
 is the filtration associated with the action of~$\CS$ on~$\CB_{0,6}$, with $\CF_{p,q}/\CF_{p-1,q+1}=E^{\infty}_{p,q}$. By Proposition~\ref{propos_E320}, we have $E^{\infty}_{2,0}=0$, and therefore $H_2(\CS)=\CF_{1,1}$. Hence, by Corollary~\ref{cor_F11}, the group~$H_2(\CS)$ is generated by the $9$ abelian cycles~$\omega^l\CA_m$, where $l\in\{0,1,2\}$ and $m\in\{1,2,3\}$. However, under the isomorphism $\CS\cong\I_M$, these $9$ abelian cycles exactly correspond to the $9$ abelian cycles listed in part~(1) of Theorem~\ref{theorem_3}. The required assertion follows.
\end{proof}

\section{Group~$E^2_{3,0}$}\label{section_E230}

The group $E^2_{3,0}$ is the kernel of the differential $d^1\colon E^1_{3,0}\to E^1_{2,0}$.

\begin{propos}\label{propos_lambdas}
The following three elements lie in~$E^2_{3,0}$:
 \begin{align*}
  \lambda_1 &= \bigl(1+\omega+\omega^2\bigr)\Bigl( (1+t_1+t_2+t_3-t_1t_2-t_1t_3-t_2t_3-t_1t_2t_3)(t_1-t_2t_3\nu)(w_1^+-w_1^-)\\
  &{}\hspace{3.1cm}+2(t_2+t_3-t_2t_3-t_1t_2t_3)t_1w_2+2(-1-t_1+t_1t_2+t_1t_3)t_2t_3\nu w_2\\
  &{}\hspace{3.1cm}-2(1-t_1)^2t_2t_3w_4+(1+t_1)(1-t_2)(1-t_3)\bigl(t_2t_3+t_1\nu\bigr)w_6\Bigr),\\
  \lambda_2 &= \bigl(1+\omega+\omega^2\bigr)\Bigl( \bigl(t_1-t_1^{-1}+t_2-t_2^{-1}+t_3-t_3^{-1}\bigr)(t_1-t_2t_3\nu)(w_1^+-w_1^-)\\
  &{}\hspace{3.1cm}+\bigl(-2+t_1-t_1^{-1}+2t_2+2t_3-t_2t_3^{-1}-t_2^{-1}t_3\bigr)t_1w_2\\
  &{}\hspace{3.1cm}+\bigl(-2+t_1^{-1}-t_1+2t_2^{-1}+2t_3^{-1}-t_2t_3^{-1}-t_2^{-1}t_3\bigr)t_2t_3\nu w_2\\
  &{}\hspace{3.1cm}+\bigl(t_1+t_1^{-1}-2\bigr)(t_1+t_2t_3)w_4\\
  &{}\hspace{3.1cm}-(1-t_2)(1-t_3)\bigl(t_2^{-1}+t_3^{-1}\bigr)(t_2t_3+t_1\nu)w_6\Bigr),\\
  \lambda_3 &=\bigl(1+\omega+\omega^2\bigr)(1-t_1)(1-t_2)(1-t_3)\\
  &{}\hspace{1cm}\Bigl(
  \bigl(t_1+t_1^{-1}-t_2t_3^{-1}-t_2^{-1}t_3\bigr)(t_1+t_2t_3\nu)(w_2-\nu w_6)+\bigl(t_1-t_1^{-1}\bigr)(t_1-t_2t_3)w_4\Bigr).
 \end{align*}
 Moreover,
 \begin{multline}\label{eq_lambda3}
   \lambda_3=\bigl(t_1-t_1^{-1}+t_2-t_2^{-1}+t_3-t_3^{-1}\bigr)\lambda_1\\
   {}-(1+t_1+t_2+t_3-t_1t_2-t_1t_3-t_2t_3-t_1t_2t_3)\lambda_2.
 \end{multline}
\end{propos}

\begin{proof}
 The equalities $d^1\lambda_1=0$ and $d^1\lambda_2=0$ are verified by a direct computation using the formulas from Proposition~\ref{propos_d1_0row} and the relations from part~(2) of Proposition~\ref{propos_E1p0}. Let us make several comments on how to shorten this computation.

 First, the element $t_1-t_2t_3\nu$ annihilates both~$z_3^+$ and~$z_3^-$. The coefficients of~$w_1^{\pm}$ in~$\lambda_1$ and~$\lambda_2$ lie in the ideal $\Z[\Gamma](t_1-t_2t_3\nu)$, and $z_3^{\pm}$ do not occur in~$d^1w_i$ for $i=2,4,6$. Therefore the coefficients of~$z_3^{\pm}$ in both~$d^1\lambda_1$ and~$d^1\lambda_2$ vanish.

 Second, $z_1$ also does not occur in~$d^1w_i$ for $i=2,4,6$. Moreover, the coefficient of~$z_1$ in $d^1(t_1-t_2t_3\nu)(w_1^+-w_1^-)$ is equal to
 $$
 \bigl(-\omega+\omega^2+r_2-r_3\bigr)(1-t_1t_2t_3\nu).
 $$
 However,
 $$
 \bigl(1+\omega+\omega^2\bigr)\bigl(-\omega+\omega^2+r_2-r_3\bigr)=0,
 $$
 and $1+\omega+\omega^2$ commutes with every symmetric Laurent polynomial in~$t_1,t_2,t_3$. Consequently, the coefficients of~$z_1$ in both~$d^1\lambda_1$ and~$d^1\lambda_2$ vanish.

 Third, since $t_1-1$ annihilates~$z_6$, and~$z_6$ does not occur in any of~$d^1w_1^{\pm}$, $d^1w_2$, and~$d^1w_6$, it follows that the coefficients of~$z_6$ in both~$d^1\lambda_1$ and~$d^1\lambda_2$ vanish.

 Fourth, $d^1w_6$ contains $(r_1-1)z_9$, while $z_9$ does not occur in any of~$d^1w_1^{\pm}$, $d^1w_2$, and~$d^1w_4$. Since $t_2-1$ annihilates~$z_9$, we obtain that $(t_2-1)(t_3-1)$ annihilates $(r_1-1)z_9$. Consequently, the coefficients of~$z_9$ in both~$d^1\lambda_1$ and~$d^1\lambda_2$ vanish.

 Fifth, the elements~$z_4^+$ and~$z_4^-$ occur in $d^1(w_1^+-w_1^-)$ and $d^1w_2$ only through the combination $z_4^+-z_4^-$. This is not the case for~$d^1w_6$. Nevertheless, $w_6$ occurs in~$\lambda_1$ and~$\lambda_2$ with coefficients from the ideal $\Z[\Gamma](t_2t_3+t_1\nu)$, and
 $$
 d^1(t_2t_3+t_1\nu)w_6=-(t_1+t_2t_3\nu)(z_4^+-z_4^-)+(t_2t_3+t_1\nu)(r_1-1)z_9.
 $$
 It follows that $z_4^+$ and~$z_4^-$ occur in~$d^1\lambda_1$ and~$d^1\lambda_2$ only through the combination~$z_4^+-z_4^-$.

 Thus, to prove that $d^1\lambda_1=0$ and $d^1\lambda_2=0$, it remains for us to show that the coefficients of~$z_4^+-z_4^-$ and of~$z_2$ in~$d^1\lambda_1$ and~$d^1\lambda_2$ are zero. This is verified by a direct computation: in this case the author knows of no way to shorten this computation.

 Finally, the formula~\eqref{eq_lambda3} is also verified by a direct computation. Since $d^1\lambda_1=0$ and~$d^1\lambda_2=0$, it follows that $d^1\lambda_3=0$.
\end{proof}

\begin{remark}
The way the formulas in Proposition~\ref{propos_lambdas} were obtained is less mysterious than it appears. In fact, the author computed the generators of the $\Lambda$-module~$E^2_{3,0}$ explicitly, in the same way as the generators~$\rho_1$ and~$\rho_2$ of~$E^2_{2,0}$ were computed above. The result is that~$E^2_{3,0}$ is generated as a $\Lambda$-module by the two elements~$\lambda_1$ and~$\lambda_2$ in Proposition~\ref{propos_lambdas}. Thus the formulas for these elements are not the result of a guess or a computer search, but rather the result of a lengthy computation similar to the computation in Section~\ref{section_E220}. However, we will not use the fact that~$\lambda_1$ and~$\lambda_2$ generate~$E^2_{3,0}$ in what follows: it will be enough for us that they lie in~$E^2_{3,0}$, that is, $d^1\lambda_1=0$ and~$d^1\lambda_2=0$. Therefore, once the formulas for~$\lambda_1$ and~$\lambda_2$ have been written out, there is no point in presenting the aforementioned multi-page computation, since the equalities $d^1\lambda_1=0$ and~$d^1\lambda_2=0$ are easier to verify directly. When the elements~$\lambda_1$ and~$\lambda_2$ are already found, it is natural to consider their linear combination that does not contain~$w_1^{\pm}$; this is how the element~$\lambda_3$ was found.
\end{remark}

\section{Computation of $d^3\lambda_3$}\label{section_d3}

\begin{propos}\label{propos_d3lambda3}
 We have $d^2\lambda_3=0$ and
\begin{equation*}%\label{eq_d3lambda3}
d^3\lambda_3=\bigl(1+\omega+\omega^2\bigr)(1-t_2)(1-t_3)(t_2-t_3)\bigl(1-t_2^{-1}t_3^{-1}\bigr)A_2.
\end{equation*}
\end{propos}

Before proving this proposition, let us derive part~(2) of Theorem~\ref{theorem_3} from it.

\begin{proof}[Proof of part~(2) of Theorem~\ref{theorem_3}]
The element $A_2\in E^1_{0,2}$ represents the abelian cycle $$\CA_2=\CA\left(T_{\delta_1},T_{\gamma_{2\bar{3}}}T_{\tgamma_{\bar{2}3}}^{-1}\right)\in H_2(\CS).$$
From Proposition~\ref{propos_d3lambda3} it follows that
$$
\bigl(1+\omega+\omega^2\bigr)(1-t_2)(1-t_3)(t_2-t_3)\bigl(1-t_2^{-1}t_3^{-1}\bigr)\CA_2=0
$$
in $H_2(\CS)\cong H_2(\I_M)$. Multiplying this relation by~$t_1t_2t_3$ and using that $t_1\CA_2=\CA_2$, we obtain the required relation~\eqref{eq_thm3_relA2}.
\end{proof}

\begin{proof}[Proof of Proposition~\ref{propos_d3lambda3}]
Similarly to the proofs of Propositions~\ref{propos_d2rho1} and~\ref{propos_d2rho2},
the proof is based on an explicit computation in the bicomplex $C_*(\CB_{0,6})\otimes_{\CS}\CR_*$ modulo $C_*(K)\otimes_{\CS}\CR_*$, where $K$ is the subcomplex introduced in Section~\ref{section_d2rho}.

 By Lemma~\ref{lem_dr}, to show that $d^2\lambda_3=0$ and to compute~$d^3\lambda_3$ it suffices to find elements
 \begin{gather*}
 \begin{aligned}
 \theta_{3,0}&\in C_3(\CB_{0,6})\otimes_{\CS}\CR_0=C_3(\CB_{0,6}),&&&&&\theta_{2,1}&\in C_2(\CB_{0,6})\otimes_{\CS}\CR_1,\\ \theta_{1,2}&\in C_1(\CB_{0,6})\otimes_{\CS}\CR_2,&&&&& \theta_{0,3}&\in C_0(\CB_{0,6})\otimes_{\CS}\CR_3,
 \end{aligned}\\
 \psi\in C_0(\CB_{0,6})\otimes_{\CS}\CR_2
 \end{gather*}
 such that $\theta_{3,0}$ represents~$\lambda_3$ and
 \begin{gather*}
  \partial'\theta_{3,0}+\partial''\theta_{2,1}=0,\\
  \partial'\theta_{2,1}+\partial''\theta_{1,2}\Kequiv 0,\\
  \partial'\theta_{1,2}+\partial''\theta_{0,3}\Kequiv \psi,\\\
  \partial''\psi=0.
 \end{gather*}
Then we will obtain that $d^2\lambda=0$, and the element~$\psi$ will represent~$d^3\lambda_3$.

Let us introduce some notation. We consider the mapping classes
\begin{align*}
 \tau_1&=T_{\bar{2}3}^{-1},&\tau_2&=T_{\bar{3}1}^{-1},&\tau_3&=T_{\bar{1}2}^{-1},\\
 \ttau_1&=T_{2\bar{3}}^{-1},&\ttau_2&=T_{3\bar{1}}^{-1},&\ttau_3&=T_{1\bar{2}}^{-1},\\
 f_1&=T_{\varepsilon_1},&f_2&=T_{\varepsilon_2},&f_3&=T_{\varepsilon_3}.
\end{align*}
Then $\kappa(\tau_i)=\kappa\bigl(\ttau_i\bigr)=t_i$ and $\kappa(f_i)=t_it_j^{-1}t_k^{-1}$ if $i,j,k$ is a permutation of~$1,2,3$.
Throughout this section, we will multiply use that certain pairs of these mapping classes commute, specifically:
\begin{itemize}
 \item $\tau_1$, $\tau_2$, and~$\tau_3$ commute with each other,
 \item $\ttau_1$, $\ttau_2$, and~$\ttau_3$ commute with each other,
 \item $\tau_i\ttau_i=\ttau_i\tau_i$ for all~$i$,
 \item $f_i$ commutes with each of the classes~$\tau_j$ and~$\ttau_j$ whenever $i\ne j$.
\end{itemize}
Moreover,
\begin{itemize}
 \item the conjugation by~$\omega$ permutes cyclically each of the triples $(\tau_1,\tau_2,\tau_3)$, $(\ttau_1,\ttau_2,\ttau_3)$, and~$(f_1,f_2,f_3)$ so that $\omega\tau_i\omega^{-1}=\tau_{i+1}$ (where indices are taken modulo~$3$), and similarly for~$\ttau_i$ and~$f_i$,
 \item the conjugation by~$r_i$ fixes each of the mapping classes~$\tau_i$, $\ttau_i$, and~$f_i$, and swaps the mapping classes in each of the pairs~$(\tau_j,\tau_k)$, $(\ttau_j,\ttau_k)$, and~$(f_j,f_k)$, where $i,j,k$ is a permutation of~$1,2,3$.
 \item $\nu\tau_i\nu=\ttau_i^{-1}$ and $\nu f_i\nu=f_i^{-1}$ for all~$i$.
\end{itemize}
Further, for each permutation $i,j,k$ of~$1,2,3$, we set
\begin{gather*}
 p_i=\tau_i+\tau_i^{-1}-\tau_j\tau_k^{-1}-\tau_j^{-1}\tau_k,\\
 q_i=\bigl(\tau_i-\tau_i^{-1}\bigr)(\tau_i-\tau_j\tau_k).
\end{gather*}
We will need the following property of these Laurent polynomials, which can be checked by a direct computation: if $i,j,k$ is a permutation of~$1,2,3$, then
 \begin{equation}\label{eq_pq_identity}
 \tau_i\tau_jp_k+\tau_jp_j+q_i=0.
 \end{equation}

We set
$$
\theta_{3,0}=\bigl(1+\omega+\omega^2\bigr)(1-\tau_1)(1-\tau_2)(1-\tau_3)\bigl(
p_1(\tau_1+\tau_2\tau_3\nu)(W_2-\nu W_6)+q_1W_4
\bigr).
$$
Then $\theta_{3,0}$ represents~$\lambda_3$.
Let us compute~$\partial'\theta_{3,0}$ using the formulas for $\partial W_i$ from Table~\ref{table_3cells}. Recall that the differential~$\partial'$ in the $0$th  row of the bicomplex $C_*(\CB_{0,6})\otimes_{\CS}\CR_*$ coincides with the differential~$\partial$ of~$C_*(\CB_{0,6})$.
We have
\begin{gather}
\begin{aligned}\label{eq_part_l3_aux1}
\partial(\tau_1+\tau_2\tau_3\nu)(W_2-\nu W_6)&= -\bigl(\tau_1\omega^2+\tau_2\tau_3\omega\bigr)Z_2 - \bigl(\tau_1\omega+\tau_2\tau_3\omega^2\bigr)\nu Z_2\\
&\quad{}+ \left(\tau_2\tau_3f_1-\tau_1\right)Z_4^+
+
\left(\tau_1 f_1^{-1}-\tau_2\tau_3\right)\nu Z_4^+\\
&\quad{} +\tau_2\tau_3(1-r_1)Z_9+\tau_1(1-r_1)\nu Z_9,
\end{aligned}\\
\partial W_4=-(1+\nu)(Z_2+Z_6).\label{eq_part_l3_aux2}
\end{gather}
The mapping class~$\ttau_1$ fixes~$Z_6$, and $\ttau_2$ fixes~$Z_9$. Therefore the mapping classes~$\tau_1$ and~$\tau_2$ fix the cells~$\nu Z_6$ and~$\nu Z_9$, respectively. Since the product $(1-\tau_1)(1-\tau_2)(1-\tau_3)$ commutes with~$r_1$, it annihilates all summands in~\eqref{eq_part_l3_aux1} and~\eqref{eq_part_l3_aux2} involving~$\nu Z_6$ or~$\nu Z_9$. Thus we have
\begin{align*}
\partial' \theta_{3,0}&=\bigl(1+\omega+\omega^2\bigr)(1-\tau_1)(1-\tau_2)(1-\tau_3)\\
&\qquad\quad{}\Bigl(-\bigl(p_1\tau_1\omega^2+p_1\tau_2\tau_3\omega+q_1\bigr)Z_2 - \bigl(p_1\tau_1\omega+p_1\tau_2\tau_3\omega^2+q_1\bigr)\nu Z_2\\
&\qquad\qquad{}+ p_1\tau_1\left(\tau_1^{-1}\tau_2\tau_3f_1-1\right)\bigl(1-f_1^{-1}\nu\bigr)Z_4^+ -q_1Z_6+p_1\tau_2\tau_3(1-\omega^2r_3)Z_9\Bigr),
\end{align*}
where computing the coefficient of the last summand we have used that $r_1=\omega^2r_3$.
Since the conjugation by~$\omega$ permutes cyclically~$\tau_1$, $\tau_2$, and~$\tau_3$ (and hence $p_1$, $p_2$, and~$p_3$), we write
\begin{align*}
\partial' \theta_{3,0}&=\bigl(1+\omega+\omega^2\bigr)(1-\tau_1)(1-\tau_2)(1-\tau_3)\\
&\qquad\quad{}\Bigl(-(p_2\tau_2+p_3\tau_1\tau_2+q_1)Z_2 - (p_3\tau_3+p_2\tau_1\tau_3+q_1)\nu Z_2\\
&\qquad\qquad{}+ p_1\tau_1\left(\tau_1^{-1}\tau_2\tau_3f_1-1\right)\bigl(1-f_1^{-1}\nu\bigr)Z_4^+ -q_1Z_6+p_3\tau_1\tau_2\omega^2Z_9 -p_2\tau_1\tau_3r_3Z_9\Bigr)\\
{}&=\bigl(1+\omega+\omega^2\bigr)(1-\tau_1)(1-\tau_2)(1-\tau_3)\\
&\qquad\quad{}\Bigl(p_1\tau_1\left(\tau_1^{-1}\tau_2\tau_3f_1-1\right)\bigl(1-f_1^{-1}\nu\bigr)Z_4^+ -q_1Z_6+p_3\tau_1\tau_2\omega^2Z_9 -p_2\tau_1\tau_3r_3Z_9\Bigr),
\end{align*}
where the last equality holds by~\eqref{eq_pq_identity}.

We now introduce notation
$$
g_i=\ttau_i\tau_i^{-1},\qquad h_i=f_i^{-1}\tau_i\tau_j^{-1}\tau_k^{-1},
$$
where $i,j,k$ is a permutation of~$1,2,3$. Since $\kappa(\tau_i)=\kappa(\ttau_i)=t_i$ and $\kappa(f_i)=t_it_j^{-1}t_k^{-1}$, we obtain that $g_i$ and~$h_i$ belong to~$\CS$.

If $Q$ is either~$Z_4^+$ or~$f_1^{-1}\nu Z_4^+$, then by Lemma~\ref{lem_d'd''} we have
$$
\left(\tau_1^{-1}\tau_2\tau_3f_1-1\right)Q=\left(h_1^{-1}-1\right)Q=\partial''(Q\otimes[h_1]).
$$
Further, $\ttau_1$ fixes each of the cells~$Z_6$, $\omega^2Z_9$, and~$r_3Z_9$. Hence if $Q$ is any of these three cells, then
$$
(\tau_1-1)Q=\bigl(g^{-1}_1-1\bigr) Q=\partial''(Q\otimes[g_1]).
$$
Consequently,
$$\partial'\theta_{3,0}=-\partial''\theta_{2,1},$$
where
\begin{align*}
\theta_{2,1}&=\bigl(1+\omega+\omega^2\bigr)(1-\tau_2)(1-\tau_3)\bullet{}\\
&\qquad\quad{}\Bigl((\tau_1-1)\tau_1 p_1\bullet\bigl(Z_4^+\otimes[h_1]-\bigl(f_1^{-1}\nu Z_4^+\bigr)\otimes[h_1]\bigr)\\
&\qquad\qquad{}-q_1\bullet (Z_6\otimes [g_1])+p_3\tau_1\tau_2\bullet\bigl(\bigl(\omega^2Z_9\bigr)\otimes[g_1]\bigr) -p_2\tau_1\tau_3\bullet\bigl((r_3Z_9)\otimes[g_1]\bigr)\Bigr).
\end{align*}

Next, let us compute $\partial'\theta_{2,1}$ up to elements of~$C_1(K)\otimes_{\CS}\CR_1$.  By the formulas from Table~\ref{table_2cells}, we have
\begin{align*}
 \partial Z_4^+&=\bigl(r_3-\omega^2\bigr)Y_2+Y_5,\\
 \partial Z_6&=(r_1-1)Y_2+Y_7,\\
 \partial Z_9&\Kequiv(\omega-f_1r_3)Y_2.
\end{align*}
Since $f_1^{-1}\nu Y_5=Y_5$, the terms in~$\partial'\theta_{2,1}$ involving~$Y_5$ cancel. The Laurent polynomial~$q_1$ is divisible by~$\tau_1-1$. The mapping class~$\tau_1$ fixes~$Y_7$ and commutes with~$g_1$. Hence the terms in~$\partial'\theta_{2,1}$ involving~$Y_7$ also cancel. Since $\omega^2r_3=r_3\omega=r_1$, it follows that
\begin{align*}
  \partial \omega^2Z_9&\Kequiv (1-f_3r_1)Y_2,\\
  \partial r_3 Z_9&\Kequiv(r_1-f_2)Y_2.
\end{align*}
Thus
\begin{align*}
\partial' \theta_{2,1}&\Kequivs\bigl(1+\omega+\omega^2\bigr)(1-\tau_2)(1-\tau_3)\bullet{}\\
&\qquad\quad{}\Bigl((\tau_1-1)\tau_1 p_1\bullet\bigl(\bigl(r_3Y_2-\omega^2Y_2\bigr)\otimes[h_1]-\bigl(f_1^{-1}r_3\nu Y_2-f_1^{-1}\omega^2\nu Y_2\bigr)\otimes[h_1]\bigr)\\
&\qquad\qquad{}-q_1\bullet\bigl((r_1Y_2-Y_2)\otimes [g_1]\bigr)+p_3\tau_1\tau_2\bullet\bigl((Y_2-f_3r_1Y_2)\otimes[g_1]\bigr)\\
&\qquad\qquad{}-p_2\tau_1\tau_3\bullet\bigl((r_1Y_2-f_2Y_2)\otimes[g_1]\bigr)\Bigr).
\end{align*}
Since $\tau_2\bigl(f_1^{-1}r_3\nu\bigr)=\bigl(f_1^{-1}r_3\nu\bigr)\ttau_1^{-1}$ and $\ttau_1$ fixes~$Y_2$, we obtain that $\tau_2$ fixes $f_1^{-1}r_3\nu Y_2$. Moreover, $\tau_2$ commutes with~$h_1$. Hence the factor~$1-\tau_2$ annihilates $\bigl(f_1^{-1}r_3\nu Y_2\bigr)\otimes[h_1]$. Similarly, the factor~$1-\tau_3$ annihilates $\bigl(f_1^{-1}\omega^2\nu Y_2\bigr)\otimes[h_1]$. Therefore the terms in~$\partial'\theta_{2,1}$ containing~$\nu Y_2$ cancel.
Further, using the relation $r_3=\omega r_1$, moving the powers of~$\omega$ to the left as in the computation of~$\partial'\theta_{3,0}$, collecting together the terms corresponding to the same cells, and then using~\eqref{eq_pq_identity}, we obtain that
\begin{align*}
\partial' \theta_{2,1}&\Kequiv\bigl(1+\omega+\omega^2\bigr)(1-\tau_2)(1-\tau_3)\bullet{}\\
&\ {}\Bigl(
p_2\tau_2\bullet\Bigl((1-\tau_1)\bullet(Y_2\otimes[h_2])-Y_2\otimes[g_1]+\tau_1\tau_2^{-1}\tau_3\bullet\bigl((f_2Y_2)\otimes[g_1]\bigr)\Bigr)\\
&\ \ {}-p_3\tau_3\bullet\Bigl((1-\tau_1)\bullet((r_1Y_2)\otimes[h_3])-(r_1Y_2)\otimes[g_1]+\tau_1\tau_2\tau_3^{-1}\bullet\bigl((f_3r_1Y_2)\otimes[g_1]\bigr)\Bigr)\Bigr).
\end{align*}
Since $\ttau_1$ fixes~$Y_2$, $\tau_1$  commutes with~$h_2$, and $f_2$ commutes with~$g_1$, we have
\begin{align*}
 (1-\tau_1)&\bullet(Y_2\otimes[h_2])-Y_2\otimes[g_1]+\tau_1\tau_2^{-1}\tau_3\bullet\bigl((f_2Y_2)\otimes[g_1]\bigr){}\\
 {}&=Y_2\otimes[h_2]-\bigl(g_1^{-1}Y_2\bigr)\otimes[h_2]-Y_2\otimes[g_1]+\bigl(h_2^{-1}Y_2\bigr)\otimes\bigl[h_2^{-1}g_1h_2\bigr]\\
 {}&=Y_2\otimes\left(
 [h_2]-g_1[h_2]-[g_1]+h_2\bigl[h_2^{-1}g_1h_2\bigr]
 \right)\\
 {}&=-Y_2\otimes d\left([g_1|h_2]-\bigl[h_2\big|h_2^{-1}g_1h_2\bigr]\right)\\
 {}&=\partial''\Bigl(Y_2\otimes \left([g_1|h_2]-\bigl[h_2\big|h_2^{-1}g_1h_2\bigr]\right)\Bigr).
\end{align*}
The same holds for the second summand in the expression for~$\partial'\theta_{2,1}$, with $Y_2$ and $h_2$ replaced by~$r_1Y_2$ and~$h_3$, respectively. Thus
$$
\partial'\theta_{2,1}\Kequiv -\partial''\theta_{1,2},
$$
where
\begin{multline*}
\theta_{1,2}=\bigl(1+\omega+\omega^2\bigr)(1-\tau_2)(1-\tau_3)\bullet\Bigl(-p_2\tau_2\bullet\bigl(Y_2\otimes \left([g_1|h_2]-\bigl[h_2\big|h_2^{-1}g_1h_2\bigr]\right)\bigr)\\
{}+p_3\tau_3\bullet\bigl((r_1Y_2)\otimes \left([g_1|h_3]-\bigl[h_3\big|h_3^{-1}g_1h_3\bigr]\right)\bigr)\Bigr).
\end{multline*}

Finally, let us compute $\partial'\theta_{1,2}$ modulo elements of~$C_0(\CB_{0,6})\otimes_{\CS}\CR_2$. Since $\partial Y_2\Kequivs -X_2$ and $r_1X_2=f_2X_2$, we have:
\begin{multline}\label{eq_part_K12}
\partial'\theta_{1,2}\Kequivs\bigl(1+\omega+\omega^2\bigr)(1-\tau_2)(1-\tau_3)\bullet\Bigl(p_2\tau_2\bullet\bigl(X_2\otimes \left([g_1|h_2]-\bigl[h_2\big|h_2^{-1}g_1h_2\bigr]\right)\bigr)\\
{}-p_3\tau_3\bullet\bigl((f_2X_2)\otimes \left([g_1|h_3]-\bigl[h_3\big|h_3^{-1}g_1h_3\bigr]\right)\bigr)\Bigr).
\end{multline}
The mapping class~$\tau_1$ fixes each of the cells~$X_2$ and~$f_2X_2$, and commutes with each of the elements~$g_1$, $h_2$, and~$h_3$. Therefore the right-hand side of~\eqref{eq_part_K12} will not change if we replace~$p_2$ with the Laurent polynomial
$$F=\tau_2+\tau_2^{-1}-\tau_3-\tau_3^{-1}=(\tau_2-\tau_3)\bigl(1-\tau_2^{-1}\tau_3^{-1}\bigr),$$
and~$p_3$ with~$-\tau_1F$.
Thus
\begin{multline*}
\partial'\theta_{1,2}\Kequiv\bigl(1+\omega+\omega^2\bigr)(1-\tau_2)(1-\tau_3)(\tau_2-\tau_3)\bigl(1-\tau_2^{-1}\tau_3^{-1}\bigr)\tau_2\bullet{}\\
\Bigl(X_2\otimes \left([g_1|h_2]-\bigl[h_2\big|h_2^{-1}g_1h_2\bigr]\right)+\tau_1\tau_2^{-1}\tau_3\bullet\bigl((f_2X_2)\otimes \left([g_1|h_3]-\bigl[h_3\big|h_3^{-1}g_1h_3\bigr]\right)\bigr)\Bigr).
\end{multline*}
Since $f_2$ commutes with~$g_1$, we have
$$
\bigl(\tau_1\tau_2^{-1}\tau_3\bigr)g_1\bigl(\tau_1\tau_2^{-1}\tau_3\bigr)^{-1}=h_2^{-1}g_1h_2.
$$
We set
$$
\tilde{h}_3=\bigl(\tau_1\tau_2^{-1}\tau_3\bigr)h_3\bigl(\tau_1\tau_2^{-1}\tau_3\bigr)^{-1}=\tau_1^{-1}\tau_2^{-1}\tau_3f_3^{-1}.
$$
Then
\begin{align*}
\partial'\theta_{1,2}
&\Kequivs\bigl(1+\omega+\omega^2\bigr)(1-\tau_2)(1-\tau_3)(\tau_2-\tau_3)\bigl(1-\tau_2^{-1}\tau_3^{-1}\bigr)\tau_2\bullet{}\\
&\quad\ {}\Bigl(X_2\otimes \bigl([g_1|h_2]-\bigl[h_2\big|h_2^{-1}g_1h_2\bigr]+
h_2\bigl[h_2^{-1}g_1h_2|\tilde{h}_3\bigr]-h_2\bigl[\tilde{h}_3\big|\tilde{h}_3^{-1}h_2^{-1}g_1h_2\tilde{h}_3\bigr]\bigr)\Bigr).
\end{align*}
Set
\begin{align*}
\theta_{0,3}
&=\bigl(1+\omega+\omega^2\bigr)(1-\tau_2)(1-\tau_3)(\tau_2-\tau_3)\bigl(1-\tau_2^{-1}\tau_3^{-1}\bigr)\tau_2\bullet{}\\
&\qquad\qquad {}\Bigl(X_2\otimes \bigl(
\bigl[g_1\big|h_2\big|\tilde{h}_3\bigr]-
\bigl[h_2\big|h_2^{-1}g_1h_2|\tilde{h}_3\bigr]+\bigl[h_2\big|\tilde{h}_3\big|g_1\bigr]\bigr)\Bigr).
\end{align*}
Now the lantern relation~\eqref{eq_lantern_eps_jk} reads as
\begin{equation*}
h_2\tilde{h}_3=T_{\delta_1'}g_1,
\end{equation*}
which, in particular, implies that $h_2\tilde{h}_3$ commutes with~$g_1$. Using this,  one can easily check that
\begin{multline}\label{eq_K12K03}
\partial'\theta_{1,2}+\partial''\theta_{0,3}\Kequivs
\bigl(1+\omega+\omega^2\bigr)(1-\tau_2)(1-\tau_3)(\tau_2-\tau_3)\bigl(1-\tau_2^{-1}\tau_3^{-1}\bigr)\tau_2\bullet{}\\
\qquad\qquad {}\Bigl(X_2\otimes \bigl(
\bigl[g_1\big|T_{\delta_1'}g_1\bigr]-\bigl[T_{\delta_1'}g_1\big|g_1\bigr]\bigr)\Bigr).
\end{multline}
 Since the mapping classes~$T_{\delta_1'}$ and~$g_1$ belong to~$\CS_{X_2}$ and commute, we see that the element $X_2\otimes \bigl(
\bigl[g_1\big|T_{\delta_1'}g_1\bigr]-\bigl[T_{\delta_1'}g_1\big|g_1\bigr]\bigr)$ lies in the kernel of~$\partial''$ and represents in~$E^1_{0,2}$ the class
$$
\iota_{X_2}\left(\CA\bigl(g_1,T_{\delta_1'}g_1\bigr)\right)=\iota_{X_2}\left(\CA\bigl(T_{\delta_1'},g_1^{-1}\bigr)\right)=
\iota_{X_2}\left(\CA\bigl(T_{\delta_1'},T_{2\bar{3}}T_{\bar{2}3}^{-1}\bigr)\right)=t_2^{-1}A_2.
$$
Therefore the right-hand side of~\eqref{eq_K12K03} lies in the kernel of~$\partial''$ and represents the class
$$
\bigl(1+\omega+\omega^2\bigr)(1-t_2)(1-t_3)(t_2-t_3)\bigl(1-t_2^{-1}t_3^{-1}\bigr)A_2.
$$
The proposition follows.
\end{proof}

\section{Relations in $H_2(\CS)$ arising from~$\lambda_1$ and~$\lambda_2$.}\label{section_end_thm3}

By Propositions~\ref{propos_A1A2} and~\ref{propos_A3}, the abelian cycles~$\CA_1$ and~$\CA_2$ lie in~$\CF_{0,2}=E^{\infty}_{0,2}$ and are represented by $A_1$ and~$A_2$, respectively, and the abelian cycle~$\CA_3$ lies in~$\CF_{1,1}$ and goes to~$\eta_4$ under the projection $\CF_{1,1}\twoheadrightarrow E^{\infty}_{1,1}$. Moreover, by Corollaries~\ref{cor_F02} and~\ref{cor_F11}, the group~$\CF_{0,2}$ is generated as a $\Z[\Gamma]$-module by~$\CA_1$ and~$\CA_2$, while $\CF_{1,1}$ is generated as a $\Z[\Gamma]$-module by~$\CA_1$, $\CA_2$, and~$\CA_3$. (In fact, $\CF_{1,1}=H_2(\CS)$ by Proposition~\ref{propos_E320}, but this is currently not important for us.) For any $\lambda\in E^3_{0,3}$, the element~$d^3\lambda$ lies in the kernel of the map $E^3_{0,2}\twoheadrightarrow E^{\infty}_{0,2}$, and therefore yields a relation between~$\CA_1$ and~$\CA_2$. In the previous section, we derived such a relation corresponding to the class~$\lambda_3$. We could try to write analogous relations for~$\lambda_1$ and~$\lambda_2$. However, the problem is that $\lambda_1$ and~$\lambda_2$ do not belong to~$E^3_{0,3}$. Of course, we could compute $d^2\lambda_i$ for $i=1,2$, but this would give us only the following partial result: if the computation gives $d^2\lambda_i=b_i\eta_4$, where $b_i\in\Z[\Gamma]$, then we obtain that $b_i\eta_4=0$ in~$E^{\infty}_{1,1}$, and hence the class~$b_i\CA_3$ belongs to $\CF_{0,2}=\Z[\Gamma]\CA_1+\Z[\Gamma]\CA_2$, i.e.,
\begin{equation}\label{eq_cc'c''}
b_i\CA_3=b_i'\CA_1+b_i''\CA_2,\qquad b_i',b_i''\in\Z[\Gamma].
\end{equation}
However, the coefficients~$b_i'$ and~$b_i''$ would remain unknown. Nevertheless, the bicomplex $C_*(\CB_{0,6})\otimes_{\CS}\CR_*$ carries information not only about the differentials of the spectral sequence~$E^*_{*,*}$ but also about how the associated factors~$E^{\infty}_{0,2}$ and~$E^{\infty}_{1,1}$ are assembled into the whole group~$H_2(\CS)$. Therefore, a computation  in this bicomplex similar to the one performed in the previous section can actually be used to obtain precise relations in~$H_2(\CS)$. Another difficulty, however, is that the formulas for~$\lambda_1$ and~$\lambda_2$ are rather long, so the corresponding computation would be very cumbersome. This is why we are going to perform the computation modulo a certain subcomplex $L\subseteq\CB_{0,6}$ which is substantially larger than the subcomplex~$K$  used above. As a result, our computation will allow us to find the coefficients~$b_i'$ of~$\CA_1$ in~\eqref{eq_cc'c''} but not the coefficients~$b_i''$ of~$\CA_2$. Nevertheless, this will be sufficient for our purposes.

Specifically, the main result of this section is the following proposition.

\begin{propos}\label{propos_rel_from_lambda}
 Let $\lambda\in E^1_{3,0}$ be an element of the form
 \begin{equation}\label{eq_lambda}
 \lambda=F(1+\omega+\omega^2)(t_1-t_2t_3\nu)\bigl(w_1^+ - w_1^-\bigr)+c_2w_2+c_4(t_1-1)^2w_4+c_6w_6,
 \end{equation}
 where $F\in\Lambda$ and $c_2,c_4,c_6\in\Z[\Gamma]$. Suppose that $d^1\lambda=0$, i.e., $\lambda\in E^2_{3,0}$. Then
 $$
 F(1+\omega+\omega^2)(\CA_1-\CA_3)\in\Z[\Gamma]\CA_2.
 $$
\end{propos}

Since both~$\lambda_1$ and~$\lambda_2$ have the form~\eqref{eq_lambda} and lie in~$E^2_{3,0}$, this proposition immediately implies the following corollary.

\begin{cor}\label{cor_sigmaM12+A3}
 The following two homology classes in~$H_2(\CS)$ belong to~$\Z[\Gamma]\CA_2$:
 \begin{gather*}
 (1+t_1+t_2+t_3-t_1t_2-t_1t_3-t_2t_3-t_1t_2t_3)\bigl(1+\omega+\omega^2\bigr)(\CA_1-\CA_3),\\
 \bigl(t_1-t_1^{-1}+t_2-t_2^{-1}+t_3-t_3^{-1}\bigr)\bigl(1+\omega+\omega^2\bigr)(\CA_1-\CA_3).
 \end{gather*}
\end{cor}

This corollary is enough to prove part~(3) of Theorem~\ref{theorem_3}.

\begin{proof}[Proof of part~(3) of Theorem~\ref{theorem_3}]
From relations~\eqref{eq_rel_BE1} it follows that $\Z[\Gamma]\CA_2$ is generated as a $\Lambda$-module by the $3$ abelian cycles $\CA_2$, $\omega\CA_2$, and~$\omega^2\CA_2$. Therefore, by Corollary~\ref{cor_sigmaM12+A3}, the homology classes
\begin{align*}
\sigma^M_1&=(1+t_1+t_2+t_3-t_1t_2-t_1t_3-t_2t_3-t_1t_2t_3)\bigl(1+\omega+\omega^2\bigr)\CA_1,\\
\sigma^M_2&=\bigl(t_1-t_1^{-1}+t_2-t_2^{-1}+t_3-t_3^{-1}\bigr)\bigl(1+\omega+\omega^2\bigr)\CA_1
\end{align*}
lie in the $\Lambda$-submodule of $H_2(\I_M)=H_2(\CS)$ generated by the $6$ abelian cycles~$\omega^l\CA_m$ with $l\in\{0,1,2\}$ and~$m\in\{2,3\}$. This is exactly assertion~(3) of Theorem~\ref{theorem_3}.
\end{proof}

In the rest of this section we prove Proposition~\ref{propos_rel_from_lambda}.

Let $L$ be the subcomplex of~$\CB_{0,6}$ consisting of all cells~$Q$ such that
$Q$ contains no vertex in the $\G$-orbit~$\G  X_1$, nor any edge in the $\G$-orbits~$\G Y_7$ and~$\G Y_8$. In other words, $L$ consists of all cells in the $\G$-orbits of the following cells:
$$
X_2,\ X_3,\ Y_2,\ Y_4,\ Y_5,\ Y_9,\ Z_2,\ Z_4^+,\ Z_4^-,\ Z_9,\ W_2,\ W_6.
$$

\begin{remark}
 Such a choice of the subcomplex~$L$ may seem artificial. However, it has the following justification. First of all, we do not want to include the $\G$-orbit of the vertex~$X_1$ in~$L$, since it is important for us not to lose information about the coefficient of the abelian cycle~$\CA_1$. Second, we definitely want to include the $\G$-orbits of the vertices~$X_2$ and~$X_3$, since to simplify the computations it is important for us that~$L$ be sufficiently large. We could take as~$L$ the subcomplex consisting of all cells that do not contain vertices from~$\G X_1$; however, this leads to the following problem. We would like the action of~$\CS$ on~$L$ to be as close to free as possible, so that we can use an analogue of Lemma~\ref{lem_dr}. However, there are three $\G$-orbits of cells~$Q$ not containing vertices from~$\G X_1$ whose stabilizers~$\CS_Q$ are nontrivial: namely, $\G X_2$, $\G Y_7$, and~$\G Y_8$. There is nothing we can do about the vertices in~$\G X_2$: as already stated above, we definitely want to include them in~$L$. Therefore we will have to overcome the difficulties caused by the nontriviality of the stabilizer~$\CS_{X_2}$. (In fact, it is precisely these difficulties that will be responsible for the fact that we do not compute the coefficients~$b_i''$ in relations~\eqref{eq_cc'c''} precisely.) However, in order to avoid additional difficulties, we do not want to include the orbits of the cells~$Y_7$ and~$Y_8$ in~$L$. But then we cannot include the orbits of the cells~$Z_6$, $Z_8$, and~$W_4$ in~$L$ either, since they contain edges from~$\G Y_7\cup\G Y_8$.
\end{remark}

Consider the following cycles in~$\Tot_2\bigl(C_*(\CB_{0,6})\otimes_{\CS}\CR_*\bigr)$:
 \begin{align*}
  \fA_1&=X_1\otimes\bigl(\bigr[T_{\delta_2}^{-1}\bigr|T_{\delta_3}^{-1}\bigr]-\bigl[T_{\delta_3}^{-1}\big|T_{\delta_2}^{-1}\bigr]\bigr),\\
  \fA_3&=(T_{31}X_2)\otimes\left(\bigl[T_{\talpha_3}\big|T_{\delta_1}^{-1}\bigr]-\bigl[T_{\delta_1}^{-1}\big|T_{\talpha_3}\bigr]\right)-Y_8\otimes\bigl[T_{\delta_1}^{-1}\bigr].
  \end{align*}
The cycle~$\fA_1$ lies in $C_0(\CB_{0,6})\otimes_{\CS}\CR_2$  and represents~$A_1$ in~$E^1_{0,2}$. Hence, it represents~$\CA_1$ in~$H_2(\CS)$. By Proposition~\ref{propos_A3}, the cycle~$\fA_3$ represents the class~$\CA_3$.

\begin{remark}
 The fact that we write $T_{\delta_i}^{-1}$ rather than~$T_{\delta_i}$ in the formula for~$\fA_1$  is completely immaterial. If $T_{\delta_2}^{-1}$ and~$T_{\delta_3}^{-1}$  are replaced by~$T_{\delta_2}$ and~$T_{\delta_3}$, respectively, the resulting cycle represents the same homology class. Nevertheless, it is convenient for us to define the cycle~$\fA_1$ in precisely this way, since this is exactly the cycle that will arise in the subsequent computation.
\end{remark}

\begin{lem}\label{lem_theta_for_mu}
 There exists an element
 $$
 \theta=\theta_{3,0}+\theta_{2,1}+\theta_{1,2}+\theta_{0,3}\in \Tot_3\bigl(C_*(\CB_{0,6})\otimes_{\CS}\CR_*\bigr)
 $$
 with $\theta_{p,q}\in C_p(\CB_{0,6})\otimes_{\CS}\CR_q$ such that $\theta_{3,0}$ represents the class
 $$
 \mu=\bigl(1+\omega+\omega^2\bigr)(t_1-t_2t_3\nu)\bigl(w_1^+ - w_1^-\bigr)\in E^1_{3,0},
 $$
 and
 \begin{equation}\label{eq_part_tot_theta}
 \partial_{\tot}\theta\Lequiv\bigl(1+\omega+\omega^2\bigr)\bullet(\fA_1-\fA_3).
 \end{equation}
\end{lem}

\begin{proof}
 Let us consecutively construct the elements~$\theta_{3,0}$, $\theta_{2,1}$, $\theta_{1,2}$, and $\theta_{0,3}$. We set
 $$
 \theta_{3,0}=\bigl(1+\omega+\omega^2\bigr)\bigl(T_{23}-T_{3\bar{1}}^{-1}T_{\bar{1}2}^{-1}\nu\bigr)\bigl(W_1^+ - W_1^-\bigr).
 $$
 Then $\theta_{3,0}$ represents~$\mu$. Let us compute~$\partial'\theta_{3,0}$ up to elements of~$C_2(L)\otimes_{\CS}\CR_0=C_2(L)$ using the formulas for~$\partial W_1^{\pm}$ from Table~\ref{table_3cells}. The cells~$Z_4^{\pm}$ belong to~$L$, so we neglect them. By~\eqref{eq_lantern_faces}, we have $T_{\varepsilon_1}=T_{\bar{1}2}T_{3\bar{1}}T_{23}$. Hence
 $$
 \bigl(T_{23}-T_{3\bar{1}}^{-1}T_{\bar{1}2}^{-1}\nu\bigr)=T_{23}\bigl(1-T_{\varepsilon_1}^{-1}\nu\bigr).
 $$
 We see that this element annihilates~$Z_3^{\pm}$. Therefore
 \begin{align*}
 \partial'\theta_{3,0}&\Lequiv\bigl(1+\omega+\omega^2\bigr)\bigl(T_{23}-T_{3\bar{1}}^{-1}T_{\bar{1}2}^{-1}\nu\bigr)\bigl(T_{2\bar{3}}(r_2-\omega)-T_{\bar{3}1}^{-1}T_{\bar{1}2}^{-1}\bigl(\omega^2-r_3\bigr)\nu\bigr)Z_1\\
 {}&=\bigl(1+\omega+\omega^2\bigr)\bigl(
 T_{23}T_{2\bar{3}}(r_2-\omega)
 +T_{3\bar{1}}^{-1}T_{\bar{1}2}^{-1}T_{3\bar{1}}T_{1\bar{2}}\bigl(\omega^2-r_3\bigr)\\
 &\hspace{3.5cm}{}+T_{3\bar{1}}^{-1}T_{\bar{1}2}^{-1}T_{\bar{2}3}^{-1}(\omega-r_2)\nu
 +T_{23}T_{\bar{3}1}^{-1}T_{\bar{1}2}^{-1}\bigl(r_3-\omega^2\bigr)\nu
 \bigr)Z_1.
 \end{align*}
 Using that $r_2=\omega^2r_1$ and $r_3=\omega r_1$ and moving the powers of~$\omega$ to the left, we obtain:
 \begin{multline}\label{eq_part*_aux1}
 \partial'\theta_{3,0}\Lequiv\bigl(1+\omega+\omega^2\bigr)\Bigl(
 \bigl(-T_{12}T_{1\bar{2}}+T_{1\bar{2}}^{-1}T_{\bar{2}3}^{-1}T_{1\bar{2}}T_{2\bar{3}}\bigr)Z_1
 +\bigl(T_{31}T_{3\bar{1}}-T_{2\bar{3}}^{-1}T_{\bar{3}1}^{-1}T_{2\bar{3}}T_{3\bar{1}}\bigr)r_1Z_1\\
 {}+\bigl(T_{2\bar{3}}^{-1}T_{\bar{3}1}^{-1}T_{\bar{1}2}^{-1}-T_{31}T_{\bar{1}2}^{-1}T_{\bar{2}3}^{-1}\bigr)\nu Z_1
 +\bigl(-T_{1\bar{2}}^{-1}T_{\bar{2}3}^{-1}T_{\bar{3}1}^{-1}+T_{12}T_{\bar{2}3}^{-1}T_{\bar{3}1}^{-1}\bigr)r_1\nu Z_1
 \Bigr).
 \end{multline}
 Since $\partial Z_1\Lequiv Y_1-Y_3$ and $Y_1$ is fixed by~$T_{\delta_2}T_{\varepsilon}\nu$, it is more convenient for furhter computations to write~$\partial'\theta_{3,0}$ through~$Z_1$, $r_1Z_1$, $T_{\delta_2}T_{\varepsilon}\nu Z_1$ and $r_1T_{\delta_2}T_{\varepsilon}\nu Z_1=T_{\delta_3}T_{\varepsilon}r_1\nu Z_1$ rather than through ~$Z_1$, $r_1Z_1$, $\nu Z_1$ and~$r_1\nu Z_1$.
 Combining relations~\eqref{eq_lantern_eps} and~\eqref{eq_lantern_faces}, we obtain:
 $$
 T_{\delta_2}T_{\varepsilon}=T_{31}T_{\bar{3}\bar{1}}T_{\varepsilon_2}^{-1}=T_{31}T_{\bar{3}\bar{1}}\bigl(T_{\bar{3}\bar{1}}^{-1}T_{2\bar{3}}^{-1}T_{\bar{1}2}^{-1}\bigr)=T_{31}T_{2\bar{3}}^{-1}T_{\bar{1}2}^{-1}.
 $$
 Conjugating this equality by~$r_1$, we also get
 $$
 T_{\delta_3}T_{\varepsilon}=T_{\bar{1}\bar{2}}T_{2\bar{3}}^{-1}T_{\bar{3}1}^{-1}.
 $$
 Therefore \eqref{eq_part*_aux1} reads as
 \begin{multline}\label{eq_part*_aux2}
 \partial'\theta_{3,0}\Lequivs\bigl(1+\omega+\omega^2\bigr)\Bigl(
 \bigl(-T_{12}T_{1\bar{2}}+T_{1\bar{2}}^{-1}T_{\bar{2}3}^{-1}T_{2\bar{3}}T_{1\bar{2}}\bigr)Z_1
 +\bigl(T_{31}T_{3\bar{1}}-T_{2\bar{3}}^{-1}T_{\bar{3}1}^{-1}T_{3\bar{1}}T_{2\bar{3}}\bigr)r_1Z_1\\
 {}+\bigl(T_{2\bar{3}}^{-1}T_{\bar{3}1}^{-1}T_{2\bar{3}}T_{31}^{-1}-T_{31}T_{\bar{2}3}^{-1}T_{2\bar{3}}T_{31}^{-1}\bigr)T_{\delta_2}T_{\varepsilon}\nu Z_1\\
 {}+\bigl(-T_{1\bar{2}}^{-1}T_{\bar{2}3}^{-1}T_{2\bar{3}}T_{\bar{1}\bar{2}}^{-1}+T_{12}T_{\bar{2}3}^{-1}T_{2\bar{3}}T_{\bar{1}\bar{2}}^{-1}\bigr)r_1T_{\delta_2}T_{\varepsilon}\nu Z_1
 \Bigr).
 \end{multline}
We introduce the following notation for mapping classes in~$\CS$:
\begin{align*}
 g_1&=T_{\delta_1}^{-1},&g_2&=T_{\delta_2}^{-1},&g_3&=T_{\delta_3}^{-1},\\
 h_1&=T_{2\bar{3}}^{-1}T_{23}^{-1},&h_2&=T_{3\bar{1}}^{-1}T_{31}^{-1},&h_3&=T_{1\bar{2}}^{-1}T_{12}^{-1},\\
 f_1&=T_{1\bar{2}}^{-1}T_{2\bar{3}}^{-1}T_{\bar{2}3}T_{1\bar{2}},&
 f_2&=T_{2\bar{3}}^{-1}T_{3\bar{1}}^{-1}T_{\bar{3}1}T_{2\bar{3}},&
 f_3&=T_{3\bar{1}}^{-1}T_{1\bar{2}}^{-1}T_{\bar{1}2}T_{3\bar{1}}.
\end{align*}
By the lantern relation (see \eqref{eq_lantern_faces}), the mapping class~$T_{1\bar{2}}T_{31}$ commutes with~$T_{\bar{2}3}$, and therefore with~$T_{\bar{2}3}T_{2\bar{3}}^{-1}$. Hence
$$
T_{31}T_{\bar{2}3}^{-1}T_{2\bar{3}}T_{31}^{-1}=T_{1\bar{2}}^{-1}T_{\bar{2}3}^{-1}T_{2\bar{3}}T_{1\bar{2}}=f_1^{-1}.
$$
Also using~\eqref{eq_lantern_faces}, we obtain
$$
T_{12}T_{\bar{2}3}^{-1}T_{2\bar{3}}T_{\bar{1}\bar{2}}^{-1}=T_{12}T_{\bar{2}3}^{-1}T_{\bar{1}\bar{2}}^{-1}T_{2\bar{3}}=T_{12}T_{\varepsilon_3}^{-1}T_{3\bar{1}}T_{2\bar{3}}=
T_{2\bar{3}}^{-1}T_{\bar{3}1}^{-1}T_{3\bar{1}}T_{2\bar{3}}=f_2^{-1}.
$$
By~\eqref{eq_Tdeltai} we have
$$
T_{\delta_1}=T_{1\bar{2}}^{-1}T_{\bar{3}1}^{-1}T_{12}^{-1}T_{31}^{-1}=T_{3\bar{1}}^{-1}T_{\bar{1}2}^{-1}T_{\bar{3}\bar{1}}^{-1}T_{\bar{1}\bar{2}}^{-1}.
$$
Using these formulas, relations~\eqref{eq_lantern_faces}, and the fact that $T_{2\bar{3}}T_{12}$ commutes with~$T_{\bar{3}1}$ (which also follows from~\eqref{eq_lantern_faces}), we get
\begin{gather*}
T_{2\bar{3}}^{-1}T_{\bar{3}1}^{-1}T_{2\bar{3}}T_{31}^{-1}=T_{12}T_{\bar{3}1}^{-1}T_{12}^{-1}T_{31}^{-1}=T_{12}T_{1\bar{2}}T_{\delta_1}=h_3^{-1}g_1^{-1},\\
T_{1\bar{2}}^{-1}T_{\bar{2}3}^{-1}T_{2\bar{3}}T_{\bar{1}\bar{2}}^{-1}=T_{31}T_{\varepsilon_2}^{-1}T_{2\bar{3}}T_{\bar{1}\bar{2}}^{-1}=T_{31}T_{\bar{1}2}^{-1}T_{\bar{3}\bar{1}}^{-1}T_{\bar{1}\bar{2}}^{-1}=T_{31}T_{3\bar{1}}T_{\delta_1}=h_2^{-1}g_1^{-1}.
\end{gather*}

Thus \eqref{eq_part*_aux2} reads as
\begin{multline*}
 \partial'\theta_{3,0}\Lequivs\bigl(1+\omega+\omega^2\bigr)\Bigl(
 \bigl(-h_3^{-1}+f_1^{-1}\bigr)Z_1
 +\bigl(h_2^{-1}-f_2^{-1}\bigr)r_1Z_1\\
 {}+\bigl(h_3^{-1}g_1^{-1}-f_1^{-1}\bigr)T_{\delta_2}T_{\varepsilon}\nu Z_1+\bigl(-h_2^{-1}g_1^{-1}+f_2^{-1}\bigr)r_1T_{\delta_2}T_{\varepsilon}\nu Z_1
 \Bigr).
 \end{multline*}

By Lemma~\ref{lem_d'd''}, we have that
$$
\partial'\theta_{3,0}\Lequiv-\partial''\tilde{\theta}_{2,1},
$$
where
\begin{multline*}
 \tilde{\theta}_{2,1}=\bigl(1+\omega+\omega^2\bigr)\bullet\Bigl(
 Z_1\otimes\bigl([h_3]-[f_1] \bigr)+(r_1Z_1)\otimes\bigl([f_2]
 -[h_2]\bigr)\\
 {}+(T_{\delta_2}T_{\varepsilon}\nu Z_1)\otimes\bigl(
 [f_1]
 -[g_1h_3]\bigr)+(r_1T_{\delta_2}T_{\varepsilon}\nu Z_1)\otimes\bigl(
 [g_1h_2]
 -[f_2]\bigr)
 \Bigr).
 \end{multline*}
We write~$\tilde{\theta}_{2,1}$ instead of~$\theta_{2,1}$, since this is not yet the element~$\theta_{2,1}$ that we will need; we will add further terms to it later.

Let us now compute~$\partial'\tilde{\theta}_{2,1}$ up to elements of~$C_1(L)\otimes_{\CS}\CR_1$. We have $\partial Z_1\Lequivs Y_1-Y_3$. Write
$$
\partial'\tilde{\theta}_{2,1}\Lequiv S_1+S_2,
$$
where $S_1$ stands for the sum of all terms involving~$Y_1$, and~$S_2$ for the sum of all terms involving~$Y_3$. Let us consider the sum~$S_2$ first.
Since $\omega$ and~$r_1$ fix~$Y_3$, all terms in~$S_2$ involving the cells~$Y_3$ and~$r_1Y_3$ cancel. Therefore
\begin{align*}
S_2&=\bigl(1+\omega+\omega^2\bigr)\bullet\bigl((g_2^{-1}T_{\varepsilon}\nu Y_3)\otimes(
 [g_1h_3]-[f_1])-(g_3^{-1}T_{\varepsilon}\nu Y_3)\otimes(
 [g_1h_2]-[f_2])
 \bigr)\\
 {}&=\bigl(1+\omega+\omega^2\bigr)\bullet\bigl((T_{\varepsilon}\nu Y_3)\otimes(
 g_2[g_1h_3]-g_2[f_1]-
 g_3[g_1h_2]+g_3[f_2])
 \bigr)\\
 {}&=\bigl(1+\omega+\omega^2\bigr)\bullet\bigl((T_{\varepsilon}\nu Y_3)\otimes(
 g_2[g_1h_3]-
 g_1[g_2h_3])
 \bigr),
\end{align*}
where the last equality holds because $\omega$ fixes~$T_{\varepsilon}\nu Y_3$ and the conjugation by~$\omega$ permutes cyclically each of the triples~$(f_1,f_2,f_3)$, $(g_1,g_2,g_3)$, and~$(h_1,h_2,h_3)$. Since $g_1$ and~$g_2$ commute, we have
$$
d([g_2|g_1h_3]-[g_1|g_2h_3])=g_2[g_1h_3]+[g_2]-g_1[g_2h_3]-[g_1].
$$
Therefore
$$
S_2=-\partial''R_2,
$$
where
$$
R_2=\bigl(1+\omega+\omega^2\bigr)\bullet\bigl((T_{\varepsilon}\nu Y_3)\otimes(
 [g_2|g_1h_3]-
 [g_1|g_2h_3])
 \bigr).
$$

Let us now consider the sum~$S_1$. Since $T_{\delta_2}T_{\varepsilon}\nu Y_1=Y_1$, we have:
\begin{equation*}
 S_1=\bigl(1+\omega+\omega^2\bigr)\bullet\bigl(
 Y_1\otimes([h_3]-[g_1h_3])+(r_1Y_1)\otimes([g_1h_2]-[h_2])
 \bigr).
 \end{equation*}
Since $g_1$ fixes~$Y_1$, we see that
$$
\partial''\bigl(Y_1\otimes[g_1|h_3]\bigr)=Y_1\otimes(-[h_3]+[g_1h_3]-[g_1]),
$$
and similarly for $(r_1Y_1)\otimes[g_1|h_2]$. We set
$$
R_1=\bigl(1+\omega+\omega^2\bigr)\bullet\bigl(
 Y_1\otimes[g_1|h_3]-(r_1Y_1)\otimes[g_1|h_2]
 \bigr).
$$
Then
\begin{align*}
S_1+\partial''R_1&=\bigl(1+\omega+\omega^2\bigr)\bullet\bigl(
 -Y_1\otimes[g_1]+(r_1Y_1)\otimes[g_1]
 \bigr)\\
 {}&=-\partial'\Bigl(\bigl(1+\omega+\omega^2\bigr)\bullet\bigl(Z_7\otimes [g_1]\bigr)\Bigr)+\bigl(1+\omega+\omega^2\bigr)\bullet\bigl(Y_8\otimes [g_1]\bigr).
\end{align*}
We now set
$$
\theta_{2,1}=\tilde{\theta}_{2,1}+\bigl(1+\omega+\omega^2\bigr)\bullet\bigl(Z_7\otimes [g_1]\bigr).
$$
Since $g_1$ fixes~$Z_7$, we have
\begin{equation}\label{eq_dtot_K1}
 \partial''\theta_{2,1}=\partial''\tilde{\theta}_{2,1}\Lequiv-\partial'\theta_{3,0}.
\end{equation}
Moreover,
\begin{equation}\label{eq_dtot_K2}
\partial'\theta_{2,1}\Lequiv\bigl(1+\omega+\omega^2\bigr)\bullet\bigl(Y_8\otimes [g_1]\bigr)-\partial''\theta_{1,2},
\end{equation}
where
$$
\theta_{1,2}=R_1+R_2.
$$

Let us now compute $\partial'\theta_{1,2}$ up to elements of~$C_0(L)\otimes_{\CS}\CR_2$. We have:
$$
\partial Y_1\Lequiv-X_1,\qquad\partial r_1Y_1\Lequiv-X_1,\qquad \partial T_{\varepsilon}\nu Y_3\Lequiv-X_1.
$$
Therefore
\begin{align*}
\partial'\theta_{1,2}&\Lequiv\bigl(1+\omega+\omega^2\bigr)\bullet\bigl(X_1\otimes([g_1|h_2]-[g_1|h_3]+[g_1|g_2h_3]-[g_2|g_1h_3])\bigr)\\
{}&=\bigl(1+\omega+\omega^2\bigr)\bullet\bigl(X_1\otimes([g_2|h_3]-[g_1|h_3]+[g_1|g_2h_3]-[g_2|g_1h_3])\bigr).
\end{align*}
We set
$$
\theta_{0,3}=\bigl(1+\omega+\omega^2\bigr)\bullet\bigl(
X_1\otimes\bigl(
\bigl[g_1\big|g_2h_3\big|h_3^{-1}\bigr]
-\bigl[g_2\big|g_1h_3\big|h_3^{-1}\bigr]
-\bigl[g_1\big|h_3\big|h_3^{-1}\bigr]
+\bigl[g_2\big|h_3\big|h_3^{-1}\bigr]
\bigr)\bigr).
$$
Since $g_1$ and~$g_2$ fix~$X_1$ and commute, we obtain that
\begin{equation}\label{eq_dtot_K3}
 \partial'\theta_{1,2}+\partial''\theta_{0,3}\Lequiv\bigl(1+\omega+\omega^2\bigr)\bullet\bigl(
 X_1\otimes([g_1|g_2]-[g_2|g_1])
 \bigr).
\end{equation}

Finally, we set
$$
\theta=\theta_{3,0}+\theta_{2,1}+\theta_{1,2}+\theta_{0,3}.
$$
Combining~\eqref{eq_dtot_K1}, \eqref{eq_dtot_K2}, and~\eqref{eq_dtot_K3}, we obtain
\begin{align*}
\partial_{\tot}\theta&\Lequiv\bigl(1+\omega+\omega^2\bigr)\bullet\bigl(
 X_1\otimes([g_1|g_2]-[g_2|g_1])+Y_8\otimes[g_1]
 \bigr)\\
 &\Lequiv\bigl(1+\omega+\omega^2\bigr)\bullet(\fA_1-\fA_3),
\end{align*}
as required.
\end{proof}

\begin{proof}[Proof of Proposition~\ref{propos_rel_from_lambda}]
We have
$$
\lambda=F\mu+c_2w_2+c_4(t_1-1)^2w_4+c_6w_6,\qquad F\in\Lambda,\ c_2,c_4,c_6\in\Z[\Gamma].
$$
Let $\theta_{3,0}$, $\theta_{2,1}$, $\theta_{1,2}$, and~$\theta_{0,3}$ be the elements from Lemma~\ref{lem_theta_for_mu}.

Choose arbitrary elements $\widehat{F},\hat{c}_2,\hat{c}_4,\hat{c}_6\in\Z[\G]$ satisfying $\kappa\bigl(\widehat{F}\bigr)=F$ and $\kappa(\hat{c}_i)=c_i$, and set:
\begin{gather*}
\Theta_{3,0}=\widehat{F}\theta_{3,0}+\hat{c}_2W_2+\hat{c}_4\bigl(T_{2\bar{3}}^{-1}-1\bigr)\bigl(T_{\bar{2}3}^{-1}-1\bigr)W_4+\hat{c}_6W_6,\\
\Theta_{p,q}=\widehat{F}\theta_{p,q},\qquad (p,q)=(2,1),(1,2),(0,3),\\
\Theta=\Theta_{3,0}+\Theta_{2,1}+\Theta_{1,2}+\Theta_{0,3}.
\end{gather*}
Then $\Theta_{3,0}$ represents~$\lambda$. Let us prove that
\begin{equation}\label{eq_part_tot_Theta}
\partial_{\tot}\Theta\Lequiv \widehat{F}\bigl(1+\omega+\omega^2\bigr)\bullet(\fA_1-\fA_3).
\end{equation}
This would be an immediate consequence of~\eqref{eq_part_tot_theta} if we replaced~$\Theta_{3,0}$ with~$\widehat{F}\theta_{3,0}$. Hence we need only to show that
\begin{equation}\label{eq_partial'inK}
\partial\bigl(\hat{c}_2W_2+\hat{c}_4\bigl(T_{2\bar{3}}^{-1}-1\bigr)\bigl(T_{\bar{2}3}^{-1}-1\bigr)W_4+\hat{c}_6W_6\bigr)\in C_2(L).
\end{equation}
We have $\partial W_2,\partial W_6\in C_2(L)$, since the cells~$W_2$ and~$W_6$ lie in~$L$. Further, $$\partial W_4\Lequiv-(\nu+1)Z_6.$$ The factors $T_{2\bar{3}}^{-1}-1$ and~$T_{\bar{2}3}^{-1}-1$ annihilate the cells~$Z_6$ and~$\nu Z_6$, respectively, and commute. Consequently, we obtain~\eqref{eq_partial'inK}, and therefore~\eqref{eq_part_tot_Theta}.

If the action of~$\CS$ on~$L$ were free, we could add to~$\Theta_{p,q}$ certain correction terms $\varphi_{p,q}\in C_p(L)\otimes_{\CS}\CR_q$ so that \eqref{eq_part_tot_Theta} would become an exact equality, as in the proof of Lemma~\ref{lem_dr}. Although the action of~$\CS$ on~$L$ is not free, it is free away from the vertices~$gX_2$ with $g\in\G$. Our goal is to show that, by adding certain correction terms~$\varphi_{p,q}$ as in the proof of Lemma~\ref{lem_dr}, we can turn \eqref{eq_part_tot_Theta} into an equality up to a cycle representing a class from $\Z[\Gamma]\CA_2$.

We add no correction term in bidegree~$(3,0)$, i.e., $\varphi_{3,0}=0$. The element $\partial'\Theta_{3,0}+\partial''\Theta_{2,1}$ lies in~$C_2(L)\otimes_{\CS}\CR_0$ by~\eqref{eq_part_tot_Theta}, and lies in the image of~$\partial''$, since $\Theta_{3,0}$ represents~$\lambda$ and~$d^1\lambda=0$. It follows that there exists an element $\varphi_{2,1}\in C_2(L)\otimes_{\CS}\CR_1$ satisfying $\partial''\varphi_{2,1}=-\partial'\Theta_{3,0}-\partial''\Theta_{2,1}$. Set $\widetilde{\Theta}_{2,1}=\Theta_{2,1}+\varphi_{2,1}$. Then
\begin{equation}\label{eq_Theta1}
\partial'\Theta_{3,0}+\partial''\widetilde{\Theta}_{2,1}=0.
\end{equation}

Let $\fA_3^{(1,1)}$ and $\fA_3^{(0,2)}$ be the bigraded components of~$\fA_3$ of bidegrees~$(1,1)$ and~$(0,2)$, respectively. Equality~\eqref{eq_part_tot_Theta} implies that the element
$$
\psi_{1,1}=\partial'\widetilde{\Theta}_{2,1}+\partial''\Theta_{1,2}+\widehat{F}\bigl(1+\omega+\omega^2\bigr)\fA_3^{(1,1)}
$$
lies in $C_1(L)\otimes_{\CS}\CR_1$. Moreover, we have
$$
\partial''\psi_{1,1}=\partial''\partial'\widetilde{\Theta}_{2,1}=-\partial'\partial''\widetilde{\Theta}_{2,1}=\partial'\partial'\Theta_{3,0}=0.
$$
Since the action of~$\CS$ on~$L$ is free away from the vertices, the $\Z[\CS]$-module~$C_1(L)$ is free, and hence the chain complex $\bigl(C_1(L)\otimes_{\CS}\CR_*,\partial''\bigr)$ is acyclic in positive dimensions, as in Lemma~\ref{lem_acyclic}. Therefore there exists an element $\varphi_{1,2}\in C_1(L)\otimes_{\CS}\CR_2$ satisfying $\partial''\varphi_{1,2}=-\psi_{1,1}$. Set $\widetilde{\Theta}_{1,2}=\Theta_{1,2}+\varphi_{1,2}$. Then
\begin{equation}\label{eq_Theta2}
\partial'\widetilde{\Theta}_{2,1}+\partial''\widetilde{\Theta}_{1,2}=-\widehat{F}\bigl(1+\omega+\omega^2\bigr)\bullet\fA_3^{(1,1)}.
\end{equation}

Equality~\eqref{eq_part_tot_Theta} implies that the element
\begin{equation}\label{eq_Theta3}
\psi_{0,2}=\partial'\widetilde{\Theta}_{1,2}+\partial''\Theta_{0,3}-\widehat{F}\bigl(1+\omega+\omega^2\bigr)\bullet\left(\fA_1-\fA_3^{(0,2)}\right)
\end{equation}
lies in $C_0(L)\otimes_{\CS}\CR_2$. Moreover, we have
\begin{align*}
\partial''\psi_{0,2}&=\partial''\partial'\widetilde{\Theta}_{1,2}+\partial''\widehat{F}\bigl(1+\omega+\omega^2\bigr)\bullet\fA_3^{(0,2)}=-\partial'\partial''\widetilde{\Theta}_{1,2}+\partial''\widehat{F}\bigl(1+\omega+\omega^2\bigr)\bullet\fA_3^{(0,2)}\\
{}&=\partial'\partial'\widetilde{\Theta}_{2,1}+\partial'\widehat{F}\bigl(1+\omega+\omega^2\bigr)\bullet\fA_3^{(1,1)}+\partial''\widehat{F}\bigl(1+\omega+\omega^2\bigr)\bullet\fA_3^{(0,2)}=0.
\end{align*}
Since the cycle~$\psi_{0,2}$ lies in $C_0(L)\otimes_{\CS}\CR_2$, and $L$ has two $\G$-orbits of vertices with representatives~$X_2$ and~$X_3$, we see that $\psi_{0,2}$ represents in~$E^1_{0,2}$ a class  belonging to the submodule
$$
\Z[\Gamma]\iota_{X_2}\bigl(H_2(\CS_{X_2})\bigr)\oplus \Z[\Gamma]\iota_{X_3}\bigl(H_2(\CS_{X_3})\bigr)\subseteq E^1_{0,2}.
$$
But $H_2(\CS_{X_3})=0$ and $\iota_{X_2}\bigl(H_2(\CS_{X_2})\bigr)$ is generated by~$t_2^{-1}A_2$. Thus $\psi_{0,2}$ represents in~$E^1_{0,2}$ a class belonging to~$\Z[\Gamma]A_2$. Consequently, $\psi_{0,2}$ regardered as a cycle in the chain complex $\Tot_*(C_*(\CB_{0,6})\otimes_{\CS}\CR_*)$ represents a class in~$\Z[\Gamma]\CA_2\subseteq H_2(\CS)$.

Finally, set $\widetilde{\Theta}=\Theta_{3,0}+\widetilde{\Theta}_{2,1}+\widetilde{\Theta}_{1,2}+\Theta_{0,3}$. Combining~\eqref{eq_Theta1}, \eqref{eq_Theta2}, and~\eqref{eq_Theta3}, we get:
$$
\partial_{\tot}\widetilde{\Theta}=\widehat{F}\bigl(1+\omega+\omega^2\bigr)\bullet(\fA_1-\fA_3)+\psi_{0,2}.
$$
The proposition follows, since the cycle $\widehat{F}\bigl(1+\omega+\omega^2\bigr)\bullet(\fA_1-\fA_3)$ represents in~$H_2(\CS)$ the class $F\bigl(1+\omega+\omega^2\bigr)(\CA_1-\CA_3)$.
\end{proof}

\part{Reduction to the subgroup generated by standard abelian cycles}

\section{Main theorem}\label{section_modA_first}
Recall that $\Sigma=\Sigma_3$ is a closed surface of genus~$3$, $H=H_1(\Sigma)$, $\I=\I_3$ is the Torelli group of~$\Sigma$, and $\bA\subseteq H_2(\I)$ is the subgroup generated by all standard abelian cycles.

The aim of this section is to prove the following theorem.
\begin{theorem}[$=$ Theorem~\ref{thm_toAC}]\label{theorem_modA}
For each $x\in H$, we have
 \begin{equation}\label{eq_tx_tors}
 (t_x-1)^6H_2(\I)\subseteq \bA.
 \end{equation}
\end{theorem}

Since $t_{kx}-1=t_x^{k^2}-1$ is divisible by~$t_x-1$, it suffices to prove~\eqref{eq_tx_tors} for primitive homology classes~$x$.
Throughout this part of the paper, we fix a primitive class~$x\in H$, and prove~\eqref{eq_tx_tors} for this particular~$x$. To do this, we consider the action of the group~$\I$ on the Bestvina--Bux--Margalit complex of cycles $\CB=\CB(x)$. We will recall the construction of this complex below; see Section~\ref{section_BBM}. At present, we wish to give the general scheme of the proof of Theorem~\ref{theorem_modA}, which will be implemented in the rest of this part.

Let $\Mod(\Sigma,x)$ and $\Sp(H,x)$ be the stabilizers of the homology class~$x$ in the groups~$\Mod(\Sigma)$ and~$\Sp(H)$, respectively. Then we have a short exact sequence
$$
1\longrightarrow \I\longrightarrow \Mod(\Sigma,x) \longrightarrow \Sp(H,x)\longrightarrow 1.
$$
We need the following properties of the complex~$\CB$.
 \begin{itemize}
  \item $\CB$ is contractible.
  \item The Torelli group~$\I$ acts on~$\CB$ cellularly and without rotations, and moreover, this action extends to a cellular action of~$\Mod(\Sigma,x)$.
  \item The cells of~$\CB$ are indexed by certain oriented multicurves on~$\Sigma$; we denote by~$\CM_p$ the set of oriented multicurves indexing the $p$-cells of~$\CB$.
 \end{itemize}

Consider the spectral sequence from Section~\ref{section_SpecSeq} for the action of~$\I$ on~$\CB$:
\begin{equation}\label{eq_SS_I3}
E^1_{p,q}=\bigoplus_{M\in\mathfrak{M}_p}H_q(\I_M) \Longrightarrow H_{p+q}(\I),
\end{equation}
where $\mathfrak{M}_p$ is some set of representatives of $\I$-orbits in~$\CM_p$. Let
\begin{equation}\label{eq_filtration_H2I3}
 0\subseteq \CF_{0,2}\subseteq \CF_{1,1}\subseteq \CF_{2,0}=H_2(\I)
\end{equation}
be the corresponding filtration with $\CF_{p,q}/\CF_{p-1,q+1}=E^{\infty}_{p,q}$. Since the action of~$\I$ on~$\CB$ extends to~$\Mod(\Sigma,x)$, by Fact~\ref{fact_Gamma} we have that \eqref{eq_SS_I3} is a spectral sequence of $\Z\bigl[\Sp(H,x)\bigr]$-modules, and moreover, the action of~$\Sp(H,x)$ on the $E^{\infty}$-page is associated with its action on the filtration~\eqref{eq_filtration_H2I3}.

The proof of the inclusion~\eqref{eq_tx_tors} will consist of three steps.
\smallskip

\textsl{Step 1.} We prove that $(t_x-1)E^2_{2,0}=0$. This immediately implies that $(t_x-1)E^{\infty}_{2,0}=0$, and therefore $(t_x-1)\CF_{2,0}\subseteq\CF_{1,1}$.
\smallskip

\textsl{Step 2.} This step is more difficult than the previous one, since the author does not know whether some power of $t_x-1$ annihilates the group~$E^{\infty}_{1,1}$. We proceed as follows: we introduce a special subgroup $Y\subseteq E^{\infty}_{1,1}$ with the following properties. First, we show that $(t_x-1)^2E^{\infty}_{1,1}\subseteq Y$. Second, we prove that $Y$ is contained in the image of the group $\bA\cap\CF_{1,1}$ under the homomorphism $\CF_{1,1}\twoheadrightarrow E^{\infty}_{1,1}$. These two properties together imply that $(t_x-1)^2\CF_{1,1}\subseteq \CF_{0,2}+\bA$.
\smallskip

\textsl{Step 3.} To prove the required inclusion~\eqref{eq_tx_tors} it now suffices to show that $(t_x-1)^3\CF_{0,2}\subseteq\bA$. The group $\CF_{0,2}=E^{\infty}_{0,2}$ is the quotient of the group
$$
E^1_{0,2}=\bigoplus_{M\in\mathfrak{M}_0}H_2(\I_M)
$$
by the images of the differentials~$d^1$, $d^2$, and~$d^3$. The set~$\CM_0$, and hence the set~$\mathfrak{M}_0$, contains one-, two-, and three-component oriented multicurves. Moreover, each of these multicurves does not separate the surface~$\Sigma$. For each three-component multicurve $M\in\mathfrak{M}_0$, assertion~(1) of Theorem~\ref{theorem_3} implies that $H_2(\I_M)$ is generated by standard abelian cycles. Hence the image of the summand $H_2(\I_M)\subseteq E^1_{0,2}$ in $\CF_{0,2}=E^{\infty}_{0,2}$ is contained in~$\bA$. If $\alpha\in \mathfrak{M}_0$ is a one-component oriented multicurve, then $[\alpha]=x$. Since the Dehn twist~$T_{\alpha}$ commutes with all elements of the stabilizer~$\I_{\alpha}$, it follows that the element $t_x-1$ annihilates the summand $H_2(\I_{\alpha})\subseteq E^1_{0,2}$. The most difficult case is that of a two-component multicurve $M\in\mathfrak{M}_0$. We need to show that $(t_x-1)^3\I_M\subseteq \bA_M$, where $\bA_M$ is the subgroup of~$H_1(\I_M)$ generated by all standard abelian cycles. To this end, we consider the action of the stabilizer~$\I_M$ on an auxiliary complex~$\CB_M$, which serves as an analogue of the complex of cycles~$\CB$ in this setting, and repeat Steps~1 and~2 described above for the spectral sequence associated with the action of~$\I_M$ on~$\CB_M$.
\smallskip

This part of the paper is organized as follows. In Section~\ref{section_BBM}, we recall the construction of the complex of cycles~$\CB$ and its basic properties. In Section~\ref{section_important}, we introduce the key notion of an \textit{important} oriented multicurve and define a chain complex~$\CC_*$ by taking the quotient of the cellular chain complex~$E^1_{*,0}=C_*(\CB/\I)$ by the subcomplex spanned by all cells corresponding to unimportant oriented multicurves. The main result of Section~\ref{section_important} is Proposition~\ref{propos_H2C0}, which asserts that $H_2(\CC_*)=0$. This result will be used twice. First, in Section~\ref{section_toF11} we use it to carry out Step~1 of the above plan. Second, we will later use it to carry out Step~2, too.

Step~2 is the most difficult part of our proof. Here we will need a number of results on the structure of the spectral sequence~$E^*_{*,*}$ obtained by the author in~\cite{Gai21}. Sections~\ref{section_active} and~\ref{section_stabilizers} are devoted to a survey of these results and some of their consequences. Section~\ref{section_ephi} contains the key construction for Step 2: the construction of a subcomplex~$\CE_*$ of the first row of the $E^1$-page of the spectral sequence~$E^*_{*,*}$ and of a surjective chain map $\varphi\colon \CE_*\to\CC_{*+1}$. Section~\ref{section_KNDZ} is devoted to the study of the kernel of the map $\varphi\colon\CE_1\to \CC_2$. Finally, in Section~\ref{section_to_F02} we combine the results obtained and complete Step~2; here we again use the fact that $H_2(\CC_*)=0$.

In Section~\ref{section_stab2}, we carry out Step~3 of our plan; essentially, in this section we run through the entire preceding proof once more in a somewhat simpler case, namely for the stabilizer~$\I_N$ of a two-component multicurve~$N$ on~$\Sigma$. Finally, in Section~\ref{section_modA_end}, we bring together all the results obtained and complete the proof of Theorem~\ref{theorem_modA}.

\section{Complex of cycles}\label{section_BBM}

Following Bestvina, Bux, and Margalit~\cite{BBM10}, we consider the following \textit{complex of minimizing cycles} (or simple a \textit{complex of cycles})~$\CB=\CB(x)$. This construction depends on the choice of a primitive homology class $x\in H$. We fix such a class throughout the rest of this paper.

The construction of this complex is completely analogous to that of the complex of relative cycles from Section~\ref{section_RCC}. In fact, the situation is, of course, exactly the opposite: Spiridonov defined the complex of relative cycles by analogy with the Bestvina–Bux–Margalit complex of minimizing cycles; it is simply that in the present paper we need the complex of relative cycles earlier.

A simple closed curve on~$\Sigma$ is called \textit{essential} if it is homotopically nontrivial, i.e., does not bound a disc. By an \textit{oriented multicurve} we always mean a finite union of pairwise non-isotopic essential oriented simple closed curves on~$\Sigma$. We will usually not distinguish between an oriented multicurve and its isotopy class. We adopt the convention that, for oriented multicurves, $N\subseteq M$ always means that each component of~$N$ is oriented in~$N$ in the same way as in~$M$.

Let $\fC$ be the set of oriented essetial simple closed curves on~$\Sigma$. Note that the same curve with two opposite orientations gives two distinct elements of~$\fC$.  The \textit{complex of cycles}~$\CB$ is the subset of~$\R_{\ge0}^{\fC}$ consisting of all points (finite formal sums) $\bgamma=\sum_{i=1}^kn_i\gamma_i$, where $\gamma_i\in\fC$ and $n_i\in\R_{>0}$, that satisfy the following conditions:
\begin{enumerate}
 \item $\sum_{i=1}^k n_i [\gamma_i]=x$,
 \item the isotopy classes $\gamma_1,\ldots,\gamma_k$ contain pairwise disjoint representatives,  and
 \item the homology classes $[\gamma_1],\ldots,[\gamma_k]$ satisfy no one-sided relations, i.e., relations of the form $\sum_{i\in I}[\gamma_i]=0$, where $I$ is a nonempty subset of~$\{1,\ldots,k\}$.
\end{enumerate}
Note that two-sided relations, i.e., relations of the form $\sum_{i\in I}[\gamma_i]=\sum_{j\in J}[\gamma_j]$, where $I$ and~$J$ are both nonempty, are allowed. Condition~(3) in particular implies that all curves $\gamma_i$ are non-separating. The points $\bgamma\in\CB$ are called \textit{cycles} for~$x$.

The set~$\CB$ has the following structure of a regular CW complex. By definition, the \textit{support} of a cycle $\bgamma\in \CB$ is the oriented multicurve obtained as the union of all oriented simple closed curves~$\gamma_i$ that enter~$\bgamma$ with nonzero coefficients. Let $\CM$ be the set of all multicurves that occur as the supports of cycles $\bgamma\in \CB$. Then, for each $M\in\CM$, the set of all cycles $\bgamma\in \CB$ with supports contained in~$M$ is a (finite-dimensional) convex polytope $P_M\subseteq\R_{\ge0}^{\fC}$. Moreover, these convex polytopes meet face-to-face and form a regular CW-decomposition of~$\CB$, see~\cite[Section~2]{BBM10}.

The group $\Mod(\Sigma,x)$ acts cellularly on~$\CB$. Note that the whole mapping class group~$\Mod(\Sigma)$ does not act on~$\CB$ because of the condition~(1). From a result of Ivanov~\cite[Corollary~1.8]{Iva92} it follows that $\I$ acts on~$\CB$ without rotations. Bestvina, Bux, and Margalit proved that the complex $\CB$ is contractible; see Theorem~E in~\cite{BBM10}. It follows that the spectral sequence~\eqref{eq_SS_I3} does indeed converges to~$H_*(\I)$.

The vertices of~$\CB$ are exactly \textit{basic cycles} for~$x$, i.e., formal sums $\bgamma=\sum_{i=1}^kn_i\gamma_i$ (with $n_i>0$) that satisfy conditions~(1)--(3) and such that the homology classes $[\gamma_1],\ldots,[\gamma_k]$ are linearly independent. Unlike in the case of the complex of relative cycles, the number of components of the support of a basic cycle can vary; namely, it can be any number from~$1$ to the genus~$g$ of the surface (in our case, $g=3$). Since $x$ is an integral homology class, the coefficients~$n_i$ of a basic cycle for~$x$ are positive integers.

Moreover, the set~$\CM$ indexing the cells of~$\CB$ admits the following characterization.

\begin{propos}[\protect{\cite[Section~2]{BBM10}}]\label{propos_M_char}
 An oriented multicurve $M$ belongs to~$\CM$ if and only if it satisfies the following conditions:
 \begin{itemize}
  \item[(i)] each component of~$M$ appears in the support of some basic cycle for~$x$ whose support is contained in~$M$,
  \item[(ii)] the components of~$M$ satisfy no one-sided relations.
 \end{itemize}
\end{propos}

We will also need the following proposition.

\begin{propos}[\protect{\cite[Lemma~2.1]{BBM10}}]\label{propos_dimPM}
We have
$$
\dim(P_M)=F-1,
$$
where $F$ is the number of connected components of~$\Sigma\setminus M$.
\end{propos}

This proposition implies that the dimension of~$\CB$ is~$2g-3$. Below we always consider the complex of cycles for the genus~$3$ case only; then $\dim(\CB)=3$.

For each~$p$, let~$\CM_p\subseteq\CM$ be the subset consisting of all oriented multicurves~$M$ such that $\dim (P_M)=p$, that is, $\Sigma\setminus M$ has $p+1$ connected components. In particular, any multicurve $M\in\CM_0$ does not separate~$\Sigma$.  Such an oriented multicurve can consist of~$1$, or~$2$, or~$3$ components. Moreover, for any such ~$M$, the class $x$ is a linear combination of the homology  classes of the components of~$M$ with positive coefficients. It follows that, if $M=\alpha$ is a one-component oriented multicurve in~$\CM_0$, then $[\alpha]=x$.

If we ignore the orientations of the components, then there are exactly six topological types of multicurves on~$\Sigma$ that contain no separating components  and whose complement consists of exactly two connected components; see Fig.~\ref{fig_types}. Hence every oriented multicurve $M\in\CM_1$ belongs to one of these types. Let $\CM_1^{(i)}\subseteq\CM_1$ be the subset consisting of all oriented multicurves of type~$i$.

\begin{figure}[!ht]
 \begin{tikzpicture}[scale=.3,
        mid arrow/.style={
            postaction={
                decorate,
                decoration={
                    markings,
                    mark=at position 0.55 with {
                        \arrow[#1]{stealth}
                    }
                }
            }
        },
        mid arrow/.default=black,
    ]

\small
% \definecolor{blue}{rgb}{0, 0, 0.7}
% \definecolor{mygreen}{rgb}{0, 0.4, 0}
\tikzset{every path/.append style={line width=.2mm}}
\tikzset{my dash/.style={dash pattern=on 2pt off 1.5pt}}

% Border lines

\begin{scope}
\draw[color=blue, thick, my dash] (-7.03,0) .. controls (-6,.3) .. (-4.97,0);
\draw[color=blue, thick, my dash] (-3.03,0) .. controls (-2,.3) .. (-.97,0);
\draw[color=blue, thick, my dash] (.97,0) .. controls (2,.3) .. (3.03,0);
\draw[color=blue, thick, my dash] (4.97,0) .. controls (6,.3) .. (7.03,0);

\draw [very thick] (-4,0) circle (1);
\draw [very thick] (0,0) circle (1);
\draw [very thick] (4,0) circle (1);
\draw [very thick]  (-4,3)--(4,3) arc (90:-90:3) -- (-4,-3) arc (270:90:3);

% Visible lines

\draw[color=blue, thick] (-7.03,0) .. controls (-6,-.3) .. (-4.97,0);
\draw[color=blue, thick] (-3.03,0) .. controls (-2,-.3) .. (-.97,0);
\draw[color=blue, thick] (.97,0) .. controls (2,-.3) .. (3.03,0);
\draw[color=blue, thick] (4.97,0) .. controls (6,-.3) .. (7.03,0);

\draw (0, -5.5) node {type 1};
\end{scope}

\begin{scope}[xshift=17cm]
\draw[color=blue, thick, my dash] (-7.03,0) .. controls (-6,.3) .. (-4.97,0);
\draw[color=blue, thick, my dash] (-3.03,0) .. controls (-2,.3) .. (-.97,0);
\draw[color=blue, thick, my dash] (0,.97) .. controls (.3,2) .. (0,3.03);
\draw[color=blue, thick, my dash] (4.97,0) .. controls (6,.3) .. (7.03,0);

\draw [very thick] (-4,0) circle (1);
\draw [very thick] (0,0) circle (1);
\draw [very thick] (4,0) circle (1);
\draw [very thick]  (-4,3)--(4,3) arc (90:-90:3) -- (-4,-3) arc (270:90:3);

% Visible lines

\draw[color=blue, thick] (-7.03,0) .. controls (-6,-.3) .. (-4.97,0);
\draw[color=blue, thick] (-3.03,0) .. controls (-2,-.3) .. (-.97,0);
\draw[color=blue, thick] (0,.97) .. controls (-.3,2) .. (0,3.03);
\draw[color=blue, thick] (4.97,0) .. controls (6,-.3) .. (7.03,0);

\draw (0, -5.5) node {type 2};
\end{scope}

\begin{scope}[xshift=34cm]
\draw[color=blue, thick, my dash] (-7.03,0) .. controls (-6,.3) .. (-4.97,0);
\draw[color=blue, thick, my dash] (0,-.97) .. controls (.3,-2) .. (0,-3.03);
\draw[color=blue, thick, my dash] (0,.97) .. controls (.3,2) .. (0,3.03);
\draw[color=blue, thick, my dash] (4.97,0) .. controls (6,.3) .. (7.03,0);

\draw [very thick] (-4,0) circle (1);
\draw [very thick] (0,0) circle (1);
\draw [very thick] (4,0) circle (1);
\draw [very thick]  (-4,3)--(4,3) arc (90:-90:3) -- (-4,-3) arc (270:90:3);

% Visible lines

\draw[color=blue, thick] (-7.03,0) .. controls (-6,-.3) .. (-4.97,0);
\draw[color=blue, thick] (0,-.97) .. controls (-.3,-2) .. (0,-3.03);
\draw[color=blue, thick] (0,.97) .. controls (-.3,2) .. (0,3.03);
\draw[color=blue, thick] (4.97,0) .. controls (6,-.3) .. (7.03,0);

\draw (0, -5.5) node {type 3};
\end{scope}

\begin{scope}[yshift=-12cm]
\draw[color=blue, thick, my dash] (-7.03,0) .. controls (-6,.3) .. (-4.97,0);
\draw[color=blue, thick, my dash] (-3.03,0) .. controls (-2,.3) .. (-.97,0);
\draw[color=blue, thick, my dash] (0,.97) .. controls (.3,2) .. (0,3.03);

\draw [very thick] (-4,0) circle (1);
\draw [very thick] (0,0) circle (1);
\draw [very thick] (4,0) circle (1);
\draw [very thick]  (-4,3)--(4,3) arc (90:-90:3) -- (-4,-3) arc (270:90:3);

% Visible lines

\draw[color=blue, thick] (-7.03,0) .. controls (-6,-.3) .. (-4.97,0);
\draw[color=blue, thick] (-3.03,0) .. controls (-2,-.3) .. (-.97,0);
\draw[color=blue, thick] (0,.97) .. controls (-.3,2) .. (0,3.03);

\draw (0, -5.5) node {type 4};
\end{scope}

\begin{scope}[xshift=17cm,yshift=-12cm]
\draw[color=blue, thick, my dash] (-7.03,0) .. controls (-6,.3) .. (-4.97,0);
\draw[color=blue, thick, my dash] (0,-.97) .. controls (.3,-2) .. (0,-3.03);
\draw[color=blue, thick, my dash] (0,.97) .. controls (.3,2) .. (0,3.03);

\draw [very thick] (-4,0) circle (1);
\draw [very thick] (0,0) circle (1);
\draw [very thick] (4,0) circle (1);
\draw [very thick]  (-4,3)--(4,3) arc (90:-90:3) -- (-4,-3) arc (270:90:3);

% Visible lines

\draw[color=blue, thick] (-7.03,0) .. controls (-6,-.3) .. (-4.97,0);
\draw[color=blue, thick] (0,-.97) .. controls (-.3,-2) .. (0,-3.03);
\draw[color=blue, thick] (0,.97) .. controls (-.3,2) .. (0,3.03);

\draw (0, -5.5) node {type 5};
\end{scope}

\begin{scope}[xshift=34cm,yshift=-12cm]
\draw[color=blue, thick, my dash] (0,-.97) .. controls (.3,-2) .. (0,-3.03);
\draw[color=blue, thick, my dash] (0,.97) .. controls (.3,2) .. (0,3.03);

\draw [very thick] (-4,0) circle (1);
\draw [very thick] (0,0) circle (1);
\draw [very thick] (4,0) circle (1);
\draw [very thick]  (-4,3)--(4,3) arc (90:-90:3) -- (-4,-3) arc (270:90:3);

% Visible lines

\draw[color=blue, thick] (0,-.97) .. controls (-.3,-2) .. (0,-3.03);
\draw[color=blue, thick] (0,.97) .. controls (-.3,2) .. (0,3.03);

\draw (0, -5.5) node {type 6};
\end{scope}

%%%%%%%
%%%%%%%

\end{tikzpicture}
 \caption{Types of multicurves in~$\CM_1$}\label{fig_types}
 \vspace{10mm}

 \begin{tikzpicture}[scale=.3,
        mid arrow/.style={
            postaction={
                decorate,
                decoration={
                    markings,
                    mark=at position 0.55 with {
                        \arrow[#1]{stealth}
                    }
                }
            }
        },
        mid arrow/.default=black,
    ]

\small
% \definecolor{blue}{rgb}{0, 0, 0.7}
% \definecolor{mygreen}{rgb}{0, 0.4, 0}
\tikzset{every path/.append style={line width=.2mm}}
\tikzset{my dash/.style={dash pattern=on 2pt off 1.5pt}}

% Border lines

\begin{scope}
\draw[color=blue, thick, my dash] (-7.03,0) .. controls (-6,.3) .. (-4.97,0);
\draw[color=blue, thick, my dash] (-3.03,0) .. controls (-2,.3) .. (-.97,0);
\draw[color=blue, thick, my dash] (.97,0) .. controls (2,.3) .. (3.03,0);
\draw[color=blue, thick, my dash] (4.97,0) .. controls (6,.3) .. (7.03,0);
\draw[color=blue, thick, my dash] (0,.97) .. controls (.3,2) .. (0,3.03);

\draw [very thick] (-4,0) circle (1);
\draw [very thick] (0,0) circle (1);
\draw [very thick] (4,0) circle (1);
\draw [very thick]  (-4,3)--(4,3) arc (90:-90:3) -- (-4,-3) arc (270:90:3);

% Visible lines

\draw[color=blue, thick] (-7.03,0) .. controls (-6,-.3) .. (-4.97,0);
\draw[color=blue, thick] (-3.03,0) .. controls (-2,-.3) .. (-.97,0);
\draw[color=blue, thick] (.97,0) .. controls (2,-.3) .. (3.03,0);
\draw[color=blue, thick] (4.97,0) .. controls (6,-.3) .. (7.03,0);
\draw[color=blue, thick] (0,.97) .. controls (-.3,2) .. (0,3.03);

\draw (0, -5.5) node {type 1};
\end{scope}

\begin{scope}[xshift=17cm]
\draw[color=blue, thick, my dash] (-7.03,0) .. controls (-6,.3) .. (-4.97,0);
\draw[color=blue, thick, my dash] (-3.03,0) .. controls (-2,.3) .. (-.97,0);
\draw[color=blue, thick, my dash] (0,.97) .. controls (.3,2) .. (0,3.03);
\draw[color=blue, thick, my dash] (4.97,0) .. controls (6,.3) .. (7.03,0);
\draw[color=blue, thick, my dash] (0,-.97) .. controls (.3,-2) .. (0,-3.03);

\draw [very thick] (-4,0) circle (1);
\draw [very thick] (0,0) circle (1);
\draw [very thick] (4,0) circle (1);
\draw [very thick]  (-4,3)--(4,3) arc (90:-90:3) -- (-4,-3) arc (270:90:3);

% Visible lines

\draw[color=blue, thick] (-7.03,0) .. controls (-6,-.3) .. (-4.97,0);
\draw[color=blue, thick] (-3.03,0) .. controls (-2,-.3) .. (-.97,0);
\draw[color=blue, thick] (0,.97) .. controls (-.3,2) .. (0,3.03);
\draw[color=blue, thick] (4.97,0) .. controls (6,-.3) .. (7.03,0);
\draw[color=blue, thick] (0,-.97) .. controls (-.3,-2) .. (0,-3.03);

\draw (0, -5.5) node {type 2};
\end{scope}

\begin{scope}[xshift=34cm]
\draw[color=blue, thick, my dash] (-7.03,0) .. controls (-6,.3) .. (-4.97,0);
\draw[color=blue, thick, my dash] (-3.03,0) .. controls (-2,.3) .. (-.97,0);
\draw[color=blue, thick, my dash] (0,.97) .. controls (.3,2) .. (0,3.03);
% \draw[color=blue, thick, my dash] (4.97,0) .. controls (6,.3) .. (7.03,0);
\draw[color=blue, thick, my dash] (0,-.97) .. controls (.3,-2) .. (0,-3.03);

\draw [very thick] (-4,0) circle (1);
\draw [very thick] (0,0) circle (1);
\draw [very thick] (4,0) circle (1);
\draw [very thick]  (-4,3)--(4,3) arc (90:-90:3) -- (-4,-3) arc (270:90:3);

% Visible lines

\draw[color=blue, thick] (-7.03,0) .. controls (-6,-.3) .. (-4.97,0);
\draw[color=blue, thick] (-3.03,0) .. controls (-2,-.3) .. (-.97,0);
\draw[color=blue, thick] (0,.97) .. controls (-.3,2) .. (0,3.03);
% \draw[color=blue, thick] (4.97,0) .. controls (6,-.3) .. (7.03,0);
\draw[color=blue, thick] (0,-.97) .. controls (-.3,-2) .. (0,-3.03);

\draw (0, -5.5) node {type 3};
\end{scope}

%%%%%%%
%%%%%%%

\end{tikzpicture}
 \caption{Types of multicurves in~$\CM_2$}\label{fig_types_M2}
\end{figure}

Similarly, there are exactly three topological types of multicurves on~$\Sigma$ that contain no separating components  and whose complement consists of exactly three connected components; see Fig.~\ref{fig_types_M2}. Let $\CM_2^{(i)}\subseteq\CM_2$ be the subset consisting of all oriented multicurves of type~$i$.

\section{Important cells and chain complex~$\CC_*$}\label{section_important}

A multicurve without separating components on a closed surface of genus~$3$ contains at most one bounding pair.
We say that an oriented multicurve $M\in\CM$ (or the corresponding cell~$P_M$) is  \textit{important} if the following conditions hold:
\begin{itemize}
 \item $M$ contains a bounding  pair~$\{\gamma,\gamma'\}$;
 \item let $S_1$ and~$S_2$ be the connected components of~$\Sigma\setminus(\gamma\cup\gamma')$; then $x$ is contained in neither of the subgroups $H_1(S_i)\subseteq H_1(\Sigma)$ for $i=1,2$.
\end{itemize}
All other oriented multicurves $M\in\CM$, and the corresponding cells~$P_M$, are called \textit{unimportant}. The property of being important is invariant under the action of the Torelli group~$\I$. A cell of the quotient complex~$\CB/\I$ is be called \textit{important} if it corresponds to an $\I$-orbit consisting of important cells of~$\CB$.

The following proposition is a straightforward consequence of the definition.

\begin{propos}\label{propos_important_obvious}
 \textnormal{(1)} All multicurves in~$\CM_0$ are unimportant.

\textnormal{(2)} A multicurve in~$\CM_1$ is important if and only if it has type~3.

\textnormal{(3)} If $M_1\subseteq M_2$, $M_1,M_2\in\CM$, and $M_1$ is important, then $M_2$ is also important. Hence unimportant cells form a subcomplex $\CU\subseteq\CB$.
\end{propos}

Recall that we have a canonical isomorphism $E^1_{*,0}\cong C_*(\CB/\I)$, and the differential~$d^1$ coincides with the cellular boundary map~$\partial$. Let~$\CP_M$ denote the cell of~$\CB/\I$ corresponding to the $\I$-orbit of the cell~$P_M$. Then $\CP_{hM}=\CP_M$ for all $h\in\I$.

We denote by~$\CC_*$ the quotient of the chain complex~$C_*(\CB/\I)$ by the subcomplex generated by all unimportant cells. In other words, $\CC_*$ is the cellular chain complex of the quotient CW complex $(\CB/\I)/(\CU/\I)$. We denote the differential in~$\CC_*$ by~$\partial_{\CC}$.

We need the following property of the chain complex~$\CC_*$.

\begin{propos}\label{propos_H2C0}
 We have $H_2(\CC_*)=0$.
\end{propos}

In the rest of this section we prove this proposition.

Following Hatcher and Margalit~\cite{HaMa12}, we consider the \textit{weight function}
$$
n\colon\CB\to\R_{\ge 0}
$$
sending every cycle $\bgamma=\sum n_i\gamma_i$ to~$\sum n_i$. (In~\cite{HaMa12} this function was denoted by~$W$.) Since the coefficients~$n_i$ are positive, we have $n(\bgamma)>0$ for all $\bgamma\in\CB$. Moreover, if a cycle~$\bgamma$ is basic, then $n(\bgamma)$ is a positive integer. The restriction of~$n$ to every cell~$P_M$ of~$\CB$ is affine linear. We set
$$
n(M)=\max_{\bgamma\in P_M}n(\bgamma).
$$
Since the maximum is attained at a vertex of~$P_M$, we obtain that $n(M)$ is a positive integer. Moreover, if $M$ contains a component $\gamma$ with $[\gamma]\ne x$, then $n(M)\ge 2$, since there exists a basic $1$-cycle for~$\gamma$ whose support is contained in~$M$ and contains~$\gamma$, and the sum of coefficients of this $1$-cycle is at least~$2$.

\begin{propos}\label{propos_important_contain}
 Any important oriented multicurve~$M\in\CM_p$, where $p$ is either~$2$ or~$3$, contains an important oriented submulticurve $K\in\CM_1$.
\end{propos}

\begin{proof}
 Let $\{\gamma,\gamma'\}$ be the bounding pair contained in~$M$. By Proposition~\ref{propos_M_char}, there exists a basic cycle for~$x$ whose support~$N$ contains~$\gamma$ and is contained in~$M$. Then $N\in \CM_0$. The homology classes of the components of~$N$ are linearly independent, hence $\gamma'$ is not contained in~$N$. The oriented multicurve $K=N\cup\gamma'$ belongs to~$\CM_1$ and is important.
\end{proof}

We now introduce the following notation. If $M\in\CM_p$ is important, where $p$ is either~$2$ or~$3$, then we denote by~$n'(M)$ the greatest of the values~$n(K)$, where $K$ runs over all submulticurves of~$M$ that belong to~$\CM_1$ and are important. It follows immediately from the definition that $n'(M)\le n(M)$. However, this inequality can be strict. Namely, this happens when the maximum $n(M)$ is attained at a vertex of~$P_M$ that does not lie in any edge~$P_K$ such that $K\subseteq M$ and~$K$ is important.

\begin{propos}\label{propos_n_unique}
Let $M\in\CM_p$ be an important oriented multicurve, where $p$ is either~$2$ or~$3$. Then there is exactly one oriented submulticurve $K\subseteq M$ such that $K\in\CM_1$, $K$ is important, and $n(K)=n'(M)$.
\end{propos}

\begin{proof}
 Suppose to the contrary that $M$ contains two distinct oriented submulticurves~$K_1$ and~$K_2$ such that $n(K_1)=n(K_2)=n'(M)$. Let $\{\gamma,\gamma'\}$ be the bounding pair contained in~$M$. Then both~$K_1$ and~$K_2$ contain~$\gamma\cup\gamma'$. Let $N_1$ and~$N_2$ be the oriented multicurves obtained from~$K_1$ and~$K_2$, respectively, by removing the component~$\gamma'$. By Proposition~\ref{propos_M_char}, $N_1,N_2\in\CM_0$ and~$n(N_i)=n(K_i)=n'(M)$. Consider the oriented multicurve $N=N_1\cup N_2$. Then $N$ belongs to~$\CM$. Indeed, first, every component of~$N$ is contained in the support of a $1$-cycle supported in either~$N_1$ or~$N_2$, and hence in~$N$. Second, since $N\subseteq M$ and $M\in\CM$, we see that the homology classes of components of~$N$ satisfy no one-sided relations. Moreover, the multicurve~$N\cup\gamma'$ also belongs to~$\CM$. Since $\dim\CB=3$ and $\dim P_{N\cup\gamma'}>\dim P_N$, we obtain that $N\in\CM_p$ for some $p\le 2$. On the other hand, $p\ge 1$, since $P_N$ contains two distinct vertices~$P_{N_1}$ and~$P_{N_2}$. Consider two cases.
 \smallskip

 \textsl{Case 1: $\dim P_N=1$.} Then $P_N$ is the segment with endpoints~$P_{N_1}$ and~$P_{N_2}$. Let $S$ be one of the two connected components of~$\Sigma\setminus N$. Then the homology classes of the components of~$N$ satisfy the only linear relation corresponding to the boundary of~$S$. Since the values of the function~$n$ at the endpoints~$P_{N_1}$ and~$P_{N_2}$ of the segment~$P_N$ are equal, we obtain that this linear relation is balanced, i.e., contains equal numbers of terms with coefficients~$+1$ and~$-1$. Therefore $S$ is a sphere with an even number of punctures. Since $N$ contains no bounding pair, it follows that $N$ has type~1 (see Fig.~\ref{fig_types}). We obtain a contradiction because it is then impossible to add to~$N$ a component~$\gamma'$ that is disjoint from~$N$ and homologous to a component of~$N$.
\smallskip

\textsl{Case 2: $\dim P_N=2$.} The multicurve $N$ does not contain a bounding pair. Hence $N$ has type~1. This means that, if we ignore the orientations of the components, the multicurve~$N$ can be drawn as in Fig.~\ref{fig_M2}. Furthermore, we have $\gamma=\alpha_0$, since it is possible to add to~$N$ the component~$\gamma'$ that is homologous to~$\gamma$ and disjoint from~$N$. Thus both $N_1$ and~$N_2$ contain~$\gamma=\alpha_0$, $N_1\cup N_2=N$, and each of the multicurves~$N_1$ and~$N_2$ does not separate~$\Sigma$. Without loss of generality, we may assume that $N_1=\alpha_0\cup\alpha_1\cup\alpha_3$ and $N_2=\alpha_0\cup\alpha_2\cup\alpha_4$. We now return to the orientations of the components of~$N$. Let $S$ be the connected component of~$\Sigma\setminus N$ homeomorphic to a four-punctured sphere. Assume that $S$ lies either to the right of both~$\alpha_1$ and~$\alpha_2$ or to the left of both~$\alpha_1$ and~$\alpha_2$. Then there exists a simple closed curve~$\beta$ disjoint from $\alpha_0\cup\alpha_3\cup\alpha_4$ such that $[\alpha_1]\cdot[\beta]=1$ and~$[\alpha_2]\cdot [\beta]=-1$. Since $x$ is a linear combination with positive coefficients of $[\alpha_0]$, $[\alpha_1]$, and $[\alpha_3]$, and at the same time a linear combination with positive coefficients of $[\alpha_0]$, $[\alpha_2]$, and $[\alpha_4]$, we arrive at a contradiction because $x\cdot [\beta]$ must be simultaneously positive and negative. Therefore, without loss of generality, we may assume that $S$ lies to the left of~$\alpha_1$ and to the right of~$\alpha_2$, and $[\alpha_0]=[\alpha_1]-[\alpha_2]$. Let $x=n_0[\alpha_0]+n_1[\alpha_1]+n_3[\alpha_3]$; then $n_0,n_1,n_3>0$ and $n(N_1)=n_0+n_1+n_3$. Now we set $N_3=\alpha_0\cup\alpha_2\cup\alpha_3$ and $K_3=N_3\cup\gamma'$. Then
$$
x=(n_0+n_1)[\alpha_0]+n_1[\alpha_2]+n_3[\alpha_3].
$$
Hence $N_3\in\CM_0$ and $n(N_3)>n(N_1)$. Consequently, $K_3\in \CM_1$, $K_3$ is important, and $n(K_3)>n(K_1)$. Since $K_3\subseteq M$, we obtain a contradiction with the assumption that $n(K_1)=n'(M)$.
\end{proof}

\begin{figure}[!ht]
 \begin{tikzpicture}[scale=.35,
        mid arrow/.style={
            postaction={
                decorate,
                decoration={
                    markings,
                    mark=at position 0.55 with {
                        \arrow[#1]{stealth}
                    }
                }
            }
        },
        mid arrow/.default=black,
    ]

\small
% \definecolor{blue}{rgb}{0, 0, 0.7}
% \definecolor{mygreen}{rgb}{0, 0.4, 0}
\tikzset{every path/.append style={line width=.2mm}}
\tikzset{my dash/.style={dash pattern=on 2pt off 1.5pt}}

% Border lines

\begin{scope}
\draw[color=blue, thick, my dash] (-7.03,0) .. controls (-6,.3) .. (-4.97,0);
\draw[color=blue, thick, my dash] (-3.03,0) .. controls (-2,.3) .. (-.97,0);
\draw[color=blue, thick, my dash] (.97,0) .. controls (2,.3) .. (3.03,0);
\draw[color=blue, thick, my dash] (4.97,0) .. controls (6,.3) .. (7.03,0);
\draw[color=blue, thick, my dash] (0,.97) .. controls (.3,2) .. (0,3.03);

\draw [very thick] (-4,0) circle (1);
\draw [very thick] (0,0) circle (1);
\draw [very thick] (4,0) circle (1);
\draw [very thick]  (-4,3)--(4,3) arc (90:-90:3) -- (-4,-3) arc (270:90:3);

% Visible lines

\draw[color=blue, thick] (-7.03,0) .. controls (-6,-.3) .. (-4.97,0) node [pos=.5, below,xshift=.8mm] {$\alpha_1$};
\draw[color=blue, thick] (-3.03,0) .. controls (-2,-.3) .. (-.97,0) node [pos=.5, below,xshift=.8mm] {$\alpha_2$};
\draw[color=blue, thick] (.97,0) .. controls (2,-.3) .. (3.03,0) node [pos=.5, below,xshift=.8mm] {$\alpha_3$};
\draw[color=blue, thick] (4.97,0) .. controls (6,-.3) .. (7.03,0) node [pos=.5, below] {$\alpha_4$};
\draw[color=blue, thick] (0,.97) .. controls (-.3,2) .. (0,3.03) node [pos=.5, left,xshift=.8mm] {$\alpha_0$};

\draw(4.2,-2) node {$S$};

\end{scope}

%%%%%%%
%%%%%%%

\end{tikzpicture}
 \caption{A multicurve $M\in\CM_2^{(1)}$}\label{fig_M2}
\end{figure}

Now, let $K\in\CM_1$ be an important oriented multicurve. By Proposition~\ref{propos_important_obvious}, $K$ has type~3. Hence $K=\gamma_0\cup\gamma_0'\cup\gamma_1\cup\gamma_2$, where $\{\gamma_0,\gamma_0'\}$ is a bounding pair, and the components~$\gamma_1$ and~$\gamma_2$ lie in the two different connected components~$S_1$ and~$S_2$  of~$\Sigma\setminus (\gamma_0\cup\gamma_0')$. For $i=1,2$, we set $R_i=S_i\setminus\gamma_i$. Let $\Upsilon(K)$ be the set consisting of all oriented multicurves~$M\in\CM_2$ such that~$M\supseteq K$ and $n'(M)=n(K)$.  For $i=1,2$, let $\Upsilon_i(K)$ be the subset of~$ \Upsilon(K)$ consisting of all~$M$ such that  $M\setminus K$ is contained in~$R_i$.

\begin{propos}\label{propos_Gamma_prop}
An oriented multicurve $M$ belongs to~$\Upsilon(K)$ if and only if $M=K\cup\delta$, where $\delta$ is a simple closed curve in~$R_i$ (where $i$ is either~$1$ or~$2$) such that  $[\delta]$ is either $[\gamma_0]+[\gamma_i]$ or $[\gamma_0]-[\gamma_i]$. In particular,
\begin{enumerate}
 \item $\Upsilon(K)=\Upsilon_1(K)\sqcup\Upsilon_2(K)$, and
 \item both subsets~$\Upsilon_1(K)$ and~$\Upsilon_2(K)$ are nonempty.
\end{enumerate}
Moreover,  if $M_1\in \Upsilon_1(K)$ and~$M_2\in\Upsilon_2(K)$, then the oriented multicurve $M_1\cup M_2$ belongs to~$\CM_3$, and $n'(M_1\cup M_2)=n(K)$.
\end{propos}

\begin{proof}
 Since $K\in \CM_1^{(3)}$, we have $x=n_0[\gamma_0]+n_1[\gamma_1]+n_2[\gamma_2]$ for some positive integers~$n_0,n_1,n_2$.

 For $i=1,2$, $R_i$ is a four-punctured sphere; hence any simple closed curve on~$R_i$ separates it. By Proposition~\ref{propos_dimPM}, any multicurve $M\in\CM_2$ divides~$\Sigma$ into  exactly three parts. Therefore any multicurve $M\in\Upsilon(K)$ is obtained from~$K$ by adding  exactly one component~$\delta$ lying in either~$R_1$ or~$R_2$. Suppose that $\delta\subseteq R_i$; let $j$ be the element of~$\{1,2\}$ different from~$i$. The curve~$\delta$ is not separating in~$\Sigma$. Moreover, the homology classes of the components of~$M$ satisfy no one-sided relations. It follows that $[\delta]$ is $[\gamma_0]+[\gamma_i]$, or~$[\gamma_0]-[\gamma_i]$, or~$[\gamma_i]-[\gamma_0]$.
 However, if $[\delta]=[\gamma_i]-[\gamma_0]$, then the multicurve $K_1=\gamma_0\cup\gamma_0'\cup\delta\cup\gamma_j$ belongs to~$\CM_1$, is important, and satisfies $n(K_1)=n_0+2n_i+n_j>n(K)$. Therefore $M\notin\Upsilon(K)$ in this case.

 Let us now show that $M\in\Upsilon(K)$ if $[\delta]$ is either $[\gamma_0]+[\gamma_i]$ or~$[\gamma_0]-[\gamma_i]$. If $[\delta]=[\gamma_0]-[\gamma_i]$, then $P_M$ is the triangle with the vertices~$P_{\gamma_0\cup\gamma_1\cup\gamma_2}$, $P_{\gamma_0'\cup\gamma_1\cup\gamma_2}$, and $P_{\delta\cup\gamma_1\cup\gamma_2}$. Hence $K$  is the only submulticurve of~$M$ that belongs to~$\CM_1$ and is important. Therefore $M\in\Upsilon(K)$.

 Now, suppose that $[\delta]=[\gamma_0]+[\gamma_i]$. Then there are two subcases. If $n_0\le n_i$, then $P_M$ is again a triangle whose vertices are~$P_{\gamma_0\cup\gamma_1\cup\gamma_2}$, $P_{\gamma_0'\cup\gamma_1\cup\gamma_2}$, and either $P_{\delta\cup\gamma_1\cup\gamma_2}$ (if $n_0<n_i$) or $P_{\delta\cup\gamma_j}$ (if $n_0=n_i$). Consequently, $K$  is again the only submulticurve of~$M$ that belongs to~$\CM_1$ and is important, and therefore $M\in\Upsilon(K)$. If $n_0>n_i$, then $P_M$ is the quadrangle  with the vertices~$P_{\gamma_0\cup\gamma_1\cup\gamma_2}$, $P_{\gamma_0'\cup\gamma_1\cup\gamma_2}$, $P_{\gamma_0\cup\delta\cup\gamma_j}$, and~$P_{\gamma_0'\cup\delta\cup\gamma_j}$. Therefore, $M$ contains exactly two important submulticurves belonging to~$\CM_1$, namely, $K$ and $K_1=\gamma_0\cup\gamma_0'\cup\delta\cup\gamma_j$. Since $n(K_1)=n_0+n_j<n(K)$, we have $M\in\Upsilon(K)$.

 Finally, let us prove the ``moreover'' part of the proposition. Since $M_1\in\Upsilon_1(K)$ and $M_2\in\Upsilon_2(K)$, we  have $M_1=K\cup\delta_1$ and~$M_2=K\cup\delta_2$, where $\delta_i\subseteq R_i$ and $[\delta_i]=[\gamma_0]\pm[\gamma_i]$. Then the oriented multicurve $M=M_1\cup M_2=K\cup\delta_1\cup\delta_2$ belongs to~$\CM$. Indeed, the homology classes of components of~$M$ satisfy no one-sided relations, and  every component of~$M$ is contained in the support of a basic cycle supported in~$M$. Moreover, $M$ is important, since $K$ is important, and $M\in\CM_3$, since $M$ divides~$\Sigma$ into four pieces. Let us  prove that $n'(M)=n(K)$. To do this, we need to show that $M$ contains no important submulticurve $K'\in \CM_1$ such that $n(K')>n(K)$. If such a submulticurve~$K'$ was contained in either~$M_1$ or~$M_2$, then the inequality $n(K')\le n(K)$ would hold, since $n'(M_1)=n'(M_2)=n(K)$. Hence the only candidate for an important submulticurve~$K'\subseteq M$ with $n(K')>n(K)$ is $K'=\gamma_0\cup\gamma_0'\cup\delta_1\cup\delta_2$. The multicurve $K'$ belongs to~$\CM_1$ only if $[\delta_i]=[\gamma_0]+[\gamma_i]$ for $i=1,2$, and $n_0>n_1+n_2$, since otherwise $x$ is not a linear combination of the classes~$[\gamma_0]$, $[\delta_1]$, and~$[\delta_2]$ with positive coefficients. Nevertheless, even if these conditions are  satisfied, then we still have $n(K')=n_0<n(K)$. The proposition follows.
\end{proof}

\begin{proof}[Proof of Proposition~\ref{propos_H2C0}]
Let $\zeta\in\CC_2$ be an arbitrary chain such that $\partial_{\CC}\zeta=0$.
We write
\begin{equation*}
\zeta=\sum k_M\CP_M,
\end{equation*}
where the sum is finite, $M\in\CM_2$ are important oriented multicurves from pairwise distinct $\I$-orbits, and $k_M\in\Z\setminus\{0\}$. If $\zeta\ne 0$, then we denote by $n(\zeta)$ the maximum of the numbers~$n'(M)$ over all $2$-cells~$\CP_M$ that enter~$\zeta$, and we denote by~$r(\zeta)$ the number of $1$-cells~$\CP_K$ such that
\begin{itemize}
 \item $\CP_K$ is contained in at least one of the $2$-cells~$\CP_M$ that enter the chain~$\zeta$,
 \item $K$ is important,
 \item $n(K)=n(\zeta)$.
\end{itemize}
From the definition of the number~$n(\zeta)$ it follows that there exists at least one $1$-cell~$\CP_K$ that satisfies these conditions, so $r(\zeta)\ge 1$. (Recall that $\CP_{hK}=\CP_K$ whenever $h\in\I$. We stress that $r(\zeta)$ is the number of distinct $1$-cells~$\CP_K$ that satisfy the above conditions rather than the number of distict oriented multicurves~$K$.)

Let us prove that $\zeta$ belongs to the image of $\partial_{\CC}\colon\CC_3\to\CC_2$ by the induction on~$n=n(\zeta)$ and, for a given~$n$, by the induction on~$r=r(\zeta)$. The base case for $\zeta=0$ is trivial.

Assume that $\zeta\ne 0$; then $n(\zeta)\ge 1$ and $r(\zeta)\ge 1$. Let $\CP_K$ be an important $1$-cell such that $\CP_K$ is contained in a $2$-cell that enters~$\zeta$, and $n(K)=n(\zeta)$. We write
\begin{equation}\label{eq_sum_beps}
\zeta=\sum_{i=1}^mk_i\CP_{M_i}+\ldots,
\end{equation}
where every $2$-cell~$\CP_{M_i}$ contains~$\CP_K$, and the dots stand for a linear combination of $2$-cells that do not contain~$\CP_K$. We choose the orientations of the cells~$\CP_K$ and~$\CP_{M_i}$ so that the incidence numbers~$\left[\CP_{M_i}:\CP_K\right]$ equal~$1$ for all~$i$. Since $\partial_{\CC}\zeta=0$, we have
\begin{equation}\label{eq_m_sum}
k_1+\cdots+k_m=0.
\end{equation}
Since $\CP_{hM}=\CP_M$ for all $h\in\I$, we may assume that $M_i\supseteq K$ for $i=1,\ldots,m$. Moreover, we see that $n'(M_i)=n(K)$ for all~$i$, since otherwise we would have $n(\zeta)\ge n'(M_i)>n(K)$. Hence $M_i\in\Upsilon(K)$ for all~$i$. We may assume that $M_i\in \Upsilon_1(K)$ for $1\le i\le m_1$ and $M_i\in \Upsilon_2(K)$ for $m_1+1\le i\le m$.

Moreover, we may assume that both numbers~$m_1$ and $m-m_1$ are nonzero. Indeed, if, say, $m-m_1=0$, that is, all multicurves $M_1,\ldots,M_m$ belong to~$\Upsilon_1(K)$, then we can choose an oriented multicurve $\widetilde{M}\in \Upsilon_2(K)$ (since $ \Upsilon_2(K)\ne \varnothing$ by Proposition~\ref{propos_Gamma_prop}) and add to the sum~\eqref{eq_sum_beps} a ``phantom'' summand $0\cdot\CP_{\widetilde{M}}$.

By Proposition~\ref{propos_Gamma_prop}, we have that $M_i\cup M_j\in\CM_3$ and $n'(M_i\cup M_j)=n(K)$ whenever $1\le i\le m_1$ and~$m_1+1\le j\le m$. We orient the $3$-cell~$\CP_{M_i\cup M_j}$ so that to obtain that $\bigl[\CP_{M_i\cup M_j}:\CP_{M_i}\bigr]=1$; then $\bigl[\CP_{M_i\cup M_j}:\CP_{M_j}\bigr]=-1$.

Now, choose integers $l_{ij}$, where $1\le i\le m_1$ and $m_1+1\le j\le m$, so that to obtain that
\begin{align*}
\sum_{j=m_1+1}^ml_{ij}&=k_i,& i&=1,\ldots,m_1,\\
\sum_{i=1}^{m_1}l_{ij}&=-k_j,&j&=m_1+1,\ldots,m.
\end{align*}
The existence of such numbers~$l_{ij}$ is provided by~\eqref{eq_m_sum}. Consider the $3$-chain
$$
\xi=\sum_{i=1}^{m_1}\sum_{j=m_1+1}^ml_{ij}\CP_{M_i\cup M_j}.
$$
Then
$$
\partial_{\CC}\xi=\sum_{i=1}^mk_i\CP_{M_i}+\ldots,
$$
where the dots stand for a linear combination of those important $2$-faces~$\CP_N$ of the $3$-cells~$\CP_{M_i\cup M_j}$ that do not contain the edge~$\CP_K$. From Proposition~\ref{propos_n_unique} it follows that $n'(N)<n(K)=n(\zeta)$ for each of such faces~$\CP_N$. Therefore the chain
$$
\zeta'=\zeta-\partial_{\CC}\xi
$$
satisfies either $n(\zeta')<n(\zeta)$ or $n(\zeta')=n(\zeta)$ and $r(\zeta')<r(\zeta)$. By the induction hypothesis, $\zeta'$ lies in the image of~$\partial_{\CC}$. Therefore $\zeta$ lies in the image of~$\partial_{\CC}$ as well. The proposition follows.
\end{proof}

\section{Reduction to~$\CF_{1,1}$}\label{section_toF11}

In this section, we prove the following proposition. Recall that $x$ is the primitive homology class that was used in the construction of the complex of cycles~$\CB=\CB(x)$.

\begin{propos}\label{propos_txE20}
We have, $(t_x-1)E^2_{2,0}=0$, and hence $(t_x-1)H_2(\I)\subseteq\CF_{1,1}$.
\end{propos}

To derive this proposition from Proposition~\ref{propos_H2C0}, we need the following assertion.

\begin{propos}\label{propos_unimp_stab}
 If an oriented multicurve $M\in\CM_2$ is unimportant, then $t_x\CP_M=\CP_M$.
\end{propos}

\begin{proof}
Since $M$ is unimportant, we have the following two possibilities.
\smallskip

\textsl{Case 1: $M$ contains no bounding pair.} Then $M\in\CM_2^{(1)}$. Hence, if we ignore the orientations of the components of~$M$, then $M$ can be drawn as in Fig.~\ref{fig_M2}. Let $S$ be the connected component of~$\Sigma\setminus M$ homeomorphic to a four-punctured sphere.
Let $a_1,a_2,a_3$ be the homology classes of the components $\alpha_1,\alpha_2,\alpha_3$, respectively, oriented so that $S$ lies to their left. (These orientations may or may not coincide with the original orientations of these components in~$M$.) With suitable orientations, the homology classes of~$\alpha_0$ and~$\alpha_4$  are~$a_1+a_2$ and  $a_1+a_2+a_3$, respectively. Moreover, $S$ contains a simple closed curve~$\gamma$ with $[\gamma]=a_1+a_3$. Since the Dehn twists $T_{\alpha_0},\ldots,T_{\alpha_4}$, and~$T_{\gamma}$ stabilize~$M$, we obtain that each of the symplectic transvections~$t_{a_1}$, $t_{a_2}$, $t_{a_3}$, $t_{a_1+a_2}$, $t_{a_1+a_3}$, and~$t_{a_1+a_2+a_3}$ stabilizes~$\CP_M$.
Finally, $x$ is a linear combination of the homology classes of~$M$; hence $$x=m_1a_1+m_2a_2+m_3a_3$$ for some $m_i\in\Z$. (The coefficients~$m_i$ can be negative, since we have changed the orientations of the components of~$M$.) From Lemma~\ref{lem_transvec3} it follows that $t_x$ lies in the abelian subgroup of~$\Sp(H,x)$ generated by~$t_{a_1}$, $t_{a_2}$, $t_{a_3}$, $t_{a_1+a_2}$, $t_{a_1+a_3}$, and~$t_{a_1+a_2+a_3}$; specifically,
$$
t_x=t_{a_1+a_2+a_3}^{m_2m_3}t_{a_1+a_2}^{m_2(m_1-m_3)}t_{a_1+a_3}^{m_3(m_1-m_2)}t_{a_1}^{m_1^2-m_1m_2-m_1m_3+m_2m_3}t_{a_2}^{m_2^2-m_1m_2}t_{a_3}^{m_3^2-m_1m_3}.
$$
Therefore $t_x\CP_M=\CP_M$.
\smallskip

\textsl{Case 2: $M$ contains a bounding pair~$\{\gamma,\gamma'\}$, but $x$ lies in the subgroup~$H_1(S_1)\subseteq H_1(\Sigma)$, where $S_1$ is one of the two connected components of~$\Sigma\setminus(\gamma\cup\gamma')$.} The surface~$S_1$ can be divided into two pairs of pants by two simple closed curves~$\alpha_1$ and~$\alpha_2$ each of which is non-separating in~$\Sigma$ and either a component of~$M$ or disjoint from~$M$. Let $a_1$ and~$a_2$ be the homology classes of~$\alpha_1$ and~$\alpha_2$, respectively. The orientations of~$\alpha_1$ and~$\alpha_2$, and hence the signs of~$a_1$ and~$a_2$, can be chosen so that $[\gamma]=a_1+a_2$. (As in the previous case, these orientations may or may not coincide with the original orientations of these components in~$M$.) Since the Dehn twists $T_{\alpha_1}$, $T_{\alpha_2}$, and~$T_{\gamma}$ stabilize~$M$, we obtain that each of the symplectic transvections~$t_{a_1}$, $t_{a_2}$, and~$t_{a_1+a_2}$ stabilizes~$\CP_M$. Since $x\in H_1(S_1)$ and $x$ has zero algebraic intersection numbers with~$a_1$ and~$a_2$, it follows that  $x=m_1a_1+m_2a_2$ for some $m_1,m_2\in \Z$. By Lemma~\ref{lem_transvec3}, $t_x$ lies in the abelian subgroup of~$\Sp(H,x)$ generated by~$t_{a_1}$, $t_{a_2}$, and~$t_{a_1+a_2}$.
Therefore $t_x\CP_M=\CP_M$.
\end{proof}

\begin{proof}[Proof of Proposition~\ref{propos_txE20}]
Let $z\in E^2_{2,0}$, and let $\zeta\in E^1_{2,0}=C_2(\CB/\I)$ be a cycle representing~$z$. Then $d^1\zeta=0$. From Proposition~\ref{propos_H2C0} it follows that there exists a chain $\eta\in E^1_{3,0}=C_2(\CB/\I)$ such that $\zeta+d^1\eta$ is a linear combination of unimportant $2$-cells~$\CP_M$ only. Then, by Proposition~\ref{propos_unimp_stab}, we obtain that $(t_x-1)(\zeta+d^1\eta)=0$, and therefore $(t_x-1)z=0$.
\end{proof}

\section{Active components and homomorphisms~$\nu_{\gamma,\CW}$}\label{section_active}

For an oriented multicurve~$M$, we denote by~$[M]$ the set of homology classes of the components of~$M$, and by~$\langle M\rangle$ the subgroup of~$H$ generated by~$[M]$. This subgroup is always an isotropic direct summand of~$H$. The rank of~$\langle M\rangle$ will be called the \textit{rank} of~$M$. In particular, oriented multicurves of rank~$3$ will be called multicurves of \textit{full rank}. If $M$ has full rank, then $\langle M\rangle$ is a Lagrangian submodule of~$H$.

A multicurve $M\in\CM_0$ has full rank if and only if it consists of three components. A multicurve $M\in\CM_1$ has full rank if and only if it has type~$1$, $2$, or~$3$. A multicurve $M\in\CM_2$ has full rank if and only if it has type either~$1$ or~$2$.

Now, let $M\in\CM_p$ be an oriented multicurve containing no bounding pair. We say that a component~$\gamma$ of~$M$ is \textit{active} if there exists a simple closed curve $\gamma'\subseteq \Sigma\setminus M$ such that the following conditions hold:
\begin{itemize}
 \item $\gamma'$ is not isotopic to a component of~$M$,
 \item $\{\gamma,\gamma'\}$ is a bounding pair,
 \item the multicurve $M\cup\gamma'$ is important.
\end{itemize}
We denote by $\Act(M)$ the set of all active components of~$M$.
We say that~$M$ is \textit{negligible} if it has no active components.   The following is a direct consequence of the definitions.

\begin{propos}\label{propos_active_obvious}
\textnormal{(1)} All oriented multicurves of not full rank are negligible.
\smallskip

\textnormal{(2)} All components of each full rank oriented multicurve in~$\CM_0$ are active.
\smallskip

\textnormal{(3)} If $M\in\CM_1^{(2)}$, then the three components of~$M$ bounding a pair of pants are active, and the fourth component is not. All multicurves in~$\CM_1^{(1)}$ are negligible.
\smallskip

\textnormal{(4)} If $M\in\CM_2^{(1)}$ is an oriented multicurve arranged as in Fig.~\ref{fig_M2} (with  some orientations of components), then the components $\alpha_1$, $\alpha_2$, $\alpha_3$, and~$\alpha_4$ are not active. Moreover, the component~$\alpha_0$ is active if and only if the homology class~$x$ lies in neither of the subgroups~$\langle[\alpha_1],[\alpha_2]\rangle$ and~$\langle[\alpha_3],[\alpha_4]\rangle$ of~$H$.
\smallskip

\textnormal{(5)} All oriented multicurves in~$\CM_3$ are negligible.
\end{propos}

We need some constructions and results obtained by the author in~\cite{Gai21}. Let $\gamma$ be a non-separating simple closed curve on~$\Sigma$. Consider the group
$$
H_{\gamma}=\langle [\gamma]\rangle^{\perp}/\langle [\gamma]\rangle,
$$
where $\perp$ denotes the orthogonal complement in~$H$ with respect to the intersection form. (Note that $H_{\gamma}$ depends only on the homology class~$[\gamma]$ rather than on~$\gamma$ itself. Nevertheless, it is more convenient to write~$H_{\gamma}$ rather than~$H_{[\gamma]}$.) Then $H_{\gamma}$ is a free abelian group of  rank~$4$, equipped with the unimodular skew-symmetric bilinear form induced by the intersection form on~$H$. Now let $\CW=\{W_1,W_2\}$ be an \textit{orthogonal splitting} of~$H_{\gamma}$, i.e., a splitting $H_{\gamma}=W_1\oplus W_2$ such that $W_1$ and~$W_2$ have rank~$2$ and are orthogonal with respect to the intersection form. We denote by $\Spl(\gamma)$ the set of all such splittings~$\CW$. Note that we consider a splitting as an unordered pair~$\{W_1,W_2\}$, that is, not distinguish between the splittings $W_1\oplus W_2$ and~$W_2\oplus W_1$.

If $\gamma'$ is a simple closed curve on~$\Sigma$ such that $\{\gamma,\gamma'\}$ is a bounding pair, then this bounding pair induces a splitting $\CW(\gamma')\in\Spl(\gamma)$. Similarly, if $\delta$ is a separating simple closed curve on~$\Sigma$ disjoint from~$\gamma$, then $\delta$ also induces a splitting $\CW(\delta)\in\Spl(\gamma)$. In~\cite[~\S 6]{Gai21}, the author constructed homomorphisms
$$
\nu_{\gamma,\CW}\colon \I_{\gamma}\to\Z,
$$
indexed by orthogonal splittings~$\CW\in\Spl(\gamma)$, with the following properties:
\begin{enumerate}
 \item if $\delta$ is a separating simple closed curve disjoint from~$\gamma$, then
 \begin{equation}\label{eq_nuWdelta}
 \nu_{\gamma,\CW}(T_{\delta})=
 \left\{
 \begin{aligned}
  &1&&\text{if }\CW=\CW(\delta),\\
  &0&&\text{if }\CW\ne\CW(\delta);
 \end{aligned}
 \right.
 \end{equation}
 \item if $\{\gamma,\gamma'\}$ is a bounding pair, then
 \begin{equation}\label{eq_nuWgamma}
 \nu_{\gamma,\CW}\bigl(T_{\gamma'}T_{\gamma}^{-1}\bigr)=
 \left\{
 \begin{aligned}
  &1&&\text{if }\CW=\CW(\gamma'),\\
  &0&&\text{if }\CW\ne\CW(\gamma');
 \end{aligned}
 \right.
 \end{equation}
 \item if $\{\varepsilon,\varepsilon'\}$ is a bounding pair such that both $\varepsilon$ and $\varepsilon'$ are disjoint from~$\gamma$, and neither of them is isotopic to~$\gamma$, then
 \begin{equation}\label{eq_nuWeps}
  \nu_{\gamma,\CW}\bigl(T_{\varepsilon}T_{\varepsilon'}^{-1}\bigr)=0;
 \end{equation}
 \item for each $h\in\I_{\gamma}$, we have $\nu_{\gamma,\CW}(h)=0$ for all but a finite number of splittings $\CW\in\Spl(\gamma)$.
\end{enumerate}

\begin{remark}
 For properties~(1)--(3), see Proposition~6.4 in~\cite{Gai21}. Property~(4) is not stated explicitly in~\cite{Gai21}. Nevertheless, it follows immediately from the same property for the homomorphisms~$d_{\CW}$ constructed in~\cite{Gai21} just before Proposition~6.4.
\end{remark}

Further, let $M\in\CM_p$ be an oriented multicurve containing no bounding pair, and let $\gamma$ be an active component of~$M$. We say that a splitting $\CW=\{W_1,W_2\}\in\Spl(\gamma)$ is \textit{compatible} with~$M$ if $\pi([M])\subseteq W_1\cup W_2$, where $\pi\colon \langle[\gamma]\rangle^{\perp}\twoheadrightarrow H_{\gamma}$ is the projection. (Recall that $[M]$ is the set of homology classes of the components of~$M$.) We denote by $\Spl(\gamma,M)$ the set of all orthogonal splittings of~$H_{\gamma}$ compatible with~$M$. Since any multicurve~$M$ containing an active component has full rank, we obtain that $M$ contains components~$\alpha_1$ and~$\alpha_2$ such that the homology classes $a_0=[\gamma]$, $a_1=[\alpha_1]$, and~$a_2=[\alpha_2]$ form a basis of a Lagrangian submodule of~$H$.

\begin{propos}\label{propos_Wk}
 Let $M$, $\gamma$, $\alpha_1$, $\alpha_2$, $a_0$, $a_1$ and~$a_2$ be  as  above. Then
 \smallskip

 \textnormal{(1)} The set $\Spl(\gamma,M)$ is nonempty.
 \smallskip

 \textnormal{(2)} If\/ $\CW_0\in\Spl(\gamma,M)$, then $\Spl(\gamma,M)$ is exactly the set of all  splittings $\CW_k=t_{a_1+a_2}^k\CW_0$, where $k\in\Z$. Moreover, the splittings~$\CW_k$ are pairwise distinct.
 \smallskip

 \textnormal{(3)} For any splitting $\CW\in\Spl(\gamma,M)$, there exists a simple closed curve~$\gamma'\subseteq \Sigma\setminus M$ such that $\{\gamma,\gamma'\}$ is a bounding pair and $\CW(\gamma')=\CW$. Moreover, the multicurve~$M\cup\gamma'$ belongs to~$\CM_{p+1}$ and is important.

 \textnormal{(4)} Suppose that $\gamma',\gamma''\subseteq \Sigma\setminus M$ are simple closed curves such that $\{\gamma,\gamma'\}$ and $\{\gamma,\gamma''\}$ are bounding pairs and  $\CW(\gamma')=\CW(\gamma'')$. Then the oriented multicurves~$M\cup\gamma'$ and~$M\cup\gamma''$ lie in the same $\I$-orbit. Moreover, if $p\ge 1$, then the curves~$\gamma'$ and $\gamma''$ lie in the same $\I_M$-orbit.
\end{propos}

\begin{proof}
Since $\gamma$ is an active component of~$M$, there exists a curve~$\gamma'\subseteq \Sigma\setminus M$ such that $\{\gamma,\gamma'\}$ is a bounding pair,  $M\cup\gamma'\in\CM_{p+1}$, and $M\cup\gamma'$ is important. Then $\CW(\gamma')\in\Spl(\gamma,M)$. Hence $\Spl(\gamma,M)\ne\varnothing$.

Since the homology classes~$a_0$, $a_1$, $a_2$ are linearly independent, it follows that $\alpha_1$ and~$\alpha_2$ lie in the distinct connected components~$S_1$ and~$S_2$  of~$\Sigma\setminus(\gamma\cup\gamma')$, respectively. For $i=1,2$, choose a homology class $b_i\in H_1(S_i)\subseteq H$ such that $a_i\cdot b_i=1$. Then $H_1(S_i)=\langle a_0,a_i,b_i\rangle$, and there exists a class~$b_0\in H$ such that $a_0,a_1,a_2,b_0,b_1,b_2$ form a symplectic basis of~$H$. Set $\bar{a}_i=\pi(a_i)$ and $\bar{b}_i=\pi(b_i)$ for  $i=1,2$. Then $\langle a_0\rangle^{\perp}=\langle a_0,a_1,b_1,a_2,b_2\rangle$ and $\bar{a}_1,\bar{a}_2,\bar{b}_1,\bar{b}_2$ form a symplectic basis of~$H_{\gamma}$. An orthogonal splitting~$\CW$ of~$H_{\gamma}$ belongs to~$\Spl(\gamma,M)$ if and only if the elements~$\bar{a}_1$ and~$\bar{a}_2$ lie in the two summands of~$\CW$. Indeed, this follows immediately from two facts: first, $M$ contains the components~$\alpha_1$ and~$\alpha_2$ with $[\alpha_i]=a_i$; second, any component of~$M$ other than $\gamma$, $\alpha_1$, and~$\alpha_2$ (if such a component exists) lies in one of the two subsurfaces~$S_i$, and therefore has one of the homology classes $\pm a_0\pm a_i$, where $i=1,2$. Consequently, $\Spl(\gamma,M)$ is precisely the set consisting of all splittings
\begin{equation*}
\CW_k\colon\quad H_{\gamma}=\bigl\langle \bar{a}_1,\bar{b}_1-k\bar{a}_2\bigr\rangle \oplus\bigl\langle \bar{a}_2,\bar{b}_2-k\bar{a}_1\bigr\rangle,\qquad k\in\Z.
\end{equation*}
The symplectic transvection $t_{a_1+a_2}$ acts on~$H_{\gamma}$ by fixing~$\bar{a}_1$ and~$\bar{a}_2$, while sending~$\bar{b}_1$ and~$\bar{b}_2$ to~$\bar{b}_1-\bar{a}_1-\bar{a}_2$ and  $\bar{b}_2-\bar{a}_1-\bar{a}_2$, respectively. Therefore, $t_{a_1+a_2}^k\CW_0=\CW_k$. We obtain assertion~(2) for the specific choice~$\CW_0=\CW(\gamma')$. Assertion~(2) for an arbitrary $\CW_0$ follows immediately, since another choice of~$\CW_0$ results only in a shift of the  numbering.

From Proposition~\ref{propos_active_obvious} it follows that $M$ is either a three-component multicurve that does not separate~$\Sigma$, or a four-component multicurve belonging to~$\CM_1^{(2)}$, or a five-component multicurve belonging to~$\CM_2^{(1)}$. In each of these cases,  there exists a simple closed curve~$\varepsilon$ in~$\Sigma\setminus M$ whose homology class is either $a_1+a_2$ or~$a_1-a_2$. Moreover, for each $i=1,2$, there exists a simple closed curve~$\delta_i$  in~$S_i$, disjoint from~$M$, such that $[\delta_i]$ is either~$a_i+a_0$ or~$a_i-a_0$. (The curve~$\delta_i$ is allowed to be isotopic to a component of~$M$.)

If $[\varepsilon]=a_1+a_2$, then the Dehn twist~$T_{\varepsilon}$ acts on~$H$ via the transvection~$t_{a_1+a_2}$. Therefore $\CW_k=\CW\bigl(T_{\varepsilon}^k(\gamma')\bigr)$. If $[\varepsilon]=a_1-a_2$, then the Dehn twist~$T_{\varepsilon}$ acts on~$H$ via the transvection
$$
t_{a_1-a_2}=t_{a_1+a_2}^{-1}t_{a_1}^2t_{a_2}^2.
$$
But the tranvections~$t_{a_1}$ and~$t_{a_2}$ act trivially on~$\Spl(\gamma,M)$, since they are induced by the Dehn twists~$T_{\alpha_1}$ and~$T_{\alpha_2}$ about the components of~$M$.
Hence $\CW_k=\CW\bigl(T_{\varepsilon}^{-k}(\gamma')\bigr)$.

Finally, let us prove assertion~(4). It suffices to prove it for~$\CW(\gamma')=\CW(\gamma'')=\CW_0$, since the general case reduces to this one by the action of the group~$\langle T_{\varepsilon}\rangle$. First, suppose that $p\ge 1$. Let $S$ be the connected component of~$\Sigma\setminus(\gamma\cup\gamma')$ that is not homeomorphic to a pair of pants. Then $S$ is a $5$-punctured sphere if $p=1$ and a $4$-punctured sphere if $p=2$. In either of the cases, there is exactly one $\PMod(S)$-orbits of simple closed curves on~$S$ that are homologous but not isotopic to~$\gamma$ in~$\Sigma$. Therefore there exists a mapping class $f\in \Mod(\Sigma)$ such that $f$ stabilizes every component of~$M$ (and preserves the orientation of each component), and $f(\gamma')=\gamma''$. Then $f$ fixes each of the homology classes~$a_0$, $a_1$, and~$a_2$. Moreover, since $\CW(\gamma')=\CW(\gamma'')$, we obtain $f(b_i)=b_i+l_ia_i+m_ia_0$ for some $l_i,m_i\in\Z$. Now for $i=1,2$ there exist integers $p_i,q_i$ such that the mapping class $g_i=T_{\alpha_i}^{p_i}T_{\delta_i}^{q_i}$ sends~$b_i$ to $b_i-l_ia_i-m_ia_0$. Indeed, if $[\delta_i]=a_i+a_0$, then $g_i=T_{\alpha_i}^{l_i-m_i}T_{\delta_i}^{m_i}$, and if $[\delta_i]=a_i-a_0$, then $g_i=T_{\alpha_i}^{l_i+m_i}T_{\delta_i}^{-m_i}$. Set $h=fg_1g_2$. Then $h$ also stabilizes every component of~$M$ (and preserves the orientation of each component), and satisfies $h(\gamma')=\gamma''$. Moreover, $h$ fixes each of the homology classes $a_0$, $a_1$, $a_2$, $b_1$, and~$b_2$. Therefore $h(b_0)=b_0+na_0$ for some $n\in\Z$. Hence the mapping class~$hT_{\gamma}^n$ belongs to~$\I_M$ and takes~$\gamma'$ to~$\gamma''$, as required.

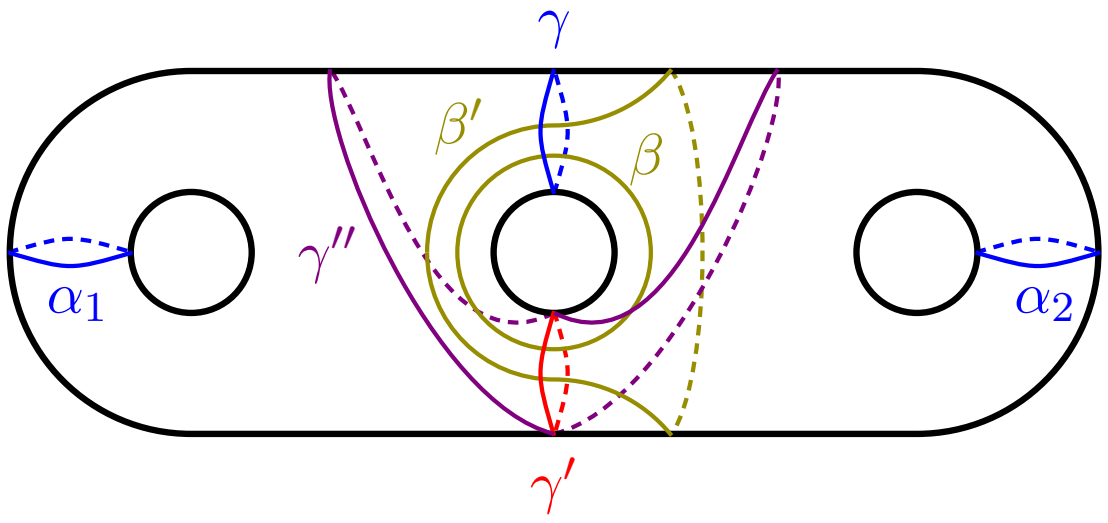
\begin{figure}[!ht]
 \begin{tikzpicture}[scale=.4]

\small
% \definecolor{blue}{rgb}{0, 0, 0.7}
 \definecolor{mygreen}{rgb}{0, 0.4, 0}
\tikzset{every path/.append style={line width=.2mm}}
\tikzset{my dash/.style={dash pattern=on 2pt off 1.5pt}}

\draw[color=olive, thick,  my dash] (1.93,3.03) .. controls (2.63,2) and (2.63,-2) .. (1.93,-3.03);
\draw[color=violet,thick,my dash] (0,-1) .. controls (-2,-2) and (-3,2) ..   (-3.7,3.03);
\draw[color=violet,thick, my dash] (0,-3) .. controls (2,-2.5) and (3.9,2) ..   (3.7,3.03);

% Border lines

\draw[color=blue, thick, my dash] (-9.03,0) .. controls (-8,.3) .. (-6.97,0);
\draw[color=red, thick, my dash] (0,-.97) .. controls (.3,-2) .. (0,-3.03);
\draw[color=blue, thick, my dash] (0,.97) .. controls (.3,2) .. (0,3.03);
\draw[color=blue, thick, my dash] (6.97,0) .. controls (8,.3) .. (9.03,0);
%\draw[color=blue] (1.93,3.03) -- (1.93,-3.03);

\draw [very thick] (-6,0) circle (1);
\draw [very thick] (0,0) circle (1);
\draw [very thick] (6,0) circle (1);
\draw [very thick]  (-6,3)--(6,3) arc (90:-90:3) -- (-6,-3) arc (270:90:3);

% Visible lines

\draw [olive,thick] (0,0) circle (1.6) node [above right=3mm,xshift=.5mm] {$\beta$};
\draw [olive,thick] (-2.1,0) arc (180:90:2.1) node [pos=.6,above left=-5pt] {$\beta'$} arc (-90:-39:2.5) ;
\draw [olive,thick] (-2.1,0) arc (180:270:2.1) arc (90:39:2.5);

\draw[color=blue, thick] (-9.03,0) .. controls (-8,-.3) .. (-6.97,0) node [pos=.5,below,xshift=.5mm] {$\alpha_1$};
\draw[color=red, thick] (0,-.97) .. controls (-.3,-2) .. (0,-3.03) node [below] {$\gamma'$};
\draw[color=blue, thick] (0,.97) .. controls (-.3,2) .. (0,3.03) node [above] {$\gamma$};
\draw[color=blue, thick] (6.97,0) .. controls (8,-.3) .. (9.03,0) node [pos=.5,below,xshift=.5mm] {$\alpha_2$};

\draw[color=violet,thick] (0,-1) .. controls (2,-2) and (3,2) ..   (3.7,3.03);
\draw[color=violet,thick] (0,-3) .. controls (-2,-2.5) and (-3.9,2) ..   (-3.7,3.03) node [pos=.6,below left=-2pt] {$\gamma''$};

\end{tikzpicture}
 \caption{Curves~$\gamma'$ and~$\gamma''$ providing the same splitting in homology}\label{fig_gamma'gamma''}
\end{figure}

Second, suppose that $p=0$. Then $S=\Sigma\setminus M$ is a $6$-punctured sphere. If $\gamma'$ and~$\gamma''$ lie in the same $\PMod(S)$-orbit, then we can repeat literally the above argument and obtain that $\gamma'$ and~$\gamma''$ belong to the same $\I_M$-orbit. However, now there are two $\PMod(S)$-orbits of simple closed curves on~$S$ that are homologous but not isotopic to~$\gamma$ in~$\Sigma$; representatives of these two orbits are the curves~$\gamma'$ and~$\gamma''$ depicted in Fig.~\ref{fig_gamma'gamma''}. Consequently, it suffices to show that the multicurves~$M\cup\gamma'$ and~$M\cup\gamma''$ lie in the same $\I$-orbit for this specific choice of~$\gamma'$ and~$\gamma''$. Consider the auxiliary simple closed curves~$\beta$ and~$\beta'$  shown in Fig.~\ref{fig_gamma'gamma''}. Then the BP-map $T_{\beta}T_{\beta'}^{-1}$ sends $\gamma'$ to~$\gamma$ and~$\gamma$ to~$\gamma''$. Thus it sends~$M\cup\gamma'$ to~$M\cup\gamma''$. The proposition follows, since $T_{\beta}T_{\beta'}^{-1}\in\I$.
\end{proof}

\begin{cor}\label{cor_ActM1M2}
 Let $M_1\subseteq M_2$ be two oriented multicurves in~$\CM$ containing no bounding pairs, and let $\gamma$ be a component of~$M_1$ such that $\gamma\in \Act(M_2)$. Then $\gamma\in\Act(M_1)$ and $\Spl(\gamma,M_1)=\Spl(\gamma,M_2)$.
\end{cor}

Another consequence of Proposition~\ref{propos_Wk}  is as follows. Let $(M,\gamma,\CW)$ be a triple such that $M$ is an unimportant multicurve from~$\CM_p$, $\gamma\in\Act(M)$, and $\CW\in\Spl(\gamma,M)$. By assertion~(3) of Proposition~\ref{propos_Wk}, we have $\CW=\CW(\gamma')$ for some important multicurve $K=M\cup\gamma'$. Moreover, by assertion~(4), all such multicurves~$K$ lie in the same $\I$-orbit. We set
$$
\CP_{M,\gamma,\CW}=\CP_K.
$$
Then this cell is well defined, that is, independent of the choice of~$\gamma'$. Since $M\cup\gamma'$ is important, we obtain that $\CP_{M,\gamma,\CW}$ becomes a well-defined element of~$\CC_{p+1}$ once we choose an orientation of this cell.

\section{Stabilizers of multicurves in~$\CM_1^{(2)}$ and~$\CM_2^{(1)}$}\label{section_stabilizers}
In this section, we study the stabilizers~$\I_M$, where $M$ belongs to either~$\CM_2^{(1)}$ or~$\CM_1^{(2)}$ and has type~2. The group~$H_1(\I_M)$ in these two cases was computed explicitly by the author in Subsections~7.2 and~7.3, respectively, of~\cite{Gai21}. We now recall these results in a convenient for us form.

For a set~$X$, we denote by~$\Z^X$ the free abelian group with basis~$X$.

The following proposition is the combination of Propositions~7.11 and~7.14 in~\cite{Gai21}.

\begin{propos}\label{propos_M2_ab}
Let $M\in \CM_2^{(1)}$, and let $\gamma$ be the active component of~$M$. Then the homomorphisms $\nu_{\gamma,\CW}\colon \I_M\to\Z$ indexed by all splittings $\CW\in\Spl(\gamma,M)$ are linearly independent, and give an isomorphism
$$
H_1(\I_M)\stackrel{\cong}{\longrightarrow}\Z^{\Spl(\gamma,M)}.
$$
In other words, $H_1(\I_M)$ is a free abelian group admitting a basis $\{e_{\CW}\}$, where $\CW\in\Spl(\gamma,M)$, such that
$$
\nu_{\gamma,\CW'}\bigl(e_{\CW}\bigr)=\left\{
\begin{aligned}
 &1&&\text{if\/ }\CW=\CW',\\
 &0,&&\text{if\/ }\CW\ne\CW'.
\end{aligned}
\right.
$$
Moreover, if $\gamma'\subseteq \Sigma\setminus M$ is a multicurve such that $\{\gamma,\gamma'\}$ is a bounding pair, then $$\bigl[T_{\gamma'}T_{\gamma}^{-1}\bigr]=e_{\CW(\gamma')}.$$
\end{propos}

Let us now give an analogue of this proposition for a type~2 multicurve from~$\CM_1$.

\begin{propos}\label{propos_type2_ab}
Let $M\in \CM_1^{(2)}$, and let $\gamma_0$, $\gamma_1$, and~$\gamma_2$ be the active components of~$M$, i.e., the three components of~$M$ that bound a pair of pants. Then the homomorphisms $\nu_{\gamma_i,\CW}\colon \I_M\to\Z$ indexed by all pairs $(\gamma_i,\CW)$ with $i=0,1,2$ and $\CW\in\Spl(\gamma_i,M)$ are linearly independent, and give an isomorphism
$$
H_1(\I_M)\stackrel{\cong}{\longrightarrow}\Z^{\Spl(\gamma_0,M)}\oplus\Z^{\Spl(\gamma_1,M)}\oplus \Z^{\Spl(\gamma_2,M)}.
$$
In other words, $H_1(\I_M)$ is a free abelian group admitting a basis $\{e_{\gamma_i,\CW}\}$, where $i=0,1,2$ and $\CW\in\Spl(\gamma_i,M)$, such that
$$
\nu_{\gamma_j,\CW'}\bigl(e_{\gamma_i,\CW}\bigr)=\left\{
\begin{aligned}
 &1&&\text{if }i=j\text{ and }\CW=\CW',\\
 &0,&&\text{otherwise.}
\end{aligned}
\right.
$$
Moreover, if $\gamma'\subseteq \Sigma\setminus M$ is a multicurve such that $\{\gamma_i,\gamma'\}$ is a bounding pair, then $$\bigl[T_{\gamma'}T_{\gamma_i}^{-1}\bigr]=e_{\gamma_i,\CW(\gamma')}.$$
\end{propos}

\begin{proof}
 Let $\gamma_3$ be the non-active component of~$M$. We set $c_i=[\gamma_i]$ for $i=0,\ldots,3$. For each $i=0,1,2$, choose a simple closed curve $\varepsilon_i$ disjoint from~$M$ such that $[\varepsilon_i]=c_i+c_3$. From Propositions~\ref{propos_M_char} and~\ref{propos_dimPM} it follows that the multicurve $M_i=M\cup\varepsilon_i$ belongs to~$\CM_2^{(1)}$. By Corollary~7.16 and Proposition~7.22 in~\cite{Gai21}, the inclusions $\I_{M_i}\le \I_M$ induce injections $H_1(\I_{M_i})\hookrightarrow H_1(\I_{M})$, and
 \begin{equation}\label{eq_M1M2M3}
 H_1(\I_{M})=H_1(\I_{M_0})\oplus H_1(\I_{M_1})\oplus H_1(\I_{M_2}).
 \end{equation}
 Moreover, if $i\ne j$, then all homomorphisms $\nu_{\gamma_j,\CW}$ vanish on~$\I_{M_i}$. This follows immediately from~\eqref{eq_nuWeps}, since by Proposition~\ref{propos_M2_ab} the abelianization~$H_1(\I_M)$ is generated by elements $\bigl[T_{\gamma'}T_{\gamma_i}^{-1}\bigr]$ such that  $\{\gamma_i,\gamma'\}$ is a bounding pair (see also the proof of Proposition~7.22 in~\cite{Gai21}). The proposition now follows from Proposition~\ref{propos_M2_ab}.
 \end{proof}

 \begin{propos}\label{propos_e_annih}
  Let $M=\gamma_0\cup\gamma_1\cup\gamma_2\cup\gamma_3$ be as in Proposition~\ref{propos_type2_ab}, with active components~$\gamma_0$, $\gamma_1$, $\gamma_2$, and non-active  component~$\gamma_3$, let $c_i=[\gamma_i]$, and let $\{e_{\gamma_i,\CW}\}$ be the basis from Proposition~\ref{propos_type2_ab}.  Then
  $$
  (t_{c_i+\sigma c_3}-1)e_{\gamma_i,\CW}=0
  $$
  for any $i\in\{0,1,2\}$, $\CW\in\Spl(\gamma_i,M)$, and~$\sigma\in\{-1,1\}$.
 \end{propos}

  \begin{proof}
   Let $\varepsilon_0$, $\varepsilon_1$, $\varepsilon_2$, $M_0$, $M_1$, $M_2$ be as in the proof of Proposition~\ref{propos_type2_ab}. Then $e_{\gamma_i,\CW}$ lies in the summand $H_1(\I_{M_i})$ of the decomposition~\eqref{eq_M1M2M3}. Since $T_{\varepsilon_i}$ acts trivially on~$H_1(\I_{M_i})$, and $[\varepsilon_i]=c_i+c_3$, we obtain that $t_{c_i+c_3}$ fixes~$ e_{\gamma_i,\CW}$. Further, since $t_{c_i-c_3}=t_{c_i+c_3}^{-1}t_{c_i}^2t_{c_3}^2$ and the twists~$T_{\gamma_i}$ and~$T_{\gamma_3}$ act trivially on~$H_1(\I_M)$, we obtain that $t_{c_i-c_3}$ fixes~$ e_{\gamma_i,\CW}$ as well.
  \end{proof}

\begin{figure}[!ht]
 \begin{tikzpicture}[scale=.3,
        mid arrow/.style={
            postaction={
                decorate,
                decoration={
                    markings,
                    mark=at position 0.55 with {
                        \arrow[#1]{stealth}
                    }
                }
            }
        },
        mid arrow/.default=black,
    ]

\small
% \definecolor{blue}{rgb}{0, 0, 0.7}
% \definecolor{mygreen}{rgb}{0, 0.4, 0}
\tikzset{every path/.append style={line width=.2mm}}
\tikzset{my dash/.style={dash pattern=on 2pt off 1.5pt}}

% Border lines

\draw[color=blue, thick, my dash] (-7.03,0) .. controls (-6,.3) .. (-4.97,0);
\draw[color=blue, thick, my dash] (-3.03,0) .. controls (-2,.3) .. (-.97,0);
\draw[color=blue, thick, my dash] (0,.97) .. controls (.3,2) .. (0,3.03);
\draw[color=blue, thick, my dash] (6.97,0) .. controls (8,.3) .. (9.03,0);
\draw[color=violet, thick, my dash] (0,-.97) .. controls (.3,-2) .. (0,-3.03);
\draw[color=red, thick, my dash] (3,3.03) .. controls (3.7,2) and (3.7,-2) .. (3,-3.03);

\draw [very thick] (-4,0) circle (1);
\draw [very thick] (0,0) circle (1);
\draw [very thick] (6,0) circle (1);
\draw [very thick]  (-4,3)--(6,3) arc (90:-90:3) -- (-4,-3) arc (270:90:3);

% Visible lines

\draw[color=blue, thick] (-7.03,0) .. controls (-6,-.3) .. (-4.97,0) node[below=-1pt,pos=.5,xshift=.5mm] {$\gamma_j$};
\draw[color=blue, thick] (-3.03,0) .. controls (-2,-.3) .. (-.97,0) node[below=-1pt,pos=.5,xshift=.5mm] {$\gamma_k$};
\draw[color=blue, thick] (0,.97) .. controls (-.3,2) .. (0,3.03) node[left=-2pt,pos=.5] {$\gamma_i$};
\draw[color=blue, thick] (6.97,0) .. controls (8,-.3) .. (9.03,0) node[below=-1pt,pos=.5] {$\gamma_3$};
\draw[color=violet, thick] (0,-.97) .. controls (-.3,-2) .. (0,-3.03) node[left=-3pt,pos=.5] {$\gamma'$};

\draw[color=red, thick] (3,3.03) .. controls (2.3,2) and (2.3,-2) .. (3,-3.03) node [pos=0.5,left=-2pt] {$\delta$};

%%%%%%%
%%%%%%%

\end{tikzpicture}
 \caption{The configuration of curves in Proposition~\ref{propos_BP_delta}}\label{fig_BP_delta}
\end{figure}
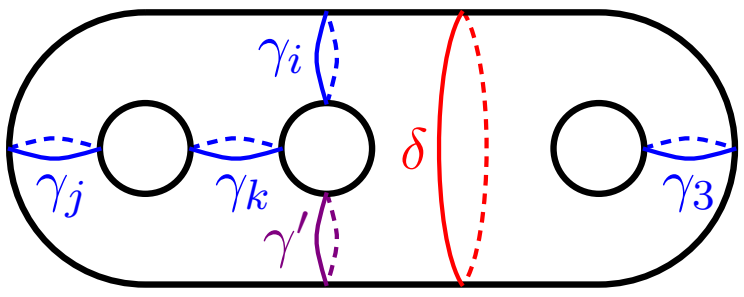

\begin{propos}\label{propos_BP_delta}
 Let $M=\gamma_0\cup\gamma_1\cup\gamma_2\cup\gamma_3$ be as in the previous two propositions. Let $i,j,k$ be a permutation of~$0,1,2$, let~$\gamma'$ be a simple closed curve disjoint from~$M$ and such that $\{\gamma_i,\gamma'\}$ is a bounding pair, and let $\delta$ be a separating simple closed curve disjoint from~$M$ and~$\gamma'$; see Fig.~\ref{fig_BP_delta}. Then, for any $\sigma_1,\sigma_2\in\{-1,1\}$, we have the following equality in~$H_1(\I_M)$:
 \begin{equation}\label{eq_tt_anih}
 (t_{c_j+\sigma_1c_3}-1)(t_{c_k+\sigma_2c_3}-1)\bigl[T_{\gamma'}T_{\gamma_i}^{-1}\bigr]=(t_{c_j+\sigma_1c_3}-1)(t_{c_k+\sigma_2c_3}-1)[T_{\delta}].
 \end{equation}
\end{propos}

\begin{proof}
  For each~$l=0,1,2$, let $\CW_l(\delta)$ be the splitting in~$\Spl(\gamma_l,M)$ induced by~$\delta$. By~\eqref{eq_nuWdelta}, we have
 $$
 \nu_{\gamma_l,\CW_l(\delta)}(T_{\delta})=1,\qquad l=0,1,2,
 $$
 and $\nu_{\gamma_l,\CW}(T_{\delta})=0$ for $\CW\ne\CW_l(\delta)$. Therefore
 \begin{equation}\label{eq_delta_3e}
 [T_{\delta}]=e_{\gamma_0,\CW_0(\delta)}+e_{\gamma_1,\CW_1(\delta)}+e_{\gamma_2,\CW_2(\delta)},
 \end{equation}
 where $\bigl\{e_{\gamma_l,\CW}\bigr\}$ is the basis from Proposition~\ref{propos_type2_ab}. On the other hand, we have $$\CW(\gamma')=\CW_i(\delta)\in\Spl(\gamma_i,M),$$ and therefore
 $$
 \bigl[T_{\gamma'}T_{\gamma_i}^{-1}\bigr]=e_{\gamma_i,\CW_i(\delta)}.
 $$
 Hence, multiplying~\eqref{eq_delta_3e} by $(t_{c_j+\sigma_1 c_3}-1)(t_{c_k+\sigma_2 c_3}-1)$ and using Proposition~\ref{propos_e_annih} and the fact that~$t_{c_j+\sigma_1 c_3}$ and~$t_{c_k+\sigma_2 c_3}$ commute, we obtain the required equality.
\end{proof}

Further we will need the following consequence of Propositions~\ref{propos_M2_ab} and~\ref{propos_type2_ab}.

\begin{propos}\label{propos_nu_trivial}
 Let $M,N\in\CM$ be oriented multicurves containing no bounding pairs and such that $N\subseteq M$.  Let $\gamma$ be a component of~$N$ such that $\gamma\in\Act(N)$ but $\gamma\notin \Act(M)$, and let $\CW\in\Spl(\gamma)$. Then the restriction of~$\nu_{\gamma,\CW}$ to the stabilizer~$\I_M$ is trivial.
\end{propos}

\begin{proof}
 Suppose that $M\in\CM_p$; then $p\in\{1,2,3\}$. If $p=3$, then $M$ is a pants decomposition without homologous components. Therefore $\I_M=1$. Assume that $p$ is either~$1$ or~$2$. Since~$N$ contains an active component, it has full rank. Hence $M$ also has full rank. Consequently, $M$ belongs to either~$\CM_1^{(1)}$, or~$\CM_1^{(2)}$, or~$\CM_2^{(1)}$.
 Then the group~$H_1(\I_M)$ is generated by BP-maps~$\bigl[T_{\varepsilon}T_{\varepsilon'}^{-1}\bigr]$ such that $\varepsilon$ and~$\varepsilon'$ are not homologous to~$\gamma$. Indeed, if $M$ belongs to either~$\CM_2^{(1)}$ or~$\CM_1^{(2)}$, this follows from Propositions~\ref{propos_M2_ab} and~\ref{propos_type2_ab}, respectively, while for $M\in\CM_1^{(1)}$ this holds by Proposition~7.1 in~\cite{Gai21}.
Therefore, by~\eqref{eq_nuWeps}, the restriction of~$\nu_{\gamma,\CW}$ to~$\I_M$ is trivial.
\end{proof}

\section{Chain complex~$\CE_*$ and chain map~$\varphi$}\label{section_ephi}

Consider the chain complex~$E^1_{*,1}$ with differential~$d^1$. We have
\begin{equation}\label{eq_E1p1}
E^1_{p,1}=\bigoplus_{M\in\mathfrak{M}_p}H_1(\I_M),
\end{equation}
where $\mathfrak{M}_p$ is a set of representatives of $\I$-orbits in~$\CM_p$. Let $\CE_p\subseteq E^1_{p,1}$ be the sum of all summands in~\eqref{eq_E1p1} that correspond to multicurves~$M$ containing no bounding pairs. We have $d^1\CE_{p}\subseteq\CE_{p-1}$. Hence $\CE_*$ is a chain subcomplex of~$E^1_{*,1}$.

Note that by Fact~\ref{fact_E1} the summands in decomposition~\eqref{eq_E1p1}, and therefore the subgroup~$\CE_p$, do not depend on the choice of the set of representatives~$\mathfrak{M}_p$. For each cell~$P_M$  (with~$M$ not necessarily belonging to~$\mathfrak{M}_p$), we denote by $\iota_{P_M}$ the embedding of $H_1(\I_M)$ onto the corresponding summand of~$E^1_{p,1}$; these embeddings are related by commutative diagrams~\eqref{eq_cd_iota} from  Fact~\ref{fact_E1}. Moreover, the embedding~$\iota_{P_M}$ reverses sign whenever we reverse the orientation of~$P_M$.

We would like to construct a $\Z$-module homomorphism $\varphi\colon\CE_p\to\CC_{p+1}$ by setting
\begin{equation}\label{eq_varphi}
 \varphi\bigl(\iota_{P_M}([h])\bigr)=\sum_{\gamma\in\Act(M)}\,\sum_{\CW\in\Spl(\gamma,M)}\nu_{\gamma,\CW}(h)\CP_{M,\gamma,\CW}
\end{equation}
for all $h\in\I_M$. Here $\CP_{M,\gamma,\CW}$ are the cells defined in the end of Section~\ref{section_active}.
The choice of the orientations of the cells~$\CP_{M,\gamma,\CW}$ in this formula is the most delicate point in the definition of the homomorphism~$\varphi$. This choice is not quite canonical. Nevertheless, it can be arranged to depend only on certain data that can be fixed once and for all. Specifically, the data are as follows. Consider all Lagrangian submodules $L\subseteq H$ such that $x\in L$. For each such $L$, choose an orientation of~$L$, i.e., a generator $\omega_L\in\wedge^3L$. Note that the orientations~$\omega_L$ for different~$L$ are not required to agree with each other in any way. We fix this choice of orientations once and for all. Now, for any basis of any Lagrangian submodule $L$ containing~$x$, we can determine whether this basis is \textit{positive} or \textit{negative}.

Now let $\{\gamma,\gamma'\}$ be a bounding pair on~$\Sigma$ such that $[\gamma]\cdot x=0$. Assume that $x$ does not belong to either of the subgroups $H_1(S_i)\subseteq H$, where $S_1$ and~$S_2$ are the connected components of~$\Sigma\setminus(\gamma\cup\gamma')$. A priori, the pair $\{\gamma,\gamma'\}$ is unordered, so the curves $\gamma$ and~$\gamma'$ play symmetric roles. We assert, however, that there is a canonical way to decide which of the two curves in the pair is \textit{first} and which is \textit{second}. To do this, we proceed as follows. Choose  orientations of~$\gamma$ and~$\gamma'$ such that $[\gamma]=[\gamma']$, and set $c=[\gamma]=[\gamma']$. Write $x=x_1+x_2$, where $x_i\in H_1(S_i)$. The summands $x_1$ and~$x_2$ are defined up to adding $kc$ to one of them and simultaneously subtracting~$kc$ from the other. Let $L$ be the Lagrangian submodule of~$H$ containing~$c$, $x_1$, and~$x_2$. Then the element $c\wedge x_1\wedge x_2\in\wedge^3L$ is well-defined and nonzero. So far we have been free to choose which of the two components of~$\Sigma\setminus(\gamma\cup\gamma')$ is~$S_1$ and which is~$S_2$. We now fix this choice by requiring that $c\wedge x_1\wedge x_2$ be a positive multiple of~$\omega_L$. Finally,  of the curves~$\gamma$ and~$\gamma'$, we call \textit{first} the one such that $S_1$ lies to its right, and \textit{second} the other one. Note that the ordering of~$\gamma$ and~$\gamma'$ does not depend on the choice of their orientations. Indeed, if we  reverse the orientations of both~$\gamma$ and~$\gamma'$ simultaneously, then the sign of~$c$ reverses, and consequently~$S_1$ and~$S_2$ are interchanged. However, the side---right or left---on which a connected component lies relative to~$\gamma$ also reverses, so we obtain the same first and second curves.

Let $P_M$ be an oriented $p$-cell such that $M$  contains no bounding pair, and let $h\in\I_M$. We are now ready to give a rule for choosing the orientations of the cells~$\CP_{M,\gamma,\CW}$ on the right-hand side of~\eqref{eq_varphi}. Recall that, by the definition in Section~\ref{section_active}, we have $\CP_{M,\gamma,\CW}=\CP_{K}$, where $K=M\cup\gamma'$ for some simple closed curve $\gamma'\subseteq \Sigma\setminus M$ such that $\{\gamma,\gamma'\}$ is a bounding pair and $\CW(\gamma')=\CW$. We orient~$P_K$ so that $[P_K:P_M]=1$ if $\gamma$ is the first component of the bounding pair $\{\gamma,\gamma'\}$ and $[P_K:P_M]=-1$ if $\gamma$ is the second component of this bounding pair. Then the orientation of~$\CP_{M,\gamma,\CW}=\CP_{K}$ is induced by this orientation of~$P_K$.

\begin{propos}\label{propos_Pwelldef}
The orientation of\/~$\CP_{M,\gamma,\CW}$ is well-defined, that is, independent of the choice of~$K$. Moreover, this orientation reverses whenever we reverse the orientation of~$P_M$. Thus~\eqref{eq_varphi} provides a well-defined $\Z$-module homomorphism~$\varphi$.
\end{propos}

\begin{proof}
First, assume that $p\ge 1$. Then by part~(4) of Proposition~\ref{propos_Wk}, the multicurve~$K=M\cup\gamma'$ is uniquely determined up to the action of~$\I_M$. If $f(\gamma')=\gamma''$ for an element $f\in\I_M$, then $\gamma$ is the first component of~$\{\gamma,\gamma'\}$ if and only if $\gamma$ is the first component of~$\{\gamma,\gamma''\}$. Hence our rule gives orientations of~$P_{M\cup\gamma'}$ and~$P_{M\cup\gamma''}$ that induce the same orientations on~$\CP_{M\cup\gamma'}=\CP_{M\cup\gamma''}=\CP_{M,\gamma,\CW}$.

Second, assume that $p=0$. Let $K=M\cup\gamma'$, and let $\varepsilon$ and~$\varepsilon'$ be the first and second components of the bounding pair $\{\gamma,\gamma'\}$, respectively.  Then $P_K$ is the segment with vertices $P_N$ and~$P_{N'}$, where $N=K\setminus\varepsilon'$ and $N'=K\setminus \varepsilon$. According to our rule, $P_K$ is oriented from $P_{N'}$ to~$P_N$, independently of which of the two curves~$\varepsilon$ and~$\varepsilon'$ is~$\gamma$ and which is~$\gamma'$. By part~(4) of Proposition~\ref{propos_Wk}, the cell~$P_K$ is uniquely determined up to the action of an element $f\in\I$. However, such an element always takes the first and second components of the bounding pair contained in~$K$ to the first and second components of the bounding pair contained in~$f(K)$, respectively. Hence the resulting orientation of $\CP_K=\CP_{M,\gamma,\CW}$ is well-defined.

The fact that the orientation of~$\CP_{M,\gamma,\CW}$ reverses whenever we reverse the orientation of~$P_M$ follows directly from the construction.
\end{proof}

\begin{remark}\label{remark_not_equiv}
Since the orientations~$\omega_L$ for different~$L$ do not agree with each other in any way, $\varphi$ is \textit{not} a homomorphism of $\Z\bigl[\Sp(H,x)\bigr]$-modules. Nevertheless, since the transvection~$t_x$ preserves the homology class of every component of any multicurve in~$\CM$, it is still true that $\varphi$ is a homomorphism of~$\Z\bigl[t_x,t_x^{-1}\bigr]$-modules.
\end{remark}

 \begin{remark}
  One can wonder whether it is possible not to deal with the orientations~$\omega_L$ and simply always choose the orientation of~$P_K$ so that $[P_K:P_M]=1$. The answer is as follows. If $p>0$, then such simple rule will indeed yield a well-defined homomorphism~$\varphi$ but this will not be the homomorphism we need further. However, for $p=0$ the homomorphism will not be well-defined at all.
 \end{remark}

Let us prove several properties of~$\varphi$.

\begin{propos}\label{propos_phi_chain}
 The homomorphism $\varphi\colon\CE_*\to\CC_{*+1}$ is a degree $1$ chain map, that is,
 \begin{equation}\label{eq_phid=dphi}
 \varphi d^1=-\partial_{\CC}\varphi.
 \end{equation}
\end{propos}

\begin{remark}
 The sign minus in~\eqref{eq_phid=dphi} agrees with the usual sign convention, since $\varphi$ has degree~$1$ and the differentials~$d^1$ and~$\partial_{\CC}$ have degree~$-1$.
\end{remark}

\begin{proof}[Proof of Proposition~\ref{propos_phi_chain}]
 Let $M$ be an oriented multicurve in~$\CM_p$ containing no bounding pair, and let
$$
\partial P_M=\sum_{i=1}^k[P_M:P_{N_i}]P_{N_i}.
$$
For each $\gamma\in\Act(M)$ and each $\CW\in\Spl(\gamma,M)$, the important subcells of~$\CP_{M,\gamma,\CW}$ are precisely the cells~$\CP_{N_i,\gamma,\CW}$ corresponding to those~$N_i$ that contain~$\gamma$. Note that $\Spl(\gamma,N_i)=\Spl(\gamma,M)$ for any such~$N_i$ by Corollary~\ref{cor_ActM1M2}. Moreover, $\CP_{M,\gamma,\CW}=\CP_{M\cup\gamma'}$ and $\CP_{N_i,\gamma,\CW}=\CP_{N_i\cup\gamma'}$ where $\gamma'\subseteq \Sigma\setminus M$ is a simple closed curve such that $\{\gamma,\gamma'\}$ is a bounding pair and $\CW(\gamma')=\CW$. Whether~$\gamma$ is the first or second component of the pair $\{\gamma,\gamma'\}$ does not depend on the other components of the multicurve. Therefore, according to our rule for orientations, $[P_{M\cup\gamma'}:P_M]=[P_{N_i\cup\gamma'}:P_{N_i}]$. Consequently, $[P_{M\cup\gamma'}:P_{N_i\cup\gamma'}]=-[P_M:P_{N_i}]$, and hence
$$
[\CP_{M,\gamma,\CW}:\CP_{N_i,\gamma,\CW}]=-[P_M:P_{N_i}].
$$
Therefore, for each $h\in\I_M$, we have:
\begin{equation}\label{eq_dphi}
 -\partial_{\CC}\varphi\bigl(\iota_{P_M}([h])\bigr)=\sum_{\gamma\in\Act(M)}\,\sum_{\CW\in\Spl(\gamma,M)}\,\sum_{i:N_i\supseteq\gamma}[P_M:P_{N_i}]\nu_{\gamma,\CW}(h)\CP_{N_i,\gamma,\CW},
\end{equation}
On the other hand,
\begin{equation}\label{eq_phid}
 \varphi d^1\bigl(\iota_{P_M}([h])\bigr)=\sum_{i=1}^k\sum_{\gamma\in\Act(N_i)}\,\sum_{\CW\in\Spl(\gamma,M)}[P_M:P_{N_i}]\nu_{\gamma,\CW}(h)\CP_{N_i,\gamma,\CW}.
\end{equation}
By Corollary~\ref{cor_ActM1M2}, any component~$\gamma$ of~$N_i$ that is active in~$M$ is also active in~$N_i$. Hence all summands on the right-hand side of~\eqref{eq_dphi} are also present on the right-hand side of~\eqref{eq_phid}. However, the sum on the right-hand side of~\eqref{eq_phid} may contain additional summands that correspond to those components $\gamma\in\Act(N_i)$ that are not active in~$M$. Nevertheless, by Proposition~\ref{propos_nu_trivial}, all such summands vanish. Thus $\varphi d^1=-\partial_{\CC}\varphi$.
\end{proof}

\begin{propos}\label{propos_phi_K}
 Let $K$ be an important oriented multicurve in~$\CM_{p+1}$ with bounding pair $\{\gamma_0,\gamma_0'\}$, and let $M=K\setminus\gamma'_0$. Then
 $$
\varphi\left(\iota_{P_M}\bigl(\bigl[T_{\gamma_0'}T_{\gamma_0}^{-1}\bigr]\bigr)\right)=\CP_K,
$$
where the orientations of~$P_M$ and~$P_K$ are chosen so that $[P_K:P_M]=1$ if $\gamma_0$ is the first component of $\{\gamma_0,\gamma_0'\}$ and $[P_K:P_M]=-1$ if $\gamma_0$ is the second component of $\{\gamma_0,\gamma_0'\}$, and the orientation of~$\CP_K$ is induced by the orientation of~$P_K$.
\end{propos}

\begin{proof}
 From formulas~\eqref{eq_nuWgamma} and~\eqref{eq_nuWeps} it follows that $\nu_{\gamma_0,\CW(\gamma_0')}\bigl(T_{\gamma_0'}T_{\gamma_0}^{-1}\bigr)=1$ and $\nu_{\gamma,\CW}\bigl(T_{\gamma_0'}T_{\gamma_0}^{-1}\bigr)=0$ for all other pairs $(\gamma,\CW)\ne (\gamma_0,\CW(\gamma_0'))$ such that $\gamma\in\Act(M)$ and~$\CW\in\Spl(\gamma,M)$. The proposition now follows from~\eqref{eq_varphi}.
\end{proof}

\begin{cor}\label{cor_phi_surj}
  For each $p$, the homomorphism $\varphi\colon\CE_p\to\CC_{p+1}$ is surjective.
\end{cor}

\section{Subgroups~$\K$, $\CN$, $\D$, and~$\CZ$}
\label{section_KNDZ}
The goal of this section is to study the kernel of the map $\varphi\colon \CE_1\to\CC_2$; we denote this kernel by~$\K$. To this end, we define three auxiliary $\Z\bigl[\Sp(H,x)\bigr]$-submodules~$\CN$, $\D$ and~$\CZ$ of~$\CE_1$, and establish connections between them and the group~$\K$. As a result, among other things, we will obtain that $\K$ is also a $\Z\bigl[\Sp(H,x)\bigr]$-submodule of~$\CE_1$; note that this is by no means obvious, since $\varphi$ is not a $\Z\bigl[\Sp(H,x)\bigr]$-module homomorphism.

First, let $\CN\subseteq\CE_1$ be the sum of all summands in~\eqref{eq_E1p1} (for $p=1$) that correspond to negligible multicurves~$M$.

Second, consider an important oriented multicurve~$K\in\CM_2$ with bounding pair $\{\gamma,\gamma'\}$. Choose an orientation of~$P_K$. The oriented multicurves $M_1=K\setminus\gamma'$ and~$M_2=K\setminus\gamma$ belong to~$\CM_1^{(2)}$. We orient the cells~$P_{M_1}$ and~$P_{M_2}$ so that $[P_K:P_{M_1}]=[P_K:P_{M_2}]=1$. Assign to~$K$ the following element of~$\CE_1$:
\begin{equation}\label{eq_Delta_K}
\Delta_K=\iota_{P_{M_1}}\bigl(\bigl[T_{\gamma'}T_{\gamma}^{-1}\bigr]\bigr)+\iota_{P_{M_2}}\bigl(\bigl[T_{\gamma}T_{\gamma'}^{-1}\bigr]\bigr).
\end{equation}
Further, if $\delta$ is a separating simple closed curve on~$\Sigma$ disjoint from~$K$, then we set
\begin{equation}\label{eq_z_K}
z_{K,\delta}=\iota_{P_{M_1}}\bigl([T_{\delta}]\bigr)+\iota_{P_{M_2}}\bigl([T_{\delta}]\bigr).
\end{equation}
By~\eqref{eq_cd_iota}, we have $\Delta_{hK}=\Delta_K$ and $z_{hK,h\delta}=z_{K,\delta}$ for all $h\in\I$.

Let $\D\subseteq\CE_1$ be the subgroup generated by the elements~$\Delta_K$, where $K$ runs over all important oriented multicurves in~$\CM_2$, and let $\CZ\subseteq\CE_1$ be the subgroup generated by the elements~$z_{K,\delta}$, where $K$ runs over all important oriented multicurves in~$\CM_2$ and $\delta$ runs over all separating simple closed curves disjoint from~$K$.

It follows immediately from these constructions that $\CN$, $\D$, and~$\CZ$ are $\Z\bigl[\Sp(H,x)\bigr]$-submodules of~$\CE_1$.
Since all multicurves in~$\CM_1^{(2)}$ are not negligible, we have $\D\cap\CN=\{0\}$ and $\CZ\cap\CN=\{0\}$.

\begin{remark}
There exists a natural bijection between the components of~$M_1$ and the components of~$M_2$ that sends~$\gamma$ to~$\gamma'$  and sends every other component to itself. This bijection induces an affine isomorphism $P_{M_1}\to P_{M_2}$. Note that the chosen orientations of~$P_{M_1}$ and~$P_{M_2}$ disagree under this isomorphism.
\end{remark}

\begin{remark}
The signs of~$\Delta_K$ and~$z_{K,\delta}$ reverse whenever we reverse the orientation of~$P_K$. Hence the correspondence~$\CP_K\mapsto\Delta_K$ yields a well-defined homomorphism $\psi\colon\CC_2\to\CE_1$. In fact, one can similarly define the homomorphisms $\psi\colon\CC_{p+1}\to\CE_p$ and show that they give a chain map $\CC_{*+1}\to\CE_*$. However, this fact will not be used below.
\end{remark}

Let us now list three properties of the subgroups~$\K$, $\CN$, $\D$, $\CZ$ of~$\CE_1$ that will be important for us:
\begin{enumerate}
 \item $\K=\D+\CN$,
 \item $(t_x-1)^2\D\subseteq\CZ$,
 \item $(t_x-1)^2\CN=0$.
\end{enumerate}
In fact, the sum~$\D+\CN$ is direct, since $\D\cap\CN=\{0\}$.
Together facts~(1)--(3) yield the following proposition.

\begin{propos}\label{propos_K_to_Z}
 We have $(t_x-1)^2\K\subseteq\CZ$.
\end{propos}

In the rest of this section we prove in turn the three above listed facts~(1)--(3), and thus obtain the proof of Proposition~\ref{propos_K_to_Z}.

\begin{propos}\label{propos_phi_ker}
 We have  $\K=\D+\CN$.
\end{propos}

\begin{proof}
If $M$ is negligible, then the sum on the right-hand side of~\eqref{eq_varphi} is empty. Hence $\varphi|_{\CN}=0$.  By Proposition~\ref{propos_phi_K}, $\varphi$ sends the two summands on the right-hand side of~\eqref{eq_Delta_K} to the same cell~$\CP_K$ with opposite orientations. Therefore $\varphi(\Delta_K)=0$ for all important oriented multicurves $K\in\CM_2$. Thus $\K\supseteq\D+\CN$. Hence we need only to prove the inverse inclusion.

 The group~$\CE_1$ is the direct sum of the groups~$H_1(\I_M)$ over some set of representatives of $\I$-orbits of oriented multicurves in~$\CM_1$  containing no bounding pairs, i.e.,  of types $1$, $2$, and~$4$. By Proposition~\ref{propos_active_obvious}, negligible multicurves are precisely all multicurves of types~$1$ and~$4$. Therefore the sum of all summands corresponding to such multicurves coincides with the group~$\CN$, and hence is contained in the kernel of~$\varphi$.

 Let $\CE^{(2)}_1\subseteq\CE_1$ be the sum of all summands~$H_1(\I_M)$ corresponding to $\I$-orbits of type~2 oriented multicurves. Let us prove that $\ker\bigl(\varphi\big|_{\CE^{(2)}_1}\bigr)\subseteq\D$. We have
 \begin{equation}\label{eq_CE12}
 \CE^{(2)}_1=\bigoplus_A\CE^{(2)}_1(A),
 \end{equation}
 where $\CE^{(2)}_1(A)$ is the sum of the  summands~$H_1(\I_M)$ corresponding to $\I$-orbits of oriented multicurves~$M\in\CM_1^{(2)}$ with the prescribed set $[M]=A$, and the sum~\eqref{eq_CE12} is taken over all~$A$ that occur as~$[M]$ for such multicurves~$M$. For each $M\in\CM_1^{(2)}$, the right-hand side of~\eqref{eq_varphi} is a linear combination of cells~$\CP_K$ for multicurves~$K$ that have the same set of homology classes~$[K]=[M]$. Since there are no linear dependencies among the cells~$\CP_K$ corresponding to multicurves~$K$ with different sets~$[K]$, we obtain
 \begin{equation*}
  \ker\left(\varphi\big|_{\CE^{(2)}_1}\right)=\bigoplus_A\ker\left(\varphi\big|_{\CE^{(2)}_1(A)}\right).
 \end{equation*}

 Let us now compute $\ker\Bigl(\varphi\big|_{\CE^{(2)}_1(A)}\Bigr)$. A set~$A$ occurs as~$[M]$ for a type~2 multicurve if and only if $A=\{a_1,a_2,a_1+a_2,a_3\}$, where $a_1,a_2,a_3$ form a basis of a Lagrangian submodule of~$H$ such that $x=n_1a_1+n_2a_2+n_3a_3$ for some $n_1,n_2,n_3>0$. By Proposition~5.2 in~\cite{Gai21}, there are exactly two $\I$-orbits of oriented multicurves~$M$ with $[M]=A$; let $M_1$ and $M_2$ be representatives of these orbits. The two orbits are distinguished by the following property: the pair-of-pants component of~$\Sigma\setminus M_1$ lies to the left of the components of~$M_1$ with homology classes~$a_1$ and~$a_2$, and to the right of the component with homology class~$a_1+a_2$, while the pair-of-pants component of~$\Sigma\setminus M_2$ lies to the right of the components of~$M_2$ with homology classes~$a_1$ and~$a_2$, and to the left of the component with the homology class~$a_1+a_2$. Then
 $$
 \CE^{(2)}_1(A)\cong H_1(\I_{M_1})\oplus H_1(\I_{M_2}).
 $$
This isomorphism depends on the choice of orientations of the cells~$P_{M_1}$ and~$P_{M_2}$, which determine the signs of the embeddings $\iota_{P_{M_i}}\colon H_1(\I_{M_i})\hookrightarrow\CE^{(2)}_1(A)$. There is a natural bijection between the components of~$M_1$ and the components of~$M_2$ that sends every component of~$M_1$ to the component of~$M_2$ in the same homology class. This bijection induces an affine isomorphism $P_{M_1}\to P_{M_2}$. We choose the orientations of~$P_{M_1}$ and~$P_{M_2}$ so that they \textit{disagree} under this isomorphism.

 Let $\gamma_0$, $\gamma_1$, and~$\gamma_2$ be the active components of~$M_1$ with the homology classes~$a_1$, $a_2$, and $a_1+a_2$, respectively, and let $\widehat{\gamma}_0$, $\widehat{\gamma}_1$, and~$\widehat{\gamma}_2$ be the active components of~$M_2$ with the same homology classes. Let $\bigl\{e_{\gamma_i,\CW}\bigr\}$ and $\bigl\{\widehat{e}_{\widehat{\gamma}_i,\CW}\bigr\}$ be the bases of the groups~$H_1(\I_{M_1})$ and~$H_1(\I_{M_2})$, respectively, given by Proposition~\ref{propos_type2_ab}. Note that $\Spl(\gamma_i,M_1)=\Spl(\widehat{\gamma}_i,M_2)$, since the definitions of these sets depend only on the homology classes of the curves. We set
 \begin{align*}
  e_{i,\CW}&=\iota_{P_{M_1}}\bigl(e_{\gamma_i,\CW}\bigr),&
  \widehat{e}_{i,\CW}&=\iota_{P_{M_2}}\bigl(\widehat{e}_{\widehat{\gamma}_i,\CW}\bigr).
 \end{align*}
These elements together form a basis of~$\CE^{(2)}_1(A)$.
By Fact~\ref{fact_E1}, this basis does not depend on the choice of the representatives~$M_1$ and~$M_2$ of the two $\I$-orbits of multicurves. Hence we can change these representatives whenever it is convenient for us.

 By~\eqref{eq_varphi}, we have:
 \begin{align*}
  \varphi(e_{i,\CW})&=\CP_{M_1,\gamma_i,\CW},\\
  \varphi(\widehat{e}_{i,\CW})&=\CP_{M_2,\widehat{\gamma}_i,\CW}
 \end{align*}
 The cells $\CP_{M_1,\gamma_i,\CW}$ with $i=0,1,2$ and $\CW\in\Spl(\gamma_i,M_1)$ are pairwise distinct, since the bounding pairs contained in them give pairwise distinct splittings in homology. Therefore the elements~$\varphi(e_{i,\CW})$ are linearly independent in~$\CC_2$. Let us prove that $e_{i,\CW}+\widehat{e}_{i,\CW}\in\D$ for all~$i$ and~$\CW$. Then we will obtain that $\ker\Bigl(\varphi\big|_{\CE^{(2)}_1(A)}\Bigr)$ is generated by these elements, and the proposition will follow.

 Fix a pair $(i,\CW)$. As we have mentioned above, we have freedom to choose convenient representatives~$M_1$ and~$M_2$. We choose $M_1$ arbitrarily. By Proposition~\ref{propos_Wk}, there exists a simple closed curve~$\gamma'\subseteq \Sigma\setminus M_1$ such that $\{\gamma_i,\gamma'\}$ is a bounding pair and $\CW$ is the splitting  induced by this bounding pair. Then we can take $M_2=(M_1\setminus\gamma_i)\cup\gamma'$. With this choice of~$M_2$, we have $\widehat{\gamma}_i=\gamma'$ and $\widehat{\gamma}_j=\gamma_j$ for $j\ne i$. Consequently,
 \begin{align*}
  e_{i,\CW}&=\iota_{P_{M_1}}\left(\bigl[T_{\gamma'}T_{\gamma_i}^{-1}\bigr]\right),&
  \widehat{e}_{i,\CW}&=\iota_{P_{M_2}}\left(\bigl[T_{\gamma_i}T_{\gamma'}^{-1}\bigr]\right).
 \end{align*}
 Set
 $
 K=M_1\cup\gamma'=M_2\cup\gamma_i.
 $
 The orientations of~$P_{M_1}$ and~$P_{M_2}$ disagree under the affine map $P_{M_1}\to P_{M_2}$ that sends $\gamma_i$ to~$\gamma'$, and sends every other component to itself. Hence  $[P_K:P_{M_1}]=[P_K:P_{M_2}]=1$ for a suitable orientation of~$P_K$. Therefore $e_{i,\CW}+\widehat{e}_{i,\CW}=\Delta_K\in\D$.
 \end{proof}

 \begin{cor}
  The group~$\K$ is a $\Z\bigl[\Sp(H,x)\bigr]$-submodule of~$\CE_1$.
 \end{cor}

 \begin{propos}\label{propos_D_to_Z}
 We have $(t_x-1)^2\D\subseteq \CZ$.
\end{propos}

\begin{proof}
 Consider an arbitrary generator $\Delta_K$ of~$\D$ given by~\eqref{eq_Delta_K}. Set
 $$
 K=\gamma\cup\gamma'\cup\gamma_1\cup\gamma_2\cup\gamma_3,
 $$
 where $\gamma_1$ and~$\gamma_2$ lie in the same connected component of~$M\setminus(\gamma\cup\gamma')$, and $\gamma_3$ lies in the other connected component. Then $\gamma$, $\gamma_1$ and~$\gamma_2$ are the active components of~$M_1=K\setminus \gamma'$, while  $\gamma'$, $\gamma_1$ and~$\gamma_2$ are the active components of~$M_2=K\setminus \gamma$. Ignore the original orientations of the components of~$K$, and orient them so that to obtain that $[\gamma]=[\gamma']=[\gamma_1]+[\gamma_2]$. Set  $c_i=[\gamma_i]$ for $i=1,2,3$. Choose a separating simple closed curve~$\delta$ disjoint from~$K$, and consider the element~$z_{K,\delta}$ given by~\eqref{eq_z_K}.  Proposition~\ref{propos_BP_delta}  implies that
 $$
 (t_{c_1+c_3}-1)(t_{c_2+c_3}-1)\Delta_K=(t_{c_1+c_3}-1)(t_{c_2+c_3}-1)z_{K,\delta}.
 $$
 It is important here that equality~\eqref{eq_tt_anih} in Proposition~\ref{propos_BP_delta} holds for all choices of signs $\sigma_1,\sigma_2\in\{-1,1\}$; hence the above equality holds even though we may have changed the signs of the classes~$c_i=[\gamma_i]$.
 Every of the transvections~$t_{c_1}$, $t_{c_2}$, $t_{c_1+c_2}$, $t_{c_3}$ fixes~$\Delta_K$, since they are induced by twists about the components of either of the multicurves~$M_i$. Moreover, $t_{c_1+c_2+c_3}$ also fixes~$\Delta_K$ by Proposition~\ref{propos_e_annih}. By Lemma~\ref{lem_transvec3}, we have
 $$
 t_{c_1+c_2+c_3}=t_{c_1+c_2}t_{c_1+c_3}t_{c_2+c_3}t_{c_1}^{-1}t_{c_2}^{-1}t_{c_3}^{-1}.
 $$
 Since all symplectic transvections entering this equality mutually commute and each of them, except~$t_{c_1+c_3}$ and~$t_{c_2+c_3}$, fixes~$\Delta_K$, we obtain that $t_{c_2+c_3}\Delta_K=t_{c_1+c_3}^{-1}\Delta_K$, and hence
 $$
 (t_{c_1+c_3}-1)(t_{c_2+c_3}-1)\Delta_K=-t_{c_1+c_3}^{-1}(t_{c_1+c_3}-1)^2\Delta_K.
 $$
 Therefore $(t_{c_1+c_3}-1)^2\Delta_K$ belongs to~$\CZ$.

Further, $x$ lies in~$\langle K\rangle$, hence $x=m_1c_1+m_2c_2+m_3c_3$ for some $m_1,m_2,m_3\in\Z$. By Lemma~\ref{lem_transvec3}, we have:
$$
t_x\Delta_K=t_{c_1+c_3}^{m_1m_3}t_{c_2+c_3}^{m_2m_3}\Delta_K=t_{c_1+c_3}^{m_1m_3-m_2m_3}\Delta_K.
$$
Consequently, the element
$$
(t_x-1)^2\Delta_K=\bigl(t_{c_1+c_3}^{m_1m_3-m_2m_3}-1\bigr)^2\Delta_K
$$
belongs to~$\CZ$.
\end{proof}

 \begin{propos}\label{propos_N_anih}
 We have $(t_x-1)^2\CN=0$.
\end{propos}

The group~$\CN$ is the sum of the groups~$H_1(\I_M)$ corresponding to oriented multicurves of types~1 and~4. We treat the two types of multicurves separately.

\begin{propos}\label{propos_type4_anih}
 If $M\in\CM_1^{(4)}$, then $(t_x-1)H_1(\I_M)=0$.
\end{propos}

\begin{proof}
 Removing from~$M$ one of the two components in the bounding pair, we  obtain a three-component multicurve $N$ that bounds a pair of pants. Moreover, the homology class of one of the three components of~$N$ is equal to the sum of the homology classes of the other two components. Let $a_1$ and~$a_2$ be the homology classes of these two remaining components. Since the Dehn twists about the components of~$N$ commute with all elements of~$\I_M$, it follows that the transvections~$t_{a_1}$, $t_{a_2}$, and~$t_{a_1+a_2}$ act trivially on~$H_1(\I_M)$. Finally, since $M\in\CM_1$, it follows that $x=n_1a_1+n_2a_2$ for some positive integers~$n_1$ and~$n_2$. It now follows from Lemma~\ref{lem_transvec3} that $t_x$ acts trivially on~$H_1(\I_M)$.
\end{proof}

The case of a type~1 multicurve is more complicated. Recall that there exists a surjective homomorphism
$$
\tau\colon \I\to\wedge^3H/H
$$
called the \textit{Johnson homomorphism}, see~\cite{Joh80b}. Here $\wedge^3H$ is the third exterior power of~$H$, and $H$ embeds into~$\wedge^3H$ by $v\mapsto v\wedge\Omega$, where $\Omega\in\wedge^2H$ is the dual of the intersection form on~$H$.
Moreover, $\tau$ is equivariant in the following sense: for any $h\in\I$ and $f\in\Mod(\Sigma)$, we have
\begin{equation}\label{eq_tau_equiv}
\tau\bigl(fhf^{-1}\bigr)=\kappa(f)\tau(h),
\end{equation}
where $\kappa\colon \Mod(\Sigma)\to\Sp(H)$ is the standard homomorphism. Hence $\tau$ induces an $\Sp(H)$-equivariant homomorphism
$$
\tau_*\colon H_1(\I)\to\wedge^3H/H
$$

A standard formula for the value of~$\tau$ on a BP-map (see~\cite{Joh80b}) in the genus~$3$ case takes the following form. Let $\{\gamma,\gamma'\}$ be a bounding pair on~$\Sigma$, and let $c=[\gamma]=[\gamma']$ for certain orientations of~$\gamma$ and~$\gamma'$. Let $S$ be the connected component of~$\Sigma\setminus(\gamma\cup\gamma')$ lying  to the right of~$\gamma$. Suppose that $H_1(S)=\langle a,b,c\rangle$ and $a\cdot b=1$. Then
\begin{equation}\label{eq_tau_value}
 \tau\bigl(T_{\gamma}T_{\gamma'}^{-1}\bigr)=a\wedge b \wedge c.
\end{equation}

\begin{propos}\label{propos_type1_Johnson}
 If $M\in\CM_1^{(1)}$, then the homomorphism
 $$
 \left(\tau|_{\I_M}\right)_*\colon H_1(\I_M)\to\wedge^3H/H
 $$
 is injective.
\end{propos}

\begin{proof}
If we ignore the original orientations of components of~$M$, then we can orient these components so that their homology classes are $a_1$, $a_2$, $a_3$, and~$a_1+a_2+a_3$ for some basis $a_1$, $a_2$, $a_3$ of a Lagrangian  submodule of~$H$. For each $i=1,2,3$, we denote by~$\alpha_i$ the component of~$M$ with the homology class~$a_i$, and we denote by~$\alpha_0$ the component with the homology class~$a_1+a_2+a_3$. We arrange the multicurve~$M$ as in Fig.~\ref{fig_M_type1}, and  consider the simple closed curves~$\gamma_1$, $\gamma_2$, $\delta_1$, $\delta_2$, $\beta_1$, $\beta_2$, and~$\beta_3$ depicted in this figure. Set $b_i=[\beta_i]$. Then $a_1,a_2,a_3,b_1,b_2,b_3$ form a symplectic basis of~$H$.

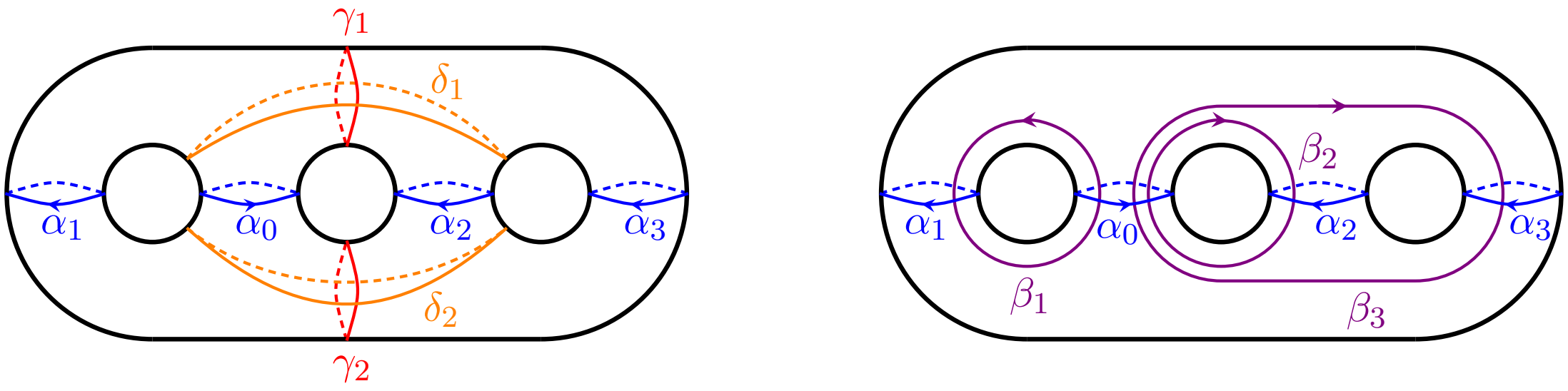
\begin{figure}[!ht]
 \begin{tikzpicture}[scale=.45]
\small
\definecolor{myblue}{rgb}{0, 0, 0.7}
\definecolor{mygreen}{rgb}{0, 0.4, 0}
\tikzset{every path/.append style={line width=.2mm}}
\tikzset{my dash/.style={dash pattern=on 2pt off 1.5pt}}

\tikzset{
  mid stealth/.style={
    postaction={
      decorate,
      decoration={
        markings,
        mark=at position 0.58 with {\arrow{stealth}}
      }
    }
  }
}

\tikzset{
  mmid stealth/.style={
    postaction={
      decorate,
      decoration={
        markings,
        mark=at position 0.5 with {\arrow{stealth}}
      }
    }
  }
}

\begin{scope}[shift={(-9,0)}]
% Invisible lines

\draw[red,my dash,thick] (0,-3.03) .. controls (-0.3,-2) .. (0,-.97); % e_0
\draw[red,my dash,thick] (0,3.03) .. controls (-0.3,2) .. (0,.97); % e_0'
\draw[color=blue,my dash, thick] (-7.03,0) .. controls (-6,0.3) .. (-4.97,0); %e_1
\draw[color=blue,my dash, thick] (7.03,0) .. controls (6,0.3) .. (4.97,0); %e_2
\draw[color=blue,my dash, thick] (-3.03,0) .. controls (-2,0.3) .. (-.97,0); % e_0+e_1
\draw[color=blue,my dash, thick] (3.03,0) .. controls (2,0.3) .. (.97,0); % e_0+e_2
\draw[color=orange,my dash,thick] (-3.293,-.707) .. controls (-1.5,-2.2) and (1.5,-2.2) .. (3.293,-.707);
\draw[color=orange,my dash,thick] (-3.293,.707) .. controls (-1.5,2.8) and (1.5,2.8) .. (3.293,.707) node [pos=.8,above=-1pt] {$\delta_1$};

% Border lines

\draw [very thick] (0,0) circle (1);
\draw [very thick] (4,0) circle (1);
\draw [very thick] (-4,0) circle (1);
\draw [very thick]  (-4,3)--(4,3);
\draw [very thick]  (-4,-3)--(4,-3);
\draw [very thick] (-4,3) arc (90:270:3); 
\draw [very thick] (4,3) arc (90:-90:3);

% Visible lines 

\draw[color=red,thick] (0,-3.03) node [below,,xshift=.5mm] {$\gamma_2$} .. controls (0.3,-2) .. (0,-.97) ;  %e_0
\draw[color=red,thick] (0,3.03) node [above=-1pt,xshift=.5mm] {$\gamma_1$} .. controls (0.3,2) .. (0,.97) ;  %e_0'

\draw[color=blue, thick, mid stealth]  (-4.97,0) .. controls (-6,-0.3) ..  (-7.03,0) node [pos=0.4,below=-1pt] {$\alpha_1$};
\draw[color=blue, thick,mid stealth] (-3.03,0) .. controls (-2,-0.3) .. (-.97,0) node [pos=.6,below=-1pt] {$\alpha_0$};
\draw[color=blue, thick, mid stealth] (3.03,0) .. controls (2,-0.3) ..  (.97,0)  node [pos=.4,below=-1pt] {$\alpha_2$};
\draw[color=blue, thick,mid stealth] (7.03,0) .. controls (6,-0.3) .. (4.97,0) node [pos=.4,below=-1pt] {$\alpha_3$};

\draw[color=orange,thick] (-3.293,-.707) .. controls (-1,-2.8) and (1,-2.8) .. (3.293,-.707) node [pos=.8,below=0pt] {$\delta_2$};
\draw[color=orange, thick] (-3.293,.707) .. controls (-1,2.2) and (1,2.2) .. (3.293,.707);

%\node[] at (0,-5) {\textit{à}};
\end{scope}

\begin{scope}[shift={(9,0)}]
% Invisible lines

\draw[color=blue,my dash, thick] (-7.03,0) .. controls (-6,0.3) .. (-4.97,0); %e_1
\draw[color=blue,my dash, thick] (7.03,0) .. controls (6,0.3) .. (4.97,0); %e_2
\draw[color=blue,my dash, thick] (-3.03,0) .. controls (-2,0.3) .. (-.97,0); % e_0+e_1
\draw[color=blue,my dash, thick] (3.03,0) .. controls (2,0.3) .. (.97,0); % e_0+e_2

% Border lines

\draw [very thick] (0,0) circle (1);
\draw [very thick] (4,0) circle (1);
\draw [very thick] (-4,0) circle (1);
\draw [very thick]  (-4,3)--(4,3);
\draw [very thick]  (-4,-3)--(4,-3);
\draw [very thick] (-4,3) arc (90:270:3); 
\draw [very thick] (4,3) arc (90:-90:3); 

% Visible lines 

%\draw [color=violet,thick] (-2.5,0) arc (0:360:1.5) ;
\draw [color=violet,thick,-stealth] (-4,1.5)  arc (-270:95:1.5) node [pos=.5,below=-1pt] {$\beta_1$};

\draw [color=violet,thick,-stealth] (0,1.5) arc (90:-275:1.5) node [pos=.15,xshift=3.5mm] {$\beta_2$} ;
\draw [color=violet,thick] (0,1.8) arc (90:270:1.8) -- (4,-1.8) node [pos=.75,below=-1pt] {$\beta_3$} arc (-90:90:1.8) -- (0,1.8) ;
\draw [color=violet,thick,-stealth] (2,1.8)--(2.6,1.8);

\draw[color=blue, thick, mid stealth] (-4.97,0)  .. controls (-6,-0.3) .. (-7.03,0)  node [pos=0.4,below=-1pt,xshift=-1mm] {$\alpha_1$};
\draw[color=blue, thick,mid stealth] (-3.03,0) .. controls (-2,-0.3) .. (-.97,0)  node [pos=.6,below=-1pt,xshift=-1.25mm,yshift=-.4mm] {$\alpha_0$};
\draw[color=blue, thick, mid stealth] (3.03,0) .. controls (2,-0.3) ..  (.97,0)  node [pos=.4,below=-1pt,xshift=1mm] {$\alpha_2$};
\draw[color=blue, thick,mid stealth] (7.03,0) .. controls (6,-0.3) .. (4.97,0) node [pos=.4,below=-1pt,xshift=1mm] {$\alpha_3$};

\end{scope}

\end{tikzpicture}
 \caption{A multicurve $M\in\CM_1^{(1)}$ and curves $\gamma_1$, $\gamma_2$, $\delta_1$, $\delta_2$, $\beta_1$, $\beta_2$, and~$\beta_3$}\label{fig_M_type1}
\end{figure}

By Proposition~7.1 in~\cite{Gai21}, the group~$\I_M$ is generated by the three elements
$$
h_1=T_{\gamma_1}T_{\gamma_2}^{-1},\qquad h_2=T_{\delta_1}T_{\delta_2}^{-1},\qquad
h_3=T_{\gamma_1}^{-1}T_{\delta_1}^{-1}T_{\gamma_1}T_{\delta_1}.
$$
Therefore $H_1(\I_M)$ is generated by the three elements~$[h_1]$, $[h_2]$, and~$[h_3]$. By~\eqref{eq_tau_value}, we have
\begin{align*}
 \tau(h_1)&=a_1\wedge b_1\wedge(a_2+a_3),\\
 \tau(h_2)&=a_2\wedge b_2\wedge(a_1+a_3).
\end{align*}
Moreover, by~\eqref{eq_tau_equiv}, we obtain
$$
\tau\bigl(T_{\delta_1}^{-1}h_1T_{\delta_1}\bigr)=t_{a_1+a_3}^{-1}\cdot\bigl(a_1\wedge b_1\wedge(a_2+a_3)\bigr)=a_1\wedge b_1\wedge(a_2+a_3)-a_1\wedge a_2\wedge a_3,
$$
and hence
\begin{equation*}
\tau(h_3)=\tau\bigl(h_1^{-1}T_{\delta_1}^{-1}h_1T_{\delta_1}\bigr)=-a_1\wedge a_2\wedge a_3.
\end{equation*}
(The last formula is, in fact, a partial case of a general formula for the value of~$\tau$ on the commutator of a simply intersecting pair due to Childers~\cite{Chi12}, Church~\cite{Chu14}, and Putman~\cite{Put09b}\footnote{This result is not included in the final journal version of the preprint~\cite{Put09b}.}.) It can be checked directly that the elements $a_1\wedge b_1\wedge(a_2+a_3)$, $a_2\wedge b_2\wedge(a_1+a_3)$, and~$a_1\wedge a_2\wedge a_3$  are linearly independent in $\wedge^3H/H$. Thus $\left(\tau|_{\I_M}\right)_*$ is injective.
\end{proof}

It is a standard fact that $(t_x-1)^2{\wedge^3H}=0$. Consequently, Proposition~\ref{propos_type1_Johnson} implies the following corollary.

\begin{cor}\label{cor_type1_anih}
 If $M\in\CM_1^{(1)}$, then $(t_x-1)^2H_1(\I_M)=0$.
\end{cor}

Proposition~\ref{propos_N_anih} follows from Proposition~\ref{propos_type4_anih} and Corollary~\ref{cor_type1_anih}.

\section{Reduction to~$\CF_{0,2}$}\label{section_to_F02}

In this section, we prove the following proposition.

\begin{propos}\label{propos_txE11}
We have $(t_x-1)^2\CF_{1,1}\subseteq\CF_{0,2}+\bA$.
\end{propos}

Let $\CY\subseteq E^1_{1,1}$ be the subgroup generated by all elements $$y_{N,\delta}=\iota_{P_N}([T_{\delta}]),$$ where
\begin{itemize}
 \item $N\in\CM_1$ is a multicurve containing a bounding pair,
 \item $\delta$ is a separating simple closed curve on~$\Sigma$ disjoint from~$N$.
\end{itemize}

Proposition~\ref{propos_txE11} will follow from the following two propositions.

\begin{propos}\label{propos_annih_E111}
 We have the inclusion
 $$
 (t_x-1)^2\ker\bigl[d^1\colon E^1_{1,1}\to E^1_{0,1}\bigr]\subseteq d^1E^1_{2,1}+\CY.
 $$
\end{propos}

\begin{propos}\label{propos_Y_Ab}
 For each generator~$y_{N,\delta}$ of~$\CY$, there exists a standard abelian cycle $\CA\in\bA$ such that $\CA\in\CF_{1,1}$ and the image of~$\CA$ under the projection $\CF_{1,1}\twoheadrightarrow E^{\infty}_{1,1}$ is represented by~$y_{N,\delta}$.
\end{propos}

Let us first derive Proposition~\ref{propos_txE11} from Propositions~\ref{propos_annih_E111} and~\ref{propos_Y_Ab}, and then prove those two propositions in turn.

\begin{proof}[Proof of Proposition~\ref{propos_txE11}]
Proposition~\ref{propos_Y_Ab} implies that $\CY$ is contained in the kernel of~$d^1$. Let $Y$ be the image of~$\CY$ in~$E^{\infty}_{1,1}$. Then Proposition~\ref{propos_annih_E111} implies that $(t_x-1)^2E^{\infty}_{1,1}\subseteq Y$.

Consider an arbitrary element $b\in\CF_{1,1}$. Let $\pi\colon \CF_{1,1}\twoheadrightarrow E^{\infty}_{1,1}$ be the projection. Then $(t_x-1)^2\pi(b)\in Y$. Therefore, by Proposition~\ref{propos_Y_Ab}, there exists $a\in\bA$ such that $\pi(a)=(t_x-1)^2\pi(b)$. Consequently, the element $(t_x-1)^2b-a$ lies in the kernel of~$\pi$, i.e., in~$\CF_{0,2}$. The proposition follows.
\end{proof}

The proof of Proposition~\ref{propos_annih_E111} will be divided into several lemmas.

\begin{lem}\label{lem_toCE}
 Suppose that $\zeta\in E^1_{1,1}$. Then there exists $\eta\in E^1_{2,1}$ such that $\zeta+d^1\eta\in\CE_1$.
\end{lem}

\begin{proof}
 It suffices to prove the lemma for $\zeta=\iota_{P_M}([h])$, where $M\in\CM_1$ is an oriented multicurve containing a bounding pair~$\{\gamma_0,\gamma_0'\}$, and $h\in\I_M$. Then $M$ is a multicurve of type~3, 5, or~6. Let $S_1$ and~$S_2$ be the connected components of~$\Sigma\setminus M$. For $i=1,2$, let $\gamma_i$ be the component of~$M$ lying in the closure $\overline{S}_i$ of~$S_i$ and different from~$\gamma_0$ and~$\gamma_0'$ (if such a component exists). For type~3, both components~$\gamma_1$ and~$\gamma_2$ are present, while for types~5 and~6, one or both of them are absent.

  We can represent $h$ by a  homeomorphism that fixes  the curves~$\gamma_0$ and~$\gamma_0'$ pointwise. Hence $h=h_1h_2$, where $h_1,h_2\in\I_M$, and $h_i$ is represented by a homeomorphism that is the identity outside~$\overline{S}_i$. Then $\zeta=\zeta_1+\zeta_2$, where $\zeta_i=\iota_{P_M}([h_i])$, and it suffices to prove the lemma for~$\zeta_1$ and~$\zeta_2$ separately. These two cases are completely analogous, so we consider the case of~$\zeta_2$.

  First, assume that $\gamma_1$ is present in~$M$. Then $S_1$ is a four-punctured sphere. Choose an oriented simple closed curve~$\varepsilon$ in~$S_1$ such that $[\varepsilon]=[\gamma_0]-[\gamma_1]$, and set $K=M\cup\varepsilon$. We have $x=n_0[\gamma_0]+n_1[\gamma_1]+n_2[\gamma_2]$ for $n_0,n_1,n_2>0$ if $\gamma_2$ is also present in~$M$, and $x=n_0[\gamma_0]+n_1[\gamma_1]$ for $n_0,n_1>0$ if $\gamma_2$ is absent in~$M$. In both cases, we see that $K\in\CM_2$ and $P_K$ is a triangle with vertices~$P_{\gamma_0\cup\gamma_1(\cup\gamma_2)}$, $P_{\gamma_0'\cup\gamma_1(\cup\gamma_2)}$, and $P_{\varepsilon\cup\gamma_1(\cup\gamma_2)}$, where $\gamma_2$ is present in all these vertices if and only if it is  present in~$M$. Hence $\partial P_K=P_M+P_{M_1}+P_{M_2}$ (for appropriate orientations of the cells), where the multicurves $M_1=K\setminus\gamma_0'$ and  $M_2=K\setminus\gamma_0$ do not contain bounding pairs. Since $h_2$ is supported in~$\overline{S}_2$, we have $h_2\in\I_K$. Then
  \begin{equation}\label{eq_toCE}
   \iota_{P_M}([h_2])-d^1\iota_{P_K}([h_2])=-\iota_{P_{M_1}}([h_2])-\iota_{P_{M_2}}([h_2])\in\CE_1,
  \end{equation}
 as required.

 Second, assume that $\gamma_1$ is absent in~$M$. Then $S_1$ is a twice-punctured torus. Divide~$S_1$ into two pairs of pants by two disjoint non-separating oriented simple closed curves~$\varepsilon_1$ and~$\varepsilon_2$ such that $[\varepsilon_1]+[\varepsilon_2]=[\gamma_0]$, and set $K=M\cup\varepsilon_1\cup\varepsilon_2$. We have $x=n_0[\gamma_0]+n_2[\gamma_2]$ for $n_0,n_2>0$ if $\gamma_2$ is present in~$M$, and $x=[\gamma_0]$ if $\gamma_2$ is absent in~$M$. It follows that $K\in\CM_2$ and $P_K$ is a triangle with vertices~$P_{\gamma_0(\cup\gamma_2)}$, $P_{\gamma_0'(\cup\gamma_2)}$, and $P_{\varepsilon_1\cup\varepsilon_2(\cup\gamma_2)}$. Hence we again obtain that $\partial P_K=P_M+P_{M_1}+P_{M_2}$, where the multicurves $M_1=K\setminus\gamma_0'$ and  $M_2=K\setminus\gamma_0$ do not contain bounding pairs, and thus again get~\eqref{eq_toCE}.
\end{proof}

Recall that we denote by~$\K$ the kernel of the homomorphism $\varphi\colon \CE_1\to\CC_2$.

\begin{lem}\label{lem_annih_E111_1}
 Suppose that $\zeta\in E^1_{1,1}$ and $d^1\zeta=0$. Then there exists $\eta\in E^1_{2,1}$ such that $\zeta+d^1\eta\in\K$.
\end{lem}

\begin{proof}
 By Lemma~\ref{lem_toCE}, we may assume that $\zeta\in\CE_1$. By Proposition~\ref{propos_phi_chain}, we have a commutative diagram
 \begin{equation*}
 \begin{tikzcd}
    \CE_0 \arrow[d,twoheadrightarrow,"\varphi"] &\CE_1 \arrow[d,twoheadrightarrow,"\varphi"]\arrow[l,"{{}\,d^1}"'] & \CE_2 \arrow[d,twoheadrightarrow,"\varphi"] \arrow[l,"{{}\,d^1}"']\\
  \CC_1 &\CC_2 \arrow[l,"{-\partial_{\CC}}"']& \CC_3 \arrow[l,"{-\partial_{\CC}}"']
 \end{tikzcd}
\end{equation*}
The vertical arrows are surjective by Corollary~\ref{cor_phi_surj}, and the lower row is exact by Proposition~\ref{propos_H2C0}. It follows easily that
$$
\ker\left[d^1\colon\CE_1\to\CE_0\right]\subseteq d^1\CE_2+ \K.
$$
The proposition follows.
\end{proof}

\begin{lem}\label{lem_annih_E111_2}
 We have $\CZ\subseteq \CY+d^1E^1_{2,1}$.
\end{lem}

\begin{proof}
 Consider an arbitrary generator~$z_{K,\delta}$ of~$\CZ$. As in Section~\ref{section_KNDZ}, let $\{\gamma,\gamma'\}$ be the bounding pair contained in~$K$, and let $M_1=K\setminus\gamma'$ and $M_2=K\setminus\gamma$, so that $z_{K,\delta}$ is given by~\eqref{eq_z_K}. Recall that, by our convention, the orientations of~$P_{M_1}$ and~$P_{M_2}$ are chosen so that $[P_K:P_{M_1}]=[P_K:P_{M_2}]=1$. We write
 $$
 \partial P_K=P_{M_1}+P_{M_2}+P_{N_1}+\cdots+P_{N_k},
 $$
 where $N_1,\ldots,N_k$ are all oriented submulticurves of~$K$ that belong to~$\CM_1$ and are different from~$M_1$ and~$M_2$. (In fact, one can easily show that $k$ is either~$1$ or~$2$, but this will not be used below.) Each multicurve~$N_i$ is contained in neither~$M_1$ nor~$M_2$. Therefore, $N_i$ contains the bounding pair~$\{\gamma,\gamma'\}$. We now have
 \begin{align*}
 d^1\bigl(\iota_{P_K}([T_{\delta}])\bigr)&=\iota_{P_{M_1}}([T_{\delta}])+\iota_{P_{M_2}}([T_{\delta}])+\iota_{P_{N_1}}([T_{\delta}])+\cdots+\iota_{P_{N_k}}([T_{\delta}])\\
 {}&=z_{K,\delta}+y_{N_1,\delta}+\cdots+y_{N_k,\delta}.
 \end{align*}
 Therefore $z_{K,\delta}\in \CY+d^1E^1_{2,1}$.
\end{proof}

Proposition~\ref{propos_annih_E111} follows from Lemmas~\ref{lem_annih_E111_1} and~\ref{lem_annih_E111_2}, and Proposition~\ref{propos_K_to_Z}.

\begin{proof}[Proof of Proposition~\ref{propos_Y_Ab}]
The multicurve~$N$ has type~3, 5, or~6. Let $\{\gamma_0,\gamma_0'\}$ be the bounding pair contained in~$N$. If $N$ has type~3, then $N=\gamma_0\cup\gamma_0'\cup\gamma_1\cup\gamma_2$ and~$\delta$ can be arranged as in Fig.~\ref{fig_yNd}. If $N$ has type~5 or~6, then the only difference is that one or both of the components~$\gamma_1$ and~$\gamma_2$ may be absent. The cell~$P_N$ is the segment with endpoints~$P_{N_1}$ and~$P_{N_2}$, where $N_1=N\setminus\gamma_0'$ and $N_2=N\setminus\gamma_0$. We may assume that $P$ is oriented from~$P_{N_1}$ to~$P_{N_2}$.

\begin{figure}[!ht]
 \begin{tikzpicture}[scale=.4,
        mid arrow/.style={
            postaction={
                decorate,
                decoration={
                    markings,
                    mark=at position 0.55 with {
                        \arrow[#1]{stealth}
                    }
                }
            }
        },
        mid arrow/.default=black,
    ]

\small
% \definecolor{blue}{rgb}{0, 0, 0.7}
% \definecolor{mygreen}{rgb}{0, 0.4, 0}
\tikzset{every path/.append style={line width=.2mm}}
\tikzset{my dash/.style={dash pattern=on 2pt off 1.5pt}}

% Border lines

\draw[color=violet,  my dash] (1.93,3.03) .. controls (2.63,2) and (2.63,-2) .. (1.93,-3.03);

\draw[color=red, thick, my dash] (4,3.03) .. controls (4.7,2) and (4.7,-2) .. (4,-3.03);

\draw[color=blue, thick, my dash] (-9.03,0) .. controls (-8,.3) .. (-6.97,0);
\draw[color=blue, thick, my dash] (0,-.97) .. controls (.3,-2) .. (0,-3.03);
\draw[color=blue, thick, my dash] (0,.97) .. controls (.3,2) .. (0,3.03);
\draw[color=blue, thick, my dash] (6.97,0) .. controls (8,.3) .. (9.03,0);
%\draw[color=blue] (1.93,3.03) -- (1.93,-3.03);

\draw [very thick] (-6,0) circle (1);
\draw [very thick] (0,0) circle (1);
\draw [very thick] (6,0) circle (1);
\draw [very thick]  (-6,3)--(6,3) arc (90:-90:3) -- (-6,-3) arc (270:90:3);

% Visible lines

\draw[color=blue, thick] (-9.03,0) .. controls (-8,-.3) .. (-6.97,0) node [pos=.5,below,xshift=.5mm] {$\gamma_1$};
\draw[color=blue, thick] (0,-.97) .. controls (-.3,-2) .. (0,-3.03) node [below,xshift=.5mm] {$\gamma_0'$};
\draw[color=blue, thick] (0,.97) .. controls (-.3,2) .. (0,3.03) node [above,xshift=.5mm] {$\gamma_0$};
\draw[color=blue, thick] (6.97,0) .. controls (8,-.3) .. (9.03,0) node [pos=.5,below,xshift=.5mm] {$\gamma_2$};;

\draw [violet] (0,0) circle (1.6) node [below right=3mm,xshift=1mm,yshift=1mm] {$\beta$};
\draw [violet] (-2.1,0) arc (180:90:2.1) arc (-90:-39:2.5);
\draw [violet] (-2.1,0) arc (180:270:2.1) node [pos=.35,below left=-2pt] {$\beta'$} arc (90:39:2.5);

\draw[color=red, thick] (4,3.03) .. controls (3.3,2) and (3.3,-2) .. (4,-3.03) node [pos=0.8,left=-1pt] {$\delta$};

\end{tikzpicture}
 \caption{Multicurve~$N$ and curves~$\delta$, $\beta$, and~$\beta'$}\label{fig_yNd}
\end{figure}
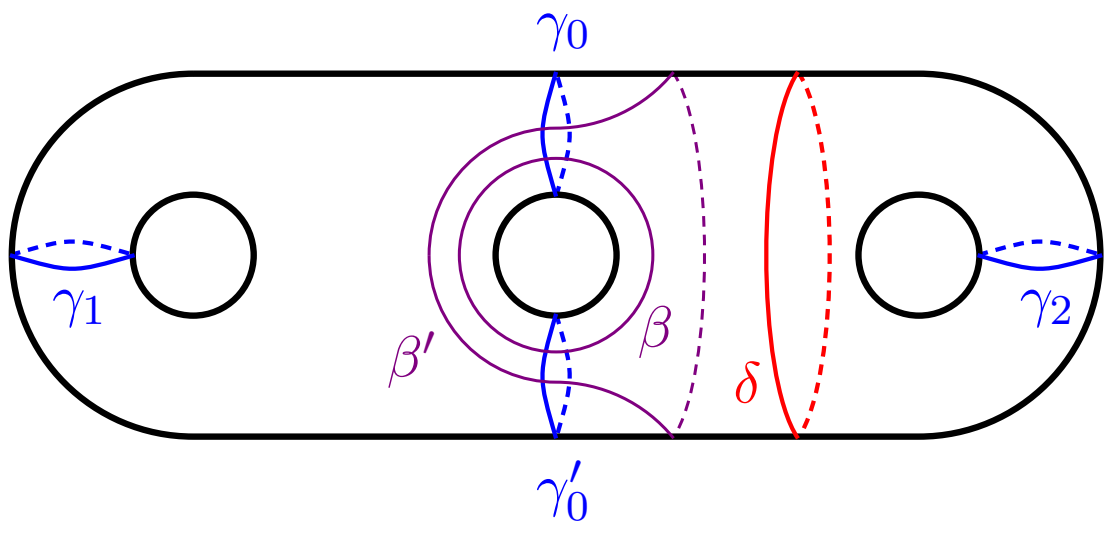

Consider the bounding pair~$\{\beta,\beta'\}$ depicted in Fig.~\ref{fig_yNd}, and set $\CA=\CA\bigl(T_{\delta},T_{\beta}T_{\beta'}^{-1}\bigr)$. Then $\CA\in\bA$. The Dehn twist~$T_{\delta}$ stabilizes~$N$, and hence fixes~$P_N$ pointwise. The BP-map~$T_{\beta}T_{\beta'}^{-1}$ sends~$\gamma_0'$ to~$\gamma_0$, and therefore sends~$P_{N_2}$ to~$P_{N_1}$. From Lemma~\ref{lem_ab_F11} it now follows that $\CA\in\CF_{1,1}$ and the image of~$\CA$ in~$E^{\infty}_{1,1}$ is represented by~$y_{N,\delta}$.
\end{proof}

\section{Stabilizer of a two-component multicurve}\label{section_stab2}

The goal of this section is to prove the following proposition.

\begin{propos}\label{propos_anih_2comp}
 Let $N=\alpha_1\cup\alpha_2$ be a two-component oriented multicurve on~$\Sigma$ that does not separate~$\Sigma$, and let $a_1=[\alpha_1]$ and $a_2=[\alpha_2]$. Then
 $$
 (t_{a_1+a_2}-1)^3H_2(\I_N)\subseteq \bA_N.
 $$
\end{propos}

Here we denote by~$\bA_N$ the subgroup of~$H_2(\I_N)$ generated by all standard abelian cycles of the form~$\CA(T_{\delta_1},T_{\delta_2})$ (where $\delta_1$ and~$\delta_2$ are disjoint and separating) and of the form~$\CA\bigl(T_{\delta},T_{\gamma_1}T_{\gamma_2}^{-1}\bigr)$ (where $\delta$ is a separating simple closed curve and $\{\gamma_1,\gamma_2\}$ is a bounding pair disjoint from it) such that the curves~$\delta_1$, $\delta_2$, $\delta$, $\gamma_1$, and~$\gamma_2$ are disjoint from~$N$.

The proof of Proposition~\ref{propos_anih_2comp} repeats, with some modifications, the content of Sections~\ref{section_BBM}--\ref{section_to_F02}. We will denote the objects that are analogues of the objects from these sections by the same letters with tildes.

Throughout this section, we fix a multicurve $N=\alpha_1\cup\alpha_2$ as in Proposition~\ref{propos_anih_2comp}.

\subsection{Complex~$\tCB$}

In this subsection we define a complex~$\tCB$ which serves as an analogue of the complex of cycles~$\CB$ for the stabilizer~$\I_N$. Let $a_1=[\alpha_1]$ and $a_2=[\alpha_2]$ be the homology classes of the components of~$N$. We choose and fix a primitive class $v\in H$ such that $a_1,a_2,v$ form a basis of a Lagrangian submodule of~$H$.

Let $\tCB=\tCB(N,v)$ be a simplicial complex defined as follows:
\begin{itemize}
 \item the vertices of~$\tCB$ are the isotopy classes of oriented simple closed curves~$\gamma$ on~$\Sigma$ such that $$[\gamma]=v+k_1a_1+k_2a_2$$ for some $k_1,k_2\in\Z$, and $\gamma$ is disjoint from~$N$;
 \item a finite set of vertices spans a simplex if and only if the corresponding isotopy classes admit pairwise disjoint representatives.
\end{itemize}

Analogues of the complex of cycles for surfaces with punctures (or boundary components) have appeared in various works. Already in the original paper~\cite{BBM10}, Bestvina, Bux, and Margalit studied the complex of cycles for the twice-punctured torus. Spiridonov~\cite[Section~7.1]{Spi23} considered the complex of cycles for the twice-punctured surface of genus~$2$. Our complex~$\tCB$ is a suitable analogue of the complex of cycles for the four-punctured torus~$\Sigma\setminus N$. The general construction of the complex of cycles for surfaces with punctures and boundary components, in sufficient generality for our purposes, to the best of the author's knowledge, first appeared in Patsankov's paper~\cite{Pat26}. Thus, Proposition~\ref{propos_tB_contr} below is a special case of Theorem~2.7 from that paper. Nevertheless, for the reader's convenience, we will give its proof.

\begin{propos}\label{propos_tB_contr}
The complex~$\tCB$ is contractible.
\end{propos}

\begin{proof}
 Let $\CB(v)$ be the usual Bestvina--Bux--Margalit complex of cycles for the homology class~$v$ on the closed surface~$\Sigma$, and let $\CB(v,N)$ be the subcomplex of~$\CB(v)$ consisting of all cells~$P_M$ such that $M$ admits a representative disjoint from~$N$. Note that $M$ may have components isotopic to the components of~$N$; moreover, the orientations of these components in~$M$ may or may not coincide with their orientations in~$N$. Then the complex~$\CB(v,N)$ is contractible. Indeed, the standard surgical proof of the contractibility of~$\CB(v)$ (see~\cite[Section~5]{BBM10}) consists in constructing a strong deformation retraction $H\colon \CB(v)\times[0,1]\to\CB(v)$ of~$\CB(v)$ onto a point~$\bgamma_0\in \CB(v)$, given by
 $$
 H(\bgamma,s)=\mathrm{Drain}(\mathrm{Surger}(s\bgamma_0+(1-s)\bgamma)),
 $$
 where $\mathrm{Surger}$ is a certain canonical surgery procedure for resolving intersections, and $\mathrm{Drain}$ is a canonical procedure for eliminating one-sided relations.  From the constructions of~$\mathrm{Surger}$ and~$\mathrm{Drain}$ in~\cite{BBM10} it follows immediately that, if both $\bgamma_0$ and~$\bgamma$ have supports disjoint from~$N$, then the same is true throughout the deformation retraction $H(\bgamma,s)$. Thus the same formula provides a strong deformation retraction of~$\CB(v,N)$ onto a point, and therefore $\CB(v,N)$ is contractible.

 Let us now show that $\tCB$ is homotopy equivalent to~$\CB(v,N)$. First, we construct an embedding $\varphi\colon \tCB\hookrightarrow\CB(v,N)$ by sending every vertex $\gamma$ of~$\tCB$ with $[\gamma]=v+k_1a_1+k_2a_2$ to the basic cycle for~$v$ given by $\gamma-k_1\alpha_1-k_2\alpha_2$, and then extending this map linearly to simplices of~$\tCB$. Note that only nonnegative coefficients are allowed in a basic cycle; hence the summands~$-k_i\alpha_i$ should be interpreted as follows: if  $k_i>0$, then $-k_i\alpha_i$ should be replaced by~$k_i\overline{\alpha}_i$, where $\overline{\alpha}_i$ is the curve~$\alpha_i$ with the orientation reversed.

 Second, let us construct a strong deformation retraction of~$\CB(v,N)$ onto~$\varphi\bigl(\tCB\bigr)$. Consider an arbitrary point $\bgamma\in\CB(v,N)$. The homology class of any non-separating simple closed curve on~$\Sigma$ that is disjoint from~$N$ and whose algebraic intersection number with~$v$ is zero belongs to one of the sets
 $$c+\langle a_1,a_2\rangle,\qquad -c+\langle a_1,a_2\rangle,\qquad \langle a_1,a_2\rangle.$$ We write
 $$
 \bgamma=m_1\gamma_1+\cdots+m_k\gamma_k+n_1\delta_1+\cdots+n_l\delta_l+q_1\varepsilon_1+\cdots+q_p\varepsilon_p,
 $$
 where $[\gamma_i]\in c+\langle a_1,a_2\rangle$, $[\delta_i]\in -c+\langle a_1,a_2\rangle$, and $[\varepsilon_i]\in \langle a_1,a_2\rangle$ for all~$i$, while the coefficients~$m_i$, $n_i$, and~$q_i$ are positive reals. Since $\bgamma$ represents~$v$, we have
 $$
 m_1+\cdots+m_k=n_1+\cdots+n_l+1\ge 1.
 $$
 Also,
 $$
 [m_1\gamma_1+\cdots+m_k\gamma_k]=(m_1+\cdots+m_k)v+r_1a_1+r_2a_2
 $$
 for some $r_1,r_2\in\R$.
 A strong deformation retraction of~$\CB(v,N)$ onto~$\varphi\bigl(\tCB\bigr)$ is given by
 $$
 \Psi(\bgamma,s)= s\cdot\frac{m_1\gamma_1+\cdots+m_k\gamma_k-r_1\alpha_1-r_2\alpha_2}{m_1+\cdots+m_k}+(1-s)\bgamma,
 $$
 where the terms~$\alpha_i$ with negative coefficients are replaced by~$\overline{\alpha}_i$ with the corresponding positive coefficients, as above.
 Thus $\tCB$ is contractible.
\end{proof}

Let $\tCM$ denote the set of all oriented multicurves~$M=\gamma_0\cup\cdots\cup\gamma_p$ on~$\Sigma$ such that $\gamma_0,\ldots,\gamma_p$ span a simplex of~$\tCB$. We denote by~$\tP_M$ the simplex of~$\tCB$ corresponding to the oriented multicurve~$M$, and by~$\tCP_M$ the corresponding simplex of the quotient complex~$\tCB/\I_N$. For each~$p$, let $\tCM_p\subseteq\tCM$ be the subset consisting of all~$M$ such that $\dim\bigl(\tP_M\bigr)=p$. Since a $p$-simplex has $p+1$ vertices, every multicurve in~$\tCM_p$ has $p+1$ components. If we ignore the orientations of the components, then there are three topological types of multicurves in~$\tCM_1$ (see Fig.~\ref{fig_types_tM1}) and two topological types of multicurves in~$\tCM_2$ (see Fig.~\ref{fig_types_tM2}). We denote by $\tCM_p^{(i)}$ the subset of~$\tCM_p$ consisting of all oriented multicurves of type~$i$. In Figs.~\ref{fig_types_tM1} and~\ref{fig_types_tM2} the components of~$N$ are shown in red, and the components of~$M$ are shown in blue.

\begin{figure}[!ht]
 \begin{tikzpicture}[scale=.3,
        mid arrow/.style={
            postaction={
                decorate,
                decoration={
                    markings,
                    mark=at position 0.55 with {
                        \arrow[#1]{stealth}
                    }
                }
            }
        },
        mid arrow/.default=black,
    ]

\small
% \definecolor{blue}{rgb}{0, 0, 0.7}
% \definecolor{mygreen}{rgb}{0, 0.4, 0}
\tikzset{every path/.append style={line width=.2mm}}
\tikzset{my dash/.style={dash pattern=on 2pt off 1.5pt}}

% Border lines

\begin{scope}
\draw[color=red, thick, my dash] (-7.03,0) .. controls (-6,.3) .. (-4.97,0);
\draw[color=blue, thick, my dash] (-3.03,0) .. controls (-2,.3) .. (-.97,0);
\draw[color=blue, thick, my dash] (.97,0) .. controls (2,.3) .. (3.03,0);
\draw[color=red, thick, my dash] (4.97,0) .. controls (6,.3) .. (7.03,0);

\draw [very thick] (-4,0) circle (1);
\draw [very thick] (0,0) circle (1);
\draw [very thick] (4,0) circle (1);
\draw [very thick]  (-4,3)--(4,3) arc (90:-90:3) -- (-4,-3) arc (270:90:3);

% Visible lines

\draw[color=red, thick] (-7.03,0) .. controls (-6,-.3) .. (-4.97,0);
\draw[color=blue, thick] (-3.03,0) .. controls (-2,-.3) .. (-.97,0);
\draw[color=blue, thick] (.97,0) .. controls (2,-.3) .. (3.03,0);
\draw[color=red, thick] (4.97,0) .. controls (6,-.3) .. (7.03,0);

\draw (0, -5.5) node {type 1};
\end{scope}

\begin{scope}[xshift=17cm]
\draw[color=red, thick, my dash] (-7.03,0) .. controls (-6,.3) .. (-4.97,0);
\draw[color=blue, thick, my dash] (-3.03,0) .. controls (-2,.3) .. (-.97,0);
\draw[color=blue, thick, my dash] (0,.97) .. controls (.3,2) .. (0,3.03);
\draw[color=red, thick, my dash] (4.97,0) .. controls (6,.3) .. (7.03,0);

\draw [very thick] (-4,0) circle (1);
\draw [very thick] (0,0) circle (1);
\draw [very thick] (4,0) circle (1);
\draw [very thick]  (-4,3)--(4,3) arc (90:-90:3) -- (-4,-3) arc (270:90:3);

% Visible lines

\draw[color=red, thick] (-7.03,0) .. controls (-6,-.3) .. (-4.97,0);
\draw[color=blue, thick] (-3.03,0) .. controls (-2,-.3) .. (-.97,0);
\draw[color=blue, thick] (0,.97) .. controls (-.3,2) .. (0,3.03);
\draw[color=red, thick] (4.97,0) .. controls (6,-.3) .. (7.03,0);

\draw (0, -5.5) node {type 2};
\end{scope}

\begin{scope}[xshift=34cm]
\draw[color=red, thick, my dash] (-7.03,0) .. controls (-6,.3) .. (-4.97,0);
\draw[color=blue, thick, my dash] (0,-.97) .. controls (.3,-2) .. (0,-3.03);
\draw[color=blue, thick, my dash] (0,.97) .. controls (.3,2) .. (0,3.03);
\draw[color=red, thick, my dash] (4.97,0) .. controls (6,.3) .. (7.03,0);

\draw [very thick] (-4,0) circle (1);
\draw [very thick] (0,0) circle (1);
\draw [very thick] (4,0) circle (1);
\draw [very thick]  (-4,3)--(4,3) arc (90:-90:3) -- (-4,-3) arc (270:90:3);

% Visible lines

\draw[color=red, thick] (-7.03,0) .. controls (-6,-.3) .. (-4.97,0);
\draw[color=blue, thick] (0,-.97) .. controls (-.3,-2) .. (0,-3.03);
\draw[color=blue, thick] (0,.97) .. controls (-.3,2) .. (0,3.03);
\draw[color=red, thick] (4.97,0) .. controls (6,-.3) .. (7.03,0);

\draw (0, -5.5) node {type 3};
\end{scope}

\end{tikzpicture}
 \caption{Types of multicurves in~$\tCM_1$}\label{fig_types_tM1}
 \vspace{10mm}

 \begin{tikzpicture}[scale=.3,
        mid arrow/.style={
            postaction={
                decorate,
                decoration={
                    markings,
                    mark=at position 0.55 with {
                        \arrow[#1]{stealth}
                    }
                }
            }
        },
        mid arrow/.default=black,
    ]

\small
% \definecolor{blue}{rgb}{0, 0, 0.7}
% \definecolor{mygreen}{rgb}{0, 0.4, 0}
\tikzset{every path/.append style={line width=.2mm}}
\tikzset{my dash/.style={dash pattern=on 2pt off 1.5pt}}

% Border lines

\begin{scope}
\draw[color=red, thick, my dash] (-7.03,0) .. controls (-6,.3) .. (-4.97,0);
\draw[color=blue, thick, my dash] (-3.03,0) .. controls (-2,.3) .. (-.97,0);
\draw[color=blue, thick, my dash] (.97,0) .. controls (2,.3) .. (3.03,0);
\draw[color=red, thick, my dash] (4.97,0) .. controls (6,.3) .. (7.03,0);
\draw[color=blue, thick, my dash] (0,.97) .. controls (.3,2) .. (0,3.03);

\draw [very thick] (-4,0) circle (1);
\draw [very thick] (0,0) circle (1);
\draw [very thick] (4,0) circle (1);
\draw [very thick]  (-4,3)--(4,3) arc (90:-90:3) -- (-4,-3) arc (270:90:3);

% Visible lines

\draw[color=red, thick] (-7.03,0) .. controls (-6,-.3) .. (-4.97,0);
\draw[color=blue, thick] (-3.03,0) .. controls (-2,-.3) .. (-.97,0);
\draw[color=blue, thick] (.97,0) .. controls (2,-.3) .. (3.03,0);
\draw[color=red, thick] (4.97,0) .. controls (6,-.3) .. (7.03,0);
\draw[color=blue, thick] (0,.97) .. controls (-.3,2) .. (0,3.03);

\draw (0, -5.5) node {type 1};
\end{scope}

\begin{scope}[xshift=20cm]
\draw[color=red, thick, my dash] (-7.03,0) .. controls (-6,.3) .. (-4.97,0);
\draw[color=blue, thick, my dash] (-3.03,0) .. controls (-2,.3) .. (-.97,0);
\draw[color=blue, thick, my dash] (0,.97) .. controls (.3,2) .. (0,3.03);
\draw[color=red, thick, my dash] (4.97,0) .. controls (6,.3) .. (7.03,0);
\draw[color=blue, thick, my dash] (0,-.97) .. controls (.3,-2) .. (0,-3.03);

\draw [very thick] (-4,0) circle (1);
\draw [very thick] (0,0) circle (1);
\draw [very thick] (4,0) circle (1);
\draw [very thick]  (-4,3)--(4,3) arc (90:-90:3) -- (-4,-3) arc (270:90:3);

% Visible lines

\draw[color=red, thick] (-7.03,0) .. controls (-6,-.3) .. (-4.97,0);
\draw[color=blue, thick] (-3.03,0) .. controls (-2,-.3) .. (-.97,0);
\draw[color=blue, thick] (0,.97) .. controls (-.3,2) .. (0,3.03);
\draw[color=red, thick] (4.97,0) .. controls (6,-.3) .. (7.03,0);
\draw[color=blue, thick] (0,-.97) .. controls (-.3,-2) .. (0,-3.03);

\draw (0, -5.5) node {type 2};
\end{scope}

\end{tikzpicture}
 \caption{Types of multicurves in~$\tCM_2$}\label{fig_types_tM2}
\end{figure}

The group~$\I_N$ acts cellularly on~$\tCB$. The fact that~$\I_N$ acts on~$\tCB$ without rotations follows immediately from the fact that $\I$ acts on~$\CB(v)$ without rotations.
Moreover, the stabilizer of a cell~$\tP_M$ in~$\I_N$ is exactly the group~$\I_{M\cup N}$. Since $\tCB$ is contractible, we have the spectral sequence
\begin{equation}\label{eq_tSS}
\tE^1_{p,q}=\bigoplus_{M\in\widetilde{\mathfrak{M}}_p}H_q(\I_{M\cup N}) \  \Longrightarrow \  H_{p+q}(\I_N),
\end{equation}
where $\widetilde{\mathfrak{M}}_p$ is some set of representatives of $\I_N$-orbits in~$\tCM_p$; the differentials of this spectral sequence will be denoted by~$\td^r$. Let
$$
0\subseteq \tCF_{0,2}\subseteq\tCF_{1,1}\subseteq\tCF_{2,0}=H_2(\I_N)
$$
be the filtration associated with this spectral sequence, i.e., such that $$\tCF_{p,q}/\tCF_{p-1,q+1}=\tE^{\infty}_{p,q}.$$

Note that the simplicial complex~$\tCB$ does not change if we replace~$v$ by any class $v+k_1a_1+k_2a_2$ with $k_1,k_2\in\Z$. Let $\Theta$ be the subgroup of~$\Sp(H)$ consisting of all transformations that fix each of the classes~$a_1$ and~$a_2$ and send $v$ to a class in~$v+\langle a_1,a_2\rangle$. Let $\Mod(\Sigma,N,\Theta)$ be the subgroup of~$\Mod(\Sigma)$ consisting of all mapping classes that stabilize~$N$ and act on~$H$ via a transformation from~$\Theta$. Then the action of~$\I_N$ on~$\tCB$ naturally extends to a cellular action of~$\Mod(\Sigma,N,\Theta)$. Therefore \eqref{eq_tSS} is a spectral sequence of $\Z[\Theta]$-modules.

\subsection{Reduction to~$\tCF_{1,1}$}

We say that a multicurve $M\in\tCM$ (or the corresponding simplex~$\tP_M$ of~$\tCB$, or the corresponding simplex~$\tCP_M$ of~$\tCB/\I_N$) is \textit{important} if it contains a bounding pair, and \textit{unimportant} otherwise. Let $\tCC_*$ be the quotient of the chain complex $\tE^1_{*,0}=C_*\bigl(\tCB/\I_N\bigr)$ by the subcomplex generated by all unimportant simplices, and let $\partial_{\tCC}$ be the differential of~$\tCC_*$. Let us prove the following analogue of Proposition~\ref{propos_H2C0}.

\begin{propos}\label{propos_H2tC0}
We have $H_2\bigl(\tCC_*\bigr)=0$.
\end{propos}
\begin{proof}
The proof mimics the proof of Proposition~\ref{propos_H2C0}, but is actually much easier, since now we do not need to perform any induction on the function~$n$. A simplification comes from the fact that an important multicurve in~$\tCM_1$ is just a bounding pair~$\{\gamma,\gamma'\}$ with $[\gamma]=[\gamma']\in v+\langle a_1,a_2\rangle$, and therefore each important simplex of~$\tCB$ contains exactly one important one-dimensional simplex. Consequently, the chain complex $\tCC_*$ is the direct sum of subcomplexes $\tCC_*(K)$, where $K=\gamma\cup\gamma'$ runs over some set of representatives of $\I_N$-orbits in~$\tCM_1$, and $\tCC_*(K)$ is generated by all simplices~$\tCP_M$ that contain~$\tCP_K$. Thus it is sufficient to prove that $H_2\bigl(\tCC_*(K)\bigr)=0$ for every~$K$.

Let $\zeta\in\tCC_2(K)$ be a chain such that $\partial_{\tCC}\zeta=0$.  Let $S_1$ and~$S_2$ be the connected components of~$\Sigma\setminus K$ containing~$\alpha_1$ and~$\alpha_2$, respectively. We write
\begin{equation*}%\label{eq_sum_tbeps}
\zeta=\sum_{i=1}^mk_i\tCP_{M_i},
\end{equation*}
where every $2$-cell~$\tCP_{M_i}$ contains~$\tCP_K$. We choose the orientations of the cells~$\tCP_K$ and~$\tCP_{M_i}$ so that $\bigl[\tCP_{M_i}:\tCP_K\bigr]=1$ for all~$i$. Since $\partial_{\tCC}\zeta=0$, we have
\begin{equation}\label{eq_m_sum_t}
k_1+\cdots+k_m=0.
\end{equation}
Since $\tCP_{hM}=\tCP_M$ for all $h\in\I_N$, we may assume that $M_i\supseteq K$ for $i=1,\ldots,m$. Every oriented multicurve~$M_i$ has exactly one component~$\gamma_i$ that is not contained in~$K$. Assume that $\gamma_i\subseteq S_1$ for $1\le i\le m_1$ and $\gamma_i\subseteq S_2$ for $m_1+1\le i\le m$. As in the proof of Proposition~\ref{propos_H2C0}, after adding a ``phantom'' summand if necessary, we may assume that both numbers~$m_1$ and $m-m_1$ are nonzero. For each pair $(i,j)$ with $1\le i\le m_1$ and~$m_1+1\le j\le m$, the oriented multicurve $M_i\cup M_j=K\cup\gamma_i\cup\gamma_j$ belongs to~$\tCM_3$. We orient the $3$-simplex~$\tCP_{M_i\cup M_j}$ so that to obtain that $\bigl[\tCP_{M_i\cup M_j}:\tCP_{M_i}\bigr]=1$; then $\bigl[\tCP_{M_i\cup M_j}:\tCP_{M_j}\bigr]=-1$.

Now, choose integers $l_{ij}$, where $1\le i\le m_1$ and $m_1+1\le j\le m$, so that to obtain that
\begin{align*}
\sum_{j=m_1+1}^ml_{ij}&=k_i,& i&=1,\ldots,m_1,\\
\sum_{i=1}^{m_1}l_{ij}&=-k_j,&j&=m_1+1,\ldots,m.
\end{align*}
The existence of such numbers~$l_{ij}$ is provided by~\eqref{eq_m_sum_t}. Consider the $3$-chain
$$
\xi=\sum_{i=1}^{m_1}\sum_{j=m_1+1}^ml_{ij}\tCP_{M_i\cup M_j}.
$$
Then $\partial_{\tCC}\xi=\zeta$.
The proposition follows.
\end{proof}

\begin{propos}\label{propos_txtE20}
We have $(t_{a_1+a_2}-1)\tE^2_{2,0}=0$, and hence $(t_{a_1+a_2}-1)H_2(\I_N)\subseteq\tCF_{1,1}$.
\end{propos}

\begin{proof}
Let  $\zeta\in \tE^1_{2,0}=C_2\bigl(\CB_N/\I_N\bigr)$ be a chain satisfying~$\td^1\zeta=0$.  From Proposition~\ref{propos_H2tC0} it follows that there exists $\eta\in \tE^1_{3,0}$ such that $\zeta+\td^1\eta$ is the sum of cells~$\tCP_M$ corresponding to unimportant multicurves $M\in\tCM_2$ only. Hence it suffices to show that $(t_{a_1+a_2}-1)\tCP_M=0$ for any such cell. For each unimportant multicurve $M\in\tCM_2$, there exists a connected component of $\Sigma\setminus (M\cup N)$ that is adjacent to both~$\alpha_1$ and~$\alpha_2$. Therefore there exists a simple closed curve $\varepsilon$ that is disjoint from~$M\cup N$ and whose homology class is either $a_1+a_2$ or~$a_1-a_2$. Then $T_{\varepsilon}$ fixes~$\tCP_M$. Since $t_{a_1-a_2}=t_{a_1+a_2}t_{a_1}^{-2}t_{a_2}^{-2}$ and the twists~$T_{\alpha_1}$ and~$T_{\alpha_2}$ also fix~$\tCP_M$, we obtain that $t_{a_1+a_2}\tCP_M=\tCP_M$ in both cases.
\end{proof}

\subsection{Reduction to~$\tCF_{0,2}$}
Let us define the analogues of the chain complex~$\CE_*$ and the chain map~$\varphi$ in our situation. We have
\begin{equation}\label{eq_tE1p1}
\tE^1_{p,1}=\bigoplus_{M\in\widetilde{\mathfrak{M}}_p}H_1(\I_{M\cup N}).
\end{equation}
Let $\tCE_p\subseteq \tE^1_{p,1}$ be the sum of all summands in~\eqref{eq_tE1p1} that correspond to unimportant multicurves~$M$, i.e., multicurves~$M$ that contain no bounding pairs. We have $\td^1\tCE_{p}\subseteq\tCE_{p-1}$. Hence $\tCE_*$ is a chain subcomplex of~$\tE^1_{*,1}$.

As in Section~\ref{section_active}, we say that a component~$\gamma$ of an unimportant oriented multicurve $M\in\tCM_p$ is active if there exists a simple closed curve $\gamma'\subseteq \Sigma\setminus(M\cup N)$  such that $\{\gamma,\gamma'\}$ is a bounding pair and $M\cup\gamma'\in\tCM_{p+1}$. We denote by~$\Act(M)$ the set of active components of~$M$. The cell $\tCP_{M,\gamma,\CW}$, where $\gamma\in\Act(M)$ and $\CW\in\Spl(\gamma,M\cup N)$, is defined in exactly the same way as the cell $\CP_{M,\gamma,\CW}$ was defined in Section~\ref{section_active}. Namely, $\tCP_{M,\gamma,\CW}=\tCP_{M\cup\gamma'}$, where $\gamma'$ is as above and $\CW(\gamma')=\CW$.

Now let $\{\gamma,\gamma'\}$ be a bounding pair on~$\Sigma$ disjoint from $N$. Assume that the curves~$\gamma$ and~$\gamma'$ are oriented so that $[\gamma]=[\gamma']$. Then we denote by~$S_1$ the connected component of~$\Sigma\setminus(\gamma\cup\gamma')$ that contains~$\alpha_1$, and by~$S_2$ the connected component that contains~$\alpha_1$. Now, of the curves~$\gamma$ and~$\gamma'$, we call \textit{first} the one such that $S_1$ lies to its right, and \textit{second} the other one. This convention is simpler than the one in Section~\ref{section_ephi}, but it depends on the orientations of the curves.

We now define a homomorphism $\tphi\colon\tCE_p\to\tCC_{p+1}$ by setting
\begin{equation}\label{eq_tphi}
 \tphi\bigl(\iota_{\tP_M}([h])\bigr)=\sum_{\gamma\in\Act(M)}\,\sum_{\CW\in\Spl(\gamma,M\cup N)}\nu_{\gamma,\CW}(h)\tCP_{M,\gamma,\CW}
\end{equation}
for all $h\in\I_{M\cup N}$. Here the orientation of each cell $\tCP_{M,\gamma,\CW}=\tCP_{M\cup\gamma'}$ is chosen as follows. We orient~$\tP_{M\cup\gamma'}$ so that $\bigl[\tP_{M\cup\gamma'}:\tP_M\bigr]=1$ if $\gamma$ is the first component of the bounding pair $\{\gamma,\gamma'\}$ and $\bigl[\tP_{M\cup\gamma'}:\tP_M\bigr]=-1$ if $\gamma$ is the second component of this bounding pair. Then the orientation of~$\tCP_{M,\gamma,\CW}=\tCP_{M\cup\gamma'}$ is induced by this orientation of~$\tP_{M\cup\gamma'}$. An analogue of Proposition~\ref{propos_Pwelldef} is proved  in exactly the same way, hence the orientation of~$\tCP_{M,\gamma,\CW}$ is well-defined. Unlike in the construction of~$\varphi$ in Section~\ref{section_ephi}, the construction of~$\tphi$ contains no uncanonical choices. Therefore, $\tphi$ is a $\Z[\Theta]$-module homomorphism. Further, Proposition~\ref{propos_phi_chain} and Corollary~\ref{cor_phi_surj}, together with their proofs, carry over verbatim to the case under consideration. Thus we obtain the following.

\begin{propos}\label{propos_tphi_chain}
 The homomorphism $\tphi\colon\tCE_*\to\tCC_{*+1}$ is a surjective degree $1$ chain map, that is, $\tphi \td^1=-\partial_{\tCC}\tphi$.
\end{propos}

Let us study the kernel~$\tCK$ of the map $\tphi\colon\tCE_1\to\tCC_2$. As in Section~\ref{section_KNDZ}, if $K\in\tCM_2$ is an important multicurve with a bounding pair~$\{\gamma,\gamma'\}$ and $\delta$ is a separating simple closed curve disjoint from~$K$, then we set
\begin{gather}
\label{eq_tDelta_K}
\tDelta_K=\iota_{\tP_{M_1}}\bigl(\bigl[T_{\gamma'}T_{\gamma}^{-1}\bigr]\bigr)+\iota_{\tP_{M_2}}\bigl(\bigl[T_{\gamma}T_{\gamma'}^{-1}\bigr]\bigr),\\
\label{eq_tz_K}
\tz_{K,\delta}=\iota_{\tP_{M_1}}\bigl([T_{\delta}]\bigr)+\iota_{\tP_{M_2}}\bigl([T_{\delta}]\bigr),
\end{gather}
where $M_1=K\setminus\gamma'$, $M_2=K\setminus\gamma$, and the simplices~$\tP_{M_1}$ and~$\tP_{M_2}$ are oriented so that $\bigl[\tP_K:\tP_{M_1}\bigr]=\bigl[\tP_K:\tP_{M_2}\bigr]=1$.
Let $\tD$ be the subgroup of~$\tCE_1$ generated by all elements~$\tDelta_K$, and let $\tCZ$ be the subgroup of~$\tCE_1$ generated by all elements~$\tz_{K,\delta}$. It follows immediately from the constructions that $\tCK$, $\tD$, and~$\tCZ$ are $\Z[\Theta]$-submodules of~$\tCE_1$.

Let $\tCN\subseteq\tCE_1$ be the sum of all summands in~\eqref{eq_tE1p1} corresponding to the neglectible oriented multicurves~$M$, i.e., the multicurves that contain no bounding pairs and  no active components. Note  that a multicurve in~$\tCM_1$ is neglectible if and only if it has type~1. Hence, Proposition~\ref{propos_type1_Johnson} immediately implies the following analogue of Proposition~\ref{propos_N_anih}.

\begin{propos}\label{propos_tN_anih}
 We have $(t_{a_1+a_2}-1)^2\tCN=0$.
\end{propos}

However, a direct analogue of Proposition~\ref{propos_phi_ker} no longer holds: namely, $\tCK$ does not coincide with the sum~$\tD+\tCN$. Instead, $\tCK$ contains one additional summand, which we now describe. This summand appear because the components of~$N$ cannot be active, and therefore a multicurve $M\in\tCM_1^{(2)}$ has only two  active components rather than three.

Let $M\in\tCM_1^{(2)}$, and let $i$ be the index for which the component~$\alpha_i$ of~$N$ is adjacent to the pair-of-pants component of~$\Sigma\setminus(M\cup N)$. Let $\alpha'$ be a simple closed curve disjoint from~$M\cup N$ such that $\{\alpha_i,\alpha'\}$ is a bounding pair. Then we set
$$
\tq_{M,\alpha'}=\iota_{\tP_M}\bigl(\bigl[T_{\alpha'}T_{\alpha_i}^{-1}\bigr]\bigr).
$$
Let $\tCQ\subseteq\tCE_1$ be the subgroup generated by all the elements~$\tq_{M,\alpha'}$.

\begin{propos}\label{propos_tphi_ker}
 We have  $\tCK=\tD+\tCQ+\tCN$.
\end{propos}

\begin{proof}
As in the proof of Proposition~\ref{propos_phi_ker}, we have $\tCK\supseteq\tD+\tCQ+\tCN$. We need to prove the inverse inclusion. Further, as in the proof of Proposition~\ref{propos_phi_ker}, it suffices to show that $$\ker\bigl(\varphi\big|_{\tCE^{(2)}_1(C)}\bigr)\subseteq\tD+\tCQ$$
for each set~$C$ that occurs as~$[M]$ for $M\in\tCM_1^{(2)}$, where $\tCE^{(2)}_1(C)$ is the sum of the  summands~$H_1(\I_{M\cup N})$ in~\eqref{eq_tE1p1} corresponding to the $\I_N$-orbits of  multicurves~$M\in\tCM_1^{(2)}$ with $[M]=C$.

A set~$C$ occurs as~$[M]$ for $M\in\tCM_1^{(2)}$ if and only if $C=\{c,c+a_i\}$ for some class $c=v+k_1a_1+k_2a_2$, where $k_1,k_2\in\Z$, and $i$ is either~$1$ or~$2$. Without loss of generality we may assume that $i=1$. From Proposition~5.2 in~\cite{Gai21} it follows that there are exactly two $\I$-orbits of four-component oriented multicurves~$K$ with $[K]=C\cup\{a_1,a_2\}$. It follows that there are exactly two $\I_N$-orbits of oriented multicurves~$M\in\tCM_1^{(2)}$ with $[M]=C$; let $M_1$ and $M_2$ be representatives of these orbits. The two orbits are distinguished by the following property: the pair-of-pants component of~$\Sigma\setminus (M_1\cup N)$ lies to the left of the components of~$M_1\cup N$ with homology classes~$a_1$ and~$c$, and to the right of the component with homology class~$c+a_1$, while the pair-of-pants component of~$\Sigma\setminus (M_2\cup N)$ lies to the right of the components of~$M_2\cup N$ with homology classes~$a_1$ and~$c$, and to the left of the component with homology class~$c+a_1$. Then
 $$
 \tCE^{(2)}_1(C)\cong H_1(\I_{M_1\cup N})\oplus H_1(\I_{M_2\cup N}).
 $$
As in the proof of Proposition~\ref{propos_phi_ker}, we choose orientations of~$P_{M_1}$ and~$P_{M_2}$ so that they \textit{disagree} under the isomorphism $P_{M_1}\to P_{M_2}$ induced by the bijection between the components of~$M_1$ and~$M_2$ that sends every component of~$M_1$ to the component of~$M_2$ in the same homology class.

 Let $\gamma_1$ and $\gamma_2$ be the components of~$M_1$ with the homology classes~$c$ and $c+a_1$, respectively, and let $\widehat{\gamma}_1$ and~$\widehat{\gamma}_2$ be the components of~$M_2$ with the same homology classes. Let $\bigl\{e_{\alpha_1,\CW},e_{\gamma_1,\CW},e_{\gamma_2,\CW}\bigr\}$ and $\bigl\{\widehat{e}_{\alpha_1,\CW},\widehat{e}_{\widehat{\gamma}_1,\CW},\widehat{e}_{\widehat{\gamma}_2,\CW}\bigr\}$ be the bases of the groups~$H_1(\I_{M_1\cup N})$ and~$H_1(\I_{M_2\cup N})$, respectively, given by Proposition~\ref{propos_type2_ab}. Note that $\Spl(\alpha_1,M_1\cup N)=\Spl(\alpha_1,M_2\cup N)$ and $\Spl(\gamma_i,M_1\cup N)=\Spl(\widehat{\gamma}_i,M_2\cup N)$ for $i=1,2$.
 We set
 \begin{align*}
  e_{0,\CW}&=\iota_{\tP_{M_1}}\bigl(e_{\alpha_1,\CW}\bigr),
  &
  \widehat{e}_{0,\CW}&=\iota_{\tP_{M_2}}\bigl(\widehat{e}_{\alpha_1,\CW}\bigr),
  && \CW\in \Spl(\alpha_1,M_1\cup N),
  \\
  e_{i,\CW}&=\iota_{\tP_{M_1}}\bigl(e_{\gamma_i,\CW}\bigr),&
  \widehat{e}_{i,\CW}&=\iota_{\tP_{M_2}}\bigl(\widehat{e}_{\widehat{\gamma}_i,\CW}\bigr),
  &&i=1,2,\ \ \CW\in \Spl(\gamma_i,M_1\cup N).
 \end{align*}
These elements together form a basis of~$\tCE^{(2)}_1(C)$.
By~\eqref{eq_tphi}, we have:
 \begin{gather*}
  \tphi(e_{0,\CW})=\tphi(\widehat{e}_{0,\CW})=0,\\
  \tphi(e_{i,\CW})=\tCP_{M_1,\gamma_i,\CW},\qquad
  \tphi(\widehat{e}_{i,\CW})=\tCP_{M_2,\widehat{\gamma}_i,\CW},\qquad i=1,2.
 \end{gather*}
 The elements~$e_{0,\CW}$ and~$\widehat{e}_{0,\CW}$ lie in~$\tCQ$. Further, the cells $\tCP_{M_1,\gamma_i,\CW}$ with $i=1,2$ and $\CW\in\Spl(\gamma_i,M_1\cup N)$ are pairwise distinct, and therefore the elements~$\tphi(e_{i,\CW})$ with $i=1,2$ are linearly independent in~$\tCC_2$. Moreover, verbatim as in the proof of Proposition~\ref{propos_phi_ker}, we obtain that the elements $e_{i,\CW}+\widehat{e}_{i,\CW}$ with $i=1,2$ lie in~$\tD$, and hence in $\ker\Bigl(\varphi\big|_{\CE^{(2)}_1(A)}\Bigr)$. It follows that $\ker\Bigl(\varphi\big|_{\CE^{(2)}_1(A)}\Bigr)$ is generated by the elements~$e_{0,\CW}$, $\widehat{e}_{0,\CW}$, and
$e_{i,\CW}+\widehat{e}_{i,\CW}$ with $i=1,2$. The proposition follows , since all these elements lie in~$\tD+\tCQ$.
 \end{proof}

 Further, let us prove the following analogue of Proposition~\ref{propos_K_to_Z}.
 \begin{propos}\label{propos_tK_to_tZ}
 We have $(t_{a_1+a_2}-1)^2\tCK\subseteq\tCZ$.
\end{propos}

\begin{proof}
 By Propositions~\ref{propos_tphi_ker} and~\ref{propos_tN_anih}, it suffices to prove that the groups $(t_{a_1+a_2}-1)^2\tD$ and $(t_{a_1+a_2}-1)^2\tCQ$ are contained in~$\tCZ$. From Proposition~\ref{propos_e_annih} it follows that $(t_{a_1+a_2}-1)\tCQ=0$.

 The proof of the inclusion $(t_{a_1+a_2}-1)^2\tD\subseteq\tCZ$ is similar to the proof of Proposition~\ref{propos_D_to_Z}. Namely, consider an arbitrary generator $\tDelta_K$ of~$\tD$ given by~\eqref{eq_tDelta_K}. Let
 $$
 K=\gamma\cup\gamma'\cup\gamma_1,
 $$
 where $\{\gamma,\gamma'\}$ is a bounding pair. Without loss of generality, we may assume that
 $\gamma_1$ and~$\alpha_1$ lie in the same connected component of~$M\setminus(\gamma\cup\gamma')$, and $\alpha_2$ lies in the other connected component. Ignore the original orientations of the components of~$K$, and orient them so that to obtain that $[\gamma]=[\gamma']=[\gamma_1]+a_1$. Set  $c=[\gamma_1]$. Choose a separating simple closed curve~$\delta$ disjoint from~$K$, and consider the element~$\tz_{K,\delta}$ given by~\eqref{eq_tz_K}.  Proposition~\ref{propos_BP_delta}  implies that
 $$
 (t_{a_1+a_2}-1)(t_{c+a_2}-1)\tDelta_K=(t_{a_1+a_2}-1)(t_{c+a_2}-1)\tz_{K,\delta}.
 $$
 Every of the transvections~$t_{a_1}$, $t_{a_2}$, $t_{c}$, $t_{c+a_1}$ fixes~$\tDelta_K$, since they are induced by twists about the components of either of the multicurves~$M_i\cup N$. Moreover, $t_{c+a_1+a_2}$ also fixes~$\tDelta_K$ by Proposition~\ref{propos_e_annih}. By Lemma~\ref{lem_transvec3}, we have
 $$
 t_{c+a_1+a_2}=t_{c+a_1}t_{c+a_2}t_{a_1+a_2}t_{c}^{-1}t_{a_1}^{-1}t_{a_2}^{-1}.
 $$
 Since all transvections entering this equality mutually commute and each of them, except~$t_{c+a_2}$ and~$t_{a_1+a_2}$, fixes~$\tDelta_K$, we obtain that $t_{c+a_2}\tDelta_K=t_{a_1+a_2}^{-1}\tDelta_K$, and hence
 $$
 (t_{a_1+a_2}-1)(t_{c+a_2}-1)\tDelta_K=-t_{a_1+a_2}^{-1}(t_{a_1+a_2}-1)^2\tDelta_K.
 $$
 Therefore $(t_{a_1+a_2}-1)^2\tDelta_K$ belongs to~$\tCZ$.
\end{proof}

Finally, the following analogue of Proposition~\ref{propos_txE11} is derived from Propositions~\ref{propos_H2tC0}, \ref{propos_tphi_chain}, and~\ref{propos_tK_to_tZ} in completely the same way as Proposition~\ref{propos_txE20} itself was derived from Propositions~\ref{propos_H2C0}, \ref{propos_phi_chain}, \ref{propos_K_to_Z}, and Corollary~\ref{cor_phi_surj} in Section~\ref{section_to_F02}. Specifically, consider the subgroup $\tCY\subseteq \tE^1_{1,1}$ generated by all elements
$$
\widetilde{y}_{\gamma\cup\gamma',\delta}=\iota_{\tP_{\gamma\cup\gamma'}}([T_{\delta}])
$$
such that $\gamma\cup\gamma'\in\tCM_1^{(3)}$ and $\delta$ is a separating simple closed curve disjoint from $\gamma\cup\gamma'\cup N$. Then Propositions~\ref{propos_annih_E111} and~\ref{propos_Y_Ab}, together with their proofs, carry over verbatim to the case under consideration. Thus we obtain the following.

\begin{propos}\label{propos_txtE11}
We have $(t_{a_1+a_2}-1)^2\tCF_{1,1}\subseteq\tCF_{0,2}+\bA_N$.
\end{propos}

\subsection{Proof of Proposition~\ref{propos_anih_2comp}}

 Each multicurve $M\in\tCM_0$ consists of a single component that together with $N$ does not separate~$\Sigma$. By assertion~(1) of Theorem~\ref{theorem_3}, the group $H_2(\I_{M\cup N})$ is generated by standard abelian cycles. Therefore, the group $\tCF_{0,2}$ is also generated by standard abelian cycles, i.e., contained in~$\bA_N$. Thus, by Propositions~\ref{propos_txtE20} and~\ref{propos_txtE11}, we obtain that $(t_{a_1+a_2}-1)^3H_2(\I_N)\subseteq\bA_N$.

\section{End of the proof of Theorem~\ref{theorem_modA}}\label{section_modA_end}

Let us return to the study of the spectral sequence~$E^*_{*,*}$ associated with the action of~$\I$ on the complex of cycles $\CB=\CB(x)$, and the corresponding filtration
$$
0\subseteq\CF_{0,2}\subseteq \CF_{1,1}\subseteq \CF_{2,0}=H_2(\I).
$$

\begin{propos}\label{propos_txE02}
 We have $(t_x-1)^3 \CF_{0,2}\subseteq\bA$.
\end{propos}

\begin{proof}
 The group $\CF_{0,2}=E^{\infty}_{0,2}$ is a quotient of the group
 $$
 E^1_{0,2}=\bigoplus_{M\in\mathfrak{M}_0}H_2(\I_M),
 $$
 where $\mathfrak{M}_0$ is a set of representatives of $\I$-orbits in~$\CM_0$. Therefore, $\CF_{0,2}$ is the sum of the images of the homomorphisms $i_M\colon H_2(\I_M)\to H_2(\I)$ induced by the  inclusions $\I_M\le \I$, where $M\in\mathfrak{M}_0$. The set $\mathfrak{M}_0$ contains one-, two-, and three-component multicurves that do not separate~$\Sigma$.

 If $M$ is three-component, then, by assertion~(1) of Theorem~\ref{theorem_3}, the group $H_2(\I_M)$ is generated by standard abelian cycles. Therefore, $i_M\bigl(H_2(\I_M)\bigr)\subseteq\bA$. It follows that $(t_x-1)^3i_M\bigl(H_2(\I_M)\bigr)\subseteq\bA$ as well, since $\bA$ is invariant under the action of~$t_x$.

 If $M=\alpha$ is one-component, then $[\alpha]=x$. Hence $t_x-1$ annihilates $H_2(\I_M)$.

 Finally, suppose that $M=\alpha_1\cup\alpha_2$ is two-component, and set $a_1=[\alpha_1]$ and $a_2=[\alpha_2]$. Then $x=n_1a_1+n_2a_2$ for some positive $n_1,n_2>0$. By Lemma~\ref{lem_transvec3}, we have
 $$
 t_x=t_{a_1+a_2}^{n_1n_2}t_{a_1}^{n_1^2-n_1n_2}t_{a_2}^{n_2^2-n_1n_2}.
 $$
 The twists~$T_{\alpha_1}$ and~$T_{\alpha_2}$, and therefore the transvections~$t_{a_1}$ and~$t_{a_2}$, act trivially on~$H_2(\I_M)$. Consequently, $t_x$ acts on~$H_2(\I_M)$ in the same way as $t_{a_1+a_2}^{n_1n_2}$ does. By Proposition~\ref{propos_anih_2comp}, we have
 $$
 (t_{a_1+a_2}-1)^3H_2(\I_M)\subseteq\bA_M.
 $$
Since $\bA_M$ is invariant under the action of~$t_{a_1+a_2}$, it follows that
 $$
 (t_x-1)^3H_2(\I_M) =\bigl(t_{a_1+a_2}^{n_1n_2}-1\bigr)^3H_2(\I_M)\subseteq\bA_M.
 $$
 Thus $(t_x-1)^3i_M\bigl(H_2(\I_M)\bigr)\subseteq\bA$. The proposition follows.
\end{proof}

Theorem~\ref{theorem_modA} follows immediately from Propositions~\ref{propos_txE20}, \ref{propos_txE11} and~\ref{propos_txE02}.

\end{document}